\documentclass{amsproc}
\usepackage{amssymb,amscd,verbatim}
\makeatletter
 
\def\diagram{\m@th\leftwidth=\z@ \rightwidth=\z@ \topheight=\z@
\botheight=\z@ \setbox\@picbox\hbox\bgroup}
 
\def\enddiagram{\egroup\wd\@picbox\rightwidth\unitlength
\ht\@picbox\topheight\unitlength \dp\@picbox\botheight\unitlength
\hskip\leftwidth\unitlength\box\@picbox}
 
\def\bfig{\begin{diagram}}
\def\efig{\end{diagram}}
\newcount\wideness \newcount\leftwidth \newcount\rightwidth
\newcount\highness \newcount\topheight \newcount\botheight
 
\def\ratchet#1#2{\ifnum#1<#2 \global #1=#2 \fi}
 
\def\putbox(#1,#2)#3{%
\horsize{\wideness}{#3} \divide\wideness by 2
{\advance\wideness by #1 \ratchet{\rightwidth}{\wideness}}
{\advance\wideness by -#1 \ratchet{\leftwidth}{\wideness}}
\vertsize{\highness}{#3} \divide\highness by 2
{\advance\highness by #2 \ratchet{\topheight}{\highness}}
{\advance\highness by -#2 \ratchet{\botheight}{\highness}}
\put(#1,#2){\makebox(0,0){$#3$}}}
 
\def\putlbox(#1,#2)#3{%
\horsize{\wideness}{#3}
{\advance\wideness by #1 \ratchet{\rightwidth}{\wideness}}
{\ratchet{\leftwidth}{-#1}}
\vertsize{\highness}{#3} \divide\highness by 2
{\advance\highness by #2 \ratchet{\topheight}{\highness}}
{\advance\highness by -#2 \ratchet{\botheight}{\highness}}
\put(#1,#2){\makebox(0,0)[l]{$#3$}}}
 
\def\putrbox(#1,#2)#3{%
\horsize{\wideness}{#3}
{\ratchet{\rightwidth}{#1}}
{\advance\wideness by -#1 \ratchet{\leftwidth}{\wideness}}
\vertsize{\highness}{#3} \divide\highness by 2
{\advance\highness by #2 \ratchet{\topheight}{\highness}}
{\advance\highness by -#2 \ratchet{\botheight}{\highness}}
\put(#1,#2){\makebox(0,0)[r]{$#3$}}}

\def\adjust[#1]{} 
 
\newcount \coefa
\newcount \coefb
\newcount \coefc
\newcount\tempcounta
\newcount\tempcountb
\newcount\tempcountc
\newcount\tempcountd
\newcount\xext
\newcount\yext
\newcount\xoff
\newcount\yoff
\newcount\gap%
\newcount\arrowtypea
\newcount\arrowtypeb
\newcount\arrowtypec
\newcount\arrowtyped
\newcount\arrowtypee
\newcount\height
\newcount\width
\newcount\xpos
\newcount\ypos
\newcount\run
\newcount\rise
\newcount\arrowlength
\newcount\halflength
\newcount\arrowtype
\newdimen\tempdimen
\newdimen\xlen
\newdimen\ylen
\newsavebox{\tempboxa}%
\newsavebox{\tempboxb}%
\newsavebox{\tempboxc}%
 
\newdimen\w@dth
 
\def\setw@dth#1#2{\setbox\z@\hbox{\m@th$#1$}\w@dth=\wd\z@
\setbox\@ne\hbox{\m@th$#2$}\ifnum\w@dth<\wd\@ne \w@dth=\wd\@ne \fi
\advance\w@dth by 1.2em}
 
\def\t@^#1_#2{\allowbreak\def\n@one{#1}\def\n@two{#2}\mathrel
{\setw@dth{#1}{#2}
\mathop{\hbox to \w@dth{\rightarrowfill}}\limits
\ifx\n@one\empty\else ^{\box\z@}\fi
\ifx\n@two\empty\else _{\box\@ne}\fi}}
\def\t@@^#1{\@ifnextchar_{\t@^{#1}}{\t@^{#1}_{}}}
\def\to{\@ifnextchar^{\t@@}{\t@@^{}}}
 
\def\t@left^#1_#2{\def\n@one{#1}\def\n@two{#2}\mathrel{\setw@dth{#1}{#2}
\mathop{\hbox to \w@dth{\leftarrowfill}}\limits
\ifx\n@one\empty\else ^{\box\z@}\fi
\ifx\n@two\empty\else _{\box\@ne}\fi}}
\def\t@@left^#1{\@ifnextchar_{\t@left^{#1}}{\t@left^{#1}_{}}}
\def\toleft{\@ifnextchar^{\t@@left}{\t@@left^{}}}
 
\def\two@^#1_#2{\allowbreak
\def\n@one{#1}\def\n@two{#2}\mathrel{\setw@dth{#1}{#2}
\mathop{\vcenter{\lineskip\z@\baselineskip\z@
                 \hbox to \w@dth{\rightarrowfill}%
                 \hbox to \w@dth{\rightarrowfill}}%
       }\limits
\ifx\n@one\empty\else ^{\box\z@}\fi
\ifx\n@two\empty\else _{\box\@ne}\fi}}
\def\tw@@^#1{\@ifnextchar _{\two@^{#1}}{\two@^{#1}_{}}}
\def\two{\@ifnextchar ^{\tw@@}{\tw@@^{}}}
 
\def\tofr@^#1_#2{\def\n@one{#1}\def\n@two{#2}\mathrel{\setw@dth{#1}{#2}
\mathop{\vcenter{\hbox to \w@dth{\rightarrowfill}\kern-1.7ex
                 \hbox to \w@dth{\leftarrowfill}}%
       }\limits
\ifx\n@one\empty\else ^{\box\z@}\fi
\ifx\n@two\empty\else _{\box\@ne}\fi}}
\def\t@fr@^#1{\@ifnextchar_ {\tofr@^{#1}}{\tofr@^{#1}_{}}}
\def\tofro{\@ifnextchar^ {\t@fr@}{\t@fr@^{}}}

\def\mon{\mathop{\m@th\hbox to
      14.6\P@{\lasyb\char'51\hskip-2.1\P@$\arrext$\hss
$\mathord\rightarrow$}}\limits} 
\def\leftmono{\mathrel{\m@th\hbox to
14.6\P@{$\mathord\leftarrow$\hss$\arrext$\hskip-2.1\P@\lasyb\char'50%
}}\limits} 
\mathchardef\arrext="0200       

\def\settypes(#1,#2,#3){\arrowtypea#1 \arrowtypeb#2 \arrowtypec#3}
\def\settoheight#1#2{\setbox\@tempboxa\hbox{#2}#1\ht\@tempboxa\relax}%
\def\settodepth#1#2{\setbox\@tempboxa\hbox{#2}#1\dp\@tempboxa\relax}%
\def\settokens`#1`#2`#3`#4`{%
     \def\tokena{#1}\def\tokenb{#2}\def\tokenc{#3}\def\tokend{#4}}
\def\setsqparms[#1`#2`#3`#4;#5`#6]{%
\arrowtypea #1
\arrowtypeb #2
\arrowtypec #3
\arrowtyped #4
\width #5
\height #6
}
\def\setpos(#1,#2){\xpos=#1 \ypos#2}

\def\settriparms[#1`#2`#3;#4]{\settripairparms[#1`#2`#3`1`1;#4]}%
 
\def\settripairparms[#1`#2`#3`#4`#5;#6]{%
\arrowtypea #1
\arrowtypeb #2
\arrowtypec #3
\arrowtyped #4
\arrowtypee #5
\width #6
\height #6
}
 
\def\resetparms{\settripairparms[1`1`1`1`1;500]\width 500}
 
\resetparms
 
\def\mvector(#1,#2)#3{
\put(0,0){\vector(#1,#2){#3}}%
\put(0,0){\vector(#1,#2){26}}%
}
\def\evector(#1,#2)#3{{
\arrowlength #3
\put(0,0){\vector(#1,#2){\arrowlength}}%
\advance \arrowlength by-30
\put(0,0){\vector(#1,#2){\arrowlength}}%
}}
 
\def\horsize#1#2{%
\settowidth{\tempdimen}{$#2$}%
#1=\tempdimen
\divide #1 by\unitlength
}
 
\def\vertsize#1#2{%
\settoheight{\tempdimen}{$#2$}%
#1=\tempdimen
\settodepth{\tempdimen}{$#2$}%
\advance #1 by\tempdimen
\divide #1 by\unitlength
}
 
\def\putvector(#1,#2)(#3,#4)#5#6{{%
\ifnum3<\arrowtype
\putdashvector(#1,#2)(#3,#4)#5\arrowtype
\else
\ifnum\arrowtype<-3
\putdashvector(#1,#2)(#3,#4)#5\arrowtype
\else
\xpos=#1
\ypos=#2
\run=#3
\rise=#4
\arrowlength=#5
\ifnum \arrowtype<0
    \ifnum \run=0
        \advance \ypos by-\arrowlength
    \else
        \tempcounta \arrowlength
        \multiply \tempcounta by\rise
        \divide \tempcounta by\run
        \ifnum\run>0
            \advance \xpos by\arrowlength
            \advance \ypos by\tempcounta
        \else
            \advance \xpos by-\arrowlength
            \advance \ypos by-\tempcounta
        \fi
    \fi
    \multiply \arrowtype by-1
    \multiply \rise by-1
    \multiply \run by-1
\fi
\ifcase \arrowtype
\or \put(\xpos,\ypos){\vector(\run,\rise){\arrowlength}}%
\or \put(\xpos,\ypos){\mvector(\run,\rise)\arrowlength}%
\or \put(\xpos,\ypos){\evector(\run,\rise){\arrowlength}}%
\fi\fi\fi
}}
 
\def\putsplitvector(#1,#2)#3#4{
\xpos #1
\ypos #2
\arrowtype #4
\halflength #3
\arrowlength #3
\gap 140
\advance \halflength by-\gap
\divide \halflength by2
\ifnum\arrowtype>0
   \ifcase \arrowtype
   \or \put(\xpos,\ypos){\line(0,-1){\halflength}}%
       \advance\ypos by-\halflength
       \advance\ypos by-\gap
       \put(\xpos,\ypos){\vector(0,-1){\halflength}}%
   \or \put(\xpos,\ypos){\line(0,-1)\halflength}%
       \put(\xpos,\ypos){\vector(0,-1)3}%
       \advance\ypos by-\halflength
       \advance\ypos by-\gap
       \put(\xpos,\ypos){\vector(0,-1){\halflength}}%
   \or \put(\xpos,\ypos){\line(0,-1)\halflength}%
       \advance\ypos by-\halflength
       \advance\ypos by-\gap
       \put(\xpos,\ypos){\evector(0,-1){\halflength}}%
   \fi
\else \arrowtype=-\arrowtype
   \ifcase\arrowtype
   \or \advance \ypos by-\arrowlength
       \put(\xpos,\ypos){\line(0,1){\halflength}}%
       \advance\ypos by\halflength
       \advance\ypos by\gap
       \put(\xpos,\ypos){\vector(0,1){\halflength}}%
   \or \advance \ypos by-\arrowlength
       \put(\xpos,\ypos){\line(0,1)\halflength}%
       \put(\xpos,\ypos){\vector(0,1)3}%
       \advance\ypos by\halflength
       \advance\ypos by\gap
       \put(\xpos,\ypos){\vector(0,1){\halflength}}%
   \or \advance \ypos by-\arrowlength
       \put(\xpos,\ypos){\line(0,1)\halflength}%
       \advance\ypos by\halflength
       \advance\ypos by\gap
       \put(\xpos,\ypos){\evector(0,1){\halflength}}%
   \fi
\fi
}
 
\def\putmorphism(#1)(#2,#3)[#4`#5`#6]#7#8#9{{%
\run #2
\rise #3
\ifnum\rise=0
  \puthmorphism(#1)[#4`#5`#6]{#7}{#8}#9%
\else\ifnum\run=0
  \putvmorphism(#1)[#4`#5`#6]{#7}{#8}#9%
\else
\setpos(#1)%
\arrowlength #7
\arrowtype #8
\ifnum\run=0
\else\ifnum\rise=0
\else
\ifnum\run>0
    \coefa=1
\else
   \coefa=-1
\fi
\ifnum\arrowtype>0
   \coefb=0
   \coefc=-1
\else
   \coefb=\coefa
   \coefc=1
   \arrowtype=-\arrowtype
\fi
\width=2
\multiply \width by\run
\divide \width by\rise
\ifnum \width<0  \width=-\width\fi
\advance\width by60
\if l#9 \width=-\width\fi
\putbox(\xpos,\ypos){#4}
{\multiply \coefa by\arrowlength
\advance\xpos by\coefa
\multiply \coefa by\rise
\divide \coefa by\run
\advance \ypos by\coefa
\putbox(\xpos,\ypos){#5} }%
{\multiply \coefa by\arrowlength
\divide \coefa by2
\advance \xpos by\coefa
\advance \xpos by\width
\multiply \coefa by\rise
\divide \coefa by\run
\advance \ypos by\coefa
\if l#9%
   \putrbox(\xpos,\ypos){#6}%
\else\if r#9%
   \putlbox(\xpos,\ypos){#6}%
\fi\fi }%
{\multiply \rise by-\coefc
\multiply \run by-\coefc
\multiply \coefb by\arrowlength
\advance \xpos by\coefb
\multiply \coefb by\rise
\divide \coefb by\run
\advance \ypos by\coefb
\multiply \coefc by70
\advance \ypos by\coefc
\multiply \coefc by\run
\divide \coefc by\rise
\advance \xpos by\coefc
\multiply \coefa by140
\multiply \coefa by\run
\divide \coefa by\rise
\advance \arrowlength by\coefa
\ifcase\arrowtype
\or \put(\xpos,\ypos){\vector(\run,\rise){\arrowlength}}%
\or \put(\xpos,\ypos){\mvector(\run,\rise){\arrowlength}}%
\or \put(\xpos,\ypos){\evector(\run,\rise){\arrowlength}}%
\fi}\fi\fi\fi\fi}}

\newcount\numbdashes \newcount\lengthdash \newcount\increment
 
\def\howmanydashes{
\numbdashes=\arrowlength \lengthdash=40
\divide\numbdashes by \lengthdash
\lengthdash=\arrowlength
\divide\lengthdash by \numbdashes
\increment=\lengthdash
\multiply\lengthdash by 3
\divide\lengthdash by 5
}
 
\def\putdashvector(#1)(#2,#3)#4#5{%
\ifnum#3=0 \putdashhvector(#1){#4}#5
\else
\ifnum#2=0
\putdashvvector(#1){#4}#5\fi\fi}
 
\def\putdashhvector(#1,#2)#3#4{{%
\arrowlength=#3 \howmanydashes
\multiput(#1,#2)(\increment,0){\numbdashes}%
{\vrule height .4pt width \lengthdash\unitlength}
\arrowtype=#4 \xpos=#1
\ifnum\arrowtype<0 \advance\arrowtype by 7 \fi
\ifcase\arrowtype
\or \advance\xpos by 10
    \put(\xpos,#2){\vector(-1,0){\lengthdash}}
    \advance\xpos by 40
    \put(\xpos,#2){\vector(-1,0){\lengthdash}}
\or \advance \xpos by 10
    \put(\xpos,#2){\vector(-1,0){\lengthdash}}
    \advance\xpos by  \arrowlength
    \advance\xpos by  -50
    \put(\xpos,#2){\vector(-1,0){\lengthdash}}
\or \advance\xpos by 10
    \put(\xpos,#2){\vector(-1,0){\lengthdash}}
\or \advance\xpos by \arrowlength
    \advance\xpos by -\lengthdash
    \put(\xpos,#2){\vector(1,0){\lengthdash}}
\or {\advance\xpos by 10
    \put(\xpos,#2){\vector(1,0){\lengthdash}}}
    \advance\xpos by \arrowlength
    \advance\xpos by -\lengthdash
    \put(\xpos,#2){\vector(1,0){\lengthdash}}
\or \advance\xpos by \arrowlength
    \advance\xpos by -\lengthdash
    \put(\xpos,#2){\vector(1,0){\lengthdash}}
    \advance\xpos by -40
    \put(\xpos,#2){\vector(1,0){\lengthdash}}
   \fi
}}
 
\def\putdashvvector(#1,#2)#3#4{{%
\arrowlength=#3 \howmanydashes
\ypos=#2 \advance\ypos by -\arrowlength
\multiput(#1,#2)(0,\increment){\numbdashes}%
    {\vrule width .4pt height \lengthdash\unitlength}
\arrowtype=#4 \ypos=#2
\ifnum\arrowtype<0 \advance\arrowtype by 7 \fi
\ifcase\arrowtype
\or \advance\ypos by \arrowlength \advance\ypos by -40
    \put(#1,\ypos){\vector(0,1){\lengthdash}}
    \advance\ypos by -40
    \put(#1,\ypos){\vector(0,1){\lengthdash}}
\or \advance\ypos by 10
    \put(#1,\ypos){\vector(0,1){\lengthdash}}
    \advance\ypos by \arrowlength \advance\ypos by -40
    \put(#1,\ypos){\vector(0,1){\lengthdash}}
\or \advance\ypos by \arrowlength \advance\ypos by -40
    \put(#1,\ypos){\vector(0,1){\lengthdash}}
\or \advance\ypos by 10
    \put(#1,\ypos){\vector(0,-1){\lengthdash}}
\or \advance\ypos by 10
    \put(#1,\ypos){\vector(0,-1){\lengthdash}}
    \advance\ypos by \arrowlength \advance\ypos by -40
    \put(#1,\ypos){\vector(0,-1){\lengthdash}}
\or \advance\ypos by 10
    \put(#1,\ypos){\vector(0,-1){\lengthdash}}
    \advance\ypos by 40
    \put(#1,\ypos){\vector(0,-1){\lengthdash}}
\fi
}}
 
\def\puthmorphism(#1,#2)[#3`#4`#5]#6#7#8{{%
\xpos #1
\ypos #2
\width #6
\arrowlength #6
\arrowtype=#7
\putbox(\xpos,\ypos){#3\vphantom{#4}}%
{\advance \xpos by\arrowlength
\putbox(\xpos,\ypos){\vphantom{#3}#4}}%
\horsize{\tempcounta}{#3}%
\horsize{\tempcountb}{#4}%
\divide \tempcounta by2
\divide \tempcountb by2
\advance \tempcounta by30
\advance \tempcountb by30
\advance \xpos by\tempcounta
\advance \arrowlength by-\tempcounta
\advance \arrowlength by-\tempcountb
\putvector(\xpos,\ypos)(1,0)\arrowlength\arrowtype
\divide \arrowlength by2
\advance \xpos by\arrowlength
\vertsize{\tempcounta}{#5}%
\divide\tempcounta by2
\advance \tempcounta by20
\if a#8 %
   \advance \ypos by\tempcounta
   \putbox(\xpos,\ypos){#5}%
\else
   \advance \ypos by-\tempcounta
   \putbox(\xpos,\ypos){#5}%
\fi}}
 
\def\putvmorphism(#1,#2)[#3`#4`#5]#6#7#8{{%
\xpos #1
\ypos #2
\arrowlength #6
\arrowtype #7
\settowidth{\xlen}{$#5$}%
\putbox(\xpos,\ypos){#3}%
{\advance \ypos by-\arrowlength
\putbox(\xpos,\ypos){#4}}%
{\advance\arrowlength by-140
\advance \ypos by-70
\ifdim\xlen>0pt
   \if m#8%
      \putsplitvector(\xpos,\ypos)\arrowlength\arrowtype
   \else
   \putvector(\xpos,\ypos)(0,-1)\arrowlength\arrowtype
   \fi
\else
   \putvector(\xpos,\ypos)(0,-1)\arrowlength\arrowtype
\fi}%
\ifdim\xlen>0pt
   \divide \arrowlength by2
   \advance\ypos by-\arrowlength
   \if l#8%
      \advance \xpos by-40
      \putrbox(\xpos,\ypos){#5}%
   \else\if r#8%
      \advance \xpos by40
      \putlbox(\xpos,\ypos){#5}%
   \else
      \putbox(\xpos,\ypos){#5}%
   \fi\fi
\fi
}}
 
\def\putsquarep<#1>(#2)[#3;#4`#5`#6`#7]{{%
\setsqparms[#1]%
\setpos(#2)%
\settokens`#3`%
\puthmorphism(\xpos,\ypos)[\tokenc`\tokend`{#7}]{\width}{\arrowtyped}b%
\advance\ypos by \height
\puthmorphism(\xpos,\ypos)[\tokena`\tokenb`{#4}]{\width}{\arrowtypea}a%
\putvmorphism(\xpos,\ypos)[``{#5}]{\height}{\arrowtypeb}l%
\advance\xpos by \width
\putvmorphism(\xpos,\ypos)[``{#6}]{\height}{\arrowtypec}r%
}}
 
\def\putsquare{\@ifnextchar <{\putsquarep}{\putsquarep%
   <\arrowtypea`\arrowtypeb`\arrowtypec`\arrowtyped;\width`\height>}}
\def\square{\@ifnextchar< {\squarep}{\squarep
   <\arrowtypea`\arrowtypeb`\arrowtypec`\arrowtyped;\width`\height>}}
\def\squarep<#1>[#2`#3`#4`#5;#6`#7`#8`#9]{{
\setsqparms[#1]
\diagram
\putsquarep<\arrowtypea`\arrowtypeb`\arrowtypec`
\arrowtyped;\width`\height>
(0,0)[#2`#3`#4`{#5};#6`#7`#8`{#9}]
\enddiagram
}}                                                 
\def\putptrianglep<#1>(#2,#3)[#4`#5`#6;#7`#8`#9]{{%
\settriparms[#1]%
\xpos=#2 \ypos=#3
\advance\ypos by \height
\puthmorphism(\xpos,\ypos)[#4`#5`{#7}]{\height}{\arrowtypea}a%
\putvmorphism(\xpos,\ypos)[`#6`{#8}]{\height}{\arrowtypeb}l%
\advance\xpos by\height
\putmorphism(\xpos,\ypos)(-1,-1)[``{#9}]{\height}{\arrowtypec}r%
}}
 
\def\putptriangle{\@ifnextchar <{\putptrianglep}{\putptrianglep
   <\arrowtypea`\arrowtypeb`\arrowtypec;\height>}}
\def\ptriangle{\@ifnextchar <{\ptrianglep}{\ptrianglep
   <\arrowtypea`\arrowtypeb`\arrowtypec;\height>}}
\def\ptrianglep<#1>[#2`#3`#4;#5`#6`#7]{{
\settriparms[#1]
\diagram
\putptrianglep<\arrowtypea`\arrowtypeb`
\arrowtypec;\height>
(0,0)[#2`#3`#4;#5`#6`{#7}]
\enddiagram
}}                                            
 
\def\putqtrianglep<#1>(#2,#3)[#4`#5`#6;#7`#8`#9]{{%
\settriparms[#1]%
\xpos=#2 \ypos=#3
\advance\ypos by\height
\puthmorphism(\xpos,\ypos)[#4`#5`{#7}]{\height}{\arrowtypea}a%
\putmorphism(\xpos,\ypos)(1,-1)[``{#8}]{\height}{\arrowtypeb}l%
\advance\xpos by\height
\putvmorphism(\xpos,\ypos)[`#6`{#9}]{\height}{\arrowtypec}r%
}}
 
\def\putqtriangle{\@ifnextchar <{\putqtrianglep}{\putqtrianglep
   <\arrowtypea`\arrowtypeb`\arrowtypec;\height>}}
\def\qtriangle{\@ifnextchar <{\qtrianglep}{\qtrianglep
   <\arrowtypea`\arrowtypeb`\arrowtypec;\height>}}
\def\qtrianglep<#1>[#2`#3`#4;#5`#6`#7]{{
\settriparms[#1]
\width=\height                                
\diagram
\putqtrianglep<\arrowtypea`\arrowtypeb`
\arrowtypec;\height>
(0,0)[#2`#3`#4;#5`#6`{#7}]
\enddiagram
}}
 
\def\putdtrianglep<#1>(#2,#3)[#4`#5`#6;#7`#8`#9]{{%
\settriparms[#1]%
\xpos=#2 \ypos=#3
\puthmorphism(\xpos,\ypos)[#5`#6`{#9}]{\height}{\arrowtypec}b%
\advance\xpos by \height \advance\ypos by\height
\putmorphism(\xpos,\ypos)(-1,-1)[``{#7}]{\height}{\arrowtypea}l%
\putvmorphism(\xpos,\ypos)[#4``{#8}]{\height}{\arrowtypeb}r%
}}
 
\def\putdtriangle{\@ifnextchar <{\putdtrianglep}{\putdtrianglep
   <\arrowtypea`\arrowtypeb`\arrowtypec;\height>}}
\def\dtriangle{\@ifnextchar <{\dtrianglep}{\dtrianglep
   <\arrowtypea`\arrowtypeb`\arrowtypec;\height>}}
\def\dtrianglep<#1>[#2`#3`#4;#5`#6`#7]{{
\settriparms[#1]
\width=\height                                
\diagram
\putdtrianglep<\arrowtypea`\arrowtypeb`
\arrowtypec;\height>
(0,0)[#2`#3`#4;#5`#6`{#7}]
\enddiagram
}}
 
\def\putbtrianglep<#1>(#2,#3)[#4`#5`#6;#7`#8`#9]{{%
\settriparms[#1]%
\xpos=#2 \ypos=#3
\puthmorphism(\xpos,\ypos)[#5`#6`{#9}]{\height}{\arrowtypec}b%
\advance\ypos by\height
\putmorphism(\xpos,\ypos)(1,-1)[``{#8}]{\height}{\arrowtypeb}r%
\putvmorphism(\xpos,\ypos)[#4``{#7}]{\height}{\arrowtypea}l%
}}
 
\def\putbtriangle{\@ifnextchar <{\putbtrianglep}{\putbtrianglep
   <\arrowtypea`\arrowtypeb`\arrowtypec;\height>}}
\def\btriangle{\@ifnextchar <{\btrianglep}{\btrianglep
   <\arrowtypea`\arrowtypeb`\arrowtypec;\height>}}
\def\btrianglep<#1>[#2`#3`#4;#5`#6`#7]{{
\settriparms[#1]
\width=\height                               
\diagram
\putbtrianglep<\arrowtypea`\arrowtypeb`
\arrowtypec;\height>
(0,0)[#2`#3`#4;#5`#6`{#7}]
\enddiagram
}}
 
\def\putAtrianglep<#1>(#2,#3)[#4`#5`#6;#7`#8`#9]{{%
\settriparms[#1]%
\xpos=#2 \ypos=#3
{\multiply \height by2
\puthmorphism(\xpos,\ypos)[#5`#6`{#9}]{\height}{\arrowtypec}b}%
\advance\xpos by\height \advance\ypos by\height
\putmorphism(\xpos,\ypos)(-1,-1)[#4``{#7}]{\height}{\arrowtypea}l%
\putmorphism(\xpos,\ypos)(1,-1)[``{#8}]{\height}{\arrowtypeb}r%
}}
 
\def\putAtriangle{\@ifnextchar <{\putAtrianglep}{\putAtrianglep
   <\arrowtypea`\arrowtypeb`\arrowtypec;\height>}}
\def\Atriangle{\@ifnextchar <{\Atrianglep}{\Atrianglep
   <\arrowtypea`\arrowtypeb`\arrowtypec;\height>}}
\def\Atrianglep<#1>[#2`#3`#4;#5`#6`#7]{{
\settriparms[#1]
\width=\height                                     
\diagram
\putAtrianglep<\arrowtypea`\arrowtypeb`
\arrowtypec;\height>
(0,0)[#2`#3`#4;#5`#6`{#7}]
\enddiagram
}}
 
\def\putAtrianglepairp<#1>(#2)[#3;#4`#5`#6`#7`#8]{{%
\settripairparms[#1]%
\setpos(#2)%
\settokens`#3`%
\puthmorphism(\xpos,\ypos)[\tokenb`\tokenc`{#7}]{\height}{\arrowtyped}b%
\advance\xpos by\height
\puthmorphism(\xpos,\ypos)[\phantom{\tokenc}`\tokend`{#8}]%
{\height}{\arrowtypee}b%
\advance\ypos by\height
\putmorphism(\xpos,\ypos)(-1,-1)[\tokena``{#4}]{\height}{\arrowtypea}l%
\putvmorphism(\xpos,\ypos)[``{#5}]{\height}{\arrowtypeb}m%
\putmorphism(\xpos,\ypos)(1,-1)[``{#6}]{\height}{\arrowtypec}r%
}}
 
\def\putAtrianglepair{\@ifnextchar <{\putAtrianglepairp}{\putAtrianglepairp%
   <\arrowtypea`\arrowtypeb`\arrowtypec`\arrowtyped`\arrowtypee;\height>}}
\def\Atrianglepair{\@ifnextchar <{\Atrianglepairp}{\Atrianglepairp%
   <\arrowtypea`\arrowtypeb`\arrowtypec`\arrowtyped`\arrowtypee;\height>}}
 
\def\Atrianglepairp<#1>[#2;#3`#4`#5`#6`#7]{{
\settripairparms[#1]
\settokens`#2`
\width=\height                                
\diagram
\putAtrianglepairp                            
<\arrowtypea`\arrowtypeb`\arrowtypec`
\arrowtyped`\arrowtypee;\height>
(0,0)[{#2};#3`#4`#5`#6`{#7}]
\enddiagram
}}
 
\def\putVtrianglep<#1>(#2,#3)[#4`#5`#6;#7`#8`#9]{{%
\settriparms[#1]%
\xpos=#2 \ypos=#3
\advance\ypos by\height
{\multiply\height by2
\puthmorphism(\xpos,\ypos)[#4`#5`{#7}]{\height}{\arrowtypea}a}%
\putmorphism(\xpos,\ypos)(1,-1)[`#6`{#8}]{\height}{\arrowtypeb}l%
\advance\xpos by\height
\advance\xpos by\height
\putmorphism(\xpos,\ypos)(-1,-1)[``{#9}]{\height}{\arrowtypec}r%
}}
 
\def\putVtriangle{\@ifnextchar <{\putVtrianglep}{\putVtrianglep
   <\arrowtypea`\arrowtypeb`\arrowtypec;\height>}}
\def\Vtriangle{\@ifnextchar <{\Vtrianglep}{\Vtrianglep
   <\arrowtypea`\arrowtypeb`\arrowtypec;\height>}}
\def\Vtrianglep<#1>[#2`#3`#4;#5`#6`#7]{{
\settriparms[#1]
\width=\height                                 
\diagram
\putVtrianglep<\arrowtypea`\arrowtypeb`
\arrowtypec;\height>
(0,0)[#2`#3`#4;#5`#6`{#7}]
\enddiagram
}}
 
\def\putVtrianglepairp<#1>(#2)[#3;#4`#5`#6`#7`#8]{{
\settripairparms[#1]%
\setpos(#2)%
\settokens`#3`%
\advance\ypos by\height
\putmorphism(\xpos,\ypos)(1,-1)[`\tokend`{#6}]{\height}{\arrowtypec}l%
\puthmorphism(\xpos,\ypos)[\tokena`\tokenb`{#4}]{\height}{\arrowtypea}a%
\advance\xpos by\height
\puthmorphism(\xpos,\ypos)[\phantom{\tokenb}`\tokenc`{#5}]%
{\height}{\arrowtypeb}a%
\putvmorphism(\xpos,\ypos)[``{#7}]{\height}{\arrowtyped}m%
\advance\xpos by\height
\putmorphism(\xpos,\ypos)(-1,-1)[``{#8}]{\height}{\arrowtypee}r%
}}
 
\def\putVtrianglepair{\@ifnextchar <{\putVtrianglepairp}{\putVtrianglepairp%
    <\arrowtypea`\arrowtypeb`\arrowtypec`\arrowtyped`\arrowtypee;\height>}}
\def\Vtrianglepair{\@ifnextchar <{\Vtrianglepairp}{\Vtrianglepairp%
    <\arrowtypea`\arrowtypeb`\arrowtypec`\arrowtyped`\arrowtypee;\height>}}
\def\Vtrianglepairp<#1>[#2;#3`#4`#5`#6`#7]{{
\settripairparms[#1]
\settokens`#2`
\diagram
\putVtrianglepairp                             
<\arrowtypea`\arrowtypeb`\arrowtypec`
\arrowtyped`\arrowtypee;\height>
(0,0)[{#2};#3`#4`#5`#6`{#7}]
\enddiagram
}}

\def\putCtrianglep<#1>(#2,#3)[#4`#5`#6;#7`#8`#9]{{%
\settriparms[#1]%
\xpos=#2 \ypos=#3
\advance\ypos by\height
\putmorphism(\xpos,\ypos)(1,-1)[``{#9}]{\height}{\arrowtypec}l%
\advance\xpos by\height
\advance\ypos by\height
\putmorphism(\xpos,\ypos)(-1,-1)[#4`#5`{#7}]{\height}{\arrowtypea}l%
{\multiply\height by 2
\putvmorphism(\xpos,\ypos)[`#6`{#8}]{\height}{\arrowtypeb}r}%
}}
 
\def\putCtriangle{\@ifnextchar <{\putCtrianglep}{\putCtrianglep
    <\arrowtypea`\arrowtypeb`\arrowtypec;\height>}}
\def\Ctriangle{\@ifnextchar <{\Ctrianglep}{\Ctrianglep
    <\arrowtypea`\arrowtypeb`\arrowtypec;\height>}}
\def\Ctrianglep<#1>[#2`#3`#4;#5`#6`#7]{{
\settriparms[#1]
\width=\height                               
\diagram
\putCtrianglep<\arrowtypea`\arrowtypeb`
\arrowtypec;\height>
(0,0)[#2`#3`#4;#5`#6`{#7}]
\enddiagram
}}                                           
\def\putDtrianglep<#1>(#2,#3)[#4`#5`#6;#7`#8`#9]{{%
\settriparms[#1]%
\xpos=#2 \ypos=#3
\advance\xpos by\height \advance\ypos by\height
\putmorphism(\xpos,\ypos)(-1,-1)[``{#9}]{\height}{\arrowtypec}r%
\advance\xpos by-\height \advance\ypos by\height
\putmorphism(\xpos,\ypos)(1,-1)[`#5`{#8}]{\height}{\arrowtypeb}r%
{\multiply\height by 2
\putvmorphism(\xpos,\ypos)[#4`#6`{#7}]{\height}{\arrowtypea}l}%
}}
 
\def\putDtriangle{\@ifnextchar <{\putDtrianglep}{\putDtrianglep
    <\arrowtypea`\arrowtypeb`\arrowtypec;\height>}}
\def\Dtriangle{\@ifnextchar <{\Dtrianglep}{\Dtrianglep
   <\arrowtypea`\arrowtypeb`\arrowtypec;\height>}}
\def\Dtrianglep<#1>[#2`#3`#4;#5`#6`#7]{{
\settriparms[#1]
\width=\height                              
\diagram
\putDtrianglep<\arrowtypea`\arrowtypeb`
\arrowtypec;\height>
(0,0)[#2`#3`#4;#5`#6`{#7}]
\enddiagram
}}                                          
\def\setrecparms[#1`#2]{\width=#1 \height=#2}%
 
\def\recursep<#1`#2>[#3;#4`#5`#6`#7`#8]{{\m@th
\width=#1 \height=#2
\settokens`#3`
\settowidth{\tempdimen}{$\tokena$}
\ifdim\tempdimen=0pt
  \savebox{\tempboxa}{\hbox{$\tokenb$}}%
  \savebox{\tempboxb}{\hbox{$\tokend$}}%
  \savebox{\tempboxc}{\hbox{$#6$}}%
\else
  \savebox{\tempboxa}{\hbox{$\hbox{$\tokena$}\times\hbox{$\tokenb$}$}}%
  \savebox{\tempboxb}{\hbox{$\hbox{$\tokena$}\times\hbox{$\tokend$}$}}%
  \savebox{\tempboxc}{\hbox{$\hbox{$\tokena$}\times\hbox{$#6$}$}}%
\fi
\ypos=\height
\divide\ypos by 2
\xpos=\ypos
\advance\xpos by \width
\bfig
\putCtrianglep<-1`1`1;\ypos>(0,0)[`\tokenc`;#5`#6`{#7}]%
\puthmorphism(\ypos,0)[\tokend`\usebox{\tempboxb}`{#8}]{\width}{-1}b%
\puthmorphism(\ypos,\height)[\tokenb`\usebox{\tempboxa}`{#4}]{\width}{-1}a%
\advance\ypos by \width
\putvmorphism(\ypos,\height)[``\usebox{\tempboxc}]{\height}1r%
\efig
}}
 
\def\recurse{\@ifnextchar <{\recursep}{\recursep<\width`\height>}}
 
\def\puttwohmorphisms(#1,#2)[#3`#4;#5`#6]#7#8#9{{%
\puthmorphism(#1,#2)[#3`#4`]{#7}0a
\ypos=#2
\advance\ypos by 20
\puthmorphism(#1,\ypos)[\phantom{#3}`\phantom{#4}`#5]{#7}{#8}a
\advance\ypos by -40
\puthmorphism(#1,\ypos)[\phantom{#3}`\phantom{#4}`#6]{#7}{#9}b
}}
 
\def\puttwovmorphisms(#1,#2)[#3`#4;#5`#6]#7#8#9{{%
\putvmorphism(#1,#2)[#3`#4`]{#7}0a
\xpos=#1
\advance\xpos by -20
\putvmorphism(\xpos,#2)[\phantom{#3}`\phantom{#4}`#5]{#7}{#8}l
\advance\xpos by 40
\putvmorphism(\xpos,#2)[\phantom{#3}`\phantom{#4}`#6]{#7}{#9}r
}}
 
\def\puthcoequalizer(#1)[#2`#3`#4;#5`#6`#7]#8#9{{%
\setpos(#1)%
\puttwohmorphisms(\xpos,\ypos)[#2`#3;#5`#6]{#8}11%
\advance\xpos by #8
\puthmorphism(\xpos,\ypos)[\phantom{#3}`#4`#7]{#8}1{#9}
}}
 
\def\putvcoequalizer(#1)[#2`#3`#4;#5`#6`#7]#8#9{{%
\setpos(#1)%
\puttwovmorphisms(\xpos,\ypos)[#2`#3;#5`#6]{#8}11%
\advance\ypos by -#8
\putvmorphism(\xpos,\ypos)[\phantom{#3}`#4`#7]{#8}1{#9}
}}
 
\def\putthreehmorphisms(#1)[#2`#3;#4`#5`#6]#7(#8)#9{{%
\setpos(#1) \settypes(#8)
\if a#9 %
     \vertsize{\tempcounta}{#5}%
     \vertsize{\tempcountb}{#6}%
     \ifnum \tempcounta<\tempcountb \tempcounta=\tempcountb \fi
\else
     \vertsize{\tempcounta}{#4}%
     \vertsize{\tempcountb}{#5}%
     \ifnum \tempcounta<\tempcountb \tempcounta=\tempcountb \fi
\fi
\advance \tempcounta by 60
\puthmorphism(\xpos,\ypos)[#2`#3`#5]{#7}{\arrowtypeb}{#9}
\advance\ypos by \tempcounta
\puthmorphism(\xpos,\ypos)[\phantom{#2}`\phantom{#3}`#4]{#7}{\arrowtypea}{#9}
\advance\ypos by -\tempcounta \advance\ypos by -\tempcounta
\puthmorphism(\xpos,\ypos)[\phantom{#2}`\phantom{#3}`#6]{#7}{\arrowtypec}{#9}
}}
 
\def\setarrowtoks[#1`#2`#3`#4`#5`#6]{%
\def\toka{#1}
\def\tokb{#2}
\def\tokc{#3}
\def\tokd{#4}
\def\toke{#5}
\def\tokf{#6}
}
\def\hex{\@ifnextchar <{\hexp}{\hexp<1000`400>}}
\def\hexp<#1`#2>[#3`#4`#5`#6`#7`#8;#9]{%
\setarrowtoks[#9]
\yext=#2 \advance \yext by #2
\xext=#1 \advance\xext by \yext
\bfig
\putCtriangle<-1`0`1;#2>(0,0)[`#5`;\tokb``\tokd]
\xext=#1 \yext=#2 \advance \yext by #2
\putsquare<1`0`0`1;\xext`\yext>(#2,0)[#3`#4`#7`#8;\toka```\tokf]
\advance \xext by #2
\putDtriangle<0`1`-1;#2>(\xext,0)[`#6`;`\tokc`\toke]
\efig
}
\makeatother

\usepackage[all]{xy}
\def \PP{{\mathbb P}}

\def \ZZ{{\mathbb Z}}

\def \NN{{\mathbb N}}

\def\Ascr{{\mathcal A}}
\def\Bscr{{\mathcal B}}
\def\Cscr{{\mathcal C}}
\def\Dscr{{\mathcal D}}
\def\Escr{{\mathcal E}}
\def\Fscr{{\mathcal F}}
\def\Gscr{{\mathcal G}}
\def\Hscr{{\mathcal H}}
\def\Iscr{{\mathcal I}}
\def\Jscr{{\mathcal J}}
\def\Kscr{{\mathcal K}}
\def\Lscr{{\mathcal L}}
\def\Mscr{{\mathcal M}}
\def\Nscr{{\mathcal N}}
\def\Oscr{{\mathcal O}}
\def\Pscr{{\mathcal P}}
\def\Qscr{{\mathcal Q}}

\def\Sscr{{\mathcal S}}
\def\Tscr{{\mathcal T}}
\def\Uscr{{\mathcal U}}
\def\Vscr{{\mathcal V}}

\def\Xscr{{\mathcal X}}

\DeclareMathOperator{\cone}{cone}

\DeclareMathOperator{\RHS}{RHS}
\DeclareMathOperator{\Dis}{Dis}
\DeclareMathOperator{\Aut}{Aut}

\def\Lotimes{\overset{L}{\otimes}}
\def\HHom{\operatorname {\Hscr \mathit{om}}}
\def\RHHom{\operatorname {R\Hscr \mathit{om}}}
\def\Lotimes{\overset{L}{\otimes}}
\def\HExt{\operatorname {\Escr \mathit{xt}}}
\def\uHHom{\underline{\HHom}}

\def\HTor{\operatorname {\Tscr \mathit{or}}}

\def\cd{\mathop{\text{cd}}}
\def\Ab{\mathbf{Ab}}

\def\Ann{\mathop{\text{Ann}}\nolimits}
\def\id{\text{id}}
\def\Id{\text{id}}

\def\ctimes{\mathbin{\hat{\otimes}}}

\def\Mod{\mathop{\text{Mod}}}
\def\mod{\mathop{\text{mod}}}
\def\MOD{\mathop{\text{MOD}}}
\def\Bimod{\mathop{\text{Bimod}}}
\def\BIMOD{\mathop{\text{BIMOD}}}
\def\cohBIMOD{\mathop{\text{cohBIMOD}}}
\def\gr{\mathop{\text{gr}}}
\def\Gr{\mathop{\text{Gr}}}
\def\GR{\mathop{\text{GR}}}
\def\Bigr{\mathop{\text{Bigr}}}
\def\BIGR{\mathop{\text{BIGR}}}

\def\Alg{\mathop{\text{Alg}}}
\def\ALG{\mathop{\text{ALG}}}
\def\Gralg{\mathop{\text{Gralg}}}
\def\QSscr{\mathop{\mathrm{Q}\mathcal{S}}}
\def\Sch{\operatorname{Sch}}
\def\iso{\operatorname{iso}}
\def\tors{\operatorname{tors}}
\def\trans{\operatorname{trans}}
\def\Div{\operatorname{Div}}
\def\Iso{\operatorname{Iso}}
\def\Tors{\operatorname{Tors}}
\def\QSch{\operatorname{Qsch}}
\def\GRALG{\mathop{\text{GRALG}}}

\def\Supp{\mathop{\text{Supp}}}

\def\rk{\operatorname{rk}}
\def\spec{\operatorname {Spec}}
\def\rad{\operatorname {rad}}
\def\gr{\operatorname {gr}}
\def\Spec{\operatorname {Spec}}
\def\diag{\operatorname {diag}}
\def\ox{{o_X}}
\def\Ext{\operatorname {Ext}}
\def\Funct{\operatorname {Funct}}
\def\Inj{\operatorname {Inj}}
\def\Hom{\operatorname {Hom}}
\def\RHom{\operatorname {RHom}}
\def\HHom{\operatorname {\Hscr \mathit{om}}}
\def\HExt{\operatorname {\Escr \mathit{xt}}}
\def\uHHom{\underline{\HHom}}

\def\HTor{\operatorname {\Tscr \mathit{or}}}

\def\im{\operatorname {im}}

\def\coker{\operatorname {coker}}
\def\ker{\operatorname {ker}}

\def\Tor{\operatorname {Tor}}
\def\End{\operatorname {End}}

\def\Pqsch{\operatorname {Pqsch}}
\def\Qsch{\operatorname {Qsch}}

\def\coh{\operatorname {mod}}
\def\gr{\operatorname {gr}}
\def\tors{\operatorname {tors}}

\def\Proj{\operatorname {Proj}}

\def\Qgr{\operatorname {Qgr}}

\def\Tors{\operatorname {Tors}}

\def\Qch{\operatorname {Mod}}

\def\gldim{\operatorname {gl\,dim}}
\def\injdim{\operatorname {inj\,dim}}

\def\r{\rightarrow}
\def\l{\leftarrow}

\let\oldtext\text
\def\text#1{\oldtext{\upshape #1}}
\let\invlim\projlim

\DeclareMathOperator{\PC}{PC}
\DeclareMathOperator{\GrPC}{GrPC}
\DeclareMathOperator{\length}{length}

\DeclareMathOperator{\GrDis}{GrDis}
\DeclareMathOperator{\QGrDis}{QGrDis}
\DeclareMathOperator{\Top}{Top}
\DeclareMathOperator{\GrTop}{GrTop}
\DeclareMathOperator{\pc}{pc}

\DeclareMathOperator{\PCFin}{PCFin}
\DeclareMathOperator{\projdim}{proj\,dim}
\DeclareMathOperator{\GKdim}{GKdim}
\DeclareMathOperator{\QGr}{QGr}
\DeclareMathOperator{\Soc}{Soc}
\newtheorem{lemma}{Lemma}[section]
\newtheorem{proposition}[lemma]{Proposition}
\newtheorem{theorem}[lemma]{Theorem}
\newtheorem{corollary}[lemma]{Corollary}

\newtheorem{lemmas}{Lemma}[subsection]
\newtheorem{propositions}[lemmas]{Proposition}
\newtheorem{theorems}[lemmas]{Theorem}
\newtheorem{corollarys}[lemmas]{Corollary}

\newtheorem{hypothesis}{Hypothesis}

\theoremstyle{definition}

\newtheorem{definition}[lemma]{Definition}

\newtheorem{examples}[lemmas]{Example}
\newtheorem{definitions}[lemmas]{Definition}
\newtheorem{conjectures}[lemmas]{Conjecture}

\newtheorem{step}{Step}
\newtheorem{case}{Case}

\theoremstyle{remark}

\newtheorem{remarks}[lemmas]{Remark}
\newtheorem{notation}{Notation}

\newdimen\uboxsep \uboxsep=1ex
\def\uboxn#1{\vtop to 0pt{\hrule height 0pt depth 0pt\vskip\uboxsep%
\hbox to 0pt{\hss #1\hss}\vss}}

\def\uboxs#1{\vbox to 0pt{\vss\hbox to 0pt{\hss #1\hss}%
\vskip\uboxsep\hrule height 0pt depth 0pt}}

\def\uboxe#1{\hbox to 0pt{\hss\vbox to 0pt{\vss #1\vss}%
\hskip\uboxsep\vrule height 0pt depth 0pt}}

\def\uboxw#1{\hbox to 0pt{\vrule height 0pt depth 0pt\hskip\uboxsep%
\vbox to 0pt{\vss #1\vss}\hss}}

\def\hoek{\hbox{\vtop{\hbox{\vrule\phantom{xx}\vrule}\hrule}\kern -0.4pt}}

\def\dirlim{\mathop{\vtop{\baselineskip -100pt\lineskip -1pt\lineskiplimit 0pt
\setbox0\hbox{\upshape lim}\copy0\hbox to \wd0{\rightarrowfill}}}\limits}
\def\invlim{\mathop{\vtop{\baselineskip -100pt\lineskip -1pt\lineskiplimit 0pt
\setbox0\hbox{\upshape lim}\copy0\hbox to \wd0{\leftarrowfill}}}\limits}

\numberwithin{equation}{section} 
\numberwithin{table}{section}
\numberwithin{figure}{section}

\def\gr{\mathop{\text{gr}}} 
 
\def\Gr{\mathop{\text{Gr}}}
 
\def\Supp{\mathop{\text{Supp}}}

\date{\today} 
\title[Non-commutative blowing up]{Blowing up 
of non-commutative smooth surfaces} 
\author{Michel Van den Bergh}
 \address{Departement WNI\\Limburgs
  Universitair Centrum\\ Universitaire Campus\\ Building D\\ 3590
  Diepenbeek\\ Belgium} 
\email{vdbergh@luc.ac.be,
  http://www.luc.ac.be/Research/Algebra/} 
\thanks{The author is
  a senior researcher at the FWO}
\begin{document}
 \begin{abstract}
   In this paper we will think of certain abelian categories with
   favorable properties as non-commutative surfaces.
   We show that under certain conditions a point on a
   non-commutative surface can be blown up. This yields a new
   non-commutative surface which is in a certain sense birational to
   the original one.  This construction is analogous to blowing up
   a Poisson surface in a point of the zero-divisor of the Poisson
   bracket.  
   
   By blowing up $\le 8$ points in the elliptic quantum plane one
   obtains global non-commutative deformations of Del-Pezzo surfaces.
   For example blowing up six points yields a non-commutative cubic
   surface.  Under a number of extra hypotheses we obtain a formula for
   the number of non-trivial simple objects on such non-commutative
   surfaces.
\end{abstract}
\maketitle
\tableofcontents
\section{Introduction}
\label{ref:1a}
Throughout this paper $k$ will be an algebraically closed field. \subsection{Motivation}
Let $X$
be a smooth connected surface over $k$ and let $q$ be a
Poisson bracket on $X$. Since we are in the dimension two, $q$ corresponds
to a section of $\omega_X^\ast$. 

Let $p\in X$ and let $\alpha:\tilde{X}\r X$  be the blowup of $X$ in $p$.
From the fact that  $\tilde{X}$ and $X$ share the same function field it is easily
seen that  $q$ extends
to $\tilde{X}$ if and only if $q$ vanishes in $p$. Denote the extended Poisson
bracket by $q'$ and let $Y$ resp. $T$ be the zero divisors of $q$ and
$q'$.
One verifies that as divisors~:
$T=\alpha^{-1}(Y)-L$, where $L=\alpha^{-1}(p)$ is the exceptional curve.
In particular  $T$ contains the strict transform $\tilde{Y}$ of $Y$, and if
$p\in Y$ is simple then actually $T=\tilde{Y}$.

Our aim in this paper is to show that there exists a non-commutative
version of this situation. That is we show that it is possible to view
the blowup of a Poisson surface as the quasi-classical analogue of a
blowup of a non-commutative surface.
Our motivation for doing this is to provide a step in the ongoing project
of classifying graded domains of low Gelfand-Kirillov dimension. Since 
the case of  dimension two was  completely solved in
\cite{Staf5} the next interesting case will very likely be
dimension three (leaving aside  rings with fractional
dimension which seem to be quite exotic).
 One may view three dimensional graded rings as non-commutative
 projective  surfaces. Motivated by some heuristic evidence Mike Artin 
 conjectures in \cite{Ar2} 
that, up to birational equivalence, there will
be only a few classes, the largest one consisting of those
algebras that are birational to a quantum $\PP^2$ (see below).

Once a birational classification exists, one might hope that there
would be some version of Zariski's theorem saying that if two
(non-commutative) surfaces are birationally equivalent then they are
related through a sequence of blowing ups and downs.  With the current
level of understanding it seems rather unlikely that Artin's
conjecture or a non-commutative version of Zariski's theorem will be
proved soon, but this paper provides at least one piece of the puzzle.

This being said, it is perhaps the right moment to point out that in this
paper we won't really define 
the notion of a non-commutative surface. Instead we first introduce
non-commutative  schemes (or quasi-schemes, to follow the
terminology of \cite{rosenberg}). These will simply be abelian categories
having sufficiently nice homological properties.  Then we will impose a few
convenient additional hypotheses which would hold for a commutative
smooth surface (see \S\ref{ref:5.1a}).

To fix ideas we will first discuss two particular cases of
quasi-schemes. If $R$ is a ring then $\Spec R$ is the category
of  right $R$-modules (the
``affine case''). 
If $A=A_0\oplus A_1\oplus\cdots$ is a graded ring then $\Proj A$ is (roughly) 
the category of graded right $A$-modules, modulo
the modules with right bounded grading (the ``projective case'').

Let us first consider the affine case. Assume that $R$ is a finitely
generated  $k$-algebra and let $C$ the commutator ideal. $C$ is the
natural analog of the zero divisor of a Poisson bracket.
Now $\Spec R/C$ is a
commutative affine scheme and a $k$-point in $\Spec R/C$ corresponds to
a maximal ideal  $m$ in $R$ with $R/m=k$. Hence a natural idea is to
define the blowup of $\Spec R$ in $p$ as  $\Proj D$  where $D$ is
the Rees  algebra  associated to $m$.
\[
D=R\oplus m\oplus m^2\oplus m^3\oplus\cdots
\]
It is easily seen however that this definition is faulty. Consider the
following example   \cite{Ar2} $R=k\langle x,y\rangle/(yx-xy-y)$,
$m=(x,y)$. Then $m^n=(x^n,y)$. Hence the analog of the
exceptional curve
\[
D/mD=R/m\oplus m/m^2\oplus m^2/m^3\oplus
\]
is isomorphic to $k[x]$. Thus  $\Proj D/m$ is a point, whereas intuitively we
would expect it to be one-dimensional in some sense.

It turns out however that in this example one can use a certain
twisting of the
Rees algebra which yields a reasonably behaved blowing
 up. This is based on the observation that the commutator ideal $C=(y)$ is
an invertible $R$-bimodule. Let $J$ be its inverse. Then we define
$I^n$ as the  image $I^{\otimes n}$ in $J^{\otimes n}$ and we define the
modified Rees algebra $D$ as 
\[
D=R\oplus I\oplus I^2\oplus I^3\oplus\cdots
\]
The blowup of $\Spec R$ in $p$ is now defined as $\Proj D$ for this
new $D$. We refer the reader to \cite{Ar2} for a detailed workout of
this example.  However we will indicate how one finds the analog of
the exceptional curve. Let $\tau$ be the automorphism of $R$ given by
$a\mapsto y^{-1}a y$.  Then $J=R_\tau$ and hence $I^n=m\tau(m)\cdots
\tau^{n-1}(m)_{\tau^n}$. Put $L=D/\tau^{-1}(m) D$. One now verifies
that $\dim L_u=u+1$. So $L$ plays the role of the exceptional curve.
Note however that $L$ is a right $D$-module but \emph{not} a left
module. In retrospect this was to be expected since, as we have said in
the first paragraph, if we blow up a Poisson surface, then the extended
Poisson bracket will in general not vanish on the exceptional curve.

This example indicates the way to go for rings whose commutator ideal is
 invertible.
The latter hypotheses is not unreasonable since if we look at the case of
a Poisson surface then we see that we expect a non-commutative smooth
surface to contain a commutative curve.
Additional motivation comes from considering the  local
rings 
$k\langle\langle x,y\rangle\rangle /(\phi)$
with $\phi$ a (non-commutative) formal power series whose lowest degree
term is a non-degenerate  quadratic tensor in $x,y$. These local rings 
are the non-commutative analogs of complete two dimensional
 regular local rings and one verifies that their commutator ideal
is indeed invertible (see e.g.\ \cite{VdBVG}).

There is one important hitch however. The commutator ideal is not
invariant under Morita equivalence! This indicates that it is important
to develop the theory in a more category-theoretical frame work. This
 will make it possible to talk about non-commutative schemes
containing a commutative curve, without refering to rings or
ideals at all.

To stress this point even more let us consider the case of graded
rings. In \cite{AS} Artin and Schelter introduced so-called regular
rings. These are basically graded rings which have the Hilbert series
of three dimensional polynomial rings, together with a few other
reasonable properties. They were classified in \cite{AS,ATV1,ATV2} and
also, with different methods, in \cite{Bondal}. Let $A$ be such a
regular ring. We view $X= \Proj A$ as a quantum $\PP^2$. Since
on $\PP^2$ the anti-canonical sheaf has degree three, the zero divisor
of a Poisson bracket will be a (possibly singular and non-reduced)
elliptic curve. Therefore we would also expect $X=\Proj A$ to contain
an elliptic curve in some reasonable sense. It turns out that this is
indeed true!  It was shown in \cite{AS,AVdB,ATV1} that $A$
contains a normal element $g$ in degree three such that $\Proj A/gA$
is equivalent with the category of quasi-coherent sheaves over an
elliptic curve $Y$. Thus if we actually identify $Y$ with $\Proj A/gA$
then $Y\hookrightarrow X$.

Now let $p\in Y$. The previous discussion suggests that it should be
possible to blow up $p$. However it is not clear how to proceed. Under
the inclusion $Y\hookrightarrow X$, $p$ corresponds to a so-called
point module \cite{ATV2} over $A$. This is by definition a graded
right module, which is generated in degree zero and which is
one-dimensional in every degree. However such a point module is only a
right module and hence it cannot be used to construct a Rees algebra.

Our solution is to construct the Rees algebra directly over $\Proj A$.
To do this we have to invoke the theory of monads \cite{ML}. Since we
only consider monads satisfying a lot of additional hypotheses we
prefer to call our monads algebras. This is at variance with
the use of ``algebra'' in the theory of categories \cite{ML} but in
our context it seems reasonable. For us an algebra over an abelian
category $\Cscr$ is in principle an algebra object in the monoidal
category of right exact functors from $\Cscr$ to itself. There are
however some technical problems with this so  we end up using a
less intuitive definition (see below).

The importance of monads in non-commutative algebraic geometry was
noticed by various people, in particular by Rosenberg. See for example
\cite{rosenberg,RL}. In the last chapter of his book Rosenberg
actually defines a blow up of an arbitrary ``closed'' subcategory of
an abelian category. While this definition is also in terms of monads,
it is as far as I can see, somewhat different from ours. To see this
let us again consider the affine case. Then Rosenberg's construction
is in terms of the functor $M\mapsto Mm$, which is not right exact. If
we replace this functor by $M\mapsto M\otimes_R m$ then one would get
the $\Proj$ of the ordinary Rees algebra of $R$, which (depending on
what one wants to achieve) might not be the right answer (as we have
shown above).
\subsection{Construction}
Following
\cite{rosenberg} we introduce the notion of a quasi-scheme. For us
this will be a Grothendieck category (that is~: an abelian category
with a generator and exact direct limits). However we tend to think of
quasi-schemes as geometric objects, so we denote them by roman
capitals $X$, $Y$, \dots.  If we really refer to the category
represented by a quasi-scheme $X$ then we write $\Qch(X)$. Note that in
fact $X=\Qch(X)$, but it is very useful to nevertheless make this
notational distinction since it allows us to introduce other notations
in a consistent way. For example we will  denote the noetherian
objects in $\Mod(X)$ by $\mod(X)$. Furthermore we can absorb 
additional structure into the symbol $X$ (such as a map to a base
quasi-scheme) which is not related to $\Mod(X)$. This would be awkward
without the two different notations $X$ and $\Mod(X)$.

A map $\alpha:X\r Y$ between quasi-schemes will be a right exact
functor $\alpha^\ast:\Mod(Y)\r \Mod(X)$ possessing a right adjoint
(denoted by $\alpha_\ast$). In this way the quasi-schemes form a
category (more precisely a two-category, see Appendix \ref{ref:Aa}).

If $X$ is a quasi-compact quasi-separated commutative scheme then
$\Qch(X)$ will be the category of
quasi-coherent sheaves on $X$. It is proved in \cite{thomasson} that
this is a Grothendieck category.  Rosenberg in \cite{rosenberg1} has
proved a reconstruction theorem which allows one to recover $X$ from
$\Qch(X)$ (generalizing work of Gabriel in the noetherian
case). He has also announced that the  functor which assigns to a
commutative scheme its associated quasi-scheme is fully
faithful if we work over $\Spec\ZZ$.

Let $X$ be a quasi-scheme. We think of objects in $\Qch(X)$ as sheaves of right
modules on $X$. However to define algebras on $X$, it is clear that we 
need
bimodules on $X$ (see \cite{VdB11} for the case where $X$ is commutative).
Let us for the moment define a bimodule on $X$ as a right exact functor
from $\Qch(X)$ to itself commuting with direct limits. Then the
corresponding category is monoidal (the tensor product
 being given by composition) and hence we can define
algebra objects. Let $\Ascr$ be such an algebra object. It is routine to
define an abelian category $\Mod(\Ascr)$ of $\Ascr$-modules. So this
seems like a reasonable starting point for the theory.

However a difficulty emerges if one wants to define Rees algebras.  As
we have seen, the main point is to take the sum of the $I^n$ for some
subbimodule $I$ of an invertible bimodule $\Lscr$. $I^n$ was defined
as the image of $I^{\otimes n}\r \Lscr^{\otimes n}$. Unfortunately to
take an image one needs an abelian category, and I don't see how to
prove that the above definition of a bimodule yields an abelian
category, even if we drop the requirement that bimodules should
commute with direct limits.  In this paper we sidestep this difficulty
by defining the category of bimodules on $X$ as the opposite category
of the category of left exact functors from $\Qch(X)$ to itself. Since
left exact functors are determined by their values on injectives, they
trivially form an abelian category. In this way one can define
Rees algebras in reasonable generality (see Definition
\ref{ref:3.5.13a}).

We will say that a quasi-scheme $X$ is noetherian if $\Mod(X)$ is
locally noetherian. That is, if $\Mod(X)$ is generated by $\mod(X)$.
As already has been indicated above, in this paper we will study a
noetherian quasi-scheme $X$ which contains a commutative curve $Y$
as a divisor.  To make this more precise we denote the identity functor
on $\Qch(X)$ by $o_X$. This is an algebra on $X$ such that
$\Mod(o_X)=\Qch(X)$.  We will assume that $o_X$ contains an invertible
subbimodule $o_X(-Y)$ such that $\Mod(o_X/o_X(-Y))$ is equivalent with
$\Qch(Y)$.

We also need some sort of smoothness condition on $X$. Since it is
obviously sufficient to impose this in a neighborhood of $Y$, we assume
that every object in $\Qch(Y)$ has finite injective dimension in
$\Qch(X)$.
This is the same setting as in \cite{VdBVG}, albeit cast in a
somewhat different language.

Now $p\in Y$ defines a subbimodule $m_p$ of $o_X$ which is the analog
of the maximal ideal corresponding to $p$. We put
$I=m_po_X(Y)$. Define
\[
 \Dscr=o_X\oplus I\oplus I^2\oplus \cdots
 \]
 This is the Rees algebra associated to $I$. We  define the
 blowup $\tilde{X}$ of $X$ in $p$ as $\Proj \Dscr$.
\subsection{General properties}
A large part of this paper is devoted to proving that $\tilde{X}$
satisfies similar properties as $X$ and furthermore that we have
obtained an analogue of the blowup of a commutative surface. However
before we start we need to have a better understanding of the formal
neighborhood of a point $p\in Y$.  This was in fact already done in
\cite{VdBVG}. The answer is given in terms of certain topological
rings (see Theorem \ref{ref:5.1.4a} for precise results).  It turns
out that the formal local structure of $p$ depends heavily on a
certain automorphism $\tau$ of $Y$ which was also introduced in
\cite{VdBVG}. To be more precise we define the normal bundle of $Y$ in
$X$ as $o_X(Y)/o_X$. This bimodule defines an autoequivalence of
$\Mod(Y)$ and by a result in \cite{AZ} such an autoequivalence must
necessarily be of the form $\tau_\ast(-\otimes_{\Oscr_Y}\Nscr)$ where
$\Nscr$ is a line bundle on $Y$ and $\tau$ is an automorphism of $Y$.
In particular if $\tau$ has infinite order (and hence $p\in Y$ is
smooth) then the formal local
structure of $p$ is given by the ring of doubly infinite lower
triangular matrices over the ring $\hat{\Oscr}_{Y,p}$. In particular
this is  independent of $p$ and $Y$ (as is the case for the
completion at a smooth point on a commutative scheme).

 In this
paper we complete the results in \cite{VdBVG} by showing that various
completion functors are exact. For this we refer the reader to
\S\ref{ref:5a}. 
 An interesting application of completion is given in section
\S\ref{ref:5.7b}. In this section we define (roughly) the
multiplicity in $p$ and the points infinitely near to $p$ of an object
in $\mod(X)$ in the case that the $\tau$-orbit of $p$ is infinite.

\smallskip

As a starting point for the study of $\tilde{X}$ we construct a
commutative diagram of quasi-schemes.
\[
\begin{CD}
\tilde{Y} @>\beta>> Y\\
@ViVV @VVjV\\
\tilde{X} @>>\alpha> X
\end{CD}
\]
where the vertical arrows are inclusions. $\tilde{X}$ is again a
noetherian quasi-scheme (Theorem \ref{ref:6.3.1a}).  $\tilde{Y}$
is a commutative curve which plays the role of the strict transform of
$Y$. $\tilde{Y}$ is again a divisor in $X$ and every object in
$\Qch(\tilde{Y})$ has finite injective dimension in $\Qch(\tilde{X})$
(see Theorem \ref{ref:6.6.3a}).

Associated to the point $p\in Y$ there is a simple object $\Oscr_p$.
We define $\Oscr_L$ as $\alpha^\ast(\Oscr_p)$ and we consider
$\Oscr_L$ as the structure sheaf on the exceptional curve in $\tilde{X}$. In fact
following a recipe given in \cite{SmithZhang} we can define a
category $\Mod(L)$. Roughly $\Mod(L)$ is generated by subquotients of
direct sums of twists of $\Oscr_L$. In this way we can speak of the
quasi-scheme $L$. It follows from Proposition \ref{ref:6.5.2a}
(together with Corollary \ref{ref:6.7.4a}) that if we view $\Mod(X)$
modulo the objects supported on the $\tau$-orbit of $p$, and
$\Mod(\tilde{X})$ modulo the objects supported on $L$ then we obtain
equivalent categories. This is the obvious analogue of the situation in
the commutive case where $X-p$ and $\tilde{X}-L$ are isomorphic (and
hence in particular $X$ and $\tilde{X}$ are birational).

In section \S\ref{ref:6.7b} we give
a precise description of $\Mod(L)$ (using results of
\cite{SmithZhang} in the case that $\tau$ has infinite order). It
will follow that $\Mod(L)$ is very closely related to the category of
quasi-coherent sheaves on $\PP^1$, illustrating again the analogy with
the commutative case.

In section \S\ref{ref:7a} we compute the derived category of
$\tilde{X}$. Our main result is that this derived category has a
semi-orthogonal decomposition given by the derived category of $X$
and the derived category of $k$. This generalizes a result by Orlov
\cite{Orlov}.
\subsection{Non-commutative  Del-Pezzo surfaces}
The results starting from Section \ref{ref:9a} are inspired by the
construction in the commutative case of (most) Del-Pezzo surfaces by
blowing up a collection of points in $\PP^2$. Let $(Y,\sigma,\Lscr)$
be an elliptic triple $Y\subset \PP^2$ as in \cite{ATV1}. We assume
that $Y$ is smooth and that $\sigma$ is a translation of infinite
order. Let $A$ be the regular algebra associated to this triple
\cite{ATV1} and let $X_0=\Proj A$. As above we consider $X_0$ as a
quantum version of $\PP^2$. The curve $Y$ is contained as a
divisor in $X_0$ and $\tau$ is equal in this case to $\sigma^3$. We
choose points $p_1,\ldots, p_n\in Y$ ($n\le 8$) and we construct a
commutative diagram of quasi-schemes
$$
\xymatrix{
&\tilde{X}_1\ar[dl]_{\alpha_1}\ar[dr]^{\delta_0}
&&\tilde{X}_2\ar[dl]_{\alpha_2}\ar[dr]^{\delta_1}
&&&&\tilde{X}_n\ar[dl]_{\alpha_n}\ar[dr]^{\delta_n}
&\\
X_1
&&X_2
&&X_3
&\cdots&X_n
&&X_{n+1}
\\
&&&&Y\ar[ullll]\ar[ull]\ar[u]\ar[urr]\ar[urrrr]
&&&&
}
$$
Here the map $\alpha_i$ is the blowup of $X_i$ in $p_i$. Morally
$X_{i+1}$ is constructed from $\tilde{X}_i$ by putting $X_{i+1}=\Proj
\left(\bigoplus_n H^0(X_i,o_{X_i}(nY))\right)$ (actually for
simplicity $X_{i+1}$ is constructed using a slightly different method
(see \S\ref{ref:9.2b}), which is easily seen to be equivalent). The
point is that in the commutative case $\delta_i$ would be an
isomorphism if the points $p_1,\ldots,p_n$ are in general position
(this follows from the Nakai criterion for ampleness, see \cite{H}). In the
non-commutative $\delta_i$ will not in general be an isomorphism. However
we show in Theorem \ref{ref:11.1.3a} that $\delta_i$ yields a derived equivalence, under suitable
general position hypotheses.
\subsection{Exceptional simple objects}
One of the aims of these notes is to classify (or rather count) the simple objects in
$\Mod(X_i)$ which are not of the form $\Oscr_q$ for some $q\in Y$. We
call such simple objects ``exceptional'' because, firstly, they do not
always exist, and secondly when they exist they are not easy to
construct or to count.

Using some results on geometry in the projective quantum plane
(\S\ref{ref:10a}) together with the above results on
derived categories we obtain in Theorem \ref{ref:11.2.1a} a formula for
the number of exceptional simple objects in $\Mod(X_{i+1})$ (if $n\le
6$, $\tau$ has infinite order and  suitable general position
hypotheses hold). It turns out that the number of exceptional simple
objects depends in a \emph{very} sensitive way upon the position of
the points $p_1,\ldots,p_n$. For example if $n=6$ (and the other
hypotheses are satisfied) then our formula yields that there may be 
between $0$ and $6$ exceptional simple objects in $\Mod(X_7)$.
\subsection{Non-commutative cubic surfaces}
Being near the end of this introduction we now indicate our original
motivation for starting this project. It concerns a problem
which is not quite completely solved but which at least has become
more tractable.

As above let $(Y,\sigma,\Lscr)$ be an elliptic triple with $Y\in
\PP^2$ and let $A$ be the associated three dimensional Artin-Schelter
regular algebra. It is easy to show that the representation theory of
$A$ is fairly trivial. At a certain point Lieven Le Bruyn (see for
example \cite{L})
suggested
that one could obtain more interesting representation theories by
considering filtered rings $C$ such that $\gr C=A$. This was motivated
by the example of enveloping algebras of Lie algebras which are
related in a similar way to polynomial algebras.

Instead of studying the filtered rings $C$, it is easier to study their
Rees rings $D=\oplus_n C_n$. These are characterized as the
$\NN$-graded rings containing a regular central element $t$ in degree
one such that $D/tD=A$. Such graded rings were studied in \cite{LSV}
but not much progress was made towards their representation theory. 

As usual a first step in the study of the representation theory of $D$
is the construction of a ``Casimir'' element. Indeed $A$ contains a
canonical central element $g$ in degree $3$ and by a computer
computation one can show that $g$ lifts to a central element $G$ in
$D$. Then instead of studying $D$ on may study the central quotients
$D_\mu=D/(G+\mu t^3)$, $\mu \in k$. The $D_\mu$ may be considered as
quantum versions of cubic surfaces in $\PP^3$.

A well known result in commutative algebraic geometry states that a
cubic surface in $\PP^3$ is obtained by blowing up six points in $\PP^2$, so one
may ask whether this is also true in the noncommutative case. I have
not been able to show this in general but a converse result is
obtained in Section \ref{ref:12a}. We show
that if we blow up six points in the elliptic quantum plane then the
resulting quasi-scheme is  contained as a cubic divisor in a
quantum $\PP^3$.

Very recently Mike Artin has explained to me that  one can probably
obtain the complete analogue of the commutative result by deformation
theory. Indeed $Z=\Proj D_\mu$ can be obtained as a deformation of a
cubic surface $Z_0$ in $\PP^3$. In  $Z_0$ we can choose six
mutually skew exceptional curves \cite{H}
 and since the structure sheaves  of these exceptional curves
have no higher $\Ext$'s they deform well. So we should find on $Z$ six
corresponding exceptional curves. Then we can blow these down, for
example using the procedure exhibited in \cite{VdB24}.
To carry out this program there are quite a few technical details that
remain to be filled in. There are some recent notes by Mike
Artin on specializing birational equivalences in the non-commutative
case. 
\subsection{Acknowledgement}
The author wishes to thank Mike Artin and Paul Smith for stimulating
discussions about the material in earlier versions of this manuscript
and about non-commutative geometry in general. I also wish to thank Paul
Smith for showing me the interesting preprint \cite{SmithZhang} which
contains related results.
\section{Preliminaries on category theory}
Our main references for categories are \cite{Gabriel,Groth1,ML,stenstrom}.
If nothing is specified then  categories will have small homsets.
This will not always be the case for some secondary categories such as
 categories of bimodules (see below). However limits and colimits are
always taken over small sets. In particular complete and cocomplete
will refer to the existence of small limits and colimits. 
When we speak of a direct limit, we mean a  colimit over a \emph{directed}
set. A similar convention applies to inverse limits.

We will often implicitly use the fact that in an abelian category 
\begin{equation}
\label{ref:2.1a}
\begin{split}
\oplus_{i\in I} Y_i&=\dirlim_{\substack{j\in J\\ J\subset I\text{
      finite}}} Y_j\\
\prod_{i\in I} Y_i&=\invlim_{\substack{j\in J\\ J\subset I\text{
      finite}}} Y_j
\end{split}
\end{equation}
Thus an additive functor commuting with direct limits
 commutes with coproducts. A similar statement applies to products.
 
We  use the following specialized version of the standard
adjoint functor theorems \cite{ML}.
\begin{theorem} 
\label{ref:2.1b}
Let $\Cscr$, $\Dscr$ be abelian categories and let $F:\Dscr\r\Cscr$,
$G:\Cscr\r\Dscr$ be respectively a right and a left exact functor.
Then
\begin{enumerate}
\item If $\Cscr$ is complete and has a
  cogenerator then $G$ has a left adjoint if and only if $G$ commutes
  with products.
\item (Dual version)  If  $\Dscr$ is cocomplete  and
has a generator then $F$ has a right adjoint if and only if $F$ commutes
with coproducts.
\end{enumerate}
\end{theorem}
Most of the abelian categories we use will be \emph{Grothendieck
  categories}. These are abelian categories which have a generator and
exact direct limits. We use the following results which are well-known
\cite{Groth1,stenstrom}.
\begin{theorem} Let $\Cscr$ be a Grothendieck category. Then
\begin{enumerate}
\item $\Cscr$ has products (not necessarily exact).
\item $\Cscr$ has enough injectives.
\item The product of the injective hulls of the quotients of a fixed
  generator is an injective cogenerator.
\end{enumerate}
\end{theorem}
A Grothendieck category which is generated by noetherian objects is
called \emph{locally noetherian}. Such categories were studied by
Gabriel in \cite{Gabriel}. One important property they have is the
following.
\begin{theorem}
  In a locally noetherian category the direct limit of injective
  objects is injective.
\end{theorem}
A \emph{noetherian} category is an abelian category in
which every object is noetherian. One has \cite[III, Th 1]{Gabriel}
\begin{theorem} Every noetherian category $\Cscr$ is equivalent
  with the category of noetherian objects in a locally noetherian
  category $\tilde{\Cscr}$.  $\tilde{\Cscr}$ is
   characterized up to equivalence by this property.
\end{theorem}
A \emph{Serre subcategory} $\Sscr$ of an abelian category $\Cscr$ is
by definition a full subcategory which is closed under subquotients
and extensions. In that case there exists a quotient category
$\Cscr/\Sscr$ which is characterized by an appropriate universal
property. If $\Cscr$ is a Grothendieck category  then we say that
$\Sscr$ is localizing if $\Sscr$ is a Serre subcategory which is
closed under direct limits. In that case $\Cscr/\Sscr$ is also a
Grothendieck category. Furthermore the quotient functor $\pi:\Cscr\r
\Cscr/\Sscr$ has a right adjoint which we will usually denote by
$\omega$. The composition $\omega\pi$ is called the \emph{localization
  functor} and will be denoted by $\tilde{(-)}$.

In this paper we basically work with the category of Grothendieck
categories, the morphisms being certain functors. 
This is an example of a
two-category. Although we usually only implicitly use this concept,
we explain some basic notions in Appendix \ref{ref:Aa}.

\section{Non-commutative geometry}
\label{ref:3a}

In this paper we will sometimes work with categories which do not come
from categories of modules over (graded) rings. Therefore in this
section we introduce some rudiments of a formalism which may be used
to imitate some of the more elementary features of commutative
algebraic geometry. This section is closely related to
\cite{rosenberg}\cite{VdB11}.

\subsection{Bimodules}
\label{ref:3.1a}
Below $\Cscr$, $\Dscr$, $\Escr$, \dots will be abelian categories. We
define the following categories.
\begin{align*}
  \Lscr(\Dscr,\Cscr)&=\{\text{left exact functors $\Dscr\r \Cscr$}\}
  \\ \BIMOD(\Cscr-\Dscr)&=\Lscr(\Dscr,\Cscr)^{\text{opp}}\\ 
  \Bimod(\Cscr-\Dscr)&=\{\Mscr\in\BIMOD(\Cscr-\Dscr)\mid \Mscr\text{
    has a left adjoint }\}
\end{align*}
Objects in $\BIMOD(\Cscr-\Dscr)$ will be called \emph{weak
  $\Cscr$-$\Dscr$ bimodules}, whereas objects in $\Bimod(\Cscr-\Dscr)$
will simply be called \emph{$\Cscr$-$\Dscr$ bimodules}.
\begin{propositions}
\label{ref:3.1.1a}
\begin{enumerate}
\item $\BIMOD(\Cscr-\Dscr)$, $\Bimod(\Cscr-\Dscr)$ have finite
  colimits, and the inclusion $\Bimod(\Cscr-\Dscr)\r
  \BIMOD(\Cscr-\Dscr)$ preserves those colimits.
\item If $\Cscr$ is complete then $\BIMOD(\Cscr-\Dscr)$ is cocomplete.
  If $\Dscr$ is in addition cocomplete then so is  $\Bimod(\Cscr-\Dscr)$
and furthermore  the inclusion $\Bimod(\Cscr-\Dscr)\r
  \BIMOD(\Cscr-\Dscr)$ preserves colimits.
\item If $\Dscr$ has enough injectives then $\BIMOD(\Cscr-\Dscr)$ is
  an abelian category.
\item If $\Dscr$ is complete and has a cogenerator then an object in
  $\BIMOD(\Cscr-\Dscr)$ is in $\Bimod(\Cscr-\Dscr)$ iff it commutes
  with products.
\item If $\Dscr$ is complete and has an injective cogenerator then an
  object in $\BIMOD(\Cscr-\Dscr)$ is in $\Bimod(\Cscr-\Dscr)$ iff it
  commutes with products when evaluated on injectives.
\item If $\Dscr$ is complete and has an injective cogenerator and
  products are exact in $\Cscr$ then $\Bimod(\Cscr-\Dscr)$ is an
  abelian subcategory of $\BIMOD(\Cscr-\Dscr)$.
\end{enumerate}
\end{propositions}
\begin{proof}
\begin{enumerate}
\item To show that $\BIMOD(\Cscr-\Dscr)$ has finite colimits we have
  to show that the opposite category $\Lscr(\Dscr,\Cscr)$ has finite
  limits. Now one easily verifies that for a functor
  $G:I\r \Lscr(\Dscr,\Cscr)$ with $I$ a finite category
\[
G'(M)=\invlim_i G(i)(M)
\]
defines the inverse limit of $G$ in $\Lscr(\Dscr,\Cscr)$.

Assume now that the $G(i)$ are in $\Bimod(\Cscr-\Dscr)$ so
they have left adjoints $F(i)$.
Define $F'$ by 
\[
F'(N)=\dirlim_i F(i)(N)
\]
One easily verifies that $F'$ is a left adjoint to $G'$ and hence $G'\in
\Bimod(\Cscr-\Dscr)$.
\item This is similar to (1).
\item Let $\Inj(\Dscr)$ denote the additive category of injectives in
  $\Dscr$. Then
$
\Lscr(\Dscr,\Cscr)
$ is equivalent to $\Funct(\Inj(\Dscr),\Cscr)$,
where ``$\Funct$'' denotes the category of additive functors. It is
now clear that $\Lscr(\Dscr,\Cscr)$ inherits the property of being an
abelian category from $\Cscr$.
\item This follows from Theorem \ref{ref:2.1b}.
\item This follows from (4),  using the fact that products are left
  exact.
\item We have to show that $\Bimod(\Cscr-\Dscr)$ is closed under
  kernels. Equivalently the subcategory of
  $\Funct(\Inj(\Dscr),\Cscr)$ of functors commuting with products has to
  be closed under cokernels.

  Let $G_1\r G_2$ be a map of functors in $\Funct(\Dscr,\Cscr)$
  commuting with products and let $G_3=\coker(G_1\r G_2)$. Then for
  $E\in \Inj(\Dscr)$ we have $G_3(E)=\coker(G_1(E)\r G_2(E))$ and
  using the fact that products are exact one easily obtains
  that $G_3$ commutes with products.\qed
\end{enumerate}
\def\qed{}
\end{proof}
We will now introduce some more suggestive notations for dealing with
$\BIMOD(\Cscr-\Dscr)$.

If $\Mscr\in \BIMOD(\Cscr-\Dscr)$ then we denote the corresponding
functor in $\Lscr(\Dscr,\Cscr)$ by
$
\HHom_\Dscr(\Mscr,-)
$.
More or less by definition we have the following facts.
\begin{propositions}
\begin{enumerate} 
\item $\HHom_\Dscr(-,-)$ is left exact in its two arguments.
\item If $\Cscr$ is complete then $\HHom_\Dscr(-,-)$ transforms
  colimits in its first argument into limits.
\item If $\Mscr\in \Bimod(\Cscr-\Dscr)$ then $\HHom_\Dscr(\Mscr,-)$
  transforms limits in its second argument into limits.
\item If $\Dscr$ has enough injectives and $E\in \Inj(\Dscr)$ then
  $\HHom_\Dscr(-,E)$ is exact.
\end{enumerate}
\end{propositions}
\begin{proof}
\begin{enumerate}
\item That $\HHom_\Dscr(-,-)$ is left exact in its second argument is by
definition. That it is left exact in its first argument follows from the
explicit construction of colimits (and hence of cokernels) in the proof of
Prop.\ \ref{ref:3.1.1a}(1).
\item
This follows as in the proof of Prop.\ \ref{ref:3.1.1a}(2).
\item If $\Mscr\in\Bimod(\Cscr-\Dscr)$ then $\HHom_\Dscr(\Mscr,-)$ has by
definition a left adjoint. Hence it commutes with limits.
\item This follows from the explicit structure of an abelian category on
$\BIMOD(\Cscr-\Dscr)$ given by the proof Prop.\ \ref{ref:3.1.1a}(3).
\qed\end{enumerate}
\def\qed{}
\end{proof}
 We write composition of functors
\[
\BIMOD(\Cscr-\Dscr)\times \BIMOD(\Dscr-\Escr)\r \BIMOD(\Cscr-\Escr)
\]
as $-\otimes_\Dscr-$. In this way we obtain for $\Mscr\in
\BIMOD(\Cscr-\Dscr)$, $\Nscr\in \BIMOD(\Dscr-\Escr)$ the satisfying
formula
\[
\HHom_\Escr(\Mscr\otimes_\Dscr\Nscr,-)=\HHom_\Dscr(\Mscr,\HHom_\Escr(\Nscr,-))
\]
Again more or less by definition we obtain the following properties.
\begin{propositions}
\begin{enumerate}
\item $-\otimes_\Dscr-$ is right exact in its two arguments.
\item If $\Cscr$ is complete then $-\otimes_\Dscr-$ preserves colimits
  in its first argument.
\item If $\Dscr$ is complete and if $\Mscr\in \Bimod(\Cscr-\Dscr)$
  then $\Mscr\otimes_\Dscr-$ preserves colimits in its second
  argument.
\item $-\otimes_\Dscr-$ sends $\Bimod(\Cscr-\Dscr)\times
  \Bimod(\Dscr-\Escr)$ to $\Bimod(\Cscr-\Escr)$.
\end{enumerate}
\end{propositions}
Now define
\[
\MOD(\Cscr)=\BIMOD(\Ab-\Cscr)
\]
The functor $M\mapsto \Hom_\Cscr(M,-)$ defines a full faithful
embedding of $\Cscr$ in $\MOD(\Cscr)$. Throughout we will identify
$\Cscr$ with its essential image under this embedding.

We now obtain the following alternative ``characterization'' of
$\Bimod(\Cscr-\Dscr)$ inside $\BIMOD(\Cscr-\Dscr)$.
\begin{propositions}
\label{ref:3.1.4a}
Let $\Mscr\in\BIMOD(\Cscr-\Dscr)$. Then $\Mscr$ is in
$\Bimod(\Cscr-\Dscr)$ if and only if the functor
\[
-\otimes_\Cscr\Mscr:\MOD(\Cscr)\r \MOD(\Dscr)
\]
sends $\Cscr$ to $\Dscr$.
\end{propositions}
A few concepts from the theory of bimodules over rings can be
generalized to our setting. We denote the derived functors of
$\HHom_\Dscr(-,-)$ in the second argument by $\HExt^i_\Dscr(-,-)$ (if
they exist).  Assume that $\Escr$ has enough injectives. We say that
$\Nscr\in \BIMOD(\Dscr-\Escr)$ is \emph{(left) flat} if
$\HHom_\Escr(\Nscr,-)$ preserves injectives. More generally for
$\Mscr\in \BIMOD(\Cscr-\Dscr)$ and $\Nscr\in\BIMOD(\Dscr-\Escr)$ we
define $\HTor^\Dscr_i(\Mscr, \Nscr)$ as the object in
$\BIMOD(\Cscr-\Dscr)$ satisfying
\begin{equation}
\label{ref:3.1b}
\HHom_\Escr(\HTor^\Dscr_i(\Mscr,\Nscr),E)\overset{}{=}
\HExt^i_\Dscr(\Mscr,\HHom_\Escr(\Nscr,E))
\end{equation}
for every injective $E$ of $\Escr$. A similar definition holds for
$\Mscr\in\Dscr$. 
\begin{propositions} 
\label{ref:3.1.5a}
Assume that $\Escr$ has enough injectives. 
\begin{enumerate}
\item $\HTor^\Dscr_i(-,-)$ is a $\delta$-functor in its two arguments.
\item $\Nscr\in \BIMOD(\Dscr-\Escr)$ is flat if and only if
  $\HTor^\Dscr_1(\Mscr,\Nscr)=0$ for all $\Mscr\in\Dscr$.
\item
If $\Nscr\in \BIMOD(\Dscr-\Escr)$ is flat then
$\HTor^\Dscr_i(\Mscr,\Nscr)=0$ for all $\Mscr\in \BIMOD(\Cscr-\Dscr)$.
\end{enumerate}
\end{propositions}
\begin{proof}
All statements follow directly from the definitions.
\end{proof}
Assume that $\Cscr$, $\Dscr$ have colimits. Then we will say that
$\Mscr\in\BIMOD(\Cscr-\Dscr)$ is coherent if $\HHom_\Dscr(\Mscr,-)$ commutes
with direct limits.
\begin{propositions}
\label{ref:3.1.6a}
Assume that $\Mscr\in\Bimod(\Cscr-\Dscr)$ and consider the following
statements.
\begin{enumerate}
\item $\Mscr$ is coherent.
\item $-\otimes_{\Cscr}\Mscr$ preserves finitely presented objects in
$\Cscr$.
\end{enumerate}
Then (1) implies (2). The converse holds if $\Cscr$ is generated by
finitely presented objects.
\end{propositions}
\begin{proof}\strut
\begin{itemize}
\item[(1)$\Rightarrow$(2)]
Assume that $T$ is a finitely presented object in $\Cscr$. I.e.
$\Hom_\Cscr(T,-)$ commutes with  direct limits. We have to show that
$\Hom_\Cscr(T\otimes_\Cscr \Mscr,-)$ commutes with  direct
limits. This follows 
from the fact that 
\[
\Hom_\Cscr(T\otimes_\Cscr\Mscr,-)=\Hom_\Cscr(T,\HHom_\Dscr(\Mscr,-))
\]
\item[(2)$\Rightarrow$(1)] We have to construct a natural isomorphism
between
$\Hom_\Cscr(U,\dirlim_i \Hom_\Dscr(\Mscr,\Nscr_i))$ and
$\Hom_{\Cscr}(U,\Hom_\Dscr(\Mscr,\dirlim \Nscr_i))$ for an arbitrary
inverse system $(\Nscr_i)_i$ in $\Dscr$ and $U\in\Cscr$. Since
$\Cscr$ is generated by finitely presented objects, it suffices to do
this for $U$ finitely presented. But this is clear by adjointness.
\qed\end{itemize}
\def\qed{}\end{proof}
\begin{propositions} 
\label{ref:3.1.7a}
Assume $\Dscr$ is a Grothendieck category. Suppose
  furthermore that $\Cscr$ has exact direct limits.
 Let $\Mscr\in \BIMOD(\Cscr-\Dscr)$. Then the following are
  equivalent
\begin{enumerate} \item $\Mscr$ is coherent.
\item $\HHom_\Dscr(\Mscr,-)$  commutes with direct limits of injectives.
\end{enumerate}
\end{propositions}
\begin{proof} 
We only have to prove 2.$\Rightarrow$1. In a Grothendieck category
embeddings into injectives can be constructed functorially
\cite{Gr}. Hence if $(T_i)_{i\in I}$ is an inverse system in
$\Dscr$ 
then there is a copresentation
\[
0\r  (T_i)_{i\in I}\r  (E_i)_{i\in I}\r  (F_i)_{i\in I}
\]
with $(E_i)_i$, $(F_i)_i$ injective. The left exactness of
$\HHom_\Dscr(\Mscr,-)$ together with the fact that direct limits are
exact in $\Cscr$ and $\Dscr$ now shows what we want.
\end{proof}
\begin{corollarys}
\label{ref:3.1.8a}
Assume that $\Cscr$ has exact direct limits and that $\Dscr$ is
locally noetherian. Then the category of
coherent  objects in $\BIMOD(\Cscr-\Dscr)$ is an abelian subcategory of
$\BIMOD(\Cscr-\Dscr)$, closed under extensions.
\end{corollarys}
\begin{proof} According to Proposition \ref{ref:3.1.7a} the fact
  whether $\Mscr\in\BIMOD(\Cscr-\Dscr)$ is coherent can be tested on
  inverse systems of injectives $(E_i)_{i\in I}$. Since $\Dscr$ is
  locally noetherian we also have that $F=\dirlim\nolimits_i E_i$ is
  injective. Hence by construction $\HHom_\Dscr(-,E_i)$ and 
$\HHom_\Dscr(-,F)$ are exact functors on $\BIMOD(\Cscr-\Dscr)$. The
corollary is now a simple application of the five-lemma.
\end{proof}

\medskip

Sometimes it is convenient  to use ``virtual'' inverse limits of
bimodules. These are defined below.
\begin{definitions}
\label{ref:3.1.9a}
Assume that $\Cscr$ is cocomplete and that $(\Nscr_n,\phi_{m,n})$ is an inverse
system in $\BIMOD(\Cscr-\Dscr)$ indexed by $\NN$. We define
$``{\invlim_n}{}"\Nscr_n$ by the rule
\begin{equation}
\label{ref:3.2a}
\HHom_\Dscr(``{\invlim_n}{}"\Nscr_n,\Mscr)=\dirlim_n
\HHom_\Dscr(\Nscr_n,\Mscr)
\end{equation}
for all $\Mscr\in \Dscr$. 
An inverse system $(\Nscr_n,\phi_{m,n})$ such that
$``{\invlim_n}{}"\Nscr_n=0$ is called a \emph{torsion inverse system}.
\end{definitions}

\begin{lemmas}
  Assume that $\Cscr$ has exact direct limits and $\Dscr$ has enough
  injectives. Then $``{\invlim_n}{}"$ is exact. In particular the
  category of torsion inverse systems is closed under subobjects,
  quotients and extensions.
\end{lemmas}
\begin{proof}
This is trivial.
\end{proof}
\begin{lemmas}
\label{ref:3.1.11a}
 Let $(\Nscr_n,\phi_{m,n})$ be an inverse system in
  $\Dscr$ (viewed as a subcategory of
  $\MOD(\Dscr)=\BIMOD(\Ab-\Dscr)$). Then the following are equivalent
\begin{enumerate} 
\item
$(\Nscr_n,\phi_{m,n})$ is torsion.
\item 
For every $m$ there exists $n\ge m$ such that
$\phi_{m,n}:\Nscr_n\r\Nscr_m$ is the zero map.
\end{enumerate}
\end{lemmas}
\begin{proof}
We prove $(1)\Rightarrow (2)$ the other direction being obvious.
We apply \eqref{ref:3.2a} with $\Mscr=\Nscr_m$. Then the
identity map $\Id_{\Nscr_m}$ must become zero in some
$\Hom_{\Dscr}(\Nscr_n,\Nscr_m)$. This is exactly (2).
\end{proof}
\begin{remarks} Note that in the previous lemma $(1)\Rightarrow (2)$
   holds more generally for inverse systems in $\BIMOD(\Cscr-\Dscr)$.
\end{remarks}

\medskip

If $\Cscr=(\Cscr,\otimes, I)$ is a monoidal category \cite{ML} then an
algebra object in $\Cscr$ is a triple $(\Ascr,\eta,\mu)$ where $\Ascr$
is an object in $\Cscr$ equipped with two maps $\eta:I\r \Ascr$ (the
unit) and $\mu:\Ascr\otimes\Ascr\r\Ascr$ (the multiplication)
satisfying the usual compatibilities.

It is  clear that $(\BIMOD(\Dscr-\Dscr),\otimes_\Dscr,\Id_\Dscr)$ and 
$(\Bimod(\Dscr-\Dscr),\otimes_\Dscr,\Id_\Dscr)$ are monoidal
categories so we denote the algebra objects in them respectively by
$\ALG(\Dscr)$ and $\Alg(\Dscr)$. The objects in $\ALG(\Dscr)$ will be
called \emph{weak algebras} and those of $\Alg(\Dscr)$ will simply be
called \emph{algebras}. Furthermore we define $\Mod(\Ascr)$
as the category consisting of pairs $(\Mscr,h)$ where $\Mscr\in \Dscr$
and $h$ is a morphism $\Mscr\otimes_\Dscr\Ascr\r\Mscr$ in
$\MOD(\Dscr)$ satisfying the usual compatibilities. Note that objects of
$\ALG(\Dscr)$ are basically monads in the sense of \cite{ML} with some
extra structure.

In order to interprete these definitions more concretely we 
recall that
$\BIMOD(\Dscr-\Dscr)=\Lscr(\Dscr,\Dscr)^{\text{opp}}$. Thus
$\ALG(\Dscr)$ is equivalent with the category of \emph{coalgebra}
objects in $\Lscr(\Dscr,\Dscr)$. Thus if
$\Ascr=(\Ascr,\eta,\mu)\in\ALG(\Dscr)$ then the unit $\eta$ is in fact
a natural transformation
\[
\eta:\HHom_\Dscr(\Ascr,-)\r \Id_\Dscr
\]
and the multiplication $\mu$ is a natural transformation
\begin{equation}
\label{ref:3.3a}
\mu:\HHom_\Dscr(\Ascr,-)\r \HHom_\Dscr(\Ascr,\HHom_\Dscr(\Ascr,-))
\end{equation}
Likewise if $\Mscr=(\Mscr,h)\in\Mod(\Ascr)$ then $h$ is a natural
transformation 
\begin{equation}
\label{ref:3.4a}
h: \Hom_\Dscr(\Mscr,-)\r \Hom_\Dscr(\Mscr,\HHom_\Dscr(\Ascr,-))
\end{equation}
Then  $\bar{h}=h_\Mscr(\Id_\Mscr)$ defines a morphism
\begin{equation}
\label{ref:3.5a}
\bar{h}:\Mscr\r \HHom_\Dscr(\Ascr,\Mscr)
\end{equation}
Conversely if one is given a morphism $\bar{h}:\Mscr\r
\HHom_\Dscr(\Ascr,\Mscr)$ as in \eqref{ref:3.5a} then the composition
\[
\Hom_\Dscr(\Mscr,-)\r \Hom_\Dscr(\HHom_\Dscr(\Ascr,\Mscr),
\HHom_\Dscr(\Ascr,-)) \xrightarrow{\bar{h}^\ast}
\Hom_\Dscr(\Mscr,\HHom_\Dscr(\Ascr,-)) 
\]
is a natural transformation as in \eqref{ref:3.4a}. 

Elaborating on this one obtains the following results~:
\begin{propositions} Let $\Ascr\in\ALG(\Dscr)$. Then $(\Mscr,h)\mapsto
  (\Mscr,\bar{h})$ defines an isomorphism between $\Mod(\Ascr)$ and
  the category of $\Ascr$-comodules where we consider $\Ascr$  as
  a coalgebra in $\Lscr(\Dscr,\Dscr)$.
\end{propositions}
\begin{lemmas}
\label{ref:3.1.14a}
Let $\Ascr\in \ALG(\Dscr)$. The forgetful functor
\begin{equation}
\label{ref:3.6a}
(-)_\Dscr:\Mod(\Ascr)\r\Dscr:(\Mscr,h)\mapsto \Mscr
\end{equation}
 has a right adjoint given by $\HHom_\Dscr(\Ascr,-)$ (with its
canonical $\Ascr$-structure given by \eqref{ref:3.3a}) and furthermore if
$\Ascr\in\Alg(\Dscr)$ then \eqref{ref:3.6a} also has a left adjoint given by
$-\otimes_\Dscr\Ascr$.
\end{lemmas}
$\Mod(\Ascr)$ inherits most of the good properties of $\Dscr$, as is
shown in the following proposition.
\begin{propositions} 
\label{ref:3.1.15a}
Let $\Ascr=(\Ascr,\eta,\mu) \in \ALG(\Dscr)$.
\begin{enumerate}
\item
$\Mod(\Ascr)$ is an abelian category.
\item
$\Mod(\Ascr)$ possesses all colimits which exist in $\Dscr$.
\item
The forgetful functor $(-)_\Dscr$ is exact, faithful  and commutes with
colimits.
\item
$\Mod(\Ascr)$ is cogenerated by objects of the form
$\HHom_\Dscr(\Ascr,\Mscr)$, $\Mscr\in\Dscr$.
\item If $E\in \Inj(\Dscr)$ then $\HHom_\Dscr(\Ascr,E)\in
  \Inj(\Mod(\Ascr))$. In particular if $\Dscr$ has enough injectives
  then so does $\Mod(\Ascr)$.
\item
If $\Dscr$ has exact direct limits then so does $\Mod(\Ascr)$.
\end{enumerate}
Assume now in addition that $\Ascr\in\Alg(\Dscr)$
\begin{enumerate}
\setcounter{enumi}{6}
\item $\Mod(\Ascr)$ possesses all limits that exist in $\Dscr$ and
  $(-)_\Dscr$ commutes with these limits.
\item
$\Mod(\Ascr)$ is generated by objects of the form
$\Mscr\otimes_\Dscr\Ascr$, $\Mscr\in\Dscr$. Hence if $\Dscr$ has a
generator then so does $\Mod(\Ascr)$.
\end{enumerate}
In particular combining (6),(8) we find that if $\Ascr\in\Alg(\Dscr)$
and $\Dscr$ is a Grothendieck category then so is $\Mod(\Ascr)$.
\end{propositions}
\begin{proof}
\begin{itemize}
\item[2.]
Let $I$ be a small category and $\Mscr:I\r \Mod(\Ascr)$ a functor. We
write
$\Mscr(i)=(\Mscr_i,h_i)$. Then one easily verifies that $\dirlim
\Mscr$ is given by the pair $(\dirlim\Mscr_i,h)$ where $\bar{h}$ is
given  by the composition
\[
\dirlim\Mscr_i\xrightarrow{\dirlim \bar{h}_i} \dirlim
\HHom_\Dscr(\Ascr,\Mscr_i)\xrightarrow{\text{can}}
\HHom_\Dscr(\Ascr,\dirlim\Mscr_i) 
\]
\item[1.] From (2)\ it follows that $\Mod(\Ascr)$ has cokernels. So we
  have to show that $\Mod(\Ascr)$ has kernels. Let $f:(\Mscr,h)\r
  (\Nscr,j)$ be a morphism in $\Mod(\Ascr)$. Then $\ker f$ is the pair
$(\Kscr,s)$ where $\Kscr=\ker(\Mscr\r\Nscr)$ and $\bar{s}$ is as in
the following commutative diagram with exact rows.
\[
\begin{CD}
0 @>>> \Kscr @>>> \Mscr @>f>>\Nscr\\
@. @V\bar{s}VV @V\bar{h}VV @V\bar{j}VV\\
0 @>>> \HHom_\Dscr(\Ascr,\Kscr) @>>> \HHom_\Dscr(\Ascr,\Mscr) @>f>>
\HHom_\Dscr(\Ascr,\Nscr)
\end{CD}
\]
\item[3.] This follows from the explicit constructions of kernels and
  colimits in (1) and (2).
\item[4.] Assume that $(\Mscr,h)\in \Mod(\Ascr)$. The composition
\[
\Mscr \xrightarrow{\bar{h}} \HHom_\Dscr(\Ascr,\Mscr)
\xrightarrow{\eta}
\Mscr
\]
is the identity so $\bar{h}$ is a monomorphism. We know already that
$\HHom_\Dscr(\Ascr,\Mscr)$ has a canonical structure as
$\Ascr$-module and it is easy to see that $\bar{h}$ is compatible with it.
\item[5.] This follows from the fact that $\HHom_\Dscr(\Ascr,-)$ has
  a left adjoint which is exact by (3).
\item[6.] This follows from the explicit construction of kernels,
  cokernels and colimits in (1)(2).
\item[7.] Let $\Mscr:I\r \Mod(\Ascr)$ be as in (2) Then $\invlim \Mscr$
  is the pair $(\invlim\Mscr_i, {h}')$ where   $\bar{h}'$ is the
  composition
\[
\invlim \Mscr_i \xrightarrow{\invlim \bar{h}_i} \invlim
\HHom_\Dscr(\Ascr,\Mscr_i) \xrightarrow {\text{can}^{-1}}
\HHom_\Dscr(\Ascr,\invlim\Mscr_i)
\]
where we have used the fact that $\Ascr$ preserves products.
\item[8.] If $(\Mscr,h)\in\Mod(\Ascr)$ then
  $\Mscr\otimes_\Dscr\Ascr\in\Mod(\Ascr)$ by lemma \ref{ref:3.1.14a}.
  Furthermore the composition
\[
\Mscr\xrightarrow{\eta} \Mscr\otimes_\Dscr \Ascr \xrightarrow{h} \Mscr
\]
is the identity so $\Mscr\otimes_\Dscr \Ascr\r \Mscr$ is an
epimorphism.\qed
\end{itemize}
\def\qed{}
\end{proof}
Assume now that $f:\Ascr\r\Bscr$ is a morphism in $\ALG(\Dscr)$. Let
$\Mscr=(\Mscr,h)\in \Mod(\Ascr)$, $\Nscr=(\Nscr,j)\in
\Mod(\Bscr)$. Then we define 
$\Nscr_\Ascr\in\Mod(\Ascr)$ as the pair $(\Nscr,j')$ were $\bar{j}'$
is the composition
\[
\Nscr \xrightarrow {\bar{j}} \HHom_\Dscr(\Bscr,\Nscr)
\xrightarrow {f^\ast} \HHom_\Dscr(\Ascr,\Nscr)
\]
$\HHom_\Ascr(\Bscr,\Mscr)\in\Mod(\Bscr)$ is the pair $(\Uscr,u)$ where
$\Uscr$ is obtained as the equalizer of 
\begin{equation}
\label{ref:3.7a}
\begin{CD}
\HHom_\Dscr(\Bscr,\Mscr) @>\text{mult}^\ast >>
\HHom_\Dscr(\Bscr\otimes_\Dscr\Bscr, \Mscr)\\
@V \bar{h} VV @V (1\otimes f)^\ast VV\\
\HHom_\Dscr(\Bscr,\HHom_\Dscr(\Ascr,\Mscr)) @=
\HHom_\Dscr(\Bscr\otimes_\Dscr\Ascr,\Mscr) 
\end{CD}
\end{equation}
and $\bar{u}$ is obtained from the fact that all objects in
\eqref{ref:3.7a} carry a $\Bscr$-structure, and the maps are compatible
with it.

Finally if $\Ascr,\Bscr\in \Alg(\Dscr)$ then $\Mscr\otimes_\Ascr\Bscr$
is defined as the coequalizer of
\[
\Mscr\otimes_\Dscr\Ascr\otimes_\Dscr\Bscr 
\mathop{\Longrightarrow}\limits^{\bar{h}\otimes 1}_{(1\otimes\text{mult})
\circ
(1\otimes
  f\otimes 1)} \Mscr\otimes_\Dscr\Bscr
\]
So, summarizing, we have constructed the standard functors
\begin{align}
(-)_\Ascr&:\Mod(\Bscr)\r \Mod(\Ascr)\label{ref:3.8a}\\
\HHom_\Ascr(\Bscr,-)&:\Mod(\Ascr)\r \Mod(\Bscr) \label{ref:3.9a}
\end{align}
and if $\Ascr, \Bscr\in \Alg(\Dscr)$
\begin{equation}
\label{ref:3.10a}
-\otimes_\Ascr\Bscr :\Mod(\Ascr) \r \Mod(\Bscr)
\end{equation}
In general \eqref{ref:3.9a} is a right adjoint to \eqref{ref:3.8a} and
if \eqref{ref:3.10a} is defined then it is a left adjoint to
\eqref{ref:3.8a}. From the constructions of limits, colimits, kernels
and cokernels in (the proof of) Proposition \ref{ref:3.1.15a} one
verifies that \eqref{ref:3.8a} \eqref{ref:3.9a}\eqref{ref:3.10a} have
the standard exactness properties and satisfy the usual
compatibilities with (co)limits.

For $\Ascr,\Bscr\in \ALG(\Dscr)$ we define
\begin{align*}
\BIMOD(\Ascr-\Bscr)&=\BIMOD(\Mod(\Ascr)-\Mod(\Bscr))\\
\Bimod(\Ascr-\Bscr)&=\Bimod(\Mod(\Ascr)-\Mod(\Bscr))
\end{align*}
We denote $-\otimes_{\Mod(\Bscr)}-$ by $-\otimes_\Bscr-$. In general
we will replace in  our notations $\Mod(\Xscr)$ by $\Xscr$ when no
confusion can arise.

Assume that $\Ascr\r \Ascr'$, $\Bscr\r \Bscr'$ are morphisms in
$\ALG(\Dscr)$. We will now define functors
\begin{align}
\Ascr'\otimes_\Ascr-\otimes_\Bscr\Bscr':\BIMOD(\Ascr-\Bscr)\r
\BIMOD(\Ascr'-\Bscr') \label{ref:3.11a}\\
{}_\Ascr(-)_\Bscr:\BIMOD(\Ascr'-\Bscr')\r \BIMOD (\Ascr-\Bscr)
\label{ref:3.12a}
\end{align}
as we did for modules.

Let $\Mscr\in\BIMOD(\Ascr-\Bscr)$,
$\Nscr\in\BIMOD(\Ascr'-\Bscr')$. Then we define
\[
\HHom_{\Bscr'}(\Ascr'\otimes_\Ascr\Mscr\otimes_\Bscr\Bscr',-)\overset{
\text{def}}{=} \HHom_\Ascr(\Ascr',\HHom_\Bscr(\Mscr,(-)_\Bscr))
\]
\begin{equation}
\label{ref:3.13a}
\HHom_\Bscr({}_\Ascr\Nscr_\Bscr,-)\overset{
\text{def}}{=}
\HHom_{\Bscr'}(\Nscr,\HHom_\Bscr(\Bscr',-))_\Ascr
\end{equation}
The functor
$-\otimes_{\Ascr'}\Ascr'\otimes_\Ascr\Mscr\otimes_\Bscr\Bscr'$ should
be given by
\[
(-)_\Ascr\otimes_\Ascr\Mscr\otimes_\Bscr\Bscr'
\]
and the functor $-\otimes_\Ascr{}_\Ascr\Nscr_\Bscr$ should be given by
\begin{equation}
\label{ref:3.14a}
((-\otimes_\Ascr\Ascr')\otimes_{\Ascr'}\Nscr)_\Bscr
\end{equation}
So we conclude that if $\Mscr\in \Bimod(\Ascr-\Bscr)$ and
$\Nscr\in\Bimod(\Ascr'-\Bscr')$ then if $\Bscr'\in \Alg(\Dscr)$ then
\eqref{ref:3.11a} respects ``Bimod'' and if $\Ascr'\in \Alg(\Dscr)$ then
\eqref{ref:3.12a} respects ``Bimod''.

Our definition of $\BIMOD(\Ascr-\Bscr)$ has the advantage that we can
directly apply Proposition \ref{ref:3.1.1a}  to obtain the properties of
this category. However we would also like to have a definition which resembles
more closely that of modules. Therefore we state the following proposition.
\begin{propositions} The following categories are equivalent. 
\begin{enumerate}
\item
$\BIMOD(\Ascr-\Bscr)$
\item
The category of triples $(\Mscr,h,h')$ where
$\Mscr\in\BIMOD(\Dscr-\Dscr)$ and $h:\Ascr\otimes_\Dscr \Mscr\r
\Mscr$, $h':\Mscr\otimes_\Dscr\Bscr\r \Mscr$ are maps in
$\BIMOD(\Dscr-\Dscr)$ satisfying the usual compatibilities.
\end{enumerate}
If $\Ascr,\Bscr\in \Alg(\Dscr)$ then in the previous statement,
``$\BIMOD$'' may be replaced by ``$\Bimod$''.
\end{propositions}
\begin{proof} This is a rather tedious verification which we leave to
the reader. Let us simply state how one associates a left exact
functor $\HHom_\Bscr(\Mscr,-):\Mod(\Bscr)\r \Mod(\Ascr)$ to a triple
$(\Mscr,h,h')$.

Let $\Nscr=(\Nscr,p)\in \Mod(\Bscr)$. Then $\HHom_\Bscr(\Mscr,\Nscr)$
will be the equalizer of 
\begin{equation}
\label{ref:3.15a}
\begin{CD}
\HHom_\Dscr(\Mscr,\Nscr) @= \HHom_\Dscr(\Mscr,\Nscr)\\
@V h^{\prime\ast}VV @V \bar{p} VV\\
\HHom_\Dscr(\Mscr\otimes_\Dscr\Bscr,\Nscr) @=\HHom_\Dscr(\Mscr,
\HHom_\Dscr(\Bscr,\Nscr))
\end{CD}
\end{equation}
The $\Ascr$-structure on $\HHom_\Bscr(\Mscr,\Nscr)$ is obtained from
the fact that all objects in \eqref{ref:3.15a} carry a canonical
$\Ascr$-structure and the maps are compatible with it. 
\end{proof}
Now let $\Ascr\r \Bscr$ be a map in $\ALG(\Dscr)$. From \eqref{ref:3.11a}
we obtain a functor
\[
-\otimes_{\Ascr}\Bscr:\MOD(\Ascr)\r \MOD(\Bscr)
\]
Assume that $\Mod(\Ascr)$ and $\Mod(\Bscr)$ have enough injectives. Then
we define associated functors
\[
\HTor_i^\Ascr(-,\Bscr):\MOD(\Ascr)\r\MOD(\Bscr)
\]
by
\[
\Hom_\Bscr(\HTor_i^\Ascr(-,\Bscr),E)
=
\Ext^i_\Ascr(-,E_\Ascr)
\]
for injectives $E\in\Mod(\Bscr)$. One easily verifies that
\[
\HTor_i^\Ascr(-,\Bscr)_\Ascr=\HTor_i^\Ascr(-,\Bscr_\Ascr)
\]
\begin{lemmas}
\label{ref:3.1.17a}
Assume that $\Dscr$ has colimits and that $\Mod(\Bscr)$ and $\Mod(\Ascr)$
are locally noetherian. Then $\HTor_i^\Ascr(-,\Bscr)$ sends coherent
objects in $\MOD(\Ascr)$ to  coherent objects in $\MOD(\Bscr)$.
\end{lemmas}
\begin{proof} This is a formal verification. We will give the proof as an
illustration of how the various hypotheses are used.

Assume that $\Mscr\in\MOD(\Ascr)$ is coherent and let $(E_i)_i$ be a
directed system of injectives in $\Mod(\Bscr)$. By Proposition
\ref{ref:3.1.7a} we have to show that
\[
\Hom_\Bscr(\HTor_i^\Ascr(\Mscr,\Bscr),\dirlim_i E_i)
=
\dirlim_i \Hom_\Bscr(\HTor_i^\Ascr(\Mscr,\Bscr),E_i)
\]
(taking into account that $\Ab$ has exact direct limits). 

From the fact that $\Mod(\Bscr)$ is locally noetherian it follows
that $\dirlim_i E_i$ is injective. Hence we have to show that 
\[
\Ext^i_\Ascr(\Mscr,(\dirlim_i E_i)_\Ascr)=
\dirlim_i \Ext^i_\Ascr(\Mscr,E_{i\Ascr})
\]
The fact that $\Mod(\Ascr)$ is locally noetherian and that $\Mscr$ is
coherent implies that the
righthand side of this equation is equal to 
$ \Ext^i_\Ascr(\Mscr,\dirlim_i E_{i\Ascr})$. So it remains to show that
$(\dirlim_i E_i)_\Ascr=\dirlim_i E_{i\Ascr}$.  This follows easily from
the explicit construction of $\dirlim$ in the proof of Proposition
\ref{ref:3.1.15a}.2 and the definition of $(-)_\Ascr$.
\end{proof}

\subsection{Graded modules, bimodules and algebras}
\label{ref:3.2b}
In this section we have stated all our definitions in the ungraded
case but we will mainly need them in the graded case. Luckily the
generalization is trivial.

If $\Dscr$ is an abelian category then we denote by $\tilde{\Dscr}$
the category of $\ZZ$-graded objects over $\Dscr$. Thus by definition
an object in $\tilde{\Dscr}$ is  a sequence of objects
$(\Mscr_n)_n$ in $\Dscr$ and
$\Hom_{\tilde{\Dscr}}((\Mscr_n)_n,(\Nscr_n)_n)=
\prod_n\Hom_\Dscr(\Mscr_n,\Nscr_n)$.
We identify an object $\Mscr\in\Dscr$ with the object of
$\tilde{\Dscr}$ which is $\Mscr$ in degree zero and zero elsewhere. As
usual $\tilde{\Dscr}$ is equipped with the shiftfunctor $\Mscr\mapsto
\Mscr(1)$ where $\Mscr(1)_n=\Mscr_{n+1}$.
Of course $\Mscr(m)$, $m\in\ZZ$ is defined similarly. To simplify the
notation we will often write $\oplus_n\Mscr_n$ for $(\Mscr_n)_n$.

By analogy with the ungraded case we define
\[
\BIGR(\Cscr-\Dscr)=\{\text{left exact functors $\tilde{\Dscr}\r\tilde{\Cscr}$
  commuting with shift}\}^{\text{opp}}
\]
\[
\Bigr(\Cscr-\Dscr)=\{\Mscr\in\BIGR(\Cscr-\Dscr)\mid\Mscr\text{ has a
  left adjoint}\}
\]
and
\[
\GR(\Cscr)=\BIGR(\Ab-\Cscr)
\]
Objects in $\BIGR(\Cscr-\Dscr)$ will be called \emph{graded weak
 $\Cscr$-$\Dscr$ bimodules} and objects in $\Bigr(\Cscr-\Dscr)$ will
 be called \emph{graded bimodules}.

We will denote the left exact functor corresponding to $\Mscr \in
\BIGR(\Cscr-\Dscr)$ with $\underline{\HHom}_\Dscr(\Mscr,-)$ in order
to avoid confusion with the ungraded case. Similar conventions will
apply to the use of $\underline{\HExt}$, $\underline{\HTor}$.

Note that if $\Mscr\in \BIGR(\Cscr-\Dscr)$ (resp.\ $\Bigr(\Cscr-\Dscr)$)
then the composition (for $n\in\ZZ$)
\[
\Dscr\hookrightarrow \tilde{\Dscr}\xrightarrow{\Mscr} \tilde{\Cscr}
\xrightarrow{ \text{degree $-n$}} \Cscr
\]
defines objects $\Mscr_{n}$ in $\BIMOD(\Cscr-\Dscr)$ (resp.
$\Bimod(\Cscr-\Dscr)$) which determine $\Mscr$. We write
$\Mscr=\oplus_n \Mscr_n$. 

Obviously $\BIGR(\Dscr-\Dscr)$, $\Bigr(\Dscr-\Dscr)$ are monoidal
categories and we will denote the algebra objects in them by
$\GRALG(\Dscr)$ resp. $\Gralg(\Dscr)$. Consistent with our earlier
coventions we speak respectively of \emph{weak algebras} and
\emph{algebra}. It is easy to see that $\Ascr\in \GRALG(\Dscr)$
is of the form $\Ascr=\oplus_n\Ascr_n$ with multiplication maps
$\Ascr_m\otimes_\Dscr\Ascr_n\r\Ascr_{m+n}$ and a unit $\Id_\Dscr\r
\Ascr_0$ satisfying the usual compatibilities. 

If $\Ascr\in \GRALG(\Dscr)$ then $\Gr(\Ascr)$ is defined as
$\Mod(\Ascr)$ but with $\Dscr$ replaced by $\tilde{\Dscr}$. Finally if  
$\Ascr,\Bscr\in\GRALG(\Dscr)$ then 
\begin{align*}
\Bigr(\Ascr-\Bscr)&=\BIGR(\Gr(\Ascr)-\Gr(\Bscr))\\
 \Bigr(\Ascr-\Bscr)&=\Bigr(\Gr(\Ascr)-\Gr(\Bscr))
\end{align*}
In the sequel we will freely use graded versions of ungraded results
stated in this section.

\subsection{Algebras which are strongly graded with respect to a Serre
subcategory}
\label{ref:3.3b}
We put this section here for lack of a better place. The results will
be used in \S\ref{ref:3.11b}.

Let $\Dscr$ be an abelian category and let $\Ascr\in\Gralg(\Dscr)$.
Put $\Qscr_n=\coker (\Ascr_n\otimes_\Dscr\Ascr_{-n}\r \Id_\Dscr)$. It is
clear that $\Qscr_n\in\Bimod(\Dscr-\Dscr)$. We call $\Ascr$
\emph{strongly graded} if $\Qscr_n=0$ for all $n$. By copying the proof
from the ring case \cite[Thm. I.3.4]{NVO} one easily shows that the
functors $(-)_0$ and $-\otimes_{\Ascr_0}\Ascr$ define inverse
equivalences between $\Gr(\Ascr)$ and $\Mod(\Ascr_0)$.

Let $\Sscr$ be as Serre subcategory of $\Mod(\Ascr_0)$. Define
$\Sscr(\Ascr)$  as the full subcategory of objects $\Tscr$ in $\Gr(\Ascr)$
such that $\Tscr_n$ is in $\Sscr$ for every $n$. It is clear that
$\Sscr(\Ascr)$ is a Serre subcategory of $\Gr(\Ascr)$.

We say that $\Ascr$ is
strongly graded with respect to $\Sscr$ if $-\otimes_{\Ascr_0}\Qscr_n$
sends $\Mod(\Ascr_0)$ to $\Sscr$ and if 
$-\otimes_{\Ascr_0}\Ascr$ sends $\Sscr$ to $\Sscr(\Ascr)$.

Then one has the following:
\begin{lemmas} 
\label{ref:3.3.1a}
Assume that $\Ascr$ is strongly graded with respect to
$\Sscr$. Then
the functors $(-)_0$ and $-\otimes_{\Ascr_0}\Ascr$ factor over
$\Sscr(\Ascr)$ and $\Sscr$ and in this way they define inverse 
equivalences between $\Gr(\Ascr)/\Sscr(\Ascr)$ and $\Mod(\Ascr_0)/\Sscr$.
\end{lemmas}
\begin{proof} For the convenience of the reader we copy the proof of 
\cite[Thm. I.3.4]{NVO} suitably modified.

If $\Mscr\in\Gr(\Ascr)$ and if $\Pscr_n$ is a submodule of $\Mscr_n$ then
we define $\Pscr_n\Ascr_m$ as the image of $\Pscr_n\otimes_\Dscr\Ascr_m$ 
in $\Mscr_{m+n}$.

We  show that $\Mscr=\Mscr_0\otimes_{\Ascr_0}\Ascr$ 
modulo $\Sscr(\Ascr)$. Working modulo $\Sscr$ we have the following.
\begin{equation}
\label{ref:3.16a}
\Mscr_{m+n}=\Mscr_{m+n}\Ascr_0=\Mscr_{m+n}\Ascr_{-m}\Ascr_{m}\subset
\Mscr_n\Ascr_m\subset\Mscr_{m+n}
 \end{equation}
 The
second equality follows from
\[
\Mscr_{m+n}/\Mscr_{m+n}\Ascr_{-m}\Ascr_{m}=\coker
(\Mscr_{m+n}\otimes_\Dscr\Ascr_{-m}\otimes_\Dscr\Ascr_{m}\r
\Mscr_{m+n})=\Mscr_{m+n}\otimes_\Dscr\Qscr_{-m}
\]
\eqref{ref:3.16a} yields $\Mscr_n\Ascr_m=\Mscr_{m+n}$ modulo
$\Sscr$. Now we claim that  the multiplication map
\[
\mu:\Mscr_0\otimes_{\Ascr_0}\Ascr\r \Mscr
\]
is in fact an isomorphism modulo $\Sscr(\Ascr)$. By what we have done
so far $\mu$ is clearly surjective modulo $\Sscr(\Ascr)$. Let $\Kscr$ be
the kernel of $\mu$. We have $\Kscr_0=0$ and hence modulo $\Sscr$,
$\Kscr_n=\Kscr_0\Ascr_n=0$. Thus $\Kscr\in\Sscr(\Ascr)$ and we are done.

We can now finish the proof, provided we can show that the functor
$-\otimes_{\Ascr_0}\Ascr$ is well defined on
$\Mod(\Ascr_0)/\Sscr$. Inspection shows that we only have to check the
following case~: assume that $\Mscr\r \Nscr$ is injective in
$\Mod(\Ascr_0)$. Then $\Kscr=\ker (\Mscr\otimes_{\Ascr_0}\Ascr\r
\Nscr\otimes_{\Ascr_0}\Ascr)$ should be in $\Sscr(\Ascr)$. Since
$\Kscr_0=0$ this follows from the previous discussion.  
\end{proof}

\subsection{Quotients of the identity functor} 
\label{ref:3.4b}
This section is related
to \cite{rosenberg}. $\Cscr$, $\Dscr$, $\Escr$ will
be abelian categories having enough injectives so that the weak
bimodule categories are abelian (cf Prop.\ \ref{ref:3.1.1a}).  Let
$\Bscr$ be a quotient of $\Id_\Dscr$ in $\ALG(\Dscr)$ (by this we mean
that the underlying $\Dscr$-bimodule map is an epimorphism). We will
denote the functor $(-)_\Dscr$ by $i_\ast$ and its right adjoint
$\HHom_\Dscr(\Bscr,-)$ by $i^!$. If $\Bscr\in\Alg(\Dscr)$ then we
denote $-\otimes_\Dscr\Bscr$ by $i^\ast$. This is the left adjoint to
$i_\ast$.
\begin{propositions}
\label{ref:3.4.1a}
\begin{enumerate}
\item $i^!i_\ast=\Id_{\Mod(\Bscr)}$
\item $i_\ast$ is fully faitful.
\item Let $\Mscr\in\Dscr$. The counitmap $i_\ast i^!\Mscr\r \Mscr$ is
  injective. Furthermore $\Mscr\in i_\ast\Mod(\Bscr)$ if and only if
  the counit map is an isomorphism.
\item $i_\ast\Mod(\Bscr)$ is  closed
under subquotients.
\item If $\Bscr\in\Alg(\Dscr)$ then $i^\ast i_\ast=\Id_{\Mod(\Bscr)}$.
\end{enumerate}
\end{propositions}
\begin{proof}
\begin{enumerate}
\item Let $\Mscr\in\Mod(\Bscr)$. The canonical map
  $\Mscr\r\HHom_\Dscr(\Bscr,\Mscr)$ (cfr.\ \eqref{ref:3.5a}) is the unit
  for the adjunction $(i_\ast,i^!)$. Thus we have to show that it is
  an isomorphism. Since $\Bscr$ is a quotient of $\Id_\Bscr$ we can
  consider the commutative diagram
\[
\begin{CD}
\Mscr @>>> \HHom_\Dscr(\Bscr,\Mscr)\\
@| @VVV\\
\Mscr @= \HHom_\Dscr(\Id_\Dscr,\Mscr)
\end{CD}
\]
where the right and upper maps are injective (see Prop.\
3.1.10.4). This implies that the upper map is
an isomorphism.
\item This is formally a consequence of (1). Indeed if
  $\Mscr,\Nscr\in\Mod(\Bscr)$ then
\[
\Hom_\Dscr(i_\ast\Mscr,i_\ast\Nscr)=\Hom_\Bscr(\Mscr,i^!i_\ast\Nscr)=
\Hom_\Bscr(\Mscr,\Nscr)
\]
\item The counit map is given by the composition
\[
i_\ast i^!\Mscr=\HHom_\Dscr(\Bscr,\Mscr)\hookrightarrow
\HHom_\Dscr(\Id_\Dscr,\Mscr)=\Mscr
\]
so it is certainly injective. If it is an isomorphism then clearly
$\Mscr=i_\ast(i^!\Mscr)\in i_\ast\Mod(\Bscr)$. Conversely if
$\Mscr=i_\ast\Nscr$ then $i_\ast i^! \Mscr =i_\ast i^!
i_\ast\Nscr=i_\ast\Nscr=\Mscr$.
\item
Since $i_\ast$ is left exact $i_\ast\Mod(\Bscr)$ is closed under
kernels.
 So it suffices to show
that $i_\ast\Mod(\Bscr)$ is closed under quotients.

Let $\Mscr\in i_\ast\Mod(\Bscr)$ and let $\Nscr$ be a quotient of $\Mscr$ in
$\Dscr$. Then we have a commutative diagram
\[
\begin{CD}
i_\ast i^!\Mscr@>\alpha>> \Mscr\\
@VVV @V\beta VV\\
i_\ast i^! \Nscr@>\gamma >>\Nscr
\end{CD}
\]
where $\alpha$ is an isomorphism, $\beta$ is surjective and $\gamma$ is
injective. This implies that $\gamma$ is an isomorphism.
\item This is proved similarly as (1).  \qed\end{enumerate}
\def\qed{}\end{proof} Let us call a full subcategory of an abelian
category \emph{closed} if it is closed under subquotients and if the
inclusion functor has a right adjoint.  In most cases this is
equivalent to the definition in \cite[IV.4]{Gabriel}. Let us call a
closed subcategory \emph{biclosed} if the inclusion functor also has a
left adjoint.

Proposition \ref{ref:3.4.1a} basically tells us that the functor
$\Bscr\r i_\ast\Mod(\Bscr)$ associates with every quotient of
$\Id_\Dscr$ a closed subcategory of $\Dscr$ and if the quotient is in
$\Alg(\Dscr)$ then the resulting category is biclosed. We will now
show that the converse to this essentially holds.

Assume that $\Escr\subset\Dscr$ is a closed subcategory. Let $i_\ast$ be the
inclusion functor and $i^!$ its right adjoint. Let $\Bscr$ be the left exact
functor $i_\ast i^!$. Then the counit $i_\ast i^!\r \Id_\Dscr$ and the
comultiplication
$i_\ast i^!=i_\ast\Id_\Dscr i^!\r i_\ast i^!i_\ast i^!$ make $\Bscr$
into a weak algebra (remember that a ``weak algebra'' is actually a
coalgebra object in the the monoidal category of left exact functors, see the
discussion after Prop.\ \ref{ref:3.1.5a}). Furthermore the counit $i_\ast i^!\r \Id_\Dscr$ makes
$\Bscr$ into a quotient of $\Id_\Dscr$ (again remember that arrows are reversed
if we pass from functors to bimodules).

Now by a routine verification (or using an appropriate version of
Beck's theorem \cite{ML}) one shows that $\Mod(\Bscr)\cong
\Escr$. Elaborating on this one may prove the following result.
\begin{propositions}
\label{ref:3.4.2a}
 With notations as above. The functors $\Bscr\mapsto i_\ast
\Mod(\Bscr)$ and $\Escr \mapsto i_\ast i^!$ induce inverse bijections between
the quotients of $\Id_\Dscr$ in $\ALG(\Dscr)$ (resp.\ $\Alg(\Dscr)$) and the
closed (resp. biclosed) subcategories in $\Dscr$.
\end{propositions}
Below we state a few elementary results concerning closed and biclosed
categories which will be useful in the sequel. Recall that a category
is said to be well-powered if the set of subobjects of an arbitrary
object is
small. In an abelian category this holds if there is a generator or a
cogenerator. 
\begin{propositions} 
\label{ref:3.4.3a}
Let $\Escr$ be a full subcategory of the abelian category
$\Dscr$.
\begin{enumerate}
\item If $\Escr$ is closed in $\Dscr$ and if $\Dscr$ is cocomplete,
  then so is $\Escr$.
\item
If $\Dscr$ is cocomplete and well-powered and if $\Escr$ is closed under
subquotients and direct sums (in $\Dscr$) then $\Escr$ is closed.
\item If $\Escr$ is biclosed and if $\Dscr$ is complete then so is 
$\Escr$.
\item
If $\Dscr$ is complete and well-powered and if $\Escr$ is closed
 and is closed under products  (in $\Dscr$) then $\Escr$ is biclosed.
\end{enumerate}
\end{propositions}
\begin{proof}
Let $i_\ast:\Escr\r\Dscr$ be the inclusion functor and let $i^!$, $i^\ast$ be
the right and left adjoint to $i_\ast$ if they exist.
\begin{enumerate}
\item
It is sufficient to show that $\Escr$ is closed under direct
sums. Let $\oplus_{j\in J} M_j$ be such a direct sum. By
construction $i^!(\oplus_{j\in J} M_j)$ is the maximal subobject of
$\oplus_{j\in J} M_j$ contained in $\Escr$. Hence for all $j$,
$M_j\subset i^!(\oplus_{j\in J} M_j)$. But this implies $i^!(\oplus_{j\in J}
M_j)=\oplus_{j\in J} M_j$ and thus $\oplus _{j\in J} M_j\in \Escr$.
\item Since the set of subobjects of an object $M$ is small and since
  $\Escr$ is closed under (small) direct unions, $M$ has a largest subobject
  $N$ which lies in $\Escr$. The assignment $M\mapsto N$ is the right
  adjoint to the inclusion functor.
\item This is proved by an argument dual to (1).
\item This is proved by an argument dual to (2). \qed
\end{enumerate}
\def\qed{}
\end{proof}
\begin{propositions} 
\label{ref:3.4.4a}
If $\Dscr$ is a closed subcategory of $\Cscr$ and $\Escr$
is a closed subcategory of $\Dscr$ then $\Escr$ is closed in $\Cscr$.
A similar statement holds if we replace ``closed'' by ``biclosed''.
\end{propositions}
\begin{proof} This follows from the fact that adjoint functors are compatible
with composition.
\end{proof}
Let us recall the definition of the \emph{Gabriel product}. If
$\Escr_1,\Escr_2$ are full subcategories of an abelian category $\Dscr$ then
$\Escr_1\cdot\Escr_2$ is the full subcategory of $\Dscr$ whose objects are
given by middle terms of exact sequences
\[
0\r M_2\r M\r M_1
\r
0
\]
with $M_1\in \Escr_1$, $M_2\in \Escr_2$. It is easy to see that if $\Escr_1$,
$\Escr_2$ are closed under subquotients then so is $\Escr_1\cdot \Escr_2$.
\begin{propositions} 
\label{ref:3.4.5a}
\cite{rosenberg} If $\Escr_1,\Escr_2\subset \Dscr$ are
(bi)closed then so is $\Escr_1\cdot\Escr_2$.
\end{propositions}
\begin{proof}
  We give the proof for closedness. Biclosedness is similar.  Let
  $i_{1\ast}$, $i_{2\ast}$, $i_{12\ast}$ be the embeddings
  $\Escr_1\subset \Dscr$, $\Escr_2\subset \Dscr$,
  $\Escr_1\cdot\Escr_2\subset \Dscr$. One has to construct the right
  adjoint to $i_{12\ast}$. Let $q:\Mscr\r \Mscr/i^!_2\Mscr$ be the
  quotient map. Then one verifies that
\[
i^!_{12}\Mscr \overset{\text{def}}{=} q^{-1} (i^!_1(\Mscr/i^!_2\Mscr))
\]
has the required properties.
\end{proof}
\subsection{Ideals in the identity functor}
\label{ref:3.5b}
In this section $\Cscr$, $\Dscr$, $\Escr$ will be abelian categories  having
enough injectives.
\begin{definitions} If $\Ascr\in \ALG(\Dscr)$ (resp.\ in
  $\Alg(\Dscr)$) then a weak ideal (resp. an
  ideal) in $\Ascr$ is a subobject of ${}_\Ascr \Ascr_\Ascr$ in
  $\BIMOD(\Ascr-\Ascr)$ (resp.\ in $\Bimod(\Ascr-\Ascr)$).
\end{definitions}
If $I$, $J$ are weak ideals in $\Ascr\in \ALG(\Dscr)$ then we define the weak
ideal
$IJ$ as the image of the composition
\[
I\otimes_\Ascr J\r \Ascr\otimes_\Ascr\Ascr=\Ascr
\]
If $f:\Ascr\r\Bscr$ is a morphism in $\ALG(\Dscr)$ then by $\ker f$ we denote
$\ker ({}_\Dscr\Ascr_\Dscr\r {}_\Dscr\Bscr_\Dscr)$. Since $\ker f$ is
canonically an object in $\BIMOD(\Ascr-\Ascr)$ we find that $\ker f$ is a weak
ideal in $\Ascr$.

We will say that  $f$ is surjective if the underlying bimodule map
${}_\Dscr\Ascr_\Dscr\r {}_\Dscr \Bscr_\Dscr$ is surjective. If this is the case
then we call the pair $(\Bscr,f)$ a quotient object of $\Ascr$ in
$\ALG(\Dscr)$. 

If $I$ is a weak ideal in $\Ascr$ then there is a  unique weak algebra
structure on $\Ascr/I$ which makes the quotient map $q:\Ascr\r \Ascr/I$ into a
morphism in $\ALG(\Dscr)$.
\begin{lemmas}
$I\mapsto (\Ascr/I,q)$ and $(\Bscr,f)\mapsto \ker f$ define inverse bijections
between
\begin{enumerate}
\item
Subobject of ${}_\Ascr\Ascr_\Ascr$ in $\BIMOD(\Ascr-\Ascr)$
\item
Isomorphism classes of quotient object of $\Ascr$.
\end{enumerate}
\end{lemmas}
\begin{proof} We leave this verification to the reader. The lemma can of course
be proved more generally in the setting of monoidal categories
$(\Cscr,\otimes,I)$ where $\Cscr$ is abelian and $\otimes$ is right exact.
\end{proof}
Recall that $\BIMOD(\Ascr-\Ascr)$ was defined as
$\BIMOD(\MOD(\Ascr)-\MOD(\Ascr))$ and similarly for ``$\Bimod$''.
Hence there is a 1-1 correspondence between (weak) ideals in $\Ascr$
and (weak) ideals in $\Id_{\Mod(\Ascr)}$ and furthermore this
correspondence is compatible with products. Hence to simplify the
notation below we will replace   $\Mod(\Ascr)$ by $\Dscr$ and $\Ascr$ by
$\Id_\Dscr$.

If $M\in \Mod(\Dscr)$ and $I$ is a weak ideal in $\Id_\Dscr$ then we define
$M_I$ as the quotient object of $M$ given by the image of $M\r
\HHom_\Dscr(I,M)$. One quickly verifies~:
\begin{lemmas}
\begin{enumerate}
\item
$M\in \Mod(\Id_\Dscr/I)$ iff $M_I=0$.
\item $(M_J)_I=M_{IJ}$.
\end{enumerate}
\end{lemmas}
\begin{propositions} 
\label{ref:3.5.4a}
Let $I,J$ be weak ideals in $\Id_\Dscr$. Then
\[
\Mod(\Id_\Dscr/IJ)=\Mod(\Id_\Dscr/I)\cdot \Mod(\Id_\Dscr/J)
\]
\end{propositions}
\begin{proof}
\begin{itemize}
\item[``$\subset$''] Let $M\in \Mod(\Id_\Dscr/IJ)$. Then
$M_{IJ}=0$ and hence $(M_J)_I=0$. Thus $M_J\in \Mod(\Id_\Dscr/I)$. Since by
construction $\ker (M\r M_J)\in \Mod(\Id_\Dscr/J)$ we are done.
\item[``$\supset$''] Let $M\in \Mod(\Id_\Dscr/I)\cdot \Mod(\Id_\Dscr/J)$. Then
there is an exact sequence
\[
0\r M_2\r M\r M_1\r 0
\]
with $M_2\in \Mod(\Id_\Dscr/J)$, $M_1\in \Mod(\Id_\Dscr/I)$. In the commutative
diagram
\[
\begin{CD}
M_2 @>>> M\\
@VVV @VVV\\
(M_2)_J @>>> M_J
\end{CD}
\]
we have $(M_2)_J=0$ and thus the composition $M_2\r M \r M_J$ is zero. In other
words, $M\r M_J$ factors through $M_1$. Thus $M_J$ is a quotient of $M_1$ and
hence  $M_{IJ}=(M_J)_I=0$. So $M\in \Mod(\Id_\Dscr/IJ)$. \qed
\end{itemize}
\def\qed{}
\end{proof}

Our next aim is to characterize the ideals in  $\Id_\Dscr$ in terms of certain
subcategories of $\Dscr$. The following is a trivial consequence of
Prop. \ref{ref:3.1.1a}(6) and Prop.\ \ref{ref:3.4.2a}.
\begin{propositions}
\label{ref:3.5.5a}
 Assume that $\Dscr$ is complete and has exact direct
products and an injective cogenerator. Then a weak ideal in $\Id_\Dscr$ is an
ideal if and only if $\Mod(\Id_\Dscr/I)$ is biclosed.
\end{propositions}
Hence our problems lie in the cases where products are not exact.

For simplicity we assume below that $\Dscr$ is complete and has an
injective cogenerator. Let us first recall the definition of the derived
functors of the product functor.  Assume that $J$ is some index
set and let $\prod_{j\in J}\Dscr$ be the category whose objects
consist of all families $(M_j)_{j\in J}$ with $M_j\in \Dscr$.  Then we
denote by $\prod$ the functor
\[
\prod :\prod_{j\in J}\Dscr \r \Dscr:(M_j)_{j\in J}\mapsto \prod_{j\in J}M_j
\]
which is obviously left exact. We denote the derived functors  of $\prod$ by
$(R^i\prod)_i$. These functors are computed in the standard way. Let
$(M_j)_{j\in J}$ be in $\prod_{j\in J}\Dscr$. One takes injective resolutions
$M_j\r E_j^\cdot$ and then one has
\[
R^i\prod (M_j)_{j\in J}=H^i\bigl(\prod_{j\in J} E_j^\cdot\bigr)
\]
\begin{definitions}
A biclosed subcategory $\Escr\subset \Dscr$ is called \emph{well-closed} if for
all families  $(E_j)_{j\in J}$ whose members are injective in $\Escr$ (but not
necessarily in $\Dscr$) one has $R^1\prod (E_j)_j=0$.
\end{definitions}
\begin{propositions}
\label{ref:3.5.7a}
Let $\Dscr$ be as above. That is, $\Dscr$ is complete and has an injective
cogenerator. Let $I\subset \Id_\Dscr$ be a weak ideal. Then $I$ is an ideal if
and only if $\Mod(\Id_\Dscr/I)$ is well-closed.
\end{propositions}
\begin{proof}
\begin{itemize}
\item [``$\Rightarrow$''] Since $I$ is an ideal, $\Id_\Dscr/I\in \Alg(\Dscr)$.
Hence $\Mod(\Id_\Dscr/I)$ is certainly biclosed. Now let $(E_j)_{j\in J}$
be a family
on injectives in $\Escr$ and let $(F_j)_{j\in J}$ be the corresponding injective
hulls in $\Dscr$, Since $F_j$ is an essential extension of $E_j$ and $E_j$ is
injective we have
\[
\HHom_\Dscr(\Id_\Dscr/I,F_j)=E_j
\]
Applying $\HHom_\Dscr(-,\prod_j F_j)$ to the exact sequence in
$\Bimod(\Dscr-\Dscr)$
\begin{equation}
\label{ref:3.17a}
0\r I\r \id_\Dscr\r \Id_\Dscr/I\r 0
\end{equation}
yields an exact sequence
\begin{equation}
\label{ref:3.18a}
0\r \prod_j E_j \r \prod_j F_j\r\prod_j \HHom_\Dscr(I,F_j)\r 0
\end{equation}
which is obviously the direct product of the short exact sequence
\begin{equation}
\label{ref:3.19a}
0\r E_j \r F_j \r \HHom_\Dscr(I,F_j)\r 0
\end{equation}
obtained from applying $\HHom(-,F_j)$ to \eqref{ref:3.17a}.

Now applying the long exact sequence for $(R^i\prod)_i$ to \eqref{ref:3.19a} yields
that $R^1\prod (E_j)_j=0$.
\item [``$\Leftarrow$''] According to Prop. \ref{ref:3.1.1a}(5) it
  suffices to show that $\HHom_\Dscr(I,-)$ commutes with products,
  when evaluated on injectives. Let $(F_j)_j$ be a family of
  injectives in $\Dscr$ and put $E_j=\HHom_\Dscr(\Id_{\Dscr}/I,F_j)$.
  We have again the short exact sequences given by \eqref{ref:3.19a} and
  using the fact that $R^1\prod (E_j)_j=0$ we obtain the short exact
  sequence \eqref{ref:3.18a} from the long exact sequence for $(R^i\prod)_i$.

On the other hand if we apply $\HHom_\Dscr(-,\prod (F_j)_j)$ to the exact
sequence \eqref{ref:3.17a} we obtain
\begin{equation}
\label{ref:3.20a}
0\r \prod_j E_j \r \prod_j F_j \r \HHom_\Dscr(I,\prod F_j)\r 0
\end{equation}
Comparing \eqref{ref:3.18a} and \eqref{ref:3.20a} yields what we want. \qed
\end{itemize}
\def\qed{}\end{proof}
Now we prove an analog for Proposition \ref{ref:3.4.4a}  for well-closed
subcategories.
\begin{propositions}
\label{ref:3.5.8a}
Assume that $\Cscr$ is complete and has an injective cogenerator. Assume that
$\Dscr$ is well-closed in $\Cscr$ and that $\Escr$ is well-closed in $\Dscr$.
Then $\Escr$ is well-closed in $\Cscr$.
\end{propositions}
\begin{proof}
  Note that by Prop. \ref{ref:3.4.3a} $\Dscr$ is complete and it is
  clear that $\Dscr$ has an injective cogenerator. So it makes sense
  to speak of a well-closed subcategory of $\Dscr$.

We already know that $\Escr$ is biclosed in $\Cscr$, so we only have to show
that $R^1\prod (E_\alpha)_\alpha=0$ for a family of injectives in $\Escr$. Let
\[
0\r E_\alpha \r F_{\alpha 0} \r F_{\alpha 1} \r F_{\alpha 2}
\]
be (truncated) injective resolutions of $E_\alpha$ in $\Dscr$. Furthermore let
\[
0\r F_{\alpha i}\r G_{\alpha i 0}\r G_{\alpha i 1} \r G_{\alpha i 2}
\]
be (truncated) Cartan Eilenberg resolutions for the complexes $F_\alpha^\cdot$.

Taking products this yields a  diagram of complexes
\[
\begin{CD}
@. @. 0\\
@. @. @VVV\\
0 @>>> \prod_\alpha E_\alpha @>>> \prod_\alpha F_\alpha^\cdot\\
@. @. @VVV\\
@. @. \prod_\alpha G_\alpha^\cdot\\
\end{CD}
\]
By hypotheses the columns of this diagram are exact and so is the first row. Then
it easily follows that
\[
0\r \prod_\alpha E_\alpha \r \prod_\alpha G_{\alpha 00} \r
\prod_\alpha G_{\alpha 10}\oplus \prod_\alpha G_{\alpha 01} \r
\prod_\alpha G_{\alpha 20}\oplus \prod_\alpha G_{\alpha 11} \oplus
\prod_\alpha G_{\alpha 0 2}
 \]
is exact. Since this is a product of truncated injective resolutions of
$E_\alpha$ in $\Cscr$ we are done.
\end{proof}
We don't know if well-closedness is compatible with the Gabriel product and
hence we don't know if the product of ideals is an ideal. In order to deal with
this problem in the sequel we introduce one more technical notion.
\begin{definitions}
\label{ref:3.5.9a}
  Assume that $\Dscr$ is complete and has an injective cogenerator. A
  biclosed subcategory $\Escr\subset \Dscr$ is \emph{very well-closed}
  if for all families $(M_j)_{j\in J}$ of objects in $\Escr$ one has
  $R^1\prod (M_j)_j=0$.
\end{definitions}
\begin{propositions} 
\label{ref:3.5.10a}
Assume that $\Cscr$ is complete and has an injective
cogenerator. Assume that $\Dscr$ is  well-closed in $\Cscr$ and $\Escr$ is
very  well-closed in $\Dscr$. Then $\Escr$ is very well-closed in $\Cscr$.
\end{propositions}
\begin{proof}
This is proved similarly as Prop. \ref{ref:3.5.8a},
\end{proof}
\begin{corollarys}
\label{ref:3.5.11a}
 Assume that $\Dscr$ is complete and has an injective
cogenerator. Then a biclosed subcategory $\Escr\subset \Dscr$ is very
well-closed if and only if
\begin{enumerate}
\item $\Escr$ is well-closed.
\item $\Escr$ has exact direct products.
\end{enumerate}
\end{corollarys}
\begin{proof} One direction is clear. To prove the other direction we use
  Proposition \ref{ref:3.5.10a} with $\Escr=\Dscr$.
\end{proof}
\begin{propositions}
\label{ref:3.5.12a}
Assume that $\Dscr$ is complete and has an injective cogenerator. Let
$\Escr_1,\Escr_2\subset \Dscr$ be very well-closed subcategories. Then
$\Escr_1\cdot \Escr_2$ is very well-closed.
\end{propositions}
\begin{proof}
Let $(M_j)_j$ be a family of objects in $\Escr_1\cdot\Escr_2$. We have exact
sequences
\[
0\r M_{j2}\r M_j \r M_{j1} \r 0
\]
with $M_{ji}\in \Escr_i$. From the long exact sequence for $(R^i\prod)_i$ and the
hypotheses we deduce that $R^1\prod (M_j)_j=0$.
\end{proof}
If $\Lscr\in\BIMOD(\Dscr-\Dscr)$ is such that $\HHom_\Dscr(\Lscr,-)$
is an equivalence of categories then we
call $\Lscr$ an \emph{invertible} bimodule. Obviously in that case $\Lscr\in
\Bimod(\Dscr-\Dscr)$, and there exist $\Lscr^{-1}\in \Bimod(\Dscr-\Dscr)$ such
that $\Lscr\otimes_\Dscr\Lscr^{-1}\cong \Lscr^{-1} \otimes_\Dscr \Lscr\cong
\Id_\Dscr$.

Assume that $\Lscr$ is invertible. We will call a weak ideal
 in $\Lscr$ a subobject $I$ of $\Lscr$ in $\BIMOD(\Dscr-\Dscr)$. If $I$
 actually lies in $\Bimod(\Dscr-\Dscr)$ then we call $I$ an ideal in $\Lscr$.
 Clearly $I\mapsto I\otimes_\Dscr\Lscr^{-1}$ and $I\mapsto \Lscr^{-1}
 \otimes_\Dscr I$ induce bijections between (weak) ideals in $\Lscr$ and (weak)
 ideals in $\Id_\Dscr$.

If $\Lscr,\Mscr$ are invertible $\Dscr-\Dscr$-bimodules and $I\subset \Lscr$,
$J\subset \Mscr$ are ideals then we define $IJ$ as the image of $I\otimes_\Dscr
J$ in $\Lscr\otimes_\Dscr\Mscr$.
\begin{definitions}
\label{ref:3.5.13a}
Let $I\subset \Lscr$ be a weak ideal in an invertible bimodule $\Lscr$. Then
the \emph{Rees algebra} $\Dscr(I)$ is the graded weak algebra given by
\[
\Id_\Dscr\oplus I\oplus I^2\oplus\cdots
\]
(with obvious multiplication).
\end{definitions}
Clearly $\Dscr(I)\in \ALG(\Dscr)$. However if $\Mod(\Id_\Dscr/I)$ is
very well-closed then by Prop. \ref{ref:3.5.4a} and \ref{ref:3.5.12a}
$\Dscr(I)$ lies in $\Alg(\Dscr)$. It would be useful if we could find
weaker conditions under which the Rees algebra of an ideal lies in
$\Alg(\Dscr)$.

We close this section by introducing some terminology which we will
use later. 

If $I\subset \Id_\Dscr$ is a weak ideal defining a (closed)
subcategory $\Escr$ of $\Dscr$ then
$I/I^2\in\BIMOD(\Escr-\Escr)$. Following \cite{rosenberg} we call
$I/I^2$ the \emph{conormal bundle} of $\Escr$ in $\Dscr$.

\subsection{Quasi-schemes}
\label{ref:3.6b}
For us a quasi-scheme $X$ will be a Grothendieck category which we
denote by $\Qch(X)$. 
A morphism $\alpha:Y\r X$ of quasi-schemes
will be an additive functor $\alpha^\ast:\Qch(X)\r \Qch(Y)$ commuting
with colimits.
 The quasi-schemes form a category which we denote by
$\Qsch$. Actually we will consider $\Qsch$ as a 2-category whose
2-cells correspond to natural isomorphisms (see Appendix \ref{ref:Aa}).

Commutative diagrams in a 2-category are
usually only assumed to be commutative up to an explicit natural isomorphism
(see Appendix \ref{ref:Aa}). Such diagrams are sometimes called
pseudo-commutative diagrams, but we will call them simply ``commutative
diagrams''. Likewise we will speak of ``functors'' when we actually mean
pseudo-functors (again see Appendix \ref{ref:Aa}).

If $\alpha:Y\r X$ is a morphism of quasi-schemes then 
 it follows from Theorem \ref{ref:2.1b}
that the adjoint to $\alpha^\ast$ exists.
This adjoint is unique up to unique isomorphism and we will denote it
by $\alpha_\ast$. The assignment $\alpha\mapsto \alpha^*$ is obviously
functorial since formally $\alpha=\alpha^*$!  However as is explained
in Appendix \ref{ref:Aa} the assignment $\alpha\mapsto
\alpha_\ast$ is also functorial if we work in the setting of
2-categories.

Denote by ``$\mathrm{Sch}$'' the category of quasi-compact, quasi-separated
schemes. If $X\in \mathrm{Sch}$ then the category of quasi-coherent sheaves on
$X$ is a Grothendieck category \cite{thomasson}.  In that case we define
$\Qch(X)$ as the category of
quasi-coherent sheaves on $X$.  Rosenberg in \cite{rosenberg1} has
proved a reconstruction theorem which allows one to recover $X$ from
$\Qch(X)$ (generalizing work of Gabriel in the noetherian case).
Furthermore he has also shown that the functor
\[
\mathrm{Sch}/\Spec \ZZ\r \QSch/\Spec \ZZ
\]
which sends a scheme to its associated quasi-scheme is fully faithful
in the sense of $2$-categories.

Now let $X$ be a quasi-schemes.  We write
$\Alg(X)=\Alg(\Qch(X))$, $\Bimod(X)=\Bimod(\Qch(X))$, etc\dots. Below by an
\emph{algebra on $X$} we will mean an object of $\Alg(X)$ unless
otherwise specified. Likewise $\Ascr$-$\Bscr$-bimodules will in
general be objects of $\Bimod(\Ascr-\Bscr)$. By $o_X$ we denote the
identity functor on $\Qch(X)$ considered as an object of $\Alg(X)$.
Obviously $\Mod(o_X)=\Qch(X)$, $\Bimod(o_X)=\Bimod(X)$.

If $\Ascr\in\Alg(X)$ then we denote by $\Spec\Ascr$ the object
$\Qsch/X$ given by the pair $(\Mod(\Ascr),-\otimes_{o_X}\Ascr)$.

If $X$ is a quasi-scheme then we define
$\Gralg(X)=\Gralg(\Qch(X))$ and also $\GRALG(X)=\GRALG(\Qch(X))$. 
 Below by a \emph{graded algebra on $X$} we will mean an object
of $\Gralg(X)$ unless otherwise specified. Likewise graded
$\Ascr$-$\Bscr$-bimodules will in general be objects of $\Bigr(\Ascr-\Bscr)$.

Related to the notion of relative categories is the notion of
\emph{enriched quasi-schemes}.  An enriched quasi-scheme will be a
pair $(X,\Oscr_X)$ where $X$ is a quasi-scheme and $\Oscr_X$ is an
arbitrary object of $\Qch(X)$. A morphism $(Y,\Oscr_Y)\r (X,\Oscr_X)$
between enriched quasi-schemes is a pair $(\alpha,s)$ where
$\alpha:Y\r X$ is a morphism between quasi-schemes and $s$ is an
isomorphism $s:\alpha^\ast(\Oscr_X)\r \Oscr_Y$. Note that if
$(Y,\alpha)\in \Qscr/X$ and $X$ is an enriched quasi-scheme then $Y$
becomes canonically an enriched quasi-scheme if we put
$\Oscr_Y=\alpha^\ast\Oscr_X$.

The prototype of an
enriched quasi-scheme is $\Spec R=(\Mod(R),R_R)$ for a ring $R$. The
following lemma will be used tacitly throughout the paper.
\begin{lemmas} Let $R$ be a ring. Then the category $\Qsch/\Spec R$ is
  equivalent (as a two-category) to the category of $R$-linear
  enriched quasi-schemes.  The equivalence is given by sending
  $(Y,\alpha)$ to $(Y,\alpha^\ast(R))$.
\end{lemmas}

If $(X,\Oscr_X)$ is an enriched quasi-scheme and $\Uscr\in \Qch(X)$
then we put $\Gamma(X,\Uscr)=\Hom_{o_X}(\Oscr_X,\Uscr)$. We say that
$\Uscr$ is generated by global sections if $\Uscr$ is a quotient of
some $\Oscr_X^{\oplus I}$.

If $\Mscr$ is a bimodule on $X$ then we put
$\Mscr_{o_X}=\Oscr_X\otimes_{o_X} \Mscr$. We think of $\Mscr_{o_X}$ as
the ``right structure'' of $\Mscr$. Care should be taken however
since $\Mscr\mapsto \Mscr_{o_X}$ is apriori not an exact
functor. This will not be a problem in our applications.

If 
 $\Mscr$ is a bimodule on $X$ then we define
$\Gamma(X,\Mscr)=\Gamma(X,\Mscr_{o_X})$. 
Again one should be careful since $\Gamma(X,-)$ is in general not left
exact. 

It is easy to see that if
$\Ascr$ is an algebra on $X$ then $A=\Gamma(X,\Ascr)$ will be a ring.
If $\Mscr$ is an $\Ascr$-module then $\Gamma(X,\Mscr)$ will be an
$A$-module.

A quasi-scheme $X$ will be called  noetherian 
if $\Qch(X)$ is locally noetherian. For an enriched quasi-scheme, we
also require  that $\Oscr_X$ is noetherian. A morphism
$\alpha:X\r Y$ will be called noetherian 
if $\alpha^\ast$ is noetherian.

An algebra $\Ascr$ on $X$ is said to be noetherian if the
functor $-\otimes_{o_X}\Ascr$ preserves noetherian objects. Clearly if
$X$ is noetherian and $\Ascr$ is noetherian then so is $\Spec\Ascr$.

A this point we will introduce a convention that will be in force
throughout this paper. 
\begin{notation} If $\Cscr=\text{Xyz}\cdots(\cdots)$ is a category
  then $\text{xyz}\cdots(\cdots)$ stands for the full subcategory of
  $\Cscr$ whose objects are the noetherian objects in $\Cscr$.
\end{notation}

\subsection{Divisors}
\label{ref:3.7b}

We will say that a  map $\alpha:Y\r X$ is a biclosed embedding if
$\alpha_\ast$ embeds $\Qch(Y)$ in $\Qch(X)$ as a biclosed subcategory.

Assume that $X,Y$ are quasi-schemes where $Y$ is embedded in
$X$ by a biclosed embedding. For simplicity we also assume that $X$ is
noetherian, although that is not strictly necessary. We denote 
\[
o_X(-Y)=\ker(o_X\r o_{Y})
\]
We say that $Y$ is a divisor in $X$ if
$o_X(-Y)$ is invertible. In the rest of this section we will assume
that $Y$ is a divisor in $X$ and we denote the inclusion mapping by
$i$.  If $\Mscr\in\Qch(X)$ then we write $\Mscr_Y$ for
$\Mscr\otimes_{o_X}o_Y$ and $\Mscr(nY)$ for $\Mscr\otimes_{o_X}
o_X(nY)$, where $o_X(nY)=o_X(Y)^{\otimes n}$.

We denote the inclusion $o_X(-Y)\r o_X$ by $t$ and we do the same
with the induced maps $\Mscr(-Y)\r \Mscr$ for $\Mscr\in \Qch(X)$.  The
normal bundle of $Y$ in $X$ is defined by $\Nscr_{Y/X}=o_X(Y)/o_X$. 

 For use below we define a few categories.
\begin{align*}
\tors_Y(X)&=\{\Mscr\in\coh(X)\mid \text{There exist $n$ such that the
   map $t^n:\Mscr(-nY)\r \Mscr$ is zero}\}\\
\iso_Y(X)&=\{\Mscr\in \coh(X)\mid \text{The  map $t:\Mscr(-1)\r
  \Mscr$ is an isomorphism}\}
\end{align*}
$\Tors_Y(X)$ and $\Iso_Y(X)$ will be the closures of $\tors_Y(X)$ and
$\iso_Y(X)$ under direct unions.
\begin{lemmas}
$\Tors_Y(X)$ and $\Iso_Y(X)$ are localizing subcategories of
$\Qch(X)$.
\end{lemmas}
\begin{proof} Since $\Qch(X)$ is a locally noetherian category, it is
  easy to see that it is sufficient to  show that $\tors_Y(X)$ and
  $\iso_Y(X)$ are Serre subcategories in $\coh(X)$. For $\tors_Y(X)$
  this is clear, so we concentrate on $\iso_Y(X)$. It is clearly
  sufficient to show that $\iso_Y(X)$ is closed under taking
  subobjects. Let $\Mscr\in \iso_Y(X)$ and let $\Nscr\subset\Mscr$. We
  have isomorphisms $t^n:\Mscr\r \Mscr(nY)$ and these yield
 an ascending chain of submodules in $\Mscr$ given by
  $t^{-n}(\Nscr(nY))$. Since $\Mscr$ is noetherian this chain must stop.
  From this we easily obtain that $\Nscr\in\iso_Y(X)$.
\end{proof}
We will say that the objects in $\Tors_Y(X)$ are \emph{supported} on
$Y$. If $\Mscr(-Y)\r\Mscr$ is injective then we will say that $\Mscr$
is $Y$-\emph{torsion free}.

 We will use the following result.
\begin{propositions} 
\label{ref:3.7.2a}
Let $\Mscr\in \coh(X)$. Then the filtration
\[
\cdots \subset t^n(\Mscr(-nY))\subset t^{n-1}(\Mscr(-(n-1)Y))\subset \cdots\subset
t(\Mscr(-Y))\subset
\Mscr
\]
satisfies the Artin-Rees condition.
\end{propositions}
\begin{proof}
Left to the reader.
\end{proof}
\begin{corollarys}
\label{ref:3.7.3a}
Let $\Mscr\in\coh(X)$. Then $\Mscr$ contains
an $Y$-torsion free submodule $\Nscr$, such that $\Mscr/\Nscr$ is
supported on $Y$.
\end{corollarys}
\begin{proof}
  Let $\Tscr$ be the maximal submodule of $\Mscr$ supported on $Y$.
  Since $\Tscr$ is noetherian we will have $t^n(\Tscr(-nY))=0$ for
  some $n$. By
  Proposition \ref{ref:3.7.2a} there will be some $m$ such that
  $t^m(\Mscr(-mY))\cap \Tscr\subset t^n(\Tscr(-nY))=0$.  Thus
  $\Nscr=t^m(\Mscr(-mY))$ is $Y$ torsion free and $\Mscr/\Nscr$ is
  supported on $Y$.
\end{proof}
If $i:(Y,\Oscr_Y)\r (X,\Oscr_X)$ is a map of enriched quasi-schemes
then we say that $i$ makes $Y$ into a divisor in $X$ if the underlying
map $Y\r X$ makes $Y$ into a divisor in $X$ in the sense of ordinary
quasi-schemes and if in addition the induced map $\Oscr_X(-Y)\r
\Oscr_X$ is injective.

\subsection{Proj}
\label{ref:3.8b}
Below $X$ will be a  noetherian quasi-scheme and
  $\Ascr=\oplus_n\Ascr_n\in \Gralg(X)$ will be noetherian.
The definition of $\Proj \Ascr$ is entirely similar to the
ring case \cite{AZ}. As before $\Gr(\Ascr)$ is the
category of $\ZZ$-graded $\Ascr$-modules, $\Tors(\Ascr)$ is the full
subcategory of $\Gr(\Ascr)$ consisting of graded modules that are unions of
right bounded modules and
\[
\Qgr(\Ascr)=\Gr(\Ascr)/\Tors(\Ascr)
\]
We use similar notations as in \cite{AZ}. Thus $\pi:\Gr(\Ascr)\r
\Qgr(\Ascr)$ is the quotient map. $\omega:\Qgr(\Ascr)\r \Gr(\Ascr)$ is
the right adjoint to $\pi$  and $\tau:\Gr(\Ascr)\r \Tors(\Ascr)$ is the
functor which associates to every object its maximal torsion
subobject. We put
$\tilde{(-)}=\omega\pi$.

Let $\Pqsch/X$ be the 2-category of triples $(Y,\alpha,s)$ where $Y$
is a quasi-scheme, $\alpha:Y\r X$  a morphism  and $s$ an
autoequivalence of $\Qch(Y)$.  Morphism between triples
$(Y,\alpha,s)\r (Z,\beta,t)$ are given by a morphism $\gamma:Y\r Z$, a
natural isomorphism $\mu:\gamma^\ast\circ\beta^\ast\r\alpha^\ast$
and a natural isomorphism $\psi:s\circ\gamma^\ast\r \gamma^\ast\circ
t$. We leave it to the reader to define natural isomorphisms between
such triples. 
We now define $\Proj \Ascr$ as the object of $\Pqsch/X$ given by the
triple $(\Qgr(\Ascr),\pi(-\otimes_\ox\Ascr),s)$ where $s$ is the shift
functor on $\Qgr(\Ascr)$ obtained from the canonical shift on
$\Gr(\Ascr)$. We will also denote by $\Proj\Ascr$ the object in $\Qsch/X$
obtained by forgetting the shift.

If we denote by $\alpha$ the structure map $\Proj\Ascr\r X$ then by definition,
$\alpha^\ast$ is given by $\pi(-\otimes_\ox\Ascr)$. Since $\omega$ is
the right adjoint to $\pi$ we deduce that $\alpha_\ast$ is given by
$\omega(-)_0$. In the sequel it will be convenient to use
$\underline{\alpha}_\ast$ as a synonym for $\omega$. Thus
\[
\underline{\alpha}_\ast(\Mscr)_\ox=\oplus_n \alpha_\ast(\Mscr(n))
\]
Below we will generalize some of the
results of \cite{AZ} to our situation since we will need
them. Usually we can simply copy the proofs in \cite{AZ}.
\begin{lemmas}
\label{ref:3.8.1a}
Let $\Ascr$, $\Bscr$ be noetherian graded algebras on $X$, $\Nscr$ a graded
$\Ascr$-module, $\Mscr_1$, $\Mscr_2$ graded $\Ascr-\Bscr$ bimodules.
\begin{enumerate}
\item Assume that $\Nscr$ is right bounded and $(\Mscr_1)_{o_X}$ is a quotient
  of $\Ascr\otimes_\ox\Mscr'_1$ where $\Mscr'_1$ is a right bounded
graded $\ox-\ox$-module. Then $\Nscr\otimes_\Ascr\Mscr_1$ is
right bounded.
\item Assume that $\Mscr_1$ is right bounded. Then
  $\Nscr\otimes_\Ascr\Mscr_1$ is torsion.
\item Assume that $\phi:\Mscr_1\r\Mscr_2$ is a morphism of graded
  $\Ascr-\Bscr$ bimodule which is an isomorphism in high degree. Then 
$\ker,\coker(\Nscr\otimes_\Ascr\Mscr_1\r\Nscr\otimes_\Ascr\Mscr_2)$ are
torsion.
\end{enumerate}
\end{lemmas}
\begin{proof} 
\begin{enumerate}
\item This is trivial since $\Nscr\otimes_\Ascr\Mscr_1$ is a quotient
  as graded
$\ox$-bimodules of
\[
\Nscr\otimes_\Ascr(\Ascr\otimes_\ox\Mscr'_1)=\Nscr\otimes_\ox\Mscr'_1
\]
\item
This is a special case of (3).
\item
Now write $\Nscr$ as a quotient of modules of the form
$\Pscr_1\r\Pscr_0$ where $\Pscr_{0,1}$ are direct sums of shifts of
modules of the form $\Nscr'\otimes_\ox\Ascr$, $\Nscr'\in\Qch(X)$ (located
in degree zero). This is possible by Proposition \ref{ref:3.1.15a}(8). 

If $\Nscr=\Nscr'\otimes_\ox\Ascr$ then the lemma is true because 
\[
(\Nscr'\otimes_\ox\Ascr)\otimes_\Ascr\Mscr=\Nscr'\otimes_\ox
{}_\ox\!\Mscr
\]
for all $\Mscr\in\Bimod(\Ascr-\Bscr)$. Hence the map
$\Nscr\otimes_\Ascr\Mscr_1\r\Nscr\otimes_\Ascr\Mscr_2$ is an
isomorphism in high degree. The general case follows from the fact
that $-\otimes_\Ascr\Mscr_{1,2}$ is compatible with colimits.
\qed\end{enumerate} \def\qed{}\end{proof} From this lemma one obtains
the following corollary (cfr.\ \cite[Prop.\ 2.5]{AZ}).
\begin{corollarys}
\label{ref:3.8.2a}
Let $\phi:\Ascr\r\Bscr$ be a morphism of noetherian graded algebras on $X$ such
that $\phi$ is an isomorphism in high degree. Then the functors 
\begin{equation}
\label{ref:3.21a}
\Gr(\Ascr)\r\Gr(\Bscr):\Nscr\r \Nscr\otimes_\Ascr\Bscr
\end{equation}
and
\begin{equation}
\label{ref:3.22a}
\Gr(\Bscr)\r\Gr(\Ascr):\Mscr\r\Mscr_\Ascr
\end{equation}
factor to give inverse equivalences between $\Qgr(\Ascr)$ and
$\Qgr(\Bscr)$. Furthermore  $\Proj(\Ascr)$ and
$\Proj(\Bscr)$ are equivalent.
\end{corollarys}
\begin{proof}
Assume that $\phi_n$ is an isomorphism for $n\ge n_0$. Then $\Bscr$ is
a quotient of $\Ascr\oplus\Ascr\otimes(\oplus_{n<n_0}\Bscr_n(-n))$. 

Assume that $\Mscr$ is a $\Bscr$-module, torsion as
$\Ascr$-module. Then $\Mscr_\Ascr$ is a quotient of $\oplus_{i\in
  I}\Mscr_i$ where the $\Mscr_i$ are right bounded
$\Ascr$-modules. Hence $\Mscr$ is a quotient of $\oplus_{i\in
  I}\Mscr_i\otimes_\Ascr\Bscr$. By lemma \ref{ref:3.8.1a}(1) all the
$\Mscr_i\otimes_\Ascr\Bscr$ are right bounded. Thus $\Mscr$ is also
torsion as $\Bscr$-module.

Thus $\Bscr$-torsion is equivalent to $\Ascr$-torsion and so we will
simply speak of torsion. \eqref{ref:3.22a} obviously preserves isomorphism
mod torsion, so we concentrate on \eqref{ref:3.21a}.

Assume that $\ker,\coker(\Mscr\xrightarrow{\theta}{\Nscr})$ are torsion for
some $\Mscr,\Nscr\in\Gr(\Ascr)$. Then we have the following
commutative diagram.
\[
\begin{CD}
\Mscr @>>> \Mscr\otimes_\Ascr\Bscr\\
@V\theta VV @V\theta \otimes 1 VV\\
\Nscr @>>> \Nscr\otimes_\Ascr\Bscr
\end{CD}
\]
According to lemma \ref{ref:3.8.1a}(3) the horizontal maps are isomorphisms modulo
torsion. Since by hypotheses this is also true for $\theta$, we obtain
that $\theta\otimes 1$ is an isomorphism modulo torsion as well.

The reader now easily verifies that \eqref{ref:3.21a}  and \eqref{ref:3.22a} are
mutual inverses modulo torsion.

If $\Ascr$, $\Bscr$ are as in the previous proposition
then \eqref{ref:3.21a} is clearly
compatible with the structure maps and the shift functors on
$\Proj\Ascr$ and $\Proj\Bscr$
whence
$\Proj\Ascr$ and $\Proj\Bscr$ are equivalent. 
\end{proof}

It  follows from the previous proposition that one may restrict
oneself to $\NN$-graded algebras since $\Proj\Ascr=\Proj\Ascr_{\ge
  0}$. In the sequel all graded algebras will be noetherian and
$\NN$-graded, unless otherwise specified. Note that just as in the
ring case $\Ascr$ noetherian implies $\Ascr_{\ge 0}$ noetherian.

This allows us to use the following technically useful result which we
state for further reference.
\begin{lemmas} \label{ref:3.8.3a}
Assume that $\Ascr$ is noetherian and
  $\NN$-graded. Then $\Tors(A)$ is closed under essential extensions
  and hence under injective hulls. Thus $\tau$ is \emph{stable} in the sense
  of \cite{stenstrom}
\end{lemmas}

Assume that $\Ascr$ is $\NN$-graded. Then if $\Mscr,\Nscr\in \Gr(\Ascr)$ 
\begin{align*}
\Hom_{\Qgr(\Ascr)}(\pi\Nscr,\pi\Mscr)&=
\dirlim_{\Nscr/\Nscr'\text{ torsion}}
\Hom_{\Gr(\Ascr)}(\Nscr',\Mscr/\tau(\Mscr))\\
&=\dirlim_{\Nscr/\Nscr'\text{ torsion}} 
\Hom_{\Gr(\Ascr)}(\Nscr',\Mscr)
\end{align*}
where the last equality follows from the fact that $\tau$ is stable
(see  \cite{stenstrom}).

If $\Nscr\in\gr(\Ascr)$ and if $\Nscr/\Nscr'$ is torsion then it is in fact
right bounded, so $\Nscr_{\ge n}\subset \Nscr'$ for $n\ge 0$. So in that case
\begin{align*}
\Hom_{\Qgr(\Ascr)}(\pi\Nscr,\pi\Mscr)&=
\dirlim_n \Hom_{\Gr(\Ascr)}(\Nscr_{\ge n},\Mscr)
&=\dirlim\Hom_{\Gr(\Ascr)}(\Nscr_{\ge n},\Mscr_{\ge n})
\end{align*}
\begin{propositions}
\label{ref:3.8.4a}
  Let $\Ascr\in \Gralg(X)$ be noetherian and $\NN$-graded.  Then for
  $\Mscr\in\Gr(\Ascr)$
\begin{align}
\tau(\Mscr)&=\dirlim\underline{\HHom}_\Ascr(\Ascr/\Ascr_{\ge
  n},\Mscr)\label{ref:3.23a}\\
\tilde{\Mscr}&=\dirlim\underline{\HHom}_\Ascr(\Ascr_{\ge n},\Mscr)
\label{ref:3.24a}
\end{align}
\end{propositions}
\begin{proof}
By localization theory
\[
\tau(\Mscr)=\ker(\Mscr\r \tilde{\Mscr})
\]
Now from the exact sequence of $\Ascr$-bimodules 
\[
0\r\Ascr_{\ge n}\r \Ascr\r \Ascr/\Ascr_{\ge n}\r 0
\]
we obtain a
left exact sequence
\[
0\r\underline{\HHom}_\Ascr(\Ascr/\Ascr_{\ge n},\Mscr)\r
\underline{\HHom}_\Ascr(\Ascr,\Mscr)\r
\underline{\HHom}_\Ascr(\Ascr_{\ge n},\Mscr)
\]
Using the fact that $\underline{\HHom}_\Ascr(\Ascr,\Mscr)=\Mscr$ and
exactness of direct limits, it follows that it suffices to prove
\eqref{ref:3.24a}.

Recall that by definition $\tilde{\Mscr}=\omega\pi\Mscr$. So to prove
\eqref{ref:3.24a} it is sufficient for every $\Nscr\in\Gr(\Ascr)$ to
construct a natural isomorphism between \begin{equation}
\label{ref:3.25a}
\Hom_{\Gr(\Ascr)}(\Nscr,\omega\pi\Mscr)
\end{equation}
and
\begin{equation}
\label{ref:3.26a}
\Hom_{\Gr(\Ascr)}(\Nscr,\dirlim\underline{\HHom}_\Ascr(\Ascr_{\ge
  n},\Mscr))
\end{equation}
and since $\Gr(\Ascr)$ is locally noetherian, it suffices to do this
in fact for $\Nscr$ noetherian. Now
\begin{equation}
\label{ref:3.27a}
\eqref{ref:3.25a}=\dirlim_n \Hom_{\Gr(\Ascr)}(\Nscr_{\ge n},\Mscr)
\end{equation}
On the other hand, again because $\Nscr$ is noetherian.
\begin{equation}
\label{ref:3.28a}
\eqref{ref:3.26a}=\dirlim\Hom_{\Gr(\Ascr)}(\Nscr\otimes_\Ascr\Ascr_{\ge
  n},\Mscr)
\end{equation}

We have to make \eqref{ref:3.27a} isomorphic to \eqref{ref:3.28a}. Now $\Nscr$ is
noetherian and hence left bounded, so by replacing $\Nscr$ by some
shift (and doing the same with $\Mscr$) we may assume that $\Nscr$ is
in fact $\NN$-graded. Using \eqref{ref:3.2a} and lemma
\ref{ref:3.1.11a}   it is now sufficient to show
that 
\begin{align*}
\Kscr_n=\ker (\Nscr\otimes_\Ascr\Ascr_{\ge n}\r \Nscr_{\ge n})\\
\Cscr_n=\coker(\Nscr\otimes_\Ascr\Ascr_{\ge n}\r \Nscr_{\ge n})
\end{align*}
are torsion inverse systems (see Definition \ref{ref:3.1.9a}).

Now we write $\Nscr$ as quotient of $\Pscr_1\r\Pscr_0$ where the
$\Pscr_i$ are finite direct sums of negative shifts of modules of the
form $\Nscr'\otimes_\ox\Ascr$. Working in the abelian category of
inverse systems modulo torsion it is sufficient to prove that
$\Kscr_n$ and $\Cscr_n$ are torsion inverse systems in the case 
$\Nscr=\Nscr'\otimes_\ox\Ascr(-m)$. In that case
\begin{align*}
\Kscr_n&=\ker (\Nscr'\otimes_\ox\Ascr_{\ge n} \r
\Nscr'\otimes_\ox\Ascr_{\ge n-m})(-m)\\
\Cscr_n&=\coker (\Nscr'\otimes_\ox\Ascr_{\ge n} \r
\Nscr'\otimes_\ox\Ascr_{\ge n-m})(-m)
\end{align*}
It is now clear that in this case actually $\Kscr_n=0$ and $\Cscr_n$
lives in degrees $[n,n+m-1]$. So $(\Cscr_n)_n$ is a torsion inverse system
as well.
\end{proof}
\subsection[Condition  ``$\chi$'' and cohomological   dimension]{Condition
{\mathversion{bold} ``$\chi$''} and cohomological   dimension} 
\label{ref:3.9b}
Below $X$ will be a noetherian quasi-scheme and $\Ascr$,
$\Bscr$ will be $\NN$-graded noetherian objects in $\Gralg(X)$.
$\alpha:\Proj \Ascr\r X$, $\beta:\Proj \Bscr\r X$ will be the
structure maps. $\pi$, $\omega$, $\tau$ refer to $\Ascr$. For $\Bscr$
we use notations such as $\pi_\Bscr$, $\omega_\Bscr$, $\tau_\Bscr$.

The following is clear.
\begin{lemmas}
\label{ref:3.9.1a}
If $\Mscr\in\Gr(\Ascr)$ then $R^i\tau\Mscr\in\Tors(\Ascr)$.
\end{lemmas}
In order to develop a non-commutative version of projective
geometry, the following definition was proposed in \cite{AZ}.
\begin{definitions}  $\Ascr$ satisfies
  $\chi$ if for every $\Mscr\in \gr(\Ascr)$ one has that 
$(R^i\tau\Mscr)_{\ge 0}\in \tors(\Ascr)$.
\end{definitions}

If $\Mscr\in \Gr(\Ascr)$ then we have a triangle 
\[
R\tau\Mscr\r \Mscr\r R\omega(\pi\Mscr)\xrightarrow{[1]}
\]
from which one deduces the
exact sequence
\begin{equation}
\label{ref:3.29a}
0\r \tau\Mscr\r\Mscr\r \tilde{\Mscr}\r R^1\tau \Mscr\r 0
\end{equation}
and isomorphisms
\begin{equation}
\label{ref:3.30a}
R^{i+1}\tau\Mscr=R^{i}\omega(\pi\Mscr)),\qquad i\ge 1
\end{equation}
So we obtain~:
\begin{lemmas}
\label{ref:3.9.3a}
$\Ascr$ satisfies $\chi$ if and only if for all $\Mscr\in\gr \Ascr$
one has
\begin{enumerate}
\item $R^i\omega (\pi\Mscr)_{\ge 0}\in \tors(\Ascr)$, $i\ge 0$.
\item $(\tilde{\Mscr})_{\ge 0}/\Mscr\in \tors(\Ascr)$
\end{enumerate}
\end{lemmas}
One should view  condition $\chi$  as a  kind of ampleness condition. 
This becomes clearer if one makes the following definitions.
\begin{definitions} 
\label{ref:3.9.4a}
A map $\alpha:Y\r X$ of noetherian
quasi-schemes is \emph{proper} if $R^i\alpha_\ast$ sends $\coh(Y)$ to
$\coh(X)$ for all $i$.
\end{definitions}
\begin{definitions} Let $(Y,\alpha,s)$ be an object in $\Pqsch/X$.
Assume that $Y$ is noetherian.
Then we say that $s$ is \emph{relatively ample} if the following two
conditions hold for $\Mscr\in\coh(Y)$.
\begin{enumerate}
\item For $i>0$ one has
\[
R^i\alpha_\ast(s^n\Mscr)=0
\]
for $n\gg 0$.
\item The adjoint map 
\[
\alpha^\ast \alpha_\ast s^n\Mscr\r s^n\Mscr
\]
is surjective for $n\gg 0$.
\end{enumerate}
\begin{remarks}
This definition of ampleness is much more restrictive than the one
used in \cite{AZ}.
\end{remarks}
\end{definitions}
 We say that $\Ascr$ is generated in
degree one if the multiplication map $\Ascr_m\otimes_{o_X}\Ascr_n\r
\Ascr_{m+n}$ is surjective for all $m,n\ge 0$.

The following is a consequence of lemma \ref{ref:3.9.3a}.
\begin{propositions}
\label{ref:3.9.7a} Let $(Y,\alpha,s)=\Proj \Ascr$. If $\Ascr$ satisfies
$\chi$ then $\alpha$ is proper and if $\Ascr$ is in addition generated
in degree one  then $s$ is relatively ample.
\end{propositions}

In the following we will often need the following technical condition
on an object $\Mscr\in\BIGR(\Ascr-o_X)$~:
\begin{itemize}
\item[(fin)] $\Mscr$ has a presentation $\Pscr_1\r \Pscr_0$ where the
  $\Pscr_i$ are finite sums of shifts of $\Ascr-o_X$-bimodules of the
  form $\Ascr\otimes_{o_X}\Pscr'$ where the $\Pscr'$ are coherent
  objects (see \ref{ref:3.1a}) in $\BIMOD(o_X-o_X)$ (located in degree
  zero).
 \end{itemize}
We should think of (fin) as a kind of ``finite presentation''
condition. (fin) is useful because of the following lemma
\begin{lemmas} 
\label{ref:3.9.8a}
Assume that $\Mscr\in\BIGR(\Ascr-o_X)$ satisfies
  (fin) and $\Nscr\in\BIGR(\Bscr-\Ascr)$ (or
  $\Nscr\in\Gr(\Ascr)$). Assume furthermore that $\Mscr$ is
  $\NN$-graded. Then the kernel and the cokernel of 
\[
\Nscr_{\ge n}\otimes_\Ascr \Mscr\r (\Nscr\otimes_\Ascr\Mscr)_{\ge n}
\]
are torsion inverse systems (cfr. Definition \ref{ref:3.1.9a}).
\end{lemmas}
Here is another obvious application.
\begin{lemmas} 
\label{ref:3.9.9a}  Let $f:\Ascr\r \Bscr$ be a morphism  in
$\Gralg(X)$ and assume that $\Bscr$ satisfies (fin). Let $\Mscr\in \Gr(\Bscr)$. Then
$\tau_{\Ascr}(\Mscr_\Ascr)=(\tau_{\Bscr}(\Mscr))_\Ascr$.
\end{lemmas}
\begin{proof}
  Using Proposition \ref{ref:3.8.4a} and adjointness this amounts to
  showing that the inverse systems $\Ascr/\Ascr_{\ge n}\otimes_\Ascr
  \Bscr$ and $\Bscr/\Bscr_{\ge n}$ are equivalent modulo torsion
  inverse systems. Since $\Ascr/\Ascr_{\ge n}\otimes_\Ascr
  \Bscr=\coker (\Ascr_{\ge n}\otimes_\Ascr\Bscr\r \Bscr)$ this follows
  directly from lemma \ref{ref:3.9.8a}.
\end{proof}

In the sequel we will need a result like \eqref{ref:3.31a}. The hypotheses
under which we can prove this are
unfortunately quite technical and almost certainly too restrictive.

\begin{propositions}
\label{ref:3.9.10a}
Let $f:\Ascr\r \Bscr$ be a morphism in $\Gralg(X)$. Assume that
$\Bscr$ satisfies (fin) and furthermore that there are graded flat
$\Ascr-o_X$-bimodules $(\Fscr_i)_i$ satisfying (fin), together with a
long exact sequence
\[
\cdots \r(\Fscr_2)_{\ge p}\r (\Fscr_1)_{\ge p}\r (\Fscr_0)_{\ge p}
\r ({}_\Ascr\Bscr_{o_X})_{\ge p}\r 0
\]
for a certain $p\in\ZZ$. Then for $\Mscr\in\Gr(\Bscr)$ one has
\begin{equation}
\label{ref:3.31a}
R^i\tau_{\Ascr}(\Mscr_\Ascr)=(R^i\tau_{\Bscr}(\Mscr))_\Ascr
\end{equation}
\end{propositions}

\begin{proof} Both sides of \eqref{ref:3.31a} are $\delta$-functors. Hence
  is suffices to show that the two sides of \eqref{ref:3.31a} are zero for
  $i>0$ and naturally isomorphic for $i=0$, when evaluated on
  injectives. 

Actually it is not necessary to take all injectives, but only a
cogenerating set.  Thus we take for our injectives all
modules of the form 
\[
F=\underline{\HHom}_\Ascr({}_\Bscr\Bscr_\Ascr,E)
\]
where $E\in \Inj(\Ascr)$. Clearly by Proposition \ref{ref:3.8.4a},
lemma \ref{ref:3.9.9a} and
adjointness
\[
\text{RHS}\eqref{ref:3.31a}=
\begin{cases}
\dirlim \underline{\HHom}_\Ascr(\Ascr/\Ascr_{\le n}\otimes_\Ascr
{}_\Ascr\Bscr_\Ascr,E)&\text {if $i=0$}\\
0&\text{otherwise}
\end{cases}
\]
On the other hand by the definition of $\underline{\HTor}$
\begin{align*}
\text{LHS}\eqref{ref:3.31a}&= \dirlim
\underline{\HExt}^i_\Ascr(\Ascr/\Ascr_{\ge n} ,F_{\Ascr})\\
&=\dirlim
\underline{\HExt}^i_\Ascr(\Ascr/\Ascr_{\ge n} ,\underline{\HHom}_\Ascr
({}_\Ascr\Bscr_\Ascr,E))\\
&=\dirlim\underline{\HHom}_\Ascr(\underline{\HTor}^{\Ascr}_i(\Ascr/\Ascr_{\le
    n},{}_\Ascr\Bscr_\Ascr),E)
\end{align*}
So the two sides of \eqref{ref:3.31a} certainly agree for $i=0$. In general
we have to show  for $i>0$ that 
\begin{equation}
\label{ref:3.32a}
\underline{\HTor}^\Ascr_i(\Ascr/\Ascr_{\ge n},{}_\Ascr\Bscr_\Ascr)
\end{equation}
is a torsion inverse system.

From fact that  the $\Fscr_i$ are flat we can deduce that,
up to inverse systems of the form
\begin{equation}
\label{ref:3.33a}
\underline{\HTor}^\Ascr_i(\Ascr/\Ascr_{\ge n},\Tscr)
\end{equation}
with $\Tscr_{o_X}$ coherent, \eqref{ref:3.32a} is given by
the middle homology of
\begin{equation}
\label{ref:3.34a}
\Ascr/\Ascr_{\ge n}\otimes_\Ascr (\Fscr_{i+1})_{\ge p}
\r
\Ascr/\Ascr_{\ge n}\otimes_\Ascr (\Fscr_{i})_{\ge p}
\r
\Ascr/\Ascr_{\ge n}\otimes_\Ascr (\Fscr_{i-1})_{\ge p}
\end{equation}
By shifting if necessary we may assume that 
$\Fscr_{i+1}$, $\Fscr_i$, $\Fscr_{i-1}$ are $\NN$-graded.

Now assume in general that $\Fscr$ is an $\NN$-graded object in $\BIGR(\Ascr-o_X)$,
satisfying (fin). 
 Then 
\[
\Ascr/\Ascr_{\ge n}\otimes_\Ascr \Fscr_{\ge p}=\coker (\Ascr_{\ge
  n}\otimes_\Ascr \Fscr_{\ge p}\r \Fscr_{\ge p})
\]
Now up to inverse systems of the form
\begin{equation}
\label{ref:3.35a}
\underline{\HTor}_i^\Ascr(\Ascr_{\ge n},\Fscr/\Fscr_{\ge p})
\end{equation}
we have
\[
\Ascr_{\ge n}\otimes_\Ascr \Fscr_{\ge p}=\Ascr_{\ge n}\otimes_\Ascr \Fscr
\]

Then using lemma \ref{ref:3.9.8a} we see that $ \Ascr_{\ge
  n}\otimes_\Ascr \Fscr $ is, modulo torsion inverse systems, equal to
$ \Fscr_{\ge n} $.  Assembling everything we find that up to torsion
inverse systems and inverse systems of the form \eqref{ref:3.33a} and
\eqref{ref:3.35a} $\Ascr/\Ascr_{\ge n}\otimes_\Ascr \Fscr_{\ge p}$
is equal to $\Fscr_{\ge p}/\Fscr_{\ge n}$ (where we let the inverse
systems start with $n=p$).

So the middle homology of \eqref{ref:3.34a} is, up to torsion inverse
systems and inverse systems of the form \eqref{ref:3.33a} and \eqref{ref:3.35a}, equal
to the middle homology of 
\[
\Fscr_{i+1}/(\Fscr_{i+1})_{\ge n}
\r
\Fscr_{i}/(\Fscr_{i})_{\ge n}
\r
\Fscr_{i-1}/(\Fscr_{i-1})_{\ge n}
\]
Since this is in fact an exact complex we are left with showing that
\eqref{ref:3.33a} and \eqref{ref:3.35a} are torsion.

Let us for example consider \eqref{ref:3.33a}. It is easy to verify
directly from the definitions that 
\[
(``{\invlim}{}"\underline{\HTor}^\Ascr_i(\Ascr/\Ascr_{\ge
  n},\Tscr))_{o_X}=
``{\invlim}{}"\underline{\HTor}^\Ascr_i(\Ascr/\Ascr_{\ge
  n},\Tscr_{o_X})=
\]
Now as in the proof of lemma \ref{ref:3.1.17a} one verifies that
$\underline{\HTor}^\Ascr_i(\Ascr/\Ascr_{\ge n},\Tscr_{o_X})$ is
coherent as weak $\Ascr-o_X$-bimodule. Hence
$\underline{\HHom}_{o_X}(\underline{\HTor}^\Ascr_i(\Ascr/\Ascr_{\ge
  n},\Tscr_{o_X}),-)$ commutes with direct sums.  Thus it suffices to
show that 
\begin{equation}
\label{ref:3.36a}
\dirlim\underline{\HHom}_{o_X}(\underline{\HTor}^\Ascr_i(\Ascr/\Ascr_{\ge
  n},\Tscr_{o_X}),E)=0
\end{equation}
for $E$ an injective object in $\Mod(X)$, viewed as a graded object
concentrated in a single degree.  Now \eqref{ref:3.36a} becomes equal to
\begin{equation}
\label{ref:3.37a}
\dirlim\underline{\HExt}^i_{\Ascr}(\Ascr/\Ascr_{\ge
  n},\underline{\HHom}_{o_X}(\Tscr_{o_X},E))
\end{equation}
Clearly $\underline{\HHom}_{o_X}(\Tscr_{o_X},E)$ is right bounded and
hence torsion. Now the fact that $\Tors(\Ascr)$ is stable (lemma
\ref{ref:3.8.3a}) together with Proposition \ref{ref:3.8.4a} easily
implies that \eqref{ref:3.37a} is zero.
\end{proof}
Now we translate \eqref{ref:3.31a} into more geometrical language. First we
indicate when a map $f:\Ascr\r \Bscr$ defines a dual map
$\bar{f}:\Proj(\Bscr)\r \Proj(\Ascr)$
\begin{propositions}
\label{ref:3.9.11a}
Let $f:\Ascr\r \Bscr$ be a morphism in $\Gralg(X)$ such that
${}_\Ascr\Bscr_{o_X}$ satisfies (fin). Then the functors
\begin{align}
\label{ref:3.38a}
-\otimes_\Ascr\Bscr&:\Gr(\Ascr)\r \Gr(\Bscr)\\
(-)_\Ascr&:\Gr(\Bscr)\r \Gr(\Ascr) \label{ref:3.39a}
\end{align}
factor through ``$\Qgr$'' and in this way define respectively $\bar{f}^\ast$
and $\bar{f}_\ast$ for a morphism $\bar{f}:\Proj(\Bscr)\r \Proj
(\Ascr)$. In this case $\bar{f}_\ast$ is exact.
\end{propositions}
\begin{proof}
That \eqref{ref:3.39a} factors through ``$\Qgr$'' as well as the
exactness of $\bar{f}_\ast$ is clear so we
concentrate on \eqref{ref:3.38a}. We have to show that if a map $\Mscr\r
\Nscr$ in $\Mod(\Ascr)$ has torsion kernel and cokernel then the same is
true for 
\[
\Mscr\otimes_\Ascr\Bscr\r\Nscr\otimes_\Ascr\Bscr
\]
Now $\Mod(\Ascr)$ is locally noetherian and $-\otimes_\Ascr\Bscr$ is
compatible with direct limits, so we may restrict ourselves to the case where
$\Mscr$, $\Nscr$ are noetherian and hence $\Mscr_n\r\Nscr_n$ is an
isomorphism for $n\gg 0$.

Now using the fact that $\Bscr$ satisfies (fin) it suffices to show that
\begin{equation}
\label{ref:3.40a}
\Mscr\otimes_\Ascr (\Ascr\otimes_{o_X}\Bscr')
\r
\Nscr\otimes_\Ascr (\Ascr\otimes_{o_X}\Bscr')
\end{equation}
has torsion kernel and cokernel for all $\Bscr'\in \BIMOD(o_X-o_X)$.
However since $-\otimes_\Ascr(\Ascr\otimes_{o_X}\Bscr')=-\otimes_{o_X}\Bscr'$ it is
clear that \eqref{ref:3.40a} is an isomorphism in high degree. 

It remains to verify the adjointness of \eqref{ref:3.38a} and
\eqref{ref:3.39a} when defined on ``$\QGr$''. To this end we have to
construct for $\Mscr\in\Gr(\Ascr)$ and $\Nscr\in \Gr(\Bscr)$ a natural
isomorphism 
\[
\Hom_{\Qgr(\Bscr)}(\pi(\Mscr\otimes_\Ascr\Bscr),\pi\Nscr)
\cong
\Hom_{\Qgr(\Ascr)}(\pi\Mscr,\pi(\Nscr_\Ascr))
\]
and as usual it  suffices to do this for  $\Mscr$ noetherian. Then
\begin{align*}
\Hom_{\Qgr(\Bscr)}(\pi(\Mscr\otimes_\Ascr\Bscr),\pi\Nscr)
&=\dirlim \Hom_{\Gr(\Bscr)}((\Mscr\otimes_\Ascr \Bscr)_{\ge n},\Nscr)\\
&=\dirlim \Hom_{\Gr(\Bscr)}(\Mscr_{\ge n}\otimes_\Ascr \Bscr,\Nscr)\\
&=\dirlim \Hom_{\Gr(\Ascr)}(\Mscr_{\ge n},\Nscr_\Ascr)\\
&=\Hom_{\Qgr(\Ascr)}(\pi\Mscr,\pi(\Nscr_\Ascr))
\end{align*}
The second equality follows from the fact that $\Bscr$ satisfies (fin) as
$\Ascr-o_X$-bimodule and thus we can apply lemma \ref{ref:3.9.8a}
\end{proof}
Thus we obtain 
\begin{propositions}
\label{ref:3.9.12a}
Let $f:\Ascr\r\Bscr$ be a morphism in $\Gralg(X)$ and assume that $\Bscr$ is
as in Proposition \ref{ref:3.9.10a}. Let $\bar{f}:\Proj\Bscr\r\Proj\Ascr$ be
defined by $\bar{f}^\ast(\pi\Mscr)=\pi(\Mscr\otimes_\Ascr\Bscr)$ for
$\Mscr\in\Gr(\Ascr)$. Then for $\Nscr\in\Qgr(\Bscr)$
\[
R^i\beta_\ast(\Nscr)=R^i\alpha_\ast(\bar{f}_\ast\Nscr)
\]
\end{propositions}
\begin{proof} Assume $\Nscr=\pi\Nscr'$, $\Nscr'\in\Gr(\Bscr)$. 
Summing over shifts, it is sufficient to show
\[
R^i\omega_\Bscr(\pi\Nscr')_\Ascr=R^i\omega_\Ascr(\pi(\Nscr'_\Ascr))
\]
If $i\ge 1$ then by \eqref{ref:3.30a} this is equivalent to
\[
R^{i+1}\tau_\Bscr(\Nscr')_\Ascr=R^{i+1}\tau_\Ascr(\Nscr'_\Ascr)
\]
which we have shown.
The case $i=0$ follows from considering \eqref{ref:3.29a} and the 5-lemma.
\end{proof}
Now we use Proposition \ref{ref:3.9.10a} to prove the following result
\begin{propositions}
\label{ref:3.9.13a}
Assume that we have a  map $f:\Ascr\r\Bscr$ of $\NN$-graded
algebras on $X$. Assume furthermore that there is an exact sequence
\[
0\r  I\r \Ascr\r{}_\Ascr\Bscr_\Ascr\r 0
\]
in $\Bigr(\Ascr-\Ascr)$ with $I$ a graded invertible $\Ascr$-bimodule,
satisfying (fin) and living in degree $\ge 1$. Then
\begin{itemize}
\item[1.] If $\Bscr$ is noetherian then so is $\Ascr$.
\end{itemize}
Assume now that  $\Bscr$ is noetherian.
\begin{itemize}
\item[2.] If $\Bscr$ satisfies $\chi$ then so does $\Ascr$.
\item[3.] Let ``$\cd$'' stand for cohomological dimension. Then 
\[
\cd \tau_\Bscr+1\ge \cd\tau_\Ascr \ge \cd \tau_\Bscr
\]
\end{itemize}
\end{propositions}
\begin{proof}
\begin{enumerate}
\item
We have to show that if $\Mscr\in\Qch(X)$ is noetherian then so is
$\Mscr\otimes_{o_X}\Ascr$.

Now if $\Nscr\in\Gr(\Ascr)$ is left bounded then a variation of the classical
argument by Hilbert shows that if $\Nscr\otimes_\Ascr\Bscr$ is noetherian
then the same holds for $\Nscr$. Since by hypotheses
\[
(\Mscr\otimes_{o_X}\Ascr)\otimes_\Ascr\Bscr=\Mscr\otimes_{o_X}\Bscr
\]
is noetherian we are done.
\item
We start by observing that $\Bscr$ clearly satisfies the hypotheses of
Proposition \ref{ref:3.9.10a}, so we can use that proposition. 

In the sequel we will denote by $G$ the functor $-\otimes_\Ascr I$ and
$t$ the natural transformation $G\r \Id_{\Gr(\Ascr)}$ coming from the
inclusion $I\r \Ascr$.  Furthermore $G^n(\Nscr)\xrightarrow {t^n}
\Nscr$ is by definition obtained from tensoring
$I^n\xrightarrow{t^n}\Ascr$ with $\Nscr$. By $\Nscr t^n$ we denote the
image of this map. We say that $\Nscr$ is annihilated by $t^n$ if
$\Nscr t^n=0$. $\Nscr$ is $t$-torsion if it is the union of subobjects
annihilated by some $t^n$. Similarly $\Nscr$ is $t$-torsion free if
multiplication by $t$ is injective.

It is easy to see that the functors
\begin{align*}
-\otimes_\Ascr\Bscr&:\Gr(\Ascr)\r \Gr(\Bscr)\\
(-)_\Ascr&:\Gr(\Bscr)\r \Gr(\Ascr)
\end{align*}
induce inverse equivalences between $\Gr(\Bscr)$ and the full subcategory of
$\Gr(\Ascr)$ consisting of objects annihilated by $t$. So we will identify
these two categories.

Let $\Mscr\in\gr(\Ascr)$. As usual $\Mscr$ is an extension
\[
0\r\Mscr_1\r \Mscr\r \Mscr_2\r 0
\]
where $\Mscr_1$ is $t$-torsion and $\Mscr_2$ is $t$-torsion free. Now
$\Mscr$ is  noetherian and hence so is $\Mscr_1$. Thus $\Mscr_1t^n=0$ for some
$n$ and in particular we can write $\Mscr_1$ as an extension of
objects annihilated by $t$.

We conclude that to verify condition $\chi$ we have to show that
$R^i\tau_\Ascr(\Nscr)_{\ge 0}$ is right bounded and noetherian for two
classes of noetherian graded $\Ascr$-modules~:
\begin{enumerate}
\item 
Those that are annihilated by $t$.
\item 
Those that are $t$-torsion free.
\end{enumerate}
Let us first treat (a). If $\Nscr t=0$ then by the previous discussion
$\Nscr=\Nscr'_\Ascr$ for some $\Nscr'\in\gr(\Bscr)$. Thus
\begin{align*}
R^i\tau_\Ascr(\Nscr)&=R^i\tau_\Ascr(\Nscr'_\Ascr)\\
&=R^i\tau_\Bscr(\Nscr')_\Ascr  \qquad(\text{Prop. \ref{ref:3.9.10a}})
\end{align*}
which shows what we want.

Now consider the case that $\Nscr$ is torsion free. We have an exact sequence
in $\Gr(\Ascr)$ 
\[
0\r G(\Nscr)\xrightarrow{t} \Nscr\r \Nscr/\Nscr t\r 0
\]
and again $\Nscr/\Nscr t=\Nscr''_\Ascr$ for some $\Nscr''\in\gr(\Bscr)$. This
gives us the following exact sequence (using Prop.\ \ref{ref:3.9.10a})
\begin{equation}
\label{ref:3.41a}
R^{i-1}\tau_\Bscr(\Nscr'')_\Ascr\r G(R^i\tau_\Ascr(\Nscr))
\xrightarrow{t} R^i\tau_\Ascr(\Nscr)
\r 
R^i\tau_\Bscr(\Nscr'')_\Ascr
\end{equation}
By lemma \ref{ref:3.9.1a} this implies $R^i\tau_A(\Nscr)\in\Tors(\Ascr)$.
$R^i\tau_\Ascr(\Nscr)$ is
$t$-torsion since $I$ lives purely in positive 
 degree. On the other hand we obtain from \eqref{ref:3.41a} and the
fact that $\Bscr$ satisfies $\chi$ that $R^i\tau_\Ascr(\Nscr)_{\ge n}$
is $t$-torsion free for $n\gg 0$. Combining this we obtain that
$R^i\tau_\Ascr(\Nscr)_{\ge n}=0$ for $n\gg 0$ which shows that
$R^i\tau_\Ascr(\Nscr)$ is right bounded.

From \eqref{ref:3.41a} we obtain exact sequences
\[
0\r \text{noetherian} \r G(R^i\tau_\Ascr(\Nscr)_{\ge m-1})\xrightarrow{t}
R^i\tau_\Ascr(\Nscr)_{\ge m} 
\]
and by descending induction on $m$ we find that
$R^i\tau_\Ascr(\Nscr)_{\ge 0}$ is noetherian.
\item 
From Proposition \ref{ref:3.9.10a} it follows immediately that $\cd \tau_\Ascr\ge
\cd \tau_\Bscr$. We therefore concentrate on the
other inequality. Assume $\Nscr\in\Gr(\Ascr)$, $p=\cd\tau_\Bscr$. We
have to show $R^q\tau_\Ascr(\Nscr)=0$, $q>p+1$.

Since $\tau_\Ascr$ commutes with direct limits and $\Gr(\Ascr)$ is
locally noetherian,  $R^q\tau_\Ascr$ also commutes with direct limits
\cite{Groth1}. Hence we may assume that $\Nscr$ is
noetherian. Using the same reduction as in (2) we may assume that
either $\Nscr t=0$ or $\Nscr$ is $t$-torsion free, The first case is
trivial by Proposition \ref{ref:3.9.10a} so we look at the last case. Using
\eqref{ref:3.41a} for $i=q$ we find that $G(R^q\tau_\Ascr(\Nscr))\r
R^q\tau_\Ascr(\Nscr)$ is an isomorphism. Since on the other hand
$R^q\tau_\Ascr(\Nscr)$ is $t$-torsion we conclude $R^q\tau_\Ascr(\Nscr)=0$.
\qed\end{enumerate}
\def\qed{}\end{proof}
Proposition \ref{ref:3.9.13a}(3) can be stated in more geometric terms. 
\begin{corollarys} Assume that $\Ascr$, $\Bscr$ are noetherian
  $\NN$-graded algebras which fit in an exact sequence as in
  Prop. \ref{ref:3.9.13a}. Then
\begin{equation}
\label{ref:3.42a}
\cd\beta_\ast\le \cd \alpha_\ast\le \cd \beta_\ast+1
\end{equation}
\end{corollarys}
\begin{proof}
It is easy to see that
\[
\cd\alpha_\ast=\cd\underline{\alpha}_\ast=\cd\omega_\Ascr
\]
and similarly for $\beta$. Furthermore from \eqref{ref:3.29a}\eqref{ref:3.30a} it follows
that \begin{equation}
\label{ref:3.43a}
\cd\omega_\Ascr=\max(\cd\tau_\Ascr-1,0)
\end{equation}
Combining this with Prop.\ \ref{ref:3.9.13a}(3) yields what we want.
\end{proof}

\begin{lemmas} 
\label{ref:3.9.15a}
Assume $\Ascr$ is noetherian and $\NN$-graded. Let
$\Tscr$   be a right bounded graded $\Ascr-\Ascr$-bimodule. Then for
  $\Mscr\in\Gr(\Ascr)$ we have
\begin{equation}
\label{ref:3.44a}
\uHHom_\Ascr(\Tscr,\Mscr)=\uHHom_\Ascr(\Tscr,\tau(\Mscr))
\end{equation}
If   $\Tscr$ is in fact  coherent as 
$o_X-o_X$-bimodule  then $\uHHom_\Ascr(\Tscr,\Mscr)$ is torsion. 
\end{lemmas}
\begin{proof}
To prove \eqref{ref:3.44a} it is sufficient to show that the right and lefthand
side of that equation represent the same functor.
This is routine, using the fact that if $\Nscr\in \Gr(\Ascr)$ then
$\Nscr\otimes_\Ascr \Tscr$ is torsion by lemma \ref{ref:3.8.1a}.2.

Now we prove the second part of the lemma. By the first part we may
clearly assume that $\Mscr$ is torsion. The coherentness of $\Tscr$ also
implies that $\Tscr$ is left bounded.

We claim that $\Tscr$ is coherent as a graded $\Ascr$-bimodule. In the
same way as in \ref{ref:3.1.6a}   it suffices to show that
$-\otimes_\Ascr\Tscr$ preserves $\gr(\Ascr)$.  Now since $\Ascr$ is
noetherian, any object in $\gr(\Ascr)$ has a presentation consisting
of objects of the form $\Pscr\otimes_{o_X}\Ascr$ where the $\Pscr$ are
noetherian $o_X$-modules. We now use the fact that
$-\otimes_\Ascr\Tscr$ is right exact.

Hence $\uHHom_\Ascr(\Tscr,-)$ commutes with direct limits. Therefore we 
may assume that $\Mscr$ is  right bounded. 
But $\uHHom_\Ascr(\Tscr,\Mscr)$ is contained  in
$\uHHom_{o_X}(\Tscr,\Mscr)$ which is now clearly also right bounded. This
proves what we want.
\end{proof}

We will need the following variant of Proposition \ref{ref:3.9.13a}
\begin{lemmas}
\label{ref:3.9.16a}
 Assume that we have a surjective map $f:\Ascr\r \Bscr$ such
that $I= \ker f$ is a coherent object in
$\Bigr(o_X-o_X)$.
Then for $\Mscr\in\Gr(\Bscr)$ we have
\begin{equation}
\label{ref:3.45a}
R^i\tau_A(\Mscr_A)\cong (R^i\tau_\Bscr(\Mscr))_\Ascr
\end{equation}
Furthermore if $\Bscr$ satisfies $\chi$ then so does $\Ascr$ and $\cd
\tau_\Bscr=\cd \tau_\Ascr$.
\end{lemmas}
\begin{proof}
  It is clear that \eqref{ref:3.45a} is a special case of 
  Proposition \ref{ref:3.9.10a}.

  Assume that $\Bscr$ satisfies $\chi$. To show that $\Ascr$ also
  satisfies $\chi$ we have to show that for $\Nscr\in\gr(\Ascr)$ one
  has $R^i\tau_\Ascr(\Nscr)_{\ge 0}\in\gr(\Ascr)$.  Let
  $\Kscr=\ker(\Nscr\r \Nscr\otimes_\Ascr\Bscr)$. Since $\Nscr$ is
  noetherian, the same is true for $\Kscr$. Furthermore by Corollary
  \ref{ref:3.8.2a} it follows that $\Kscr$ is torsion.  Hence
\[
R\tau^i_\Ascr(\Kscr)=
\begin{cases} \Kscr&\text{if $i=0$}\\
0&\text{if $i>0$}
\end{cases}
\]
Thus the long exact sequence for $R^i\tau_\Ascr$ together with
\eqref{ref:3.45a} yields an exact sequence
\[
0\r \Kscr
\r \tau_\Ascr(\Nscr)\r \tau_\Bscr(\Nscr\otimes_\Ascr\Bscr)\r 0
\]
and isomorphisms
\[
R^i\tau_\Ascr(\Nscr)\r R^i\tau_\Bscr(\Nscr\otimes_\Ascr\Bscr)
\]
for $i\ge 1$. This yields that $\Ascr$ satisfies $\chi$ and also that
$\cd\tau_\Bscr\ge \cd \tau_\Ascr$. Since \eqref{ref:3.45a} implies that
$\cd\tau_\Bscr\le \cd \tau_\Ascr$, we are through.
 \end{proof}
\subsection{Higher inverse images}
\label{ref:3.10b}
Assume that $\alpha:Y\r X$ is a map of quasi-schemes. As was observed
above the definition of $R^i\alpha_\ast$ presents no difficulty since
$\Mod(Y)$ has enough injectives. However the definition of
$L_i\alpha^\ast$ is more delicate since we have not assumed that
$\Mod(X)$ has enough flat objects. As an approximation we define
$L_i\alpha^\ast$ as the functor $\MOD(X)\r \MOD(Y)$  given by
\[
\HHom_{o_Y}(L_i\alpha^\ast \Mscr,E)=\HExt^i_{o_X}(\Mscr,\alpha_\ast E)
\]
where $E$ runs through the injectives in $\Qch(Y)$. 
Furthermore we will define $\cd \alpha^\ast$ as the maximum $i$ such
that $L_i\alpha^\ast$, restricted to $\Mod(Y)$ is non-zero.
\begin{lemmas}
Let $X$ be a quasi-scheme and let $\Ascr\in\Alg(X)$. Put
$Y=\Spec \Ascr$ and let $\alpha:Y\r X$ be the structure
map. Then one has $L_i\alpha^\ast(-)=\HTor_i^{o_X}(-,\Ascr)$.
\end{lemmas}
This lemma is proved in a similar way as the lemma below which covers
the graded case. In the graded case we need more hypotheses since we
have defined things in less generality.
\begin{lemmas}
\label{ref:3.10.2a}
Let $X$ be a noetherian quasi-scheme and let $\Ascr$ be a noetherian
$\NN$-graded algebra on $X$. Put $Y=\Proj\Ascr$. Assume that
$\Mscr\in\Mod(X)$ is such that
$\underline{\HTor}_i^{o_X}(\Mscr,\Ascr)\in \Gr(\Ascr)$. Then
\[
L_i\alpha^\ast\Mscr=\pi(\underline{\HTor}_i^{o_X}(\Mscr,\Ascr))
\]
\end{lemmas}
\begin{proof}
We have to show that 
\begin{equation}\
\label{ref:3.46a}
\Ext^i_{o_X}(\Mscr,\alpha_\ast F)= \Hom_{\QGr(\Ascr)}
( \pi(\underline{\HTor}_i^{o_X}(\Mscr,\Ascr)) ,F)
\end{equation}
where $F$ runs through the injectives in $\Mod(Y)$. Put $F=\pi E$. It
follows from lemma
\ref{ref:3.8.3a}  that   $E$ is
an injective in $\Gr(\Ascr)$ satisfying $\tilde{E}=E$. Then by
adjointness the righthand side of \eqref{ref:3.46a} becomes 
\[
\Ext^i_{o_X}(\Mscr,\HHom_{\Gr(\Ascr)}(\Ascr,E))
=
\Ext^i_{o_X}(\Mscr,E_0)
\]
Hence the assertion we have to prove  boils down to $\alpha_\ast \pi
E=E_0$. This follows from the fact that $\alpha_\ast(-)=\omega(-)_0$.
\end{proof}

\subsection{On the positive part of an algebra strongly graded
with   respect to a localizing subcategory.} 
\label{ref:3.11b}
In this section we fix the
following situation. $X$ is a noetherian quasi-scheme, $\Sscr$ is a
localizing subcategory in $\Qch(X)$. 
$\Ascr$ is a noetherian algebra on $X$,
strongly graded with respect to $\Sscr$ (see \S\ref{ref:3.3b}). With
$\Ascr_{\ge 0}$ we denote the positive part of $\Ascr$. It is easy to see that $\Ascr_{\ge 0}$ is
also noetherian.

We note the following.
\begin{lemmas} 
\label{ref:3.11.1a}
Let $\Mscr\in\Gr(\Ascr_{\ge 0})$. Then the $\Ascr_{\ge
    0}$-module structure on $\tilde{\Mscr}$ extends in a natural way
  to  an $\Ascr$-module structure.
\end{lemmas}
\begin{proof}
Multiplication defines graded $\Ascr_{\ge 0}$-bimodule maps
\[
\Ascr_{\ge m}\otimes_{\Ascr_{\ge 0}}\Ascr_{\ge n}\r \Ascr_{\ge m+n}
\]
Applying $\uHHom_{\Ascr_{\ge 0}}(-,\Mscr)$ yields a map
\[
\uHHom_{\Ascr_{\ge 0}}(\Ascr_{\ge m+n},\Mscr)\r 
\uHHom_{\Ascr_{\ge 0}}(\Ascr_{\ge m},\uHHom_{\Ascr_{\ge 0}}(\Ascr_{\ge
n},\Mscr)) \]
which by adjointness yields a map
\[
\uHHom_{\Ascr_{\ge 0}}(\Ascr_{\ge m+n},\Mscr)\otimes_{\Ascr_{\ge
    0}}\Ascr_{\ge m} 
\r \uHHom_{\Ascr_{\ge 0}}(\Ascr_{\ge n},\Mscr)
\]
Taking direct limits over $n$, using \eqref{ref:3.24a} and letting $m$ go to
$-\infty$ 
we find a map \[
\tilde{\Mscr}\otimes_{\Ascr_{\ge 0}}\Ascr\r \tilde{\Mscr}
\]
A straightforward, but mildly tedious verification shows that this is
an $\Ascr$-module structure.
\end{proof}

According to Corollary \ref{ref:3.8.2a} we have the following inverse
equivalences
\begin{equation}
\label{ref:3.47a}
\Gr(\Ascr)/\Tors(\Ascr) 
\xymatrix{
\ar@<1ex>[rr]^-{(-)_{\Ascr_{\ge 0}}}
&
&
\ar@<1ex>[ll]^-{-\otimes_{\Ascr_{\ge 0}}\Ascr}
}
 \Gr(\Ascr_{\ge 0})/\Tors(\Ascr_{\ge 0})
\end{equation}
By lemma \ref{ref:3.3.1a} we have inverse equivalences
\begin{equation}
\label{ref:3.48a}
\Mod(\Ascr_0)/\Sscr \xymatrix{
\ar@<1ex>[rr]^-{-\otimes_{\Ascr_0}\Ascr}
&
&
\ar@<1ex>[ll]^-{(-)_{{0}}}
}
 \Gr(\Ascr)/\Sscr(\Ascr)
\end{equation}
In order to combine these equivalences we 
observe that
$
\Tors(\Ascr)\subset \Sscr(\Ascr)
$.
Indeed let $\Mscr\in\Tors(\Ascr)$. Since $\Tors(\Ascr)$ and
$\Sscr(\Ascr)$ are closed under direct limits we may assume that $\Mscr$ is
right bounded. But then 
 $\Mscr(N)_0=0\in\Sscr$ for $N\gg 0$
and hence according to lemma \ref{ref:3.3.1a} we have
$\Mscr(N)\in\Sscr(\Ascr)$. Thus the same holds for $\Mscr$.

Define $\QSscr(\Ascr_{\ge
  0})$  as the image in $\QGr(\Ascr_{\ge 0})$ of 
$\Sscr(\Ascr_{\ge 0})$ under the quotient map. 
One notes that $\QSscr(\Ascr_{\ge 0})$ has the following
 alternative description.
\begin{lemmas}
$\QSscr(\Ascr_{\ge 0})$ is precisely the image of $\Sscr(\Ascr)$
under  $\pi_\Ascr(-)_{\Ascr_{\ge 0}}$.
\end{lemmas}
\begin{proof} 
If we look at the commutative diagram
\[
\begin{CD}
\Gr(\Ascr_{\ge 0})@> \pi_{\Ascr_{\ge 0}}>> \QGr(\Ascr_{\ge 0})\\
@A (-)_{\Ascr_{\ge 0}} AA @A (-)_{\Ascr_{\ge 0}} AA\\
\Gr(\Ascr) @>\pi_\Ascr>> \QGr(\Ascr)
\end{CD}
\]
then we see that it is sufficient to show that the image of
$\Sscr(\Ascr)$ under $(-)_{\Ascr_{\ge 0}}$ is $\Sscr(\Ascr_{\ge 0})$
modulo $\Tors(\Ascr_{\ge 0})$. It is clear that this image is indeed
contained in $\Sscr(\Ascr_{\ge 0})$. Conversely let
 $\Mscr\in\Sscr(\Ascr_{\ge 0})$.  According to
Corollary \ref{ref:3.8.2a} we have $\Mscr= (\Mscr\otimes_{\Ascr_{\ge
    0}}\Ascr)_{\Ascr_{\ge 0}}$ modulo $\Tors(\Ascr_{\ge 0})$. Since
$\Mscr\otimes_{\Ascr_{\ge 0}}\Ascr$ is  contained in $\Sscr(\Ascr)$
we are through.
\end{proof}

Combining \eqref{ref:3.47a},\eqref{ref:3.48a} and the previous lemma
yields equivalences.
\[
\Mod(\Ascr_0)/\Sscr \xymatrix{
\ar@<1ex>[rr]^-{-\otimes_{\Ascr_0}\Ascr}
&
&
\ar@<1ex>[ll]^-{(-)_{{0}}}
}
 \Gr(\Ascr)/\Sscr(\Ascr)  \xymatrix{
\ar@<1ex>[rr]^-{(-)_{\Ascr_{\ge 0}}}
&
&
\ar@<1ex>[ll]^-{-\otimes_{\Ascr_{\ge 0}}\Ascr}
}
\QGr(\Ascr_{\ge
  0})/\QSscr(\Ascr_{\ge 0})
\]
Looking only at the outer categories yields equivalences
\[
\Mod(\Ascr_0)/\Sscr  \xymatrix{
\ar@<1ex>[rr]^-{-\otimes_{\Ascr_0}\Ascr_{\ge 0}}
&
&
\ar@<1ex>[ll]^-{(-\otimes_{\Ascr_{\ge 0}}\Ascr)_0}
}
\QGr(\Ascr_{\ge
  0})/\QSscr(\Ascr_{\ge 0})
\]
Define $U=\Spec \Ascr_0$, $Y=\Proj \Ascr_{\ge 0}$ and let the
map $\alpha :Y\r U$ be given by
$\alpha^\ast=\pi(-\otimes_{\Ascr_0}\Ascr_{\ge 0})$. 
\begin{lemmas}
  Let $\Mscr\in\Gr(\Ascr_{\ge 0})$. Modulo $\Sscr$, $\alpha_\ast\pi\Mscr$
  is given by $(\Mscr\otimes_{\Ascr_{\ge 0}}\Ascr)_0$.
\end{lemmas}
\begin{proof}
  Without loss
  of generality we may  assume that $\Mscr=\tilde{\Mscr}$. Hence
  in particular by lemma \ref{ref:3.11.1a} $\Mscr=\Nscr_{\Ascr_{\ge
      0}}$ where $\Nscr$ is a  graded $\Ascr$-module. Furthermore
  $\alpha_\ast\pi\Mscr= \Mscr_0=\Nscr_0$.

According to Corollary \ref{ref:3.8.2a} we now have that the canonical map
\[
\Mscr\otimes_{\Ascr_{\ge 0}}\Ascr=\Nscr_{\Ascr_{\ge
    0}}\otimes_{\Ascr_{\ge 0}} \Ascr\r \Nscr
\]
is an isomorphism modulo  $\Tors(\Ascr)\subset \Sscr(\Ascr)$. Hence the
restricted map
\[
(\Mscr\otimes_{\Ascr_{\ge 0}}\Ascr)_0\r \Nscr_0
\]
 is an isomorphism modulo $\Sscr$. This is precisely what we had to prove.
\end{proof}
Let us now introduce the more suggestive notation
$\alpha^{-1}(\Sscr)$ for $\QSscr(\Ascr_{\ge 0})$. Summarizing everything,
we have shown.
\begin{propositions} \label{ref:3.11.4a}
Let $U,Y,\alpha,\Sscr$ be as above. Then the
functors $\alpha^\ast$, $\alpha_\ast$ factor through $\Sscr$ and
$\alpha^{-1}(\Sscr)$ to define inverse equivalences between
$\Qch(U)/\Sscr
$ and $\Qch(Y)/\alpha^{-1}(\Sscr)$.
\end{propositions}

\subsection{Veronese subalgebras}
\label{ref:3.12b}
Let $X$ be a  quasi-scheme and let $\Ascr\in\Gralg(X)$ be noetherian and
$\NN$-graded.
The $n$'th Veronese of $\Ascr$ is the graded subalgebra $\Ascr^{(n)}$ of
$\Ascr$ defined by  $(\Ascr^{(n)})_m=\Ascr_{nm}$. If $\Mscr\in\Gr(\Ascr)$
then $\Mscr^{(n)}$ is defined similarly.

Then we have the
following.
\begin{lemmas} 
\label{ref:3.12.1a}
Assume that $X$ and $\Ascr$ are noetherian and that
$\Ascr$ is generated in degree one (cfr \S\ref{ref:3.9b}).
 Then the functors
\begin{align*}
\Gr(\Ascr)\r \Gr(\Ascr^{(n)}):&\Mscr\r \Mscr^{(n)}\\
\Gr(\Ascr^{(n)})\r\Gr(\Ascr):\Nscr\r \Nscr\otimes_{\Ascr^{(n)}}\Ascr
\end{align*}
factor over $\Tors(\Ascr)$ and $\Tors(\Ascr^{(n)})$ and in this way define
inverse equivalences between $\Proj\Ascr$ and $\Proj\Ascr^{(n)}$.
\end{lemmas}
\begin{proof}
This is formally similar to the ring case. See for example
\cite{V}.
\end{proof}

\section{Pseudo-compact rings}
\label{ref:4a}
In the sequel we will study the formal local structure of some
specific quasi-schemes. It turns out that this is best described in
terms of pseudo-compact rings, so for the convenience of the reader we
collect some of the properties of such rings here. Most of this is
taken from \cite{VdBVG} and \cite{Gabriel}. We refer the reader to
these papers for more information.

A topological right  module $M$ over a topological ring $A$
is \emph{pseudo-compact} if it is Hausdorf, complete and its topology
is generated by right submodules of finite colength. $A$ is said to be a
\emph{pseudo-compact ring} if $A$ is pseudo-compact as a right $A$-module.
Left pseudo-compact is defined similarly. By $\PC(A)$,
we will denote the category of right pseudo-compact modules over a
right pseudo-compact ring. By \cite{Gabriel}  $\PC(A)$ is an
abelian category satisfying AB5$^\ast$ and AB3. 

Let
$A$  be a pseudo-compact ring.
The dual category of $\PC(A)$ is a  locally finite category. That is a
Grothendieck category generated by  objects of finite length.
Conversely, if $\Cscr$ is a locally finite category then $\Cscr$ can be
realized as $\PC(A)^\circ$ for some pseudo-compact ring $A$ 
\cite{Gabriel}. $A$ is constructed as follows. Let $E$ be an injective
cogenerator of $\Cscr$, containing every indecomposable injective at
least once and put $A=\End_\Cscr(E)^\circ$. If $S$ is a
finite length subobject of $E$ then
\[
\frak{l}(S)=\{f\in A\mid f(S)=0\}
\]
defines a right ideal of finite colength in $A$. We take these right
ideal as the basis for a topology on $A$. In this way $A$ becomes a
pseudo-compact ring.
If $M\in \Cscr$ then we topologize $\Hom_\Cscr(M,E)$ in a similar way.
The functor $M\mapsto \Hom_\Cscr(M,E)$ defines a duality between $\Cscr$
and $\PC(A)$.

As an application we obtain~:
\begin{lemma}
\label{ref:4.1a}
Inverse limits of projectives are projective in $\PC(A)$.
\end{lemma}
\begin{proof} This follows from the fact that $\PC(A)^\circ$ is a
  locally finite category and hence direct limits of injectives are
  injective.
\end{proof}
In nice cases there are good relations between the properties of
$\PC(A)$ and $\Mod(A)$. For example 
by \cite{Gabriel} it
follows 
that the forgetful functor $\PC(A)\r \Mod(A)$ is faithful and commutes
with kernels, cokernels and products. 
By \cite[Lemma 3.4]{VdBVG}, an object in $\PC(A)$ is simple in
$\PC(A)$ if and only it is  simple in $\Mod(A)$. A similar result holds
for the property of being noetherian \cite[Cor 3.10]{VdBVG}. 
As usual we denote by $\pc(A)$ the category of
noetherian pseudo-compact $A$-modules.

An object $M$ in $\PC(A)$ is said to be finitely generated in $\PC(A)$ if 
 $M$ is a quotient of $A^n$ in $\PC(A)$ for some $n$. If
$M,N\in\PC(A)$, $M$ finitely generated then according to \cite[Prop.\
3.5]{VdBVG} we have
\[
\Hom_{\PC(A)}(M,N)=\Hom_{\Mod(A)}(M,N)
\]
In particular, $M$ is finitely generated in $\PC(A)$ if and only if it is
finitely generated in $\Mod(A)$. From
this one can deduce~:
\begin{lemma}
\label{ref:4.2a}
\cite[Cor.\ 3.8]{VdBVG} If $M$ is a finitely generated pseudo-compact
$A$-module then $L\subset M$ is open if and only if $M/L$ is of finite
length and pseudo-compact when equipped with the discrete topology.
\end{lemma}
(Note that a linear Hausdorf topology on a module of finite length is
automatically discrete.)  If $A$ is noetherian then it follows from
\cite[Prop.\ 3.19]{VdBVG} that the forgetful functor $\pc(A)\r
\mod(A)$ is an equivalence of categories. More generally we say that
$A$ is \emph{locally noetherian} if for every primitive idempotent $e$
in $A$ we have that $eA$ is noetherian. Assume that $A$ is locally
noetherian.  Then by \cite[Cor. 3.15]{VdBVG} we have that the
forgetful functor $\pc(A)\r\mod(A)$ is fully faithful and closed under
extensions.  Let $\PCFin(A)$ denote the category of finite length
objects in $\PC(A)$.  It follows that the objects are precisely the
finite length objects in $\Mod(A)$ whose Jordan-Holder quotients are
pseudo-compact simples.  It also follows that if $M$ is a noetherian
$A$-module then the topology on $M$ is simply the cofinite topology.

Let $A$ be an arbitrary pseudo-compact ring. Then we denote by $\Dis(A)$
the category of topological $A$-modules which are discrete. It is
clear that $\Dis(A)$ is the  full subcategory of $\Mod(A)$ consisting of
modules $M$ such that for all $m\in M$, $\Ann_A(m)$ is open 
in $A$, or equivalently that $mA$ is
pseudo-compact of finite length (this follows for example from lemma
\ref{ref:4.2a}). From this we deduce
that $\Dis(A)$ is a locally finite category. Clearly $\PC(A)\cap
\Dis(A)=
\PCFin(A)$  (where the intersection is taken
inside the category of topological right $A$-modules).

It is interesting to observe that since $\Dis(A)$ is locally finite there
must necessarily exist another pseudo-compact ring $A^\ast$ such that
$\PC(A^\ast)=\Dis(A)^\circ$. In nice cases we have $A^\ast=A^\circ$ (see 
for example Proposition \ref{ref:4.10a} below).

Let $\Top(A)$ be the additive category of topological right $A$-modules.
For $M,N\in \Top(A)$ we  have functors
\begin{gather}
\label{ref:4.1b}
\Hom_{\Top(A)}(M,-): \Dis(A)\r \Ab\\
\Hom_{\Top(A)}(-,N):\PC(A)^\circ\r \Ab
\label{ref:4.2b}
\end{gather}
\begin{proposition}
\label{ref:4.3a}
Let $(M_i)_{i\in J}$, $(N_j)_{j\in J}$ be respectively an inverse
system in $\PC(A)$ and a directed system in $\Dis(A)$. Then
\begin{equation}
\label{ref:4.3b}
\Hom_{\Top(A)}(\invlim_i M_i,\dirlim_j
N_j)=\dirlim_i\dirlim_j\Hom_{\Top(A)}(M_i,N_j)
\end{equation}
\end{proposition}
\begin{proof}
  If $M\in\PC(A)$, $N\in\Dis(A)$ then every continuous morphism $f:M\r
  N$ has an open kernel $M'$. This means that $N'=\im f\cong M/M'$ has
  finite length.  Thus we have the following equalities
\begin{align}
\label{ref:4.4a}
\Hom_{\Top(A)}(M,N)&=\dirlim_{N'}\Hom_{\PC(A)}(M,N')\\
\label{ref:4.5a}
&=\dirlim_{M'}\Hom_{\Dis(A)}(M/M',N)\\
\label{ref:4.6a}
&=\dirlim_{M',N'}\Hom_A(M/M',N')
\end{align}
where $M'$ runs trough the open submodules in $M$ and $N'$ runs through
the finite length submodules of $N$. 

Now let us for example show that $\Hom_{\Top(A)}(-,N)$ sends inverse
limits  to direct limits. We have
\begin{align*}
\Hom_{\Top(A)}(\invlim_i M_i,N)&=
\dirlim_{N'} \Hom_{\PC(A)}(\invlim_i M_i,N')\qquad\text{(eq.
\eqref{ref:4.4a})}\\
&=\dirlim_{N'} \Hom_{\PC(A)^\circ}(N',\dirlim_i M_i)
\\ &=\dirlim_{i,N'} \Hom_{\PC(A)^\circ}(N',M_i)\\
&=\dirlim_{i,N'}\Hom_{\PC(A)}(M_i,N')\\
&=\dirlim_i \Hom_{\Top(A)} (M_i,N)\qquad\text{(eq.
\eqref{ref:4.4a})}
\end{align*}
The third equality follows from the fact that $N'$ is a noetherian object
in the locally noetherian category $\PC(A)^\circ$.

The proof that $\Hom_{\Top(A)}(M,-)$ commutes with direct limits is similar.
\end{proof}

Since $\Dis(A)$ has enough injectives and $\PC(A)$ has enough
projectives, we can define the derived functors of \eqref{ref:4.1b} and
\eqref{ref:4.2b}. 
Let us  temporarily denote them by $\Ext_{I}^i(M,-)$ and $\Ext_{I\!
I}^i(-,N)$. 
\begin{lemma} $\Ext_I^i$ and $\Ext_{I\!I}^i$ coincide when both are
defined. That is if $M\in\PC(A)$ and $N\in \Dis(A)$ then
\[
\Ext_I^i(M,N)=\Ext_{I\!I}^i(M,N)
\]
\end{lemma}
\begin{proof} To prove this we have to show that if $P$ is projective in
$\PC(A)$ then $\Hom_{\Top(A)}(P,-)$ is exact on $\Dis(A)$ and if $E$ is
an injective in $\Dis(A)$ then $\Hom_{\Top(A)}(-,E)$ is exact when
evaluated on $\PC(A)$. Since these are obviously dual statements we
only prove the first one.

By Proposition \ref{ref:4.3a}
\[
\Hom_{\Top(A)}(P,N)=\dirlim_{N'} \Hom_{\PC(A)}(P,N')
\]
where $N'$ runs through the finite length submodules of $N$. Since
$\Hom_{\PC(A)}(P,-)$ is exact on $\PC(A)$ and
$\dirlim$ is exact on $\Ab$ we have to show that every exact sequence
in $\Dis(A)$
\[
0\r N_1\r N\r N_2\r 0
\]
can be obtained as a direct limit of exact sequences of the form 
\[
0\r N'_1\r N'\r N_2'\r 0
\]
where $N'_1,N_2',N'$ are finite length subobjects of $N_1,N_2,N$. That
this is true follows easily from the fact that $\Dis(A)$ is locally
finite.
\end{proof}
Henceforth we neglect the distinction between $\Ext_{I}^i$ and
$\Ext_{I\!I}^i$, and we simply write $\Ext^i_{\Top(A)}$.
We then obtain
the following generalization of Proposition \ref{ref:4.3a}
\begin{proposition}
\label{ref:4.5b}
Let $(M_i)_{i\in I}$, $(N_j)_{j\in J}$ be respectively an inverse
system in $\PC(A)$ and a directed system in $\Dis(A)$. Then
\begin{equation}
\label{ref:4.7a}
\Ext^i_{\Top(A)}(\invlim_i M_i,\dirlim_j
N_j)=\dirlim_i\dirlim_j\Ext^i_{\Top(A)}(M_i,N_j)
\end{equation}
\end{proposition}
\begin{proof}
Let $N\in \Dis(A)$ and let $E^\cdot$ be an injective resolution of $N$ in
$\Dis(A)$. Then
\begin{align*}
\Ext^i_{\Top(A)}(\invlim_i M_i,N)&=H^i(\Hom_{\Top(A)}(\invlim_i
M_i,E^\cdot))\\
&=\dirlim_i H^i(\Hom_{\Top(A)}(M_i,E^\cdot))\qquad \text{(Proposition
\ref{ref:4.3a})}
\\ &=\dirlim_i \Ext^i(M_i,N)
\end{align*}
The fact that $\Ext^i_{\Top(A)}(M,-)$ is compatible with direct limits is
proved similarly.
\end{proof}
Of course the ordinary ``$\Ext$'' in the abelian categories $\PC(A)$
and $\Dis(A)$ is also defined. 
 Since $\PC(A)$ has enough projectives we
clearly have
\[
\Ext^i_{\PC(A)}(M,N)=\Ext^i_{\Top(A)}(M,N)
\]
if $M,N\in \PC(A)$. Similarly if $M,N\in \Dis(A)$ then 
\[
\Ext^i_{\Dis(A)}(M,N)=\Ext^i_{\Top(A)}(M,N)
\]
since $\Dis(A)$ has enough injectives.
\begin{proposition}
\label{ref:4.6b}
Let $M\in\PC(A)$, $N\in \Dis(A)$. We have the following formulas
\[
\label{ref:4.6c}
\projdim_{\PC(A)} M=\sup_i \{i\mid\Ext^i_{\PC(A)}(M,-)\neq 0\}
\]
\[
\label{ref:4.6d}
\injdim_{\Dis(A)} N=\sup_i \{i\mid\Ext^i_{\Dis(A)}(-,N)\neq 0\}
\]
If $A$ is locally noetherian and if $M$ is noetherian in $\PC(A)$ then
\[
\label{ref:4.6e}
\projdim_{\PC(A)} M=\sup_i\{i\mid \exists S\text{ simple }: \Ext^i_{\PC(A)}(M,S)\neq 0\}
\]
and if $A^\ast$ is locally noetherian and $N$ is artinian in $\Dis(A)$ then
\[
\label{ref:4.6f}
\injdim_{\Dis(A)} N=\sup_i\{i\mid \exists T\text{ simple }:
\Ext^i_{\Dis(A)}(T,N)\neq 0\}
\]
\end{proposition}
\begin{proof}
This is entirely classical.
 Let us for example prove the fourth equality.
By degree shifting this amounts to showing that if $\Ext^1(T,N)=0$ for
all simple $T$ and if $N$ is artinian then $N$ is injective.

Let $E$ be the injective hull of $N$ and $U=E/N$. Let $q:E\r U$ be the
quotient map.
By hypotheses $E$ and hence $U$ is artinian. Also by hypotheses, the
restriction of $q$ to 
$\Soc(E)\r \Soc(U)$ is surjective. Since socles are by definition 
semisimple this last map has
a splitting which can be lifted to a map $t:U\r E$. By hypotheses $s=qt$
is the identity on $\Soc(U)$, from which it follows that $s$ is
injective. Also $s^n(U)$ is a descending chain of submodules in $U$ and
thus $s^n(U)=s^{n+1}(U)$. The injectivity of $s$ yields that $U=sU$ and
thus $s$ is an automorphism of $U$. Now $ts^{-1}$ is a splitting of $q$
and hence $N$ is a direct summand of $E$, whence injective.
\end{proof}
As usual we define
\begin{align*}
\gldim \PC(A)&=\sup_{M\in \PC(A)} \projdim M\\
\gldim \Dis(A)&=\sup_{N\in \Dis(A)}\injdim N
\end{align*}
\begin{corollary}
\label{ref:4.7b}
The following holds.
\begin{enumerate}
\item $\gldim \PC(A)$ is equal to the supremum of the projective
dimensions of the simple pseudo-compact $A$-modules 
(as objects in $\PC(A)$).
\item $\gldim \Dis(A)$ is equal to the supremum of the injective
dimensions of the simple pseudo-compact $A$-modules 
(as objects in $\Dis(A)$).
\item If $A$ and $A^\ast$ are locally noetherian then 
$\gldim \PC(A)$ and $\gldim \Dis(A)$ are both equal to the supremum of the
$i\in\NN$ such that there exist simple pseudo-compact-$A$-modules 
$S,T$ such that $\Ext^i(S,T)\neq 0$.
\end{enumerate}
\end{corollary}
\begin{proof}
1. and 2. are true by \cite[Lemma 5.1]{VdBVG}. 3. follows from 1. and 2.
and the foregoing proposition.
\end{proof}

An object in $\PC(A)$ is said to be cosemisimple if it is a direct
product of simple modules. This is equivalent with being semisimple in
the dual category $\PC(A)^\circ$.  If $M\in \PC(A)$ then we define
$M/\rad(M)$ as  the quotient of $M$ which is the socle of $M$ in
$\PC(A)^\circ$. By construction $M/\rad(M)$ is the largest cosemisimple
object in $\PC(A)$ which is a quotient of $M$.
The \emph{radical} of $M$, denoted by $\rad(M)$, is  defined as the
kernel of $M\r M/\rad(M)$. Since  it is closed, it is pseudo-compact.
From the fact that taking socles is left exact it also follows that the
functor $M\mapsto M/\rad(M)$ is right exact. It is shown in
\cite{Gabriel} that $\rad(A)$ is a twosided ideal and coincides with
the ordinary Jacobson radical of $A$. From the fact that $\rad(A)$
annihilates all cosemisimple objects in $\PC(A)$ we obtain
\begin{equation}
\label{ref:4.8a}
M\rad(A)\subset \rad(M)
\end{equation}
\begin{lemma} 
\label{ref:4.8b}
(Nakayama's lemma) If $M\in \PC(A)$ then $M=\rad(M)$ if
and only if $M=0$.
\end{lemma}
\begin{proof} It is easy to see that a non-zero object in $\PC(A)$ maps 
onto at least one simple object. This proves the lemma.
\end{proof}
\begin{lemma}
\label{ref:4.9a}
The following are equivalent.
\begin{enumerate}
\item $M$ is finitely generated.
\item $M/M\rad(A)$ is a finitely generated $A/\rad(A)$-module.
\item $M/\rad(M)$ is a finitely generated $A/\rad(A)$-module.
\end{enumerate}
If any one of these conditions holds then $\rad(M)=M\rad(A)$.
\end{lemma}
\begin{proof}
It is clear that 1. implies 2. From \eqref{ref:4.8a} it follows 
that 2. implies 3. Hence we have to show that 3. implies 1.  By lifting
the generators of $M/\rad(M)$ we can construct a map $\theta:A^k\r M$
which becomes surjective after applying the functor $T\mapsto T/\rad(T)$.
Given the right exactness of this functor we obtain that $C=\rad(C)$
for $C=\coker \theta$. By Nakayama it follows that $C=0$. This proves
the first part of the lemma.

If condition 2. holds then $M/M\rad(A)$ is a quotient of $(A/\rad(A))^k$
for some $k$, and hence $M/M\rad(A)$ is cosemisimple. We conclude that
$\rad(M)\subset M\rad(A)$.
\end{proof}
For further reference we state the following formula
\begin{equation}
\label{ref:4.9b}
M/\rad(M)=\prod_{S\text{ simple}} S^{\alpha_{M,S}}
\end{equation} 
where $\alpha_{M,S}=\dim_{\End_{\PC(A)}(S)}\Hom_{\PC(A)}(M,S)$.
This is easily proved by looking at the dual statement in
$\PC(A)^\circ$.

The definition of a pseudo-compact ring is essentially onesided, which
is somewhat inconvenient. We now introduce a symmetric notion.

Let $k$ be a field. A pseudo-compact ring
which is a $k$ algebra is said to be \emph{cofinite} if all simple
pseudo-compact $A$-modules are finite dimensional over $k$.
\begin{proposition} \label{ref:4.10a} Assume that $A$ is cofinite. Then
\begin{enumerate}
\item $A/\rad(A)=\prod_{i\in I} M_{n_i}(D_i)$, for finite dimensional
  division algebras $(D_i)_{i\in I}$.
\item
The topology on $A$
is generated by twosided ideals of finite codimension. 
\item $A$ is left and right pseudo-compact.
\item $S\mapsto \Hom_k(S,k)$ defines a duality between left and right
pseudo-compact $A$-modules of finite length. In particular we can take
$A^\ast=A^\circ$.
\item (Matlis-duality) $M\mapsto \Hom_{\Top(k)}(M,k)$ defines a
  duality between pseudo-compact left (right) $A$-modules and discrete
  right (left)
  $A$-modules. Thus $\PC(A)^\circ=\Dis(A^\circ)$ and
  $\Dis(A)^\circ=\PC(A^\circ)$.
\end{enumerate}
\end{proposition}
\begin{proof}
\begin{enumerate}
\item According to Gabriel $A/\rad(A)$ is a product of endomorphism
  rings of vectorspaces over division algebras. The simple
  pseudo-compact modules over such a ring will be finite dimensional if
  and only if $A/\rad(A)$ has the indicated form.
\item
Let $L\subset A$ be an open ideal. Then $S=A/L$ is a pseudo-compact
right $A$-module of finite length. Put $T=\End_k(S)$ and consider $T$ as
an $A$-bimodule in the obvious way. Then as right $A$-module, we have
$T=S^{\dim_k S}$ and in particular $T$ is pseudo-compact. There is a
canonical map of $A$-bimodules $A\r T$, given by the right action of $A$
on $S$. Let $M$ be the kernel of this map. Then $M$ is an open twosided
ideal contained in $L$.
\item Since the topology on $A$ is generated by twosided ideals of finite
codimension (by 2.) we see that $A$ is also pseudo-compact on the left.
\item This follows again easily from the fact that $S$ is annihilated by
an open twosided ideal.
\item This is a consequence of 3. and the fact that objects in
  $\Dis(A)$  are direct limits of finite length pseudo-compact modules and
  objects in $\PC(A)$ are inverse limits of  finite length
  pseudo-compact modules. 
\qed
\end{enumerate}
\def\qed{}\end{proof} 
\begin{corollary}
\label{ref:4.11a}
 For a
left and right locally noetherian cofinite pseudo-compact $k$-algebra the numbers
$\gldim \PC(A)$, $\gldim \Dis(A)$, $\gldim \PC(A^\circ)$ and $\gldim
\Dis(A^\circ)$ all coincide. 
We call this common value the global
dimension of $A$ and denote it by $\gldim A$.
\end{corollary}
\begin{proof}
This follows from
Corollary \ref{ref:4.7b}.3 and Proposition \ref{ref:4.10a}.5.
\end{proof}
\begin{definition}
\label{ref:4.12a}
  Let $A,B$ be cofinite $k$-algebras. Let $M$ be a topological
  $A$-$B$-bimodule. Then $M$ is bi-pseudo-compact if the topology on
  $M$ is Hausdorf and complete and is generated by subbimodules
  $M'\subset M$ of finite codimension.
\end{definition}
It is easy to see that bi-pseudo-compact bimodules form an abelian category
satisfying AB3 and AB5${}^\ast$. We denote this category by $\PC(A-B)$.

If $A,B,C$ are three cofinite algebras and $M\in \PC(A-B)$, $N\in
\PC(B-C)$ then we define 
\[
M\ctimes_B N =\invlim_{M',N'} M/M'\otimes_B N/N'
\]
where $M',N'$ run through all open subbimodules in $M$, $N$. By
construction $M\ctimes_B N\in \PC(A-C)$. There are analogous definitions
 when $M$ or $N$ are onesided pseudo-compact modules.

One easily obtains
\begin{lemma}
  $-\ctimes_B-$ is right exact and commutes with inverse limits in
  both its factors.
\end{lemma}
The following is also standard.
 \begin{lemma}
\label{ref:4.14a}
 If $M$ is finitely presented as right $B$-module then
 \[
 M\ctimes_B N=M\otimes_B N
 \]
\end{lemma}
If $A,B$ are cofinite $k$-algebras then $A\ctimes_k B$ carries a canonical
ring structure which makes it into a cofinite $k$-algebra.
One has $\PC(A-B)=\PC(A^\circ\ctimes_k B)$. Furthermore the
pseudo-compact simple $A^\circ\ctimes B$-modules are of the form
$S\otimes_k T$ with $S$ simple in $\PC(A)$ and $T$ simple in $\PC(B)$
(this follows from the corresponding statement for finite dimensional
algebras). From this one deduces $A\ctimes_k B/\rad (A\ctimes_k
B)=A/\rad(A)\ctimes_k B/\rad(B)$.

 This observation
allows one to prove
\begin{lemma} 
\label{ref:4.15a}
The forgetful functor $\PC(A-B)\r \PC(B)$ preserves
projectives. 
\end{lemma}
\begin{proof} It suffices to prove this for the functor
  $\PC(A^\circ\ctimes B)\r \PC(B)$. By the above discussion all
  projectives in $\PC(A^\circ\ctimes B)$ are products of $ Ae\ctimes
  fB$ where $e,f$ are primitive idempotents in $A$, $B$. These are
  direct summands of $A\ctimes B$, whence it suffices to show that
  $(A\ctimes B)_B$ is projective. This follows from the fact that
\begin{align*}
A\ctimes B&
=\invlim_{I,J} (A/I)\otimes (B/J)\\
&=\invlim_I(\invlim_J (A/I)\otimes (B/J))\\
&=\invlim_I A/I\otimes B
\end{align*}
where $I,J$ run through the open twosided ideals in $A$ and $B$. Hence
$(A\ctimes B)_B$ is an inverse limit of projective $B$-modules, and is
therefore itself projective (lemma \ref{ref:4.1a})
\end{proof}
We now give a  structure theorem on cofinite $k$-algebras, which was more
or less proved 
 in \cite{VdBVG} (in slightly greater generality).
Below let $(e_i)_{i\in I}$ be a summable set
of primitive idempotents in a pseudo-compact ring $A$, having sum 1 (as in
\cite{Gabriel}). 
\begin{proposition}
\label{ref:4.16a}
  Let $A$ be cofinite, let $M$ be a pseudo-compact $A$-module and let
  $N$ be a bipseudo-compact $A$-bimodule. Put $A_{ij}=e_iAe_j$,
  $M_i=Me_i$, $N_{ij}=e_iNe_j$, equipped with the induced topology.
  Then \begin{align} A&=\prod_{i,j} A_{ij}\label{ref:4.10b}\\ 
    M&=\prod_i M_i\label{ref:4.11b}\\
N&=\prod_{i,j} e_iNe_j 
\end{align}
These equalities are in fact homeomorphisms if we equip the righthand
sides with the product topology. Furthermore the $A_i$ are cofinite
with $A_i/\rad(A_i)=D_i$ (where $D_i$ is as in Proposition
\ref{ref:4.10a}), the $M_i$ are pseudo-compact $A_i$-modules and the
$A_{ij}$ and $N_{ij}$ are bipseudo-compact $A_i-A_j$-modules. \end{proposition}
\begin{proof}
Similar to \cite[Prop.\ 4.3]{VdBVG}.
\end{proof}
In the case of a cofinite $k$-algebra, there is a simple test for a
pseudo-compact module to be finitely generated.
\begin{lemma}
\label{ref:4.17a}
Assume that $A$ is cofinite. Then $M\in \PC(A)$ is finitely generated
if and only if for every simple pseudo-compact $S$ we have that
$\dim_{\End_{\PC(A)}(S)}\Hom_{\PC(A)}(M,S)$ is finite and is
bounded independently of $S$.
\end{lemma}
\begin{proof}
This follows from Lemma \ref{ref:4.9a} and \eqref{ref:4.9b}.
\end{proof}
The following is a generalization of lemma \ref{ref:4.17a}.
\begin{lemma}
\label{ref:4.18a}
Assume that $A$ is cofinite. Let $M\in\PC(A)$. Then $M$ has a
 free resolution of length $n$ in $\PC(A)$ (or, equivalently, in
$\Mod(A)$) \[
F_n\r F_{n-1}\r \cdots \r F_1\r F_0\r M\r 0
\]
where $F_i$ is  of finite rank over $A$, if and only if
for every
$i\in\{0,\ldots,n\}
$ and for every simple $S$ the dimension of $\Ext^i_{\PC(A)}(M,S)
$ is finite and bounded independently of $S$.
\end{lemma} \begin{proof}
This is easily proved from the case $n=0$, by degree shifting.
\end{proof}

We will need the following result.
\begin{proposition}
\label{ref:4.19a}
 Assume that $A$ is cofinite, $M\in \PC(A-A)$
  and $S\in \Dis(A)$. Then 
 $\Ext^i_{\Top(A)}(M,S)\in \Dis(A)$. Here the ``$\Ext$'' is taken with
respect to the right $A$-structure on $M$.
\end{proposition}
\begin{proof} 
Let $P^\cdot$ be a projective resolution of $M$ in $\PC(A-A)$. According
to lemma \ref{ref:4.15a} the terms of $P^\cdot$ are projective in $\PC(A)$.
Hence we find
\[
\Ext^i_{\PC(A)}(M,S)=H^i(\Hom_{\PC(A)}(P^\cdot,S))
\]
Hence it suffices to show that if $P$ is a projective object in
$\PC(A-A)$  then $\Hom_{\PC(A)}(P,S)\in \Dis(A)$.  Such a projective
object is a product of direct summands of $A^\circ\ctimes A$, and since
$\Hom_{\PC(A)}(-,S)$ transforms products in the first argument into direct
sums (according to Proposition \ref{ref:4.3a}) it suffices to show that 
$\Hom_{\PC(A)}(A\ctimes A,S)$ lies in $\Dis(A)$. Now we have
\begin{align*}
\Hom_{\PC(A)}(A\ctimes A,S)&=\Hom_{\PC(A)}(\invlim_I A/I\otimes A,S)\\
&=\dirlim_I \Hom_{\PC(A)}(A/I\otimes A,S)\qquad (\text{Proposition
\ref{ref:4.3a}})\\ 
&=\dirlim_I (A/I)^\ast\otimes S
\end{align*}
Here $I$ runs through the open twosided ideals in $A$. The right
$A$-structure on $(A/I)^\ast\otimes S$ we use is the one on $(A/I)^\ast$.
Thus $(A/I)^\ast\otimes S=(A/I)^t$ for some $t$, and we are done.
\end{proof}
We can define the derived functors of $-\ctimes_{B}-$ in both 
arguments. For lack of a better notation we denote them by
$\Tor_i^{\PC(B)}(-,-)$. 
Thus if $A,B,C$ are cofinite $k$-algebras and $M\in \PC(A-B)$,
$N\in\PC(B-C)$ then we compute $\Tor^{\PC(B)}_i(M,N)$ in the  usual
way. 
For example we can start with a projective resolutions of $M$  in
$\PC(A-B)$. According to lemma  \ref{ref:4.15a}  this
yields a projective resolution of $M$ in $\PC(B)$. We can also start with
a projective resolution of $N$ and get the same result.

We will need the following.
\begin{lemma}
\label{ref:4.20a}
Let $E$ be an injective object in $\Dis(C)$. Then
\[
\Hom_{\Top(C)}(\Tor^{\PC(B)}_i(M,N),E)=\Ext^i_{\Top(B)}
(M,\Hom_{\Top(C)}(N,E))
\]
\end{lemma}
\begin{proof}
This follows easily if we replace $M$ by a projective resolution.
\end{proof}

\section{Cohen-Macaulay curves embedded in quasi-schemes}
\label{ref:5a}

\subsection{Preliminaries}
\label{ref:5.1a}
In the sequel $k$ will be an algebraically closed base field. We will
(usually tacitly) assume that all quasi-schemes are in $\QSch/\Spec
k$. Note that if $(X,\gamma)\in \QSch/\Spec k$ then $X$ contains the
canonical object $\Oscr_X=\gamma^\ast k$ (the ``structure
sheaf''). However this extra structure on $X$ will not be used until
\S\ref{ref:6.6b}.

Below $i:Y\r X$ will be a biclosed embedding of a commutative
Cohen-Macaulay curve $Y/k$ as a divisor (in the enriched sense) in a
noetherian quasi-scheme $X/k$ (\S\ref{ref:3.7b}).  In the commutative
case this hypothesis would imply that $X$ is a surface in a
neighborhood of $Y$.

Throughout this paper we will impose the following smoothness
  condition on $X$.
  {\def\thehypothesis{(*)}
\begin{hypothesis} Every object in $\Qch(Y)$ has finite injective
  dimension in $\Qch(X)$.
\end{hypothesis}
}
It is easy to see that this is equivalent with the seemingly weaker
condition.
  {\def\thehypothesis{(*')}
\begin{hypothesis} For every $p$ one has that $\Oscr_p$ has finite injective
  dimension in $\Qch(X)$.
\end{hypothesis}
}
The latter condition is sometimes automatic as can be seen from the
following lemma.
\begin{lemmas}
\label{ref:5.1.1a}
Assume that $Y$ is smooth in $p$. Then $\Oscr_p$ has finite injective
dimension in $\Qch(X)$.
\end{lemmas}
\begin{proof} It is easy to see that we have to show that there is some
$n$ such that $\Ext^n_{\Qch(X)}(\Fscr,\Oscr_p)=0$ for all $\Fscr\in\coh(X)$. Using
the long exact sequence for $\Ext$ we only have to show this in the
following two cases.
\begin{enumerate}
\item
The canonical map $\Fscr(-Y)\r \Fscr$ is injective.
\item $\Fscr\in\coh(Y)$.
\end{enumerate}
The lemma now follows from Propositions \ref{ref:5.1.2a} below.
\end{proof}
\begin{propositions} 
\label{ref:5.1.2a}
Let $\Fscr\in \Qch(X)$, $\Sscr\in \Qch(Y)$. Then
\begin{enumerate}
\item If the canonical map $\Fscr(-Y)\r \Fscr$ is an injection then
\[
\Ext^i_{\Qch(X)}(\Fscr,\Sscr) =\Ext^i_{\Qch(Y)}(\Fscr/\Fscr(-Y),\Sscr)
\]
\item If $\Fscr\in \Qch(Y)$ then there is a long exact sequence
\begin{multline*}
\r \Ext^{i-2}_{\Qch(Y)}(\Fscr,\Sscr(Y))
\r \Ext^i_{\Qch(Y)}(\Fscr,\Sscr)\r \Ext^i_{\Qch(X)}(\Fscr,\Sscr)
\\ \r \Ext^{i-1}_{\Qch(Y)}(\Fscr,\Sscr(Y))
\r 
\Ext^{i+1}_{\Qch(Y)}(\Fscr,\Sscr)\r
\end{multline*}
\end{enumerate}
\end{propositions}
\begin{proof} This is proved  in the same way as if $X$ were a
  commutative scheme. As an example let us prove 2.  Choose an
  injective resolution $E^\cdot$ for $\Sscr$ in $\Mod(X)$. Write
  $F^\cdot=\HHom_{o_X}(o_Y,E^\cdot)$. The complex $F^\cdot$ consists
  of injectives in $\Mod(Y)$. Applying the exact functors
  $\HHom_{o_X}(-,E_i)$ to
\[
0\r o_X(-Y)\r o_X\r o_Y\r 0
\]
yields an exact sequence of complexes
\begin{equation}
\label{ref:5.1b}
0\r  F^\cdot \r E^\cdot \r E^\cdot(Y)\r 0
\end{equation}
Applying the long exact sequence for homology to
\eqref{ref:5.1b} we see that the homology of
 $F^\cdot$ is $\Sscr$ and $\Sscr(Y)$ in degrees
0 and 1 respectively. This means that we have a triangle in $\Mod(Y)$
\begin{equation}
\label{ref:5.2a}
\Atriangle<1`-1`1; >[\Sscr(Y)[-1]`
\Sscr`F^\cdot;``]
\end{equation}
On the other hand we have
\begin{align*}
\RHom_{o_X}(\Fscr,\Sscr)&=\Hom_{o_X}(\Fscr,E^\cdot)\\
&=\Hom_{o_Y}(\Fscr,F^\cdot)\\
&=\RHom_{o_Y}(\Fscr,F^\cdot)
\end{align*}
Thus applying $\RHom_{o_Y}(\Fscr,-)$ to the triangle
\eqref{ref:5.2a} we obtain a new triangle
\begin{equation}
\Atriangle<1`-1`1; >[R\Hom_{o_Y}(\Fscr,\Sscr(Y)[-1])`\RHom_{o_Y}(\Fscr,\Sscr)
`\RHom_{o_X}(\Fscr,\Sscr);``]
\end{equation}
The long exact sequence for homology of this triangle is precisely 2.
\end{proof}

Since $\Nscr_{Y/X}$ is an invertible bimodule on $Y$ it follows from
\cite[Prop. 6.8]{AZ} that we have $\Nscr_{Y/X}=\Nscr_\tau$ for some line
bundle $\Nscr$ on $Y$ and an automorphism $\tau$ of $Y$. Recall that
by definition
\begin{equation}
\label{ref:5.4a}
(-\otimes_{o_Y} \Nscr_\tau)=\tau_\ast(-\otimes_{\Oscr_Y}\Nscr)
\end{equation}
 
 By $\Cscr_f$ we denote the finite length objects in $\Qch(X)$ whose
 Jordan-Holder quotients lie in $\Qch(Y)$. By $\Cscr$ we denote the
 corresponding locally finite subcategory of $\Qch(X)$.
 
 If $p\in Y$ then we denote by $\Cscr_{f,p}$ the full subcategory of
 $\Cscr_f$ consisting of objects whose Jordan-Holder quotients are among
 $(O_{\tau^np})_n$. Again $\Cscr_p$ is the corresponding locally closed
 subcategory of $\Cscr$. Clearly $\Cscr_f=\oplus_{p\in
 Y/\langle\tau\rangle}\Cscr_{f,p}$, $\Cscr=\oplus_{p\in
 Y/\langle\tau\rangle}\Cscr_p$. 
 
 From the fact that $o_X(Y)/o_X=\Nscr_\tau$ we deduce that 
 \begin{equation}
 \label{ref:5.5a}
 O_q(Y)\cong O_{\tau q}
 \end{equation}
 In particular $\Cscr_p$ is stable under $-\otimes_{o_X}o_X(Y)$. From
 this one deduces
 \begin{propositions} (\cite[Prop.\ 8.4]{VdBVG})
 \label{ref:5.1.3a}
 $\Cscr_p$ is closed under injective hulls in $\Qch(X)$.
 \end{propositions}
 We now translate (and slightly generalize) the main result of
 \cite{VdBVG} to our situation.
 \begin{theorems}
\label{ref:5.1.4a}
We have the following.
\begin{enumerate}
\item There is a category equivalence $\hat{(-)}_{p}$ between
  $\Cscr_{p}$ and the category $\Dis(C_p)$ (\S\ref{ref:4a}) for a
  certain pseudo-compact ring $C_p$. This ring $C_p$ has the following
  form~:
\begin{enumerate}
\item
If $|O_\tau(p)|=\infty$ then $C_p$ is given by the $\ZZ\times\ZZ$ lower
triangular matrices with entries in $\hat{\Oscr}_{Y,p}$. In this case $p$
is regular on $Y$ and thus
we  have $\hat{\Oscr}_{Y,p}\cong k[[x]]$.
\item If $|O_\tau(p)|=n$ then $C_p$ is given by a ring of $n\times n$
matrices of the form
\[
\begin{pmatrix}
R& RU &\cdots &RU\\
\vdots &\ddots& \ddots&\vdots\\
\vdots &&\ddots &RU\\
R&\cdots&\cdots&R
\end{pmatrix}
\]
where $R$ is a complete local ring of the form
\[
R=k\langle\langle x,y\rangle\rangle/(\psi)
\]
with
\begin{equation}
\label{ref:5.6a}
\psi=yx-qxy+\mathrm{higher\ terms}
\end{equation}
for some $q\in k^\ast$, or
\begin{equation}
\label{ref:5.7a}
\psi=yx-xy-x^2+\mathrm{higher\ terms}
\end{equation}
$U$ is a regular normalizing element in $\rad(R)$ such that
$R/(U)=\hat{\Oscr}_{Y,p}$. If $p$ is not  fixed under $\tau$ then 
$p$ is regular on $Y$
and also $U\not\in\rad^2(R)$.
\end{enumerate}
In all cases $R$ carries the usual topology and $C_p$ carries the
corresponding product topology.
\item
Let $I=\ZZ$ if $|O_\tau(p)|=\infty$ and $I=\ZZ/n\ZZ$ if
$|O_\tau(p)|=n$.  In this way the elements of $C_p$ correspond to
$I\times I$-matrices. For $i\in I$ let $e_i$ be the corresponding 
diagonal idempotent. 
Put
$S_i=e_iC_p/\rad(e_iC_p)$.
Then $
(\Oscr_{\tau^i p})^\wedge_p=S_i$. 
\item Define the following normal element $N$ of $C_p$.
\begin{enumerate} 
\item If $|O_\tau(p)|=\infty$ then $N$ is given by the matrix whose
entries are everywhere zero except on the lower subdiagonal where they
are one.
\item
If $|O_\tau(p)|<\infty$ then
\[
N=
\begin{pmatrix}
0&\cdots & 0 & U\\
1&\ddots && 0\\
\vdots &\ddots & \ddots & \vdots\\
0&\cdots &1 & 0
\end{pmatrix}
\]
\end{enumerate}
Let $\phi=N\cdot N^{-1}$ 
Then we have the following  commutative diagram
\begin{equation}
\label{ref:5.8a}
\begin{CD}
\Cscr_{p} @>-\otimes o_X(-Y)>> \Cscr_{p}\\
@VV\hat{(-)}_{p}V @VV\hat{(-)}_{p}V\\
\Dis(C_{p}) @>(-)_\phi >> \Dis(C_{p})
\end{CD}
\end{equation}
\item If $\Fscr$ is an object in $\Cscr_{p}$ then one has the following
commutative diagram.
\[
\begin{CD}
(\Fscr(-Y))^\wedge_{p} @> >> (\Fscr)^\wedge_{p}\\
@V\cong VV @|
\\
((\Fscr)^\wedge_{p})_\phi @>\cdot N>>(\Fscr)^\wedge_{p}
\end{CD}
\]
where the top arrow is obtained from the inclusion $o_X(-Y)\hookrightarrow
o_X$ and the left arrow from \eqref{ref:5.8a}. 
\item
Let $\Fscr$ be a finite length object in $\Qch(Y)$.  Then $\hat{\Fscr
}_{p}$ is a $C_p/(N)=\prod_{q\in O_\tau(p)}\hat{\Oscr}_{Y,p}$-module
and   
\begin{equation}
\label{ref:5.9a}
\hat{\Fscr
}_{p}=\prod_{q\in O_\tau(p)} \hat{\Fscr}_{Y,q}
\end{equation}
where we have written $\hat{(-)}_{Y,q}$ for the ordinary completion at
$q$ on $Y$.
\end{enumerate}
\end{theorems}
In \cite{VdBVG} this theorem was proved with
$\Cscr_p$ replaced by $\Cscr_{f,p}$. However it is easy to see that one 
obtains the current theorem by taking direct limits.

Below we write $m$ for the maximal ideal of $R$ and $m_i$ will be the
maximal ideal of $C_p$ corresponding to $S_i$. We also use $S_i$ for the
bimodule $C_p/m_i$.

$R$ is clearly noetherian. $C_p$ is noetherian if the orbit of $p$ is
finite, and locally noetherian otherwise.  Furthermore it is also
clear that $C_p$ is cofinite.

The following was proved in \cite[Thm 1.1.4]{VdBVG}.
\begin{lemmas} Every finite dimensional $C_p$-representation $F$
 is in $\PCFin(C_p)$. Hence $F=\prod_i Fe_i$.
\end{lemmas}
The ring $C_p$ has another good
property, which wasn't stated in \cite{VdBVG}.
\begin{propositions}
\label{ref:5.1.6a}
The ring $C_p$ is coherent.
\end{propositions}
\begin{proof}
This is clear if the orbit of $p$ is finite, so we assume it to be
infinite.

We first show that if $B$ is a finitely generated $N$-torsion free right
pseudo-compact $C_p$-module then $B$ is finitely presented.

 We have an exact sequence of pseudo-compact modules
\[
0\r K\r C_p^m \r B\r 0
\]
Tensoring with $C_{Y,p}\overset{\text{def}}{=}C_p/C_pN$ yields an
exact sequence
\[
0\r K/NK \r C_{Y,p}^m \r B/NB\r 0
\]
and using Lemma \ref{ref:4.9a} we see that it is sufficient
to show that $K/NK$ is finitely generated. 
Now $C_{Y,p}=\prod_i \hat{\Oscr}_{Y,\tau^ip}$, and from the theory of
discrete valuation rings we see that the number of generators of a submodule  of
$C_{Y,p}^m e_i=\hat{\Oscr}_{Y,\tau^ip}^m$ is bounded by $m$. This easily
implies what we want. 

Now we prove that $C_p$ is coherent.
 We have to show that the kernel of an arbitrary map $\alpha:
C_p^m\r C_p$ is finitely generated. Clearly $B=\im \alpha$ is
pseudo-compact and finitely generated. We can now apply  the
result in the previous paragraph.
\end{proof}
We will also need~:
\begin{propositions}
\label{ref:5.1.7a}
The global  dimension of $C_p$ is equal to two.
\end{propositions}
\begin{proof}
According to Corollaries \ref{ref:4.11a}  and  \ref{ref:4.7b} it
suffices to show that the projective dimension of each $S_i$ is equal
to two.  

 Put $P_i=e_iC_p$. One easily
checks that the minimal resolution of $S_i$ is given by
\[
0\r P_{i-1}\r P_{i-1}\oplus P_i\r P_i\r S_i\r 0
\]
This implies what we want.
\end{proof}

The following result will be used. 
\begin{lemmas} 
\label{ref:5.1.8a}
Assume that $p$ is a fixed point for $\tau$. Then the
  multiplicity of $p$ on $Y$ is equal to the largest integer $n$ such
  that $RU\subset \rad^n R$.
\end{lemmas}
\begin{proof} Put $S=R/(U)$ and let $m$ be the maximal ideal of
  $R$. By Theorem \ref{ref:5.1.4a}.5 
  $S=\hat{\Oscr}_{Y,p}$.  Assume that $U\in m^\mu-m^{\mu+1}$. 
  Equip $R$ and $S$ with the $m$-adic
  filtration.  Since $\gr R$ is a domain (direct
  verification) it is easy to see that there is an exact sequence
\begin{equation}
\label{ref:5.10a}
0\r \gr R(-\mu)\xrightarrow{\cdot U} \gr R\r \gr S\r 0
\end{equation}
 Then from \eqref{ref:5.10a} we find
that $\mu$ is equal to
   $\dim_k \rad^u S/\rad^{u+1}S$ for large $u$. Hence $\mu$ is the
  multiplicity of $p$ in $Y$.
\end{proof}

\subsection{Some computations}

\leavevmode 

Our view point is that $C_p$ encodes the local structure around a
point $p\in Y$. In order to blow up $p$ we will consequently need some
computations in $C_p$. Our aim here is to prove Corollary
\ref{ref:5.2.4a} below. This corollary is easy if $p$ is a fixed point
for $\tau$ and fairly easy if $p$ has infinite $\tau$-orbit. So the
main purpose will be to treat the case where $n=|O_\tau(p)|$ satisfies
$2\le n<\infty$. However we will develop a formalism which also
includes the case $n=\infty$. Perhaps this has some independent
interest. Our main result will be Proposition \ref{ref:5.2.2a} which
is however more elaborate than what we need for the application to
Corollary \ref{ref:5.2.4a}.

 Notations will be as in Theorem \ref{ref:5.1.4a}.
$m$ will be the maximal ideal of $R$ and $m_i$ will be the twosided
maximal ideal of $C_p$ corresponding to $S_i$. Note that usually $i$
is a taken modulo  $n$ here.
\begin{definitions}
An array $a=(a_q)_{q\in\ZZ}$ with entries in $\NN\cup \{+\infty\}$ will be 
called \emph{good} if it is non-decreasing, bounded below, and if $a_q$ is
 infinite for $t\gg 0$ and finite for $t\ll 0$. If $a,b$ are good then
 $a\ge b$ iff $a_q\ge b_q$ for all $q$.
\end{definitions}
Our aim is to use good arrays as a bookkeeping device in order to
study certain right $C_p$-modules in the case that $n\ge 2$. Let $a$
be a good array and let $I$ be as in Theorem \ref{ref:5.1.4a}. We distinguish two cases.
\begin{itemize}
\item
\textbf{The case {\mathversion{bold} $n< \infty$}. }

Fix an arbitrary $x\in m-m^2$
such that $x-U\not\in m^2$
for $R$.
For $s\in \ZZ$ define 
\[
H_{a,s}=\sum_{t\in\ZZ}  x^{a_{s-nt}}U^t R
\]
where this sum is taken inside the fraction field of $R$ and where by
convention $x^\infty=0$. The $H_{a,s}$ are clearly fractional right $R$
ideals.

Put
 \[
 P_{a}=(H_{a,0},\ldots, H_{a,n-1})
 \]
\item
\textbf{The case {\mathversion{bold} $n=\infty$}. }
We consider this as a limiting case of the previous case.
We define
$P_a$ as the row matrices $(P_{a,j})_j$ with
\[
P_{a,j}=m^{a_j}
\]
\end{itemize}
\begin{propositions}
\label{ref:5.2.2a}
 In this proposition $a,b,c$ will be good
  arrays. As above let $n=|I|\ge 2$. If $i\in I$ then
\[
\tilde{\imath}=
\begin{cases}
\text{the unique element of $\{0,\ldots,n-1\}$ congruent to $i$}
  &\text{if $n<\infty$}\\
i&\text{if $n=\infty$}
\end{cases}
\]
\begin{enumerate}
\item The $P_a$ are right $C_p$-modules (with the obvious $C_p$-action).
\item Let $K$ be the fraction field of $R$ and consider the $P_a$ as
  submodules of  $K^I$. Then 
\[
P_{a}\subset P_{b}\iff a\ge b
\]
\item One also has
\begin{align*}
P_{a}+P_{b}&=P_{\inf(a,b)}\\
\label{ref:5.14a}
P_{a}\cap P_{b}&=P_{\sup(a,b)}
\end{align*}
\item
If $b\ge a$ and
$a_i=\infty\iff b_i=\infty$ 
then  $P_{a}/P_{b}$ has finite length and the composition factors  are
given by
\[
\oplus_{i\in\ZZ} S_{\tilde\imath}^{b_i-a_i}
\]
(with multiplicity).
\item
One has
\[
P_{a}m_l=P_{c}
\]
where
\[
c_q=
\begin{cases}
a_q&\text{if $q\not \cong l$ mod $n$}\\
\min(a_{q+1},a_q+1)&\text{if $q \cong l$ mod $n$}
\end{cases}
\]
\item
One has
\[
\rad P_a=P_c
\]
where
\[
c_q=
\begin{cases}
a_q+1&\text{if $a_q\neq a_{q+1}$}\\
a_q&\text{otherwhise}
\end{cases}
\]
\item
One has
\[
P_{a}/\rad P_a= \bigoplus_{q\in\ZZ,a_q\neq a_{q+1}} S_{\tilde{q}} 
\]
\item
For $i\in I$ write $P_i=e_i C_p$. Then
\[
P_i=P_{c}
\]
where
\[
c=(\ldots,0,\ldots,0,\infty,\ldots)
\]
with the first $\infty$ occuring in position $\tilde\imath+1$.
\end{enumerate}
\end{propositions}
\begin{proof}
All this is fairly easy if $n=\infty$, so we concentrate on the case
$n<\infty$.

It is easy to see that $x^i U^j$ is a topological $k$-basis
for $R$
and from the fact that $U$ is normalizing one
obtains the following alternative form of $H_{a,s}$.
\begin{equation}
\label{ref:5.11a}
H_{a,s}=\prod_{t\in\ZZ}\prod_{i\ge a_{s-nt}} k x^i U^t
\end{equation}
From the fact that $a_q$ is ascending we obtain $H_{a,s}\subset
H_{a,s-1}$. We also have $H_{a,s-n}=H_{a,s}U^{-1}$. This easily
implies 1.

If $a,b$ are good then 
\eqref{ref:5.11a} implies that
\begin{equation}
\label{ref:5.12a}
\forall s:H_{a,s}\subset H_{b,s}\iff a\ge b
\end{equation}
and also
\begin{align}
\label{ref:5.13a}
H_{a,s}+H_{b,s}&=H_{\inf(a,b),s}\\
\label{ref:5.14a}
H_{a,s}\cap H_{b,s}&=H_{\sup(a,b),s}
\end{align}
It is clear that \eqref{ref:5.12a}\eqref{ref:5.13a} and \eqref{ref:5.14a}
imply the corresponding properties for $P_a, P_b$. This proves 2. and 3.

If $a,a'$ are such that
\[
a'_q=\begin{cases} a_q &\text{if $q\neq i$}\\
a_q+1&\text{if $q=i$, $a_i\ne \infty$}
\end{cases}
\]
then it follows from \eqref{ref:5.11a} that
\[
H_{a,s}/H_{a',s}=
\begin{cases} R/m&\text{if $s\cong i$ mod $n$}
\\
0&\text{otherwise}
\end{cases}
\]

 We deduce that 
\[
P_{a}/P_{a'}=S_{\bar\imath }
\]
This yields 4.

We now compute $P_{a}m_l$ for an arbitrary $l\in I$. We find
\begin{equation}
\label{ref:5.15a}
(P_{a}m_l)_j=
\begin{cases}
H_{a,j}&\text{if $j\neq l$}\\
H_{a,j}m+H_{a,j+1}&\text{if $j=l$}
\end{cases}
\end{equation}
To compute the righthand side of \eqref{ref:5.15a} we use the
observation that $U^{-b}xU^b\not\cong U\text{ mod } m^2$ and hence
$m=UR+U^{-b}xU^b R$. Thus
\[
x^aU^bm=x^aU^b(UR+U^{-b}xU^bR)=x^{a+1}U^bR+x^aU^{b+1}R
\]
One now easily obtains 5. Items 6. and 7. are consequences of the fact
that $\rad P_a=\bigcap_l P_a m_l$. Finally 8. is a simple verification
which we leave to the reader.
\end{proof}

\begin{examples}
\label{ref:5.2.3a}
Let $n=\infty$, $a=(\ldots,0,0,1,1,1,3,\infty,\infty,\ldots)$ where the first
$\infty$ occurs in position $2$. Then 
\begin{equation}
\label{ref:5.16a}
P_a/\rad P_a= S_{1}\oplus S_{0}\oplus S_{-3}
\end{equation}
If $n=3$  and $a$ is the same then
\begin{equation}
\label{ref:5.17a}
P_a/\rad P_a= S_{\bar{0}}^2\oplus S_{\bar{1}}
\end{equation}
Note that going from \eqref{ref:5.16a} to \eqref{ref:5.17a} amounts to
introducing periodicity modulo $3$ among the $S_i$. 
\end{examples}
\begin{corollarys}
\label{ref:5.2.4a}
$\dim_k (C_p/m_0m_{-1}m_{-2}\cdots m_{-p+1})=\frac{p(p+1)}{2}$
\end{corollarys}
\begin{proof} This is clear if $n=1$. In the case $n\ge 2$ we use the
  fact that we have $C_p=\prod_{i\in I} P_i$. The corollary now
  follows easily from 5. and 8. of the foregoing proposition.
\end{proof}

\subsection[Completion of objects in {$\coh(X)$}]{Completion of objects in {\mathversion{bold} $\coh(X)$}}
\label{ref:5.3a}
We have already defined $\hat{(-)}_{p}$ on $\Cscr_p$. In a different
direction, it is also possible to extend $\hat{(-)}_{p}$ to $\coh(X)$,
but then it loses some of its good properties. For $\Fscr\in\coh(X)$
define
\[
\hat{\Fscr}_{p}=
\invlim_{\Fscr/\Fscr'\in\Cscr_p} (\Fscr/\Fscr')^\wedge_{p}
\]
\begin{theorems}
\label{ref:5.3.1a}
$(-)^\wedge_{p}$ is an exact functor from $\coh(X)$ to $\PC(C_p)$.
The analogs of 3.,4.,5. of Theorem \ref{ref:5.1.4a} hold for
$\Fscr\in\coh(X)$.  If $\Fscr\in\coh(X)$ then $\hat{\Fscr}_{p}$ is
finitely presented.  Furthermore, $\hat{\Fscr}_{p}$ lies in $\pc(C_p)$
if and only if the intersection of the support of $\Fscr/\Fscr(-Y)\in
\Qch(Y)$ and the orbit of $p$ is finite.
\end{theorems}
\begin{proof}
Let 
\[
0\r \Fscr\xrightarrow{\phi} \Gscr\xrightarrow{\theta} \Hscr\r 0
\]
be an exact sequence in $\coh(X)$. We have to show that 
\[
0\r \invlim_{\Fscr/\Fscr'\in\Cscr_p} (\Fscr/\Fscr')^\wedge_{p}
\r \invlim_{\Gscr/\Gscr'\in\Cscr_p} (\Gscr/\Gscr')^\wedge_{p}
\r \invlim_{\Hscr/\Hscr'\in\Cscr_p} (\Hscr/\Hscr')^\wedge_{p}
\r
0
\]
is exact. Given the exactness of $(-)^\wedge_{p}$ on finite length
objects and the exactness of $\invlim$ on pseudo-compact modules this means
we have to show that \begin{enumerate}
\item
$\theta(\Gscr')_{\Gscr/\Gscr'\in\Cscr_p}$ is cofinal in 
$(\Hscr')_{\Hscr/\Hscr'\in\Cscr_p}$
\item
$\phi^{-1}(\Gscr')_{\Gscr/\Gscr'\in\Cscr_p}$ is cofinal in
$(\Fscr')_{\Fscr/\Fscr'\in\Cscr_p}$.
\end{enumerate}
The first statement is clear. The
second statement is the Artin-Rees condition, which is equivalent
with $\Cscr_p$ being closed under injective hulls in $\Qch(X)$. This is
precisely Proposition \ref{ref:5.1.3a}.
The fact that the analog of Theorem \ref{ref:5.1.4a} holds is
easy to see.

Now we prove that $\hat{\Fscr}_p$ is finitely presented. According to lemma
\ref{ref:4.18a} it is sufficient to show that for every $i$ and every 
$q\in O_\tau(p)$, $\Ext^i_{\PC(C_p)}(\hat{\Fscr}_{p},(\Oscr_q)^\wedge_{p})$ has 
finite dimension,
bounded independently of $q$.  By Proposition \ref{ref:5.3.4a} 
below we have \[
\Ext^i_{\PC(C_p)}(\hat{\Fscr}_{p},(\Oscr_q)^\wedge_{p}))=
\Ext^i_{\Qch(X)}(\Fscr,\Oscr_q)
\]
and the dimension of the righthand side of this equation  is indeed
bounded independently of $q$ by lemma \ref{ref:5.3.2a}.

Now we concentrate on the second part of the theorem.  
We have an exact sequence in $\Qch(X)$
\[
 \Fscr(-Y)\r \Fscr\r \Fscr/\Fscr(-Y)\r 0
 \]
 which yields by the analog of Theorem \ref{ref:5.1.4a}.5 an exact
 sequence in $\PC(C_p)$
\[
(\hat{\Fscr}_{p})_\phi\xrightarrow{\cdot N}  
\hat{\Fscr}_{p} \r \prod_{q\in O_\tau(p)}(\Fscr/\Fscr(-Y))^\wedge_{Y,q}\r 0
\] By \cite[Prop 3.22]{VdBVG} $\hat{\Fscr}_{p}$
will be noetherian if and only if 
$\prod_{q\in O_\tau(p)}(\Fscr/\Fscr(-Y))^\wedge_{Y,q}$ is noetherian,
that is, if and only if, for almost all $q\in O_\tau(p)$ we have 
$(\Fscr/\Fscr(-Y))^\wedge_{Y,q}=0$. 
This is equivalent with the intersection
of the support of $\Fscr/\Fscr(-Y)$ with the $\tau$-orbit of $p$
being finite. \end{proof}

We now supply  the details that were used in the proof of the
previous theorem.
\begin{lemmas} 
\label{ref:5.3.2a}
Assume that $\Fscr\in \coh(X)$, $q\in Y$. Then 
$\dim \Ext^i_{\Qch(X)}(\Fscr,\Oscr_q)$ is finite, and bounded independently
of $q$.
\end{lemmas}
\begin{proof}
As usual, one  reduces to one of the  following cases.
\begin{enumerate}
\item The canonical map $\Fscr(-Y)\r \Fscr$ is injective.
\item $\Fscr$ is in $\coh(Y)$.
\end{enumerate}
Using Proposition \ref{ref:5.1.2a}, we can then reduce to
showing that if $\Gscr\in\coh(Y)$ then $\dim
\Ext^i_{\Qch(Y)}(\Gscr,\Oscr_q)$ is finite and bounded independently of
$q$. This now follows from the fact
that $\Gscr$ has a resolution consisting of vector bundles on $Y$.
\end{proof}
The following result is useful.
\begin{lemmas}\label{ref:5.3.3a} $\hat{\Oscr}_{X,p}\cong C_p$ as right $C_p$-modules.
\end{lemmas}
\begin{proof} The analogous result for $Y$ is trivially true. We can lift
this to $X$ using Nakayama's lemma (lemma \ref{ref:4.8b}).
\end{proof}

We now have to be able to compare $\Ext$-groups in $\PC(C_p)$ and
$\Qch(X)$. In fact we have the following result.
\begin{propositions}
\label{ref:5.3.4a}
Assume that $\Fscr\in\coh(X)$ and $\Sscr\in \Cscr_{p}$. Then
\[
\Ext_{\Top(C_p)}^i(\hat{\Fscr}_{p},\hat{\Sscr}_{p})=
\Ext_{\Qch(X)}^i(\Fscr,\Sscr)
\]
\end{propositions}
\begin{proof} 
We can reduce to the case $i=0$ by replacing $\Sscr$ with an injective
resolution in $\Cscr_p$.
We have
\begin{align*}
\Hom_{\Top(A)}(\hat{\Fscr}_p,\hat{\Sscr}_p)&=
\Hom_{\Top(A)}(\invlim_{\Fscr'}(\Fscr/\Fscr')_p^\wedge,\hat{\Sscr}_p)\\
&=\dirlim_{\Fscr'}\Hom_{\Dis(A)}((\Fscr/\Fscr')_p^\wedge,\hat{\Sscr}_p)
\qquad \text{(Proposition \ref{ref:4.3a})}\\
&=\dirlim_{\Fscr'} \Hom_{\Qch(X)}(\Fscr/\Fscr',\Sscr)\\
&=\Hom_{\Qch(X)}(\Fscr,\Sscr)
\end{align*}
Here $\Fscr'$ runs through the subobjects in $\Fscr$ such that
$\Fscr/\Fscr'\in\Cscr_{f,p}$.
\end{proof}
\begin{corollarys} 
\label{ref:5.3.5a}
Assume $\Fscr\in\coh(X)$. Then $\hat{(-)}_p$
  defines a one-one correspondence between open subobjects of
  $\hat{\Fscr}_p$ and the subobjects $\Fscr'$ of $\Fscr$ such that
  $\Fscr/\Fscr'\in\Cscr_p$. 
\end{corollarys}
\begin{proof}
  Assume for example that $H\subset \hat{\Fscr}_p$ is open. Put
  $S=\hat{\Fscr}_p/ H$. Then $S\in \PCFin(C_p)$. By Theorem
  \ref{ref:5.1.4a} we have $S=\hat{\Sscr}_p$ for some $\Sscr\in
  \Cscr_{f,p}$.
 Let $q:\hat{\Fscr}_p\r
S$ be the quotient map. According to Proposition \ref{ref:5.3.4a}, $q$
corresponds to a map $p:\Fscr\r \Sscr$. Define $\Hscr=\ker p$. Then
$\hat{\Hscr}_p=H$. 
\end{proof}

\subsection{Completion of bimodules}
\label{ref:5.4b}

\leavevmode

If is also possible to define $\hat{(-)}_{p}$ on certain
$o_X-o_X$-bimodules. Unless otherwise specified, when we say ``bimodule''
we mean an object of $\Bimod(o_X-o_X)$.
\begin{definitions} We say that an $o_X-o_X$ bimodule $\Mscr$ is
$\Cscr_{p}$-preserving if both $-\otimes_{o_X}\Mscr$ and
$\HHom_{o_X}(\Mscr,-)$ preserve $\Cscr_{p}$.
\end{definitions}
In this section we will write $(-)^\sim_{p}:\Dis(C_p)\r \Cscr_{p}$
for the inverse of the  functor $(-)^\wedge_{p}$.

Assume that $\Mscr$ is a coherent $\Cscr_{p}$ preserving bimodule. Then
$-\otimes_{o_X} \Mscr$ preserves $\Cscr_{f,p}$ by \ref{ref:3.1.6a}.
We define
\begin{equation}
\label{ref:5.18a}
\hat{\Mscr}_{p}=\invlim_I ((C_p/I)^\sim_{p}\otimes_{o_X}
\Mscr)^\wedge_{p}
\end{equation}
where $I$  runs over the  open twosided ideals in $C_p$.
\begin{propositions} 
\label{ref:5.4.2a}
Assume that $\Mscr$ is a coherent $\Cscr_{p}$
preserving $o_X-o_X$-bimodule. Then
$\hat{\Mscr}_{p}$ is a bipseudo-compact $C_p$-bimodule.
\end{propositions}
\begin{proof}
  If follows from functoriality that $((C_p/I)^\sim_{p}\otimes_{o_X}
  \Mscr)^\wedge_{p}$ is annihilated by $I$ on the left.  On the right
  it is in $\PCFin(C_p)$ and hence it is also annihilated by some open
  twosided ideal $J\subset C_p$. Thus $((C_p/I)^\sim_{p}\otimes_{o_X}
  \Mscr)^\wedge_{p}$ is bipseudo-compact, and hence so is the inverse
  limit.
\end{proof}

We have the
following analog of Proposition \ref{ref:5.3.4a}.
\begin{propositions}
\label{ref:5.4.3a}
Assume that $\Mscr$ is a coherent $\Cscr_{p}$-preserving bimodule and
$\Sscr\in\Cscr_{p}$. Then 
$\Ext^i_{\Top(C_p)}(\hat{\Mscr}_{p},\hat{\Sscr}_{p})\in \Dis(C_p)$ and 
$\HExt^i_{\Qch(X)}(\Mscr,\Sscr)\in \Cscr_p$. 
Furthermore as
objects of $\Dis(C_p)$
\[
\Ext^i_{\Top(C_p)}(\hat{\Mscr}_{p},\hat{\Sscr}_{p})
=
\HExt^i_{\Qch(X)}(\Mscr,\Sscr)^\wedge_{p}
\]
\end{propositions}
\begin{proof} 
That $\HExt^i_{\Qch(X)}(\Mscr,\Sscr)\in \Cscr_p$ follows from
the fact that $\HHom(\Mscr,-)$ is $\Cscr_{p}$ preserving and
Proposition \ref{ref:5.1.3a}. The fact that 
$\Ext^i_{\Top(C_p)}(\hat{\Mscr}_{p},\hat{\Sscr}_{p})\in \Dis(C_p)$ is
precisely Proposition \ref{ref:4.19a}.

Below $\Uscr$ will run through $(C_p/I)^\sim_p$, where $I$
is an open twosided ideal in $C_p$.
We have
\begin{equation}
\label{ref:5.19a}
\begin{split}
\Ext^i_{\Top(C_p)}(\hat{\Mscr}_{p},\hat{\Sscr}_{p})&=
\Ext^i_{\Top(C_p)}(\invlim_\Uscr (\Uscr\otimes_{o_X}
\Mscr)^\wedge_{p},\hat{\Sscr}_{p})\\
&=
\dirlim_\Uscr \Ext^i_{\Dis(C_p)}((\Uscr\otimes_{o_X}
\Mscr)^\wedge_{p}, \hat{\Sscr}_{p})\qquad \text{(Proposition
\ref{ref:4.3a})}
\\ &=\dirlim_\Uscr \Ext^i_{\Cscr_p}(\Uscr\otimes_{o_X}\Mscr,\Sscr)
\qquad\text{(Theorem
\ref{ref:5.1.4a})}\\
&=\dirlim_\Uscr \Ext^i_{\Qch(X)}(\Uscr\otimes_{o_X}\Mscr,\Sscr)
\quad\text{(Proposition
\ref{ref:5.1.3a})}\end{split}
\end{equation}
Let $E^\cdot$ be an injective resolution of $\Sscr$ in $\Cscr_{p}$.
Then
\begin{align*}
\Ext^i_{\Qch(X)}(\Uscr\otimes_{o_X}\Mscr,\Sscr)&
=H^i(\Hom_{\Qch(X)}(\Uscr\otimes_{o_X} \Mscr,E^\cdot)\\
&=H^i(\Hom_{\Qch(X)}(\Uscr,\HHom_{\Qch(X)}(\Mscr,E^\cdot))\\
&=
H^i(\Hom_{\Dis(C_p)}(\hat{\Uscr}_p,\HHom_{\Qch(X)}(\Mscr,E^\cdot)^\wedge_p )
\\
&=H^i(\Hom_{\Dis(C_p)}(C_p/I,\HHom_{\Qch(X)}(\Mscr,E^\cdot)^\wedge_p )
\end{align*}
Combining this with \eqref{ref:5.19a} yields
\begin{align*}
\Ext^i_{\Top(C_p)}(\hat{\Mscr}_{p},\hat{\Sscr}_{p})
&=
\dirlim_I 
H^i(\Hom_{\Dis(C_p)}(C_p/I,
\HHom_{\Qch(X)}(\Mscr,E^\cdot)^\wedge_p )\\
&=H^i(\Hom_{\Top(C_p)}(\invlim_I(C_p/I),
\HHom_{\Qch(X)}(\Mscr,E^\cdot)^\wedge_p ))\\
&=H^i(\Hom_{\Top(C_p)}(C_p,\HHom_{\Qch(X)}(\Mscr,E^\cdot)^\wedge_p ))\\
&=\HExt^i_{\Qch(X)}(\Mscr,\Sscr)^\wedge_p 
\qed\end{align*}
\def\qed{}\end{proof}

\begin{propositions} 
\label{ref:5.4.4a}
The functor $(-)_{p}^\wedge$ preserves short
exact sequences of 
coherent $\Cscr_{p}$-preserving bimodules.
\end{propositions}
\begin{proof}
Let 
\[
0\r \Kscr\r \Mscr\r \Nscr\r 0
\]
be an exact sequence of coherent $\Cscr_{p}$ preserving bimodules. We
have to show that 
\begin{equation}
\label{ref:5.20a}
0\r \hat{\Kscr}_p\r \hat{\Mscr}_p \r \hat{\Nscr}_p\r 0
\end{equation}
is exact. Since the bimodules in \eqref{ref:5.20a} are in $\PC(C_p)$
(Proposition \ref{ref:5.4.2a}) it suffices to show that 
\begin{equation}
\label{ref:5.21a}
0\r \Hom_{\Top(C_p)}(\hat{\Nscr}_p,\hat{E}_p)\r
\Hom_{\Top(C_p)}(\hat{\Mscr}_p,\hat{E}_p)\r
\Hom_{\Top(C_p)}(\hat{\Kscr}_p,\hat{E}_p)\r
0
\end{equation}
is exact where $E$ is the sum of the injective hulls of the $\Oscr_{\tau^i
p}$ (one uses the fact that $\Hom_{\Top(C_p)}(-,\hat{E}_p)$ is exact
and faithful on $\PC(C_p)$). Now by Proposition \ref{ref:5.4.3a} it
follows that \eqref{ref:5.21a} is obtained from completing the
following exact sequence in $\Cscr_p$.
\[
0\r \HHom_{\Qch(X)}(\Nscr,E)\r 
\HHom_{\Qch(X)}(\Mscr,E)\r 
\HHom_{\Qch(X)}(\Kscr,E)\r 
0
\]
and hence we are done!
\end{proof}
\begin{propositions} 
\label{ref:5.4.5a}
Completion commutes with tensor product in the
following sense. Let $\Fscr$ be a coherent object in $\Qch(X)$ and let
$\Mscr$, $\Nscr$ be coherent $\Cscr_{p}$ preserving bimodules. Then
there are natural isomorphisms
\begin{enumerate} 
\item $(\Fscr\otimes_{o_X} \Nscr)^\wedge_{p}=\hat{\Fscr}_{p}\ctimes_{C_p}
  \hat{\Nscr}_{p}$ (note that by Theorem \ref{ref:5.3.1a} and lemma
  \ref{ref:4.14a} we may replace ``$\ctimes$'' by ``$\otimes$'').
 \item $(\Mscr\otimes_{o_X}
 \Nscr)^\wedge_{p}=\hat{\Mscr}_{p}\ctimes_{C_p} 
 \hat{\Nscr}_{p}$
 \end{enumerate}
 \end{propositions}
\begin{proof}
\begin{enumerate}
\item
By definition we have
\[
(\Fscr\otimes_{o_X} \Nscr)^\wedge_{p}=\invlim_\Tscr(
(\Fscr\otimes_{o_X} \Nscr) /\Tscr)^\wedge_p
\]
where $\Tscr$ runs through the subobjects of $\Fscr\otimes_{o_X} \Nscr$ such
that $(\Fscr\otimes_{o_X}\Nscr)\Tscr\in \Cscr_p$. Now we claim that for
every $\Tscr$ there exists a $\Fscr'\subset \Fscr$ such that
$\Fscr/\Fscr'\in\Cscr_p$ and such that the image
of $\Fscr'\otimes_{o_X} \Nscr$ in $\Fscr\otimes_{o_X} \Nscr$ is contained in
$\Tscr$.

Let $\Qscr=(\Fscr\otimes_{o_X} \Nscr) /\Tscr$ and let
$\Fscr\otimes_{o_X}\Nscr\r \Qscr$ be the quotient map. By adjointness
there is a corresponding map $ \Fscr\r \HHom(\Nscr,\Qscr)$. We define
$\Fscr'$ as the kernel of this map.  This $\Fscr'$ has the properties
we want.

We obtain
\begin{equation}
\label{ref:5.22a}
(\Fscr\otimes_{o_X} \Nscr)^\wedge_{p}=\invlim_{\Fscr'}
((\Fscr/\Fscr')\otimes_{o_X} \Nscr)^\wedge_p
\end{equation}
For $K\in \PC(C_p)$ we define the functor 
\[
F(K)=\invlim_{K'} ((K/K')^\sim_p \otimes_{o_X} \Nscr)^\wedge_p
\]
where $K'$ runs through the open subobjects of $K$. One shows  that
this functor is right exact and commutes with products. Furthermore,
according to \eqref{ref:5.18a} we
have $F(C_p)=\hat{\Nscr}_p$. By 
 \eqref{ref:5.22a} and Corollary \ref{ref:5.3.5a}
\begin{align*}
F(\hat{\Fscr}_p)&=\invlim ((\Fscr/\Fscr')\otimes_{o_X} \Nscr)^\wedge_p\\
&=(\Fscr\otimes_{o_X}\Nscr)^\wedge_p
\end{align*}
Now by Theorem \ref{ref:5.3.1a} we know that $\hat{\Fscr}_p$ is coherent. Take a
presentation
\[
C^n_p\r C_p^m\r \hat{\Fscr}_p\r 0
\]
By right exactness of tensor product we find that the cokernel of
$(C^n_p\r C_p^m)\otimes_{C_p} \hat{\Nscr}_p$ is equal to
$\hat{\Fscr}_p\otimes_{o_X} \hat{\Nscr}_p$. On the other hand by right
exactness of $F$ we find that this cokernel is equal to
$F(\hat{\Fscr}_p)=(\Fscr\otimes\Nscr)^\wedge_p$.
\item
This follows from 1. We have
\begin{align*}
(\Mscr\otimes\Nscr)^\wedge_p&=\invlim_I ((C_p/I)^\sim_p\otimes_{o_X}
\Mscr\otimes_{o_X}\Nscr)^\wedge_p\\
&=\invlim_I ((C_p/I)_p^\sim\otimes_{o_X} \Mscr)^\wedge_p\ctimes_{C_p} \hat{\Nscr}_p\\
&=\invlim_I ((C_p/I)\ctimes_{C_p} \hat{\Mscr}_p\ctimes_{C_p} \hat{\Nscr}_p\\
&= \hat{\Mscr}_p\ctimes_{C_p} \hat{\Nscr}_p
\end{align*}
where we have used the fact that $\ctimes$ commutes with inverse limits.
\qed\end{enumerate}
\def\qed{}
\end{proof}

\subsection[The category $\tilde{\Cscr}_{f,p}$]{The category 
{\boldmath $\tilde{\Cscr}_{f,p}$}}
\label{ref:5.5b}

\leavevmode

If $\Gscr\in\coh(Y)$ then we can view $\Gscr$ as a coherent
$\Oscr_Y$-bimodule, by defining $\HHom_{o_Y}(\Gscr,-)$ as
$\HHom_{\Oscr_Y}(\Gscr,-)$ where the second ``$\HHom$'' is the
ordinary $\HHom$ for sheaves.  $\HHom_{\Oscr_Y}(\Gscr,-)$ has a left adjoint
given by $-\otimes_{\Oscr_Y} \Gscr$ so
$\Gscr\in\Bimod(o_Y-o_Y)$. Hence if $q\in Y$ then we can view $o_q$
as a coherent object in $\Bimod(o_Y-o_Y)$. We also write $o_q$
for the coherent object in $\Bimod(o_X-o_X)$ given by ${}_{o_X} o_q
{}_{o_X}$ (cfr. \eqref{ref:3.12a} for notations).
\begin{lemmas} 
\label{ref:5.5.1a}
\begin{enumerate}
\item Assume that $\Fscr\in\Qch(X)$. Then 
\begin{equation}
\label{ref:5.23a}
\HHom_{\Qch(X)}(o_q,\Fscr)=\Oscr_q\otimes_k \Hom_{\Qch(X)}(\Oscr_q,\Fscr)
\end{equation}
\item
If $\Fscr\in \coh(X)$ then
\begin{equation}
\label{ref:5.24a}
\Fscr\otimes_{o_X} o_q=\Oscr_q\otimes_k \Hom_{\Qch(X)}(\Fscr,\Oscr_q)^\ast
\end{equation}
\item
We have a canonical identification
\begin{equation}
\label{ref:5.25a}
o_q\otimes_{o_X} o_X(Y)=o_X(Y)\otimes_{o_X} o_{\tau q}
\qed\end{equation}
\end{enumerate}
\end{lemmas}
\begin{proof}
\begin{enumerate}
\item Let $\Fscr\in \Qch(X)$. 
According to  \eqref{ref:3.13a} we have
\[
\HHom_{\Qch(X)}(o_q,\Fscr)=\HHom_{\Qch(Y)}(o_q,i^!(\Fscr))_{o_X}
\]
and by adjunction we also have
\[
\Hom_{\Qch(X)}(\Oscr_q,\Fscr)=\Hom_{\Qch(Y)}(\Oscr_q,i^!(\Fscr))
\]
Now we use the fact that  the analog of \eqref{ref:5.23a}
 holds on $Y$.
\item
This is also proved by reduction to $\coh(Y)$, using
\[
\Fscr\otimes_{o_X} o_q=i_\ast(i^\ast(\Fscr)\otimes_{o_Y} o_q)\qquad
\text{(see \eqref{ref:3.14a})}
\]
\item
It suffices to show that the left and righthand side of \eqref{ref:5.25a}
take the same values on arbitrary $\Fscr\in \Qch(X)$.

We have
\begin{equation}
\begin{split}
\label{ref:5.26a}
\HHom_{\Qch(X)}(o_q\otimes_{o_X} o_X(Y),\Fscr)&=
\HHom_{\Qch(X)}(o_q,\Fscr(-Y))\\
&=\Oscr_q\otimes_k \Hom(\Oscr_q,\Fscr(-Y))\\
&=\Oscr_q\otimes_k \Hom(\Oscr_q(Y),\Fscr))
\end{split}
\end{equation}
and
\begin{equation}
\label{ref:5.27a}
\begin{split}
\HHom_{\Qch(X)}(o_X(Y)\otimes_{o_X} o_{\tau q},\Fscr)
&=\HHom_{\Qch(X)}(o_X(Y),\HHom_{\Qch(X)}(o_{\tau q},\Fscr))\\
&=(\Oscr_{\tau q}\otimes_k \Hom_{\Qch(X)}(\Oscr_{\tau q},\Fscr))(-Y)\\
&=\Oscr_{\tau q}(-Y)\otimes_k \Hom_{\Qch(X)}(\Oscr_{\tau q},\Fscr)
\end{split}
\end{equation}
Now by our conventions $\Oscr_q(Y)\cong \Oscr_{\tau q}$ (see
\eqref{ref:5.5a}). From this it follows that the righthand sides in
\eqref{ref:5.26a} and \eqref{ref:5.27a} are isomorphic.  The reason
that this is canonical is that we use the (non-canonical) identification
$\Oscr_q(Y)\cong\Oscr_{\tau q}$ twice, in such a way that the
ambiguities cancel.\qed
\end{enumerate}
\def\qed{} \end{proof}
Let $\cohBIMOD(o_x-o_X)$ be the full subcategory of $\BIMOD(o_X-o_X)$
consisting of coherent objects. According to Corollary \ref{ref:3.1.8a} this
is an abelian subcategory of $\BIMOD(o_X-o_X)$, closed under extensions.
\begin{propositions} $o_q$ is a simple object in 
$\cohBIMOD(o_X-o_X)$.
\end{propositions}
\begin{proof} An object in $\cohBIMOD(o_X-o_X)$ is determined by the
values it takes on indecomposable injectives.
We have by lemma \ref{ref:5.5.1a}.1
\[
\HHom_{\Qch(X)}(o_q,E)=\begin{cases} \Oscr_q&\text{if $E$ is the
    injective hull of $\Oscr_q$}\\ 0&\text{otherwise}
\end{cases}
\]
Hence a proper left exact subfunctor of $\HHom_{\Qch(X)}(o_q,-)$
commuting with direct sums will be the zero-functor.
\end{proof}
From this proposition we easily deduce
\begin{corollarys} The objects $o_q\otimes_{o_X} o_X(nY)$ with $n\in\ZZ$ are
simple objects in $\cohBIMOD(o_X-o_X)$.
\end{corollarys}
Now we define $\tilde{\Cscr}_{f,p}$ as the full subcategory of
$\cohBIMOD(o_X-o_X)$ consisting of finite length objects whose
Jordan-Holder quotients are of the form $o_{\tau ^m p}\otimes o_X(nY)$, $m,n\in \ZZ$. 
\begin{propositions} 
\label{ref:5.5.4a}
Assuming that $\Cscr$, $\Dscr$ are abelian
categories possessing an injective cogenerator. Assume that we have an
exact sequence in $\BIMOD(\Cscr-\Dscr)$
\[
0\r \Kscr\r \Mscr\r \Nscr\r 0
\]
with $\Nscr\in\Bimod(\Cscr-\Dscr)$, such that $R^1\prod_{i\in I}
\HHom_\Dscr(\Nscr,E_i)=0$ (see the discussion after Proposition
\ref{ref:3.5.5a}) for all families of injectives $(E_i)_{i\in I}$ in
$\Dscr$.  Then $\Kscr\in \Bimod(\Cscr-\Dscr)$ if and only if
$\Mscr\in\Bimod(\Cscr-\Dscr)$.
\end{propositions}
\begin{proof}
We have to show that $\HHom_\Dscr(\Kscr,-)$ commutes products when
evaluated on injectives, if and only if the same holds for
$\HHom_\Dscr(\Mscr,-)$. Using the hypotheses this is a
direct consequence of the five-lemma.
\end{proof}
\begin{propositions}
\label{ref:5.5.5a}
Assume that $\Sscr\in \tilde{\Cscr}_{f,p}$. Then 
$R^1\prod_{i\in I}
\HHom_{\Qch(X)}(\Sscr,E_i)=0$  for all families of 
injectives $(E_i)_{i\in I}$ in $\Qch(X)$.
\end{propositions}
\begin{proof}
 Using the long exact sequence for
$R\prod_{i\in I}$ one sees that it is sufficient that  $R^1\prod_{i\in I}
\HHom_{\Qch(X)}(o_q,E_i) =0$ for $q\in Y$ and for all
families of injectives in $\Qch(X)$.

If $\Mscr\in\Qch(X)$ then $\HHom_{\Qch(X)}(o_q,\Mscr)\in \Qch(q)$.
By Definition \ref{ref:3.5.9a} it now suffices to show that $q$ is
very well closed in $X$. Since $q$ is obviously defined by an ideal in
$o_Y$, we have that $q$ is well closed in $Y$ by Proposition
\ref{ref:3.5.7a}.  $\Qch(q)$ is of course equivalent with the category
of $k$-vectorspaces.  In particular $\Qch(q)$ has exact direct
products and hence $q$ is very well closed in $Y$ by Corollary
\ref{ref:3.5.11a}.

 Furthermore $Y$ is also defined by an ideal inside $X$ and hence  by
 Proposition \ref{ref:3.5.7a} it follows that $Y$
is well closed in $X$. Now we invoke Proposition \ref{ref:3.5.10a} which
yields that $q$ is very well closed in $X$.
\end{proof}
\begin{corollarys}
\label{ref:5.5.6a}
If we have an exact sequence in $\BIMOD(o_X-o_X)$
\[
0\r \Mscr\r \Nscr\r \Sscr\r 0
\]
with $\Sscr\in\tilde{\Cscr}_{f,p}$ then $\Mscr\in \Bimod(o_X-o_X)$ if and
only if $\Nscr \in \Bimod(o_X-o_X)$.
\end{corollarys}
\begin{propositions}
$\tilde{\Cscr}_{f,p}\subset \Bimod(o_X-o_X)$.
\end{propositions}
\begin{proof}
 Let 
$\Sscr\in \tilde{\Cscr}_{f,p}
$. 
By induction there is a
short exact sequence in $\tilde{\Cscr}_{f,p}$ 
\[
0\r \Sscr_1\r \Sscr\r \Sscr_2\r 0
\]
such that $S_i\in \Bimod(o_X-o_X)$. It now suffices to apply Corollary
\ref{ref:5.5.6a} above.
\end{proof}
Let us denote by $\PCFin(C_p-C_p)$ the full subcategory of
$\PC(C_p-C_p)$ consisting of finite length objects. Thus
$\PCFin(C_p-C_p)=\PCFin(C_p^\circ\ctimes C_p)$. Recall that
$C_p/\rad(C_p)=k^I$ where $I$ is as in Theorem \ref{ref:5.1.4a}. Thus
the twosided bipseudo-compact maximal ideals in $C_p$ are naturally
indexed by $I$. The corresponding simple modules were denoted by $S_i$
in Theorem \ref{ref:5.1.4a}. Below we will view the $S_i$ as
$C_p$-bimodules.

It is clear that $(C_p)^\circ \ctimes C_p$ has a normalizing sequence
given by $N\ctimes 1$ and $1\ctimes N$ and the quotient is given by
$(C_p/(N))^\circ\ctimes (C_p/(N))$ which is locally noetherian. It
follows from \cite[Prop. 3.23]{VdBVG} that $(C_p)^\circ\ctimes C_p$ is
locally noetherian. Hence in particular the category $\PCFin(C_p-C_p)$ is
given by the category of $C_p$-bimodules which are finite extensions of
the $S_i$.

The functor $\hat{(-)}_p$ is defined on 
$\tilde{\Cscr}_{f,p}$. We have
\begin{lemmas}
\label{ref:5.5.8a}
$(o_{\tau^i p})^\wedge_p=S_i$
\end{lemmas}
\begin{proof}
We use \eqref{ref:5.24a}. Let $q=\tau^i p$.
\begin{align*}
(o_q)^\wedge_p&=\invlim_I ((C_p/I)^\sim_p\otimes_{C_p} o_q)^\wedge_p\\
&=\invlim_I((\Oscr_q)^\wedge_p\otimes_k  \Hom_{\Qch(X)}((C_p/I)^\sim_p
,\Oscr_q)^\ast)\qquad (\text{lemma \ref{ref:5.5.1a}})      \\
&=\invlim_I(S_i\otimes_k  \Hom_{\PC(C_p)}(C_p/I,S_i)^\ast)\\
&=S_i
\end{align*}
where as usual $I$ runs through the open twosided ideals in $C_p$.
\end{proof}

We have the following result.
 \begin{theorems}
\label{ref:5.5.9a}
 The functor $\hat{(-)}_p$ defines an equivalence between
$\tilde{C}_{p,f}$ and $\PCFin(C_p-C_p)$.
\end{theorems}
\begin{proof} First note that thanks to lemma \ref{ref:5.5.8a} and
  Proposition \ref{ref:5.4.4a} the
image of $\tilde{C}_{p,f}$ under $(-)^\wedge_p$ is contained in
$\PCFin(C_p-C_p)$. To show that this is actually an equivalence we will
construct an inverse.

Assume $U\in \PCFin(C_p-C_p)$ and $\Fscr\in \coh(X)$. Define
\begin{equation}
\label{ref:5.28a}
T(\Fscr)=(\hat{\Fscr}_p\otimes_{C_p} U)^\sim_p
\end{equation}
Since $\hat{\Fscr}_p$ is finitely generated (Theorem \ref{ref:5.3.1a}) we find
that $T$ defines a right exact additive functor $\coh(X)\r
\Cscr_{f,p}\subset \coh(X)$. Hence this functor extends to a functor
$\tilde{T}:\Qch(X)\r \Qch(X)$ commuting with colimits. We denote by
$\tilde{U}_p$ the right adjoint to this functor. If we view $\tilde{U}_p$
as an object in $\Bimod(o_X-o_X)$ then
 notationally we have  \begin{equation}
\label{ref:5.29a}
\Fscr\otimes_{o_X}
\tilde{U}_p=T(\Fscr)=(\hat{\Fscr}_p\otimes_{C_p} U)^\sim_p 
\end{equation}

We now show that $\tilde{(-)}_p$ is a left inverse to $\hat{(-)}_p$. Let
$\Sscr\in \Cscr_{f,p}$. To show that $(\hat{\Sscr}_p)^\sim_p=\Sscr$ it
suffices to construct a natural isomorphism
\[
\Fscr\otimes_{o_X} (\hat{\Sscr}_p)^\sim_p=\Fscr\otimes_{o_X}\Sscr
\]
where $\Fscr\in\coh(X)$. With the help of \eqref{ref:5.29a} we compute
\begin{align*}
\Fscr\otimes_{o_X} (\hat{\Sscr}_p)^\sim_p&=
(\hat{\Fscr}_p\otimes_{C_p} \hat{\Sscr}_p)^\sim_p
\\
&=((\Fscr\otimes_{o_X} \Sscr)^\wedge_p)^\sim_p\qquad\text{(Proposition
  \ref{ref:5.4.5a})} \\
&=\Fscr\otimes_{o_X} \Sscr
\end{align*}
Thus the composition
\begin{equation}
\label{ref:5.30a}
\tilde{\Cscr}_{f,p}\xrightarrow{\hat{(-)}_p}
\PCFin(C_p-C_p)\xrightarrow{\tilde{(-)}_p}
\Bimod(o_X-o_X)
\end{equation}
is the identity on $\tilde{\Cscr}_{f,p}$.  We now claim that the
essential image of $\tilde{(-)}_p$ is contained in
$\tilde{\Cscr}_{f,p}$. Since $\tilde{\Cscr}_{f,p}$ is closed under
extensions in $\Bimod(o_X-o_X)$ and since by \eqref{ref:5.29a}
$\tilde{(-)}_p$ is at least right exact (it is of course exact) it
suffices to show that $\tilde{S}_{i,p}\in \tilde{\Cscr}_{f,p}$. This
follows from the fact that $(o_{\tau^i p})^\wedge_p=S_i$ by lemma
\ref{ref:5.5.8a}, whence $\tilde{S}_{i,p}=o_{\tau^i p}$.

Now we show that $\tilde{(-)}_p$ is also a right inverse to
$\hat{(-)}_p$. Let $S\in \PCFin(C_p-C_p)$. We  show that
$(\tilde{S}_p)^\wedge_p=S$. By definition
\begin{align*}
(\tilde{S}_p)^\wedge_p&=\invlim_I
 ((C_p/I)^\sim_p\otimes_{o_X}
\tilde{S}_p)^\wedge_p\\
&=\invlim_I ((((C_p/I)^\sim_p)^\wedge_p\otimes_{C_p}
S)^\sim_p)^\wedge_p\qquad\text{(eq. \eqref{ref:5.29a})}\\
&=\invlim (C_p/I)\otimes_{C_p} S\\
&=S\qed
\end{align*}
\def\qed{}\end{proof}

To close this section we discuss $\HTor$ between objects of
$\tilde{\Cscr}_{f,p}$.
\begin{theorems}
\label{ref:5.5.10a}
Let $\Fscr\in \coh(X)$, $\Tscr\in \tilde{\Cscr}_{f,p}$ and let $\Sscr$
be a coherent $\Cscr_p$-preserving $o_X$-bimodule such that
$\hat{\Sscr}_p$ is finitely presented on the right. Then
\begin{align}
\label{ref:5.31a}
\HTor_i^{o_X}(\Fscr,\Tscr)&\in \Cscr_{f,p}\\
\label{ref:5.32a}
\HTor_i^{o_X}(\Sscr,\Tscr)&\in\tilde{\Cscr}_{f,p}\\
\label{ref:5.33a}
\Tor^{\PC(C_p)}_i(\hat{\Fscr}_p,\hat{\Tscr}_p)&\in \PCFin(C_p)\\
\label{ref:5.34a}
\Tor^{\PC(C_p)}_i(\hat{\Sscr}_p,\hat{\Tscr}_p)&\in \PCFin(C_p-C_p)
\end{align}
Furthermore we have
\begin{align}
\label{ref:5.35a}
\HTor_i^{o_X}(\Fscr,\Tscr)^\wedge_p&
=\Tor^{\PC(C_p)}_i(\hat{\Fscr}_p,\hat{\Tscr}_p)\\ 
\label{ref:5.36a}
\HTor_i^{o_X}(\Sscr,\Tscr)^\wedge_p&
=\Tor^{\PC(C_p)}_i(\hat{\Sscr}_p,\hat{\Tscr}_p)
\end{align}
\end{theorems} 
\begin{proof}
  By Proposition \ref{ref:5.1.6a} and Theorem \ref{ref:5.3.1a}
  $\hat{\Fscr}_p$ has   a resolution by free modules of finite rank. This implies
  \eqref{ref:5.33a}. From the hypotheses we also find that
  $\Tor^{\PC(C_p)}_i(\hat{\Sscr}_p,\hat{\Tscr}_p)$ lies in
  $\PCFin(C_p)$.  Since it also lies in $\PC(C_p-C_p)$, it must
  necessarily lie in $\PCFin(C_p-C_p)$. This proves \eqref{ref:5.34a}.

To prove the other statements of the theorem it suffices to prove 
\begin{align}
\label{ref:5.37a}
\HTor_i^{o_X}(\Fscr,\Tscr)
&=\Tor^{\PC(C_p)}_i(\hat{\Fscr}_p,\hat{\Tscr}_p)^\sim_p\\
\label{ref:5.38a}
\HTor_i^{o_X}(\Sscr,\Tscr)
&=\Tor^{\PC(C_p)}_i(\hat{\Sscr}_p,\hat{\Tscr}_p)^\sim_p
\end{align}
Recall that to prove these equations, we have to show that they represent
the same functor on injectives.

We'll only prove \eqref{ref:5.38a}. \eqref{ref:5.37a} is similar.

Let $E$ be an injective object in $\Qch(X)$. Let $E_1$ be the $\Cscr_p$
part of $E$. We have 
\begin{align*}
\HHom_{\Qch(X)}(\Tor^{\PC(C_p)}_i&
(\hat{\Sscr}_p,\hat{\Tscr}_p)^\sim_p,E)\\
&=\HHom_{\Qch(X)}(\Tor^{\PC(C_p)}_i(\hat{\Sscr}_p,
\hat{\Tscr}_p)^\sim_p,E_1)
\\
&=\Hom_{\Top(C_p)}
(\Tor^{\PC(C_p)}_i(\hat{\Sscr}_p,\hat{\Tscr}_p),\hat{E}_{1,p})^\sim_p
\qquad\text{(Prop. \ref{ref:5.4.3a})}
\\
&=\Ext^i_{\Top(C_p)}(\hat{\Sscr}^p,\Hom_{\Dis(C_p)}(\hat{\Tscr}_p
,\hat{E}_{1,p}))^\sim_p\qquad \text{(lemma \ref{ref:4.20a})}\\
&=\HExt^i_{\Qch(X)}(\Sscr,\HHom_{\Qch(X)}(\Tscr,E_1))\qquad
\text{(Prop. \ref{ref:5.4.3a})}\\
&=\HExt^i_{\Qch(X)}(\Sscr,\HHom_{\Qch(X)}(\Tscr,E))\\
&=\HHom_{\Qch(X)}(\HTor_i(\Sscr,\Tscr),E)\qquad 
\text{(eq. \eqref{ref:3.1b})}
\end{align*}
Done!
\end{proof}

\subsection{Completion of algebras}
In the sequel we will need the following concept.
\begin{definitions} 
\label{ref:5.6.1a}
A topological $\ZZ$-graded ring $A$ is said to be graded cofinite
if $A_0$ is cofinite and if the homogeneous parts of $A$ are
bipseudo-compact $A_0$-modules.
\end{definitions}
If $A$ is as above then we denote by $\GrTop(A)$ the category of
graded right topological $A$-modules. With $\GrDis(A)$, resp.
$\GrPC(A)$ we denote the full subcategories consisting of objects
$\oplus_m M_m$ such that $M_m$ is in $\Dis(A_0)$, resp. $\PC(A_0)$.

We say that $A$ is (right) locally noetherian if for every indecomposable
idempotent in $A_0$ we have that $e_iA$ is noetherian in $\GrPC(A)$. Left
locally noetherian is defined similarly.

If $A$ is locally
noetherian then it easily follows that $\GrDis(A)$ is locally noetherian.
 By analogy with the definition
of $\Proj$ of an ordinary graded ring we define $\Tors(A)$ as the full
subcategory of $\GrDis(A)$ consisting of objects which are direct
limits of right bounded ones. We then put
$\QGrDis(A)=\GrDis(A)/\Tors(A)$.

Now let $X,Y,p$ etc\dots be as in the previous sections. Let
$\Ascr=\oplus_n \Ascr_n$ be  in $\Gralg(X)$ and
assume that the $\Ascr_n$ are coherent and $\Cscr_p$ preserving. Then
it makes sense to define $\hat{\Ascr}_p$ as $\oplus_n
\hat{\Ascr}_p$. The compatibility of $\hat{(-)}_p$ with tensor product
implies that $\hat{\Ascr}_p$ is a cofinite algebra.

As in \S\ref{ref:3.3b} we define $\Cscr_{p}(\Ascr)$ as the graded
$\Ascr$-modules $\Mscr$ such $\Mscr_n\in\Cscr_p$ for all $n$.
Similarly $Q\Cscr_p(\Ascr)$ is the image of $\Cscr_p(\Ascr)$ in
$\Gr(\Ascr)$.

 It is easy to prove the following. 
\begin{propositions} 
\label{ref:5.6.2a}
The functor $\hat{(-)}_p$ defines an  equivalence
  between $\Cscr_p(\Ascr)$ and $\GrDis(\hat{\Ascr}_p)$ and between
  $Q\Cscr_p(\Ascr)$ and $\QGrDis(\hat{\Ascr}_p)$.
\end{propositions}
\begin{proof} This follows easily from 
  compatibility of completion with tensor product (Proposition
  \ref{ref:5.4.5a}). 
\end{proof}
 
\subsection{Multiplicities in the case that $\tau$ has infinite order} \label{ref:5.7b}
This section is somewhat more special than the previous ones. Our aim
is to use the completion functor to attach to certain
$\Fscr\in\mod(X)$ an invariant $T_p(\Fscr)$ which, in some sense, is
an analogue for the multiplicities of $\Fscr$ in $p$ and in the points
infinitely near to $p$. 
  A straightforward treatment seems
only to be possible in the case that the $\tau$-orbit of $p$ is
infinite, \emph{so we assume this throughout this section}. Note that
 this implies in particular that $p$ is
smooth on $Y$.

To simplify the notations we will write $\Fscr_Y$ for
$i^\ast\Fscr=\Fscr/\Fscr(-Y)$.  We will say that $\Fscr\in\mod(X)$ is
$Y$-transversal if $\Fscr_Y$ has finite length. It then
follows automatically that $\ker (\Fscr(-Y)\r \Fscr)$ also has finite
length (exercise).

The category of $Y$-transversal objects is denoted by $\trans_Y(X)$.
If $\Fscr\in \trans_Y(X)$ then, in a slight generalization of
\cite{Ajitabh,VdB23}, we define $\Div(\Fscr)$ as the difference of the
divisors on $Y$  associated respectively to the finite
length $\Oscr_Y$-modules given by $\Fscr_Y$ and $\ker (\Fscr(-Y)\r
\Fscr)$.  It is easy to see that $\Div$ is additive on short exact
sequences.  Furthermore if one has a surjective map $\Fscr\r \Oscr_q$,
with kernel $\Gscr$ then the following formula holds \cite{VdB23}
\begin{equation}
\label{ref:5.39a}
\Div(\Gscr)=\Div(\Fscr)-(q)+(\tau^{-1} q)
\end{equation}
The disadvantage of the invariant $\Div(\Fscr)$ is that it depends on
$\Fscr$ itself and not only on the image of $\Fscr$ in
$\mod(X)/\Cscr_{f,p}$. We now define a
better behaving invariant.

Assume that $\Fscr\in \trans_Y(X)$. We look at
 the exact sequence
\begin{equation}
\label{ref:5.40a}
0\r \Gscr\r \Fscr(-Y)\r \Fscr\r \Hscr\r 0
\end{equation}
By definition $\Gscr$, $\Hscr$ are finite length modules on $Y$.

The completion of \eqref{ref:5.40a} looks like
\begin{equation}
\label{ref:5.41a}
0\r  \hat{\Gscr}_p\r \hat{\Fscr}_{p,\phi}\xrightarrow{\times N} \hat{\Fscr}_p
\r \hat{\Hscr}_p\r 0 
\end{equation}
Now $\hat{\Gscr}_p$, $\hat{\Hscr}_p$ are finite length modules over
$C_p/(N)$. It follows from Nakayama's lemma(lemma \ref{ref:4.8b}) that
$\hat{\Fscr}_p$ is a quotient of a finite direct sum of $P_i$'s.

Using \eqref{ref:4.11b} we can decompose
\eqref{ref:5.41a} as a product of exact sequences of the form
\begin{equation}
\label{ref:5.42a}
0\r \hat{\Gscr}_pe_i\r \hat{\Fscr}_p e_{i+1}\xrightarrow{\times N}
\hat{\Fscr}_p e_i 
\r \hat{\Hscr}_p e_i\r 0
\end{equation}
Since $\hat{\Fscr}_p$ is a quotient of a finite direct sum of $P_i$'s
it follows that $\hat{\Fscr}_p e_i=0 $ will be zero if $i\gg 0$.
Therefore it follows from \eqref{ref:5.42a} that $\hat{\Fscr}_p e_i$ is a
finite dimensional $R$-module.  Furthermore since $\hat{\Gscr}_pe_i$,
$\hat{\Hscr}_pe_i$ are only non-zero for a finite number of $i$, it
follows that $ \hat{\Fscr}_p e_{i+1}\xrightarrow{\times N}
\hat{\Fscr}_p e_i $ is an isomorphism for large negative $i$. We
define $T_p(\Fscr)$ as $\invlim_{i} \hat{\Fscr}_pe_i$ (with transition
maps given by $N$). Thus $T_p$ is a functor from $\trans_Y(X)$ to
finite dimensional $R$-modules.
\begin{propositions}
$T_p$ has the following properties. 
\begin{enumerate}
\item $T_p$ is exact.
\item $T_p(\Fscr)$ does only depend on the image of $\Fscr$ in
  $\mod(X)/\Cscr_{f,p}$.
\end{enumerate}
\end{propositions}
\begin{proof}
\begin{enumerate}
\item This follows from the exactness of the completion functor.
\item
This follows from the fact that for a finite dimensional $C_p$-module
$G$ one has $Ge_i=0$ for almost all $i$. \qed 
\end{enumerate}
\def\qed{}\end{proof} Let us now indicate a direct way of computing
$T_p(\Fscr)$. We say that $\Fscr$ is \mbox{($p$-)\emph{normalized}} if
$\Fscr$ is $\Cscr_{f,p}$-torsion free and $\tau^ip$ is not in the
support of $\Fscr_Y=0$, unless $i=0$. We have the following result.
\begin{propositions}
\label{ref:5.7.2a}
\begin{enumerate}
\item Every $\Fscr\in \trans_Y(X)$ is equivalent modulo $\Cscr_{f,p}$
  to a unique $p$-normalized object in $\mod(X)$.
\item If $\Fscr$ is $p$-normalized then
  $T_p(\Fscr)=\hat{\Fscr}_{Y,p}$. 
\end{enumerate}
\end{propositions}
\begin{proof}
\begin{enumerate} 
\item Let us first consider uniqueness. If we have an inclusion of
  $p$-normalized objects $\Fscr\subset\Fscr'$ such that
  $\Fscr'/\Fscr\in\Cscr_{f,p}$ then by \eqref{ref:5.39a}
  $\Fscr=\Fscr'$. Now assume that $\Fscr$ and $\Fscr'$ are
  $p$-normalized objects which are equivalent modulo $\Cscr_{f,p}$.
  This means that there exist a $\Cscr_{f,p}$-torsion free $\Gscr$,
  together with inclusions
\[
\Fscr\r \Gscr\l \Fscr'
\]
whose  cokernel lies in $\Cscr_{f,p}$. 
We now replace $\Fscr$, $\Fscr'$ by their images in $\Gscr$ and
we put
$\Hscr=\Fscr+\Fscr'$. Since $\Hscr$ is a quotient of $\Fscr\oplus
\Fscr'$ it follows that $\Hscr$ is also $p$-normalized. But then by
the above it follows that
$\Fscr=\Hscr=\Fscr'$. 

Let us now consider existence. Assume that $\Fscr\in \trans_Y(X)$.
Without loss of generality we may assume that $\Fscr$ is
$\Cscr_{f,p}$-torsion free. We will then modify $\Div(\Fscr)$ step by
step until it satisfies the condition  for $p$-normalization.

If $\tau^i p\in\Div(\Fscr)$ with $i>0$ then we let $\Fscr'$ be the
kernel of the associated surjective map $\Fscr\r \Oscr_{\tau^i p}$. By \eqref{ref:5.39a} it follows that
$\Fscr'$ is closer to being $p$-normalized than the original $\Fscr$.

Assume now $\tau^i p\in\Div(\Fscr)$ with $i<0$. In this case we use
the associated injective map $\Oscr_{\tau^{i+1} p}\r\Fscr_Y(Y)$ and we let
$\Fscr'$ be the pullback of the maps
\[
\begin{CD}
\Fscr(Y) @>>> \Fscr_Y(Y)\\
@. @AAA\\
@. \Oscr_{\tau^{i+1} p}
\end{CD}
\]
One now easily checks 
\[
\Div(\Fscr')=\Div(\Fscr)-(\tau^i p)+(\tau^{i+1} p)
\]
so that we have again made progress.
Repeating these constructions we eventually
find a normalized object, which is equivalent to $\Fscr$.
\item
Assume that $\Fscr$ is normalized. Then 
$\hat{\Fscr}_{Y,p}e_i$ is zero, except if $i=0$. Whence we
find from the exact sequence \eqref{ref:5.42a}
$T_p(\Fscr)=\hat{\Fscr}_{Y,p}e_0=\hat{\Fscr}_{Y,p}$. \qed
\end{enumerate}
\def\qed{}
\end{proof}
 Let $N_p$ be the functor which associates to $\Fscr\in
\trans_Y(X)$ its normalization. Then we have proved the following
\begin{corollarys} $N_p$ defines an equivalence of categories between
  $\trans_Y(X)/\Cscr_{f,p}$ and the full subcategory of $\trans_Y(X)$
  consisting of $p$-normalized objects.
\end{corollarys}

If $\Fscr$ is $p$-normalized then we can reconstruct $\hat{\Fscr}_p$ from
$T_p(\Fscr)$. 
\begin{lemmas} 
\label{ref:5.7.4a}
Assume that $\Fscr\in \trans_Y(X)$ is $p$-normal. Then
  $\hat{\Fscr}_p$ is isomorphic to the row vector 
\[
(\cdots T_p(\Fscr) \cdots T_p(\Fscr)\ 0 \cdots 0 \cdots)
\]
where the right most $T_p(\Fscr)$ occurs in position $0$.
\end{lemmas}
\begin{proof} Apriori $\hat{\Fscr}_p$ is given by the row vector
  $(\hat{\Fscr}_p e_i)_i$ \eqref{ref:4.11b}. By
  the previous discussion $\hat{\Fscr}_p e_i=0$ for $i\gg 0$ and
  multiplication by $N$ is an isomorphism on $\hat{\Fscr}_pe_i$ for
  $i\neq 0$. It now follows from the definition that
\[
\hat{\Fscr}_pe_i=
\begin{cases}
T_p(\Fscr)&\text{if $i\le 0$}\\
0&\text{if $i>0$}
\end{cases}
\qed\]
\def\qed{}\end{proof}

\section[Blowing up a point on a commutative divisor]{Blowing up a point on
a commutative divisor}
\label{ref:6a}
\subsection{Some ideals}
\label{ref:6.1a}
\leavevmode

Let $X,Y,p$ be as above  
and let $q$ be an arbitrary point on $Y$. The bimodule $o_q$ defined
above is by construction
 a quotient of $o_Y$ and hence also of $o_X$. We define
\begin{align*}
m_{Y,q}&=\ker(o_Y\r o_q)\\
m_q&=\ker (o_X\r o_q)
\end{align*}
By Corollary \ref{ref:5.5.6a}, $m_{Y,q}$, $m_q\in \Bimod(o_X-o_X)$. We
also define
\begin{equation}
\label{ref:6.1b}
\begin{split}
I_Y&= m_{Y,p}\otimes_{o_X} o_X(Y)\\
&= m_{Y,p}\otimes_{o_Y}\Nscr_{Y/X}\\
{I}&= m_p\otimes_{o_X}o_X(Y)
\end{split}
\end{equation}
Clearly $I_Y\subset \Nscr_{Y/X}$, ${I}\subset \Oscr_X(Y)$. Below
we give some properties of ${I^n}\overset{\text{def}}{=}
\im(I^{\otimes n}\r o_X(nY))$. It will be clear that
suitable analogs for $I_Y^n$ hold.
\begin{propositions}
\label{ref:6.1.1a}
\begin{enumerate}
\item $o_X(nY)/{I}^n\in \tilde{\Cscr}_{f,p}$. \item ${I}^n$ is a
  $\Cscr_p$-preserving coherent bimodule contained in
  $\Bimod(o_X-o_X)$.
\end{enumerate}
\end{propositions}
\begin{proof}
\begin{enumerate}
\item
We use induction on $n$, the case $n=1$ being clear. Put
$\Sscr_n=o_X(nY)/{I^n}$. Tensoring the exact sequence
\begin{equation}
\label{ref:6.2a}
0\r {I^n}\r o_X(nY)\r \Sscr_n\r 0
\end{equation}
with ${I}$ yields an exact sequence
\[
{I}\otimes_{o_X} {I^n}\r {I} \otimes_{o_X} o_X(nY)\r {I}\otimes_{o_X} \Sscr_n\r 0
\]
and thus
\begin{equation}
\label{ref:6.3a}
I^{n+1}=\ker ({I} \otimes_{o_X} o_X(nY)\r {I}\otimes_{o_X} \Sscr_n)
\end{equation}
Viewing ${I}\otimes o_X(nY)$ and $I^{n+1}$ as subbimodules of
$o_X((n+1)Y)$ we obtain from \eqref{ref:6.3a} an exact sequence 
\[
0\r {I}\otimes_{o_X} \Sscr_n\r \Sscr_{n+1}\r \Sscr_1(nY)\r 0
\]
Hence it suffices to prove that ${I}\otimes_{o_X} \Sscr_n\in
\tilde{\Cscr}_{p,f}$.

Tensoring 
\[
0\r {I} \r o_X(Y)\r o_{\tau p}\r 0
\]
with $\Sscr_n$ yields
\[
0\r \HTor_1^{o_X}(o_{\tau p},\Sscr_n)\r {I}\otimes_{o_X} \Sscr_n \r
o_X(Y)\otimes_{o_X} \Sscr_n\r o_{\tau p}\otimes_{o_X} \Sscr_n\r 0
\]
By \eqref{ref:5.32a} it follows that indeed ${I}\otimes_{o_X} \Sscr_n\in
\tilde{\Cscr}_{f,p}$.
\item
From the exact sequence \eqref{ref:6.2a}, Proposition \ref{ref:3.1.8a}
and  Corollary
\ref{ref:5.5.6a} it
follows that indeed $I^{n}$ is coherent and is contained in
$\Bimod(o_X-o_X)$. 

Let $\Tscr\in \Cscr_p$. Applying $\Tscr\otimes_{o_X}-$ and
$\HHom_{o_X}(-,\Tscr)$ to the sequence \eqref{ref:6.2a} yields long
exact sequences
\begin{gather*}
0\r\HTor_1^{o_X}(\Tscr,\Sscr_n)\r \Tscr\otimes_{o_X} {I^n}
\r
\Tscr(nY)\r \Tscr\otimes_{o_X} \Sscr_n\r 0
\\
0\r\HHom_{o_X}(\Sscr_n,\Tscr)\r \Tscr(-nY)\r \HHom_{o_X}({I^n},\Tscr)\r
\HExt^1_{o_X}(\Sscr_n,\Tscr)\r 0
\end{gather*}
Thanks to \eqref{ref:5.31a} and Proposition \ref{ref:5.4.3a} we can
conclude that ${I^n}$ is $\Cscr_p$ preserving.
\end{enumerate}
\def\qed{}\end{proof} The following proposition gives a little
additional information on ${I^n}$ which we will need below.
\begin{propositions} 
\label{ref:6.1.2a}
${I^n}$  has the following additional properties. Assume that
$\Mscr\in \coh(X)$. Then
\[
\HTor_i^{o_X}(\Mscr,{I^n})\qquad
\begin{cases}
\in \coh(X) &\text{if $i=0$}\\
\in \Cscr_{f,p}&\text{if $i=1$}\\
=0&\text{if $i\ge 2$}
\end{cases}
\]
\end{propositions}
\begin{proof} The case $i=0$ is already covered by Proposition 
\ref{ref:6.1.1a}. Hence we consider the case $i>0$. Tensoring the exact
sequence \eqref{ref:6.2a} on the left with $\Mscr$ yields
\[
\HTor_i^{o_X}(\Mscr,{I^n})=\HTor_{i+1}^{o_X}(\Mscr,\Sscr_n)
\]
Hence by Theorem \ref{ref:5.5.10a} $\HTor_i^{o_X}(\Mscr,{I^n})\in \Cscr_{f,p}$ for $i\ge 1$.

It remains to be shown that $\HTor_j^{o_X}(\Mscr,\Sscr_n)=0$ for $j\ge 3$.
Using Lemma \ref{ref:5.5.8a} and Theorem \ref{ref:5.5.10a}, it is
sufficient to know that the left projective dimension of $S_i$ is $\le
2$. This is the version for left modules of  Proposition \ref{ref:5.1.7a}.
\end{proof}
The following proposition gives a more explicit description of
${I^n}$.
\begin{propositions} 
\label{ref:6.1.3a}
One has the following alternative expression for
${I^n}$. 
\begin{gather*}
{I^n}=m_{p}m_{\tau^{-1}p}\cdots m_{\tau^{n-1}p}\otimes_{o_X}
o_X(nY)
\end{gather*}
We have
\begin{align}
\label{ref:6.4a}
({I}^n)^\wedge_p&=\hat{I}^n_p\\
\label{ref:6.5a}
&=(m_0m_{-1}m_{-2}\cdots m_{-n+1})_{\phi^{-n+1}}
\end{align}
as subobjects of the invertible $C_p$-bimodule $(C_p)_{\phi^{-n+1}}$.
Here $\phi$ is as in Theorem \ref{ref:5.1.4a} and for $i\in\ZZ$,
$m_i$ is the twosided maximal ideal in $C_p$ corresponding to $S_i$.
\end{propositions}
\begin{proof}
The first statement follows easily 
 from \eqref{ref:5.25a} which implies that
\[
o_X(Y)\otimes_{o_X}  m_q=m_{\tau^{-1}q}\otimes_{o_X} o_X(Y)
\]
(as subobjects of $o_X(Y)$). 

Now we prove the second statement. 
Completing the exact sequence 
\[
0\r m_p \r o_X\r o_{p}\r 0
\]
and using lemma \ref{ref:5.5.8a} gives $\hat{m}_p=m_0$. Furthermore
applying definition \eqref{ref:5.18a} together with Theorem
\ref{ref:5.1.4a}.4 yields $o_X(Y)^\wedge_p=(C_p)_{\phi^{-1}}$. 
Using the compatibility of tensor product with $(-)^\wedge_p$ we then deduce
from \eqref{ref:6.1b}
\begin{align}
\label{ref:6.6a}
\hat{I}_p&=(m_0)_{\phi^{-1}}
\end{align}

\eqref{ref:6.4a} is easily proved by induction. By
\eqref{ref:6.3a}
 we have an exact sequence
\[
0\r I^{n+1}\r {I}\otimes_{o_X} o_X(nY)\r {I}\otimes_{o_X} \Sscr_n\r 0
\]
which yields  an exact sequence (using exactness of $\hat{(-)}$ and
compatibility with tensor product, see Propositions
\ref{ref:5.4.4a} and \ref{ref:5.4.5a})
\begin{equation}
\label{ref:6.7a}
0\r (I^{n+1})^\wedge_p \r (\hat{I}_p)_{\phi^{-n}}
\r  \hat{I}_p\otimes_{o_X} \hat{\Sscr}_{n,p}\r 0
\end{equation}
On the other hand completing the exact sequence
\[
0\r I^{n}\r o_X(nY)\r  \Sscr_n\r 0
\]
and tensoring with $\hat{I}_p$ yields by induction an exact sequence
\begin{equation}
\label{ref:6.8a}
0\r \hat{I}_p\cdot \hat{I}^{n}_p
\r (\hat{I}_p)_{\phi^{-n}}
\r  \hat{I}_p\otimes_{o_X} \hat{\Sscr}_{n,p}\r 0
\end{equation}
Comparing \eqref{ref:6.7a} and \eqref{ref:6.8a} yields what we want.

\eqref{ref:6.5a} now follows from \eqref{ref:6.4a},\eqref{ref:6.6a} and the
easily verified fact that  $\phi^{-1}(m_i)=m_{i-1}$.
\end{proof}

Before we continue we make a few remarks on the case that $p$ is a
fixed point for $\tau$. In that case it follows from \eqref{ref:5.25a} that
\begin{equation}
\label{ref:6.9a}
m_p\otimes_{o_X} o_X(Y)=o_X(Y)\otimes_{o_X} m_p
\end{equation}
 as subobjects of
$o_X(Y)$.

Let $\mu$ be the multiplicity of $p$ in $Y$. If $p$ is a fixed point
then it is possible that $\mu\ge 1$ by Theorem \ref{ref:5.1.4a}. We
have the following description of $\mu$ which is completely analogous
to the classical case.
\begin{lemmas} 
\label{ref:6.1.4a}
Assume that $p$ is a fixed point for $\tau$. Then $\mu$
  is the largest integer $n$ such that $o_X(-Y)\subset m_p^n$ as
  subbimodules of $o_X$.
\end{lemmas}
\begin{proof}
We have to find the largest integer such that 
\[
o_X/m^n_p\r o_X/(o_X(-Y)+m^n_p)
\]
is  an isomorphism. Using the properties of $\hat{(-)}$ and in
particular Theorem \ref{ref:5.5.9a} we see that
this is equivalent with 
\[
R/m^n\r R/(RU+m^n)
\]
being an isomorphism, where $m$ is the maximal ideal of $R$. This in
turn is the same as saying that $U\in m^n$.  Now lemma \ref{ref:5.1.8a}
yields what we want.
\end{proof}

\begin{propositions}\label{ref:6.1.5a}
Let $\mu$ be the multiplicity of $Y$ in $p$ and let $n\in\ZZ$. Then the
natural exact sequence 
\[
0\r o_X((n-1)Y)\r o_X(nY) \r o_Y(nY)\r 0
\]
restricts to an exact sequence
\begin{equation}
\label{ref:6.10a}
0\r I^{n-\mu}((\mu-1)Y) \r {I^n}\r I^n_{Y}\r 0
\end{equation}
Here $I^n$ for $n<0$ is to be interpreted as $o_X(nY)$.
\end{propositions}
\begin{proof} 
Note that by \eqref{ref:6.9a} there is no ambiguity in the notation
$I^{n-\mu}((\mu-1)Y)$. 

First we verify that the lower sequence is well defined. We have
\begin{equation}
\label{ref:6.11a}
\Oscr(- Y)\subset m_pm_{\tau^{-1} p}\cdots m_{\tau^{\mu-1}p}
\end{equation}
If $\mu=1$ then this is by definition. If $\mu>1$ then $p$ is a fixed
point and we can use lemma \ref{ref:6.1.4a}. Tensoring \eqref{ref:6.11a} with
$o_Y(\mu Y)$
yields an inclusion
\[
\Oscr((\mu-1)Y)\subset I^\mu
\]
Assume $n\ge \mu$. Multiplying with $I^{n-\mu}$ then yields an inclusion
\begin{equation}
\label{ref:6.12a}
I^{n-\mu}((\mu-1)Y)\subset {I^n}
\end{equation}
In a similar way one verifies directly from \eqref{ref:6.11a} that
\eqref{ref:6.12a} also holds if $n<\mu$. Hence
\eqref{ref:6.10a} is indeed well defined.

To prove the proposition we look at the following commutative diagram.
\[
\begin{CD}
@.0 @.0 @.0 @.\\
@. @VVV @VVV @VVV @.\\
0 @>>> I^{n-\mu}((\mu-1)Y) @>\alpha>> {I^n} @>\beta>> I^n_{Y}
@>>> 0\\ @. @VVV @VVV @VVV @.\\
0 @>>> o_X((n-1)Y) @>>> o_X(nY)  @>>> o_Y(nY) @>>> 0\\
@. @VVV @VVV @VVV @.\\
@. \Sscr_{n-\mu}((\mu-1)Y) @>\gamma>>  \Sscr_n  @>\delta >> \Sscr_{Y,n}
@>>> 0\\
 @. @VVV @VVV @VVV @.\\
@.0 @.0 @.0 @.
\end{CD}
\]
where as usual $\Sscr_n=o_X(nY)/{I^n}$, $\Sscr_{Y,n}=o_Y(nY)/I^n_{Y,p}$.
If $n<0$ then $\Sscr_n$, $\Sscr_{n,Y}$ are defined as zero.

In this diagram the middle row is exact, $\alpha$ is injective and
$\beta$ is surjective. Elementary diagram chasing shows that the lower
row is exact.

We now have to show that $\ker \beta=\im\alpha$. A diagram chase
shows that this is equivalent to $\gamma$ being injective. This is
clear if $n<\mu$ so we may assume that $n\ge \mu$. Injectivity of $\gamma$ is then
in turn equivalent to 
\begin{equation}
\label{ref:6.13a}
\length \Sscr_n=\length\Sscr_{n-\mu}+\length \Sscr_{Y,n}
\end{equation}
By Theorem
\ref{ref:5.5.9a} this can be checked on $C_p$.  By Corollary
\ref{ref:5.2.4a} we have
\[
\length {\Sscr}_{n}=\frac{n(n+1)}{2}
\]
If $|O_\tau(p)|>1$ then by Theorem \ref{ref:5.1.4a} $p$ is smooth on
$Y$ and hence \[
\length \Sscr_{Y,n}=n
\]
One sees that in this case \eqref{ref:6.13a} is satisfied.  

Assume now that $\tau p=p$. By Theorem \ref{ref:5.5.9a} we may check
the injectivity of $\gamma$ after completing. We have $C_p=R$ as in Theorem
\ref{ref:5.1.4a}.  Let $m$ be the maximal ideal of $R$. We then have
$\hat{\Sscr}_n=R/m^n$ and the map $\hat{\gamma}_p$ is given by right
multiplication by $U$.  From lemma \ref{ref:5.1.8a} it follows that $U\in
m^\mu-m^{\mu+1}$. From this and the fact that $\gr R$ is a domain (for
the $m$-adic filtration) we deduce that right multiplication
by $U$ defines an injective map $R/m^{n-\mu}\r R/m^n$.
\end{proof}

\subsection{Some Rees algebras}
\label{ref:6.2b}
We will
use $I_Y$, ${I}$ to define the following Rees algebras.
\begin{align*}
\Dscr_Y&=o_Y\oplus I_Y\oplus I^2_{Y}\oplus\cdots\\
\Dscr&=o_X\oplus {I}\oplus  I^2 \oplus \cdots
\end{align*}
As a corollary to Proposition \ref{ref:6.1.1a} we immediately deduce
\begin{corollarys}
$\Dscr_Y$, $\Dscr\in \Alg(X)$.
\end{corollarys}
Now we deduce some other good
properties of $\Dscr$. We denote shifting
in $\Gr(\Dscr)$ and $\Bigr(\Dscr)$ by $(-)$.

As usual the case where $p$ is  fixed point for $\tau$ is somewhat
peculiar.
In that case it follows easily from \eqref{ref:6.9a} that
\begin{equation}
\label{ref:6.14a}
\Dscr(nY)\overset{\text{def}}{=}\Dscr\otimes_{o_X}o_X(nY)
=o_X(nY)\otimes_{o_X}\Dscr
\end{equation}
is a twosided invertible graded bimodule over $\Dscr$. This is clearly
false if $\tau p\neq p$. 

It is sometimes convenient to define a modified Rees algebra $\tilde{\Dscr}$ by
\begin{equation}
\label{ref:6.15a}
\tilde{\Dscr}_n
=
\begin{cases}
o_X(nY) &\text{if $n<0$}\\
I^n & \text{if $n\ge 0$}
\end{cases}
\end{equation}
$\tilde{\Dscr}_Y$ is defined similarly but with $o_Y,I_Y$ replacing
$o_X,I$. 

\begin{theorems} 
\label{ref:6.2.2a}
Let $\Dscr$, $\Dscr_Y$ be as above. Let $\mu$ be
the multiplicity of $p$ on $Y$. Then 
\begin{enumerate}
\item There is an exact sequence of graded $\tilde{\Dscr}$-bimodules.
\begin{equation}
\label{ref:6.16a}
0\r  \tilde{\Dscr}(-\mu)((\mu-1)Y)\r \tilde{\Dscr}\r
\tilde{\Dscr}_Y\r 0 
\end{equation}
\item $\Dscr$ is noetherian. 
\item $\Dscr$ satisfies $\chi$.
\item $\cd\tau_{\Dscr}\le 2$.
\end{enumerate}
\end{theorems}
\begin{proof}
The fact that \eqref{ref:6.16a} exists follows directly from
Proposition \ref{ref:6.1.5a}.

Our aim is now to deduce 2., 3. and 4. from the corresponding statements
for $\Dscr_Y$. So we first have to handle this case. If $\mu=1$
then $I_Y$ is an invertible bimodule. Then we can apply Proposition
\ref{ref:3.9.13a} to the exact sequence
\[
0\r I_Y\otimes_{o_X} \Dscr_Y(-1) \r\Dscr_Y \r \Oscr_{Y}\r
0 \]
and we are done. Hence by Theorem \ref{ref:5.1.4a} we only have to treat
the case that $p$ is a fixed point. But this means that $m_{Y,p}$ is
$\tau$-invariant. 
Let $\Escr_Y$ be the ordinary rees algebra associated to $m_{Y,p}$.
$\Escr_Y$ may be viewed as a sheaf of algebras in the ordinary sense
and the reader may  verify that it satisfies properties 2., 3. and $\cd
\tau_{\Escr_Y}\le 1$.
 Furthermore there is a category equivalence
 given by
\begin{equation}
\label{ref:6.17a}
 \Gr(\Dscr_Y)\r \Gr(\Escr_Y):\oplus_n\Mscr_n\mapsto
 \oplus_n\Mscr_n(-nY)
\end{equation}
 From this it is routine to
 pull the good properties of $\Escr_Y$ back to $\Dscr_Y$.

Define now $\Dscr'_{Y}$ by
\begin{equation}
\label{ref:6.18a}
(\Dscr'_{Y})_n=
\begin{cases}
\Dscr_Y&\text{if $n\ge \mu$}\\
\Dscr&\text{otherwise}
\end{cases}
\end{equation}
Then \eqref{ref:6.16a} gives an exact sequence
\begin{equation}\label{ref:6.19a}
0\r \Dscr(-\mu)((\mu-1)Y)\r \Dscr\r
\Dscr'_{Y}\r 0 
\end{equation}
It is easy to see that the fact that $\Dscr_Y$ is noetherian,
together with the fact that ${I^n}$ is coherent for $n=0,\ldots,\mu-1$,
implies that $\Dscr'_{Y}$ is noetherian.  Furthermore $\Dscr_Y$ and
$\Dscr'_{Y}$ have the same tails. So from lemma \ref{ref:3.9.16a} we obtain
that
 $\Dscr'_{Y}$  satisfies $\chi$ and furthermore $\cd
 \tau_{\Dscr'_{Y}}=\cd \tau_{\Dscr_Y}=1$.

 Now
applying Proposition \ref{ref:3.9.13a} to \eqref{ref:6.19a} implies
2., 3. and 4.
\end{proof}
\begin{remarks}
The inequality in \ref{ref:6.2.2a}.4 is of course an equality. However
we won't need this.
\end{remarks}

\subsection{Definition of blowing up}
\label{ref:6.3b}
We define ${\tilde{Y}}=\Proj \Dscr_Y$, $\tilde{X}=\Proj
\Dscr$  and we call $\tilde{X}$ the blowing up of $X$ in $p$. We
denote by $\pi$ resp. $\pi_Y$ the quotient maps $\Gr(\Dscr)\r
\Qgr({\Dscr})$ resp. $\Gr(\Dscr_Y)\r
\Qgr({\Dscr}_{Y})$.

It follows from \eqref{ref:6.16a} that we have a commutative
diagram of quasi-schemes
\begin{equation}
\label{ref:6.20a}
\begin{CD}
{\tilde{Y}} @>\beta>> Y\\
@V j VV @V i VV \\
\tilde{X} @>\alpha >> X
\end{CD}
\end{equation}
where $i$ is as before, $j$ comes from the quotient map $\Dscr\r
\Dscr_Y$ (through Proposition \eqref{ref:3.9.11a}) and
$\alpha^\ast\Mscr=\pi(\Mscr\otimes_{o_X} \Dscr)$, with the analogous
definition for $\beta^\ast$.

The derived functors $R^i\alpha_\ast$ are defined as
usual. $L_i\alpha_\ast$ is defined as in Section
\S\ref{ref:3.10b}. From lemma \ref{ref:3.10.2a} together with
Proposition \ref{ref:6.1.2a} it follows that $L_i\alpha^\ast$ actually
defines a functor $\Mod(X)\r \Mod({\tilde{X}})$.

The following theorem summarizes some of the main properties
of blowing up.
\begin{theorems} 
\label{ref:6.3.1a}
\begin{enumerate}
\item
The pair $({\tilde{Y}},\beta)$ is isomorphic in $\Qsch/Y$ to the
ordinary commutative blowing up at $p$ of $Y$. In particular $\beta$ is
an isomorphism if $Y$ is smooth in $p$.
\item
${\tilde{X}}$ is a noetherian quasi-scheme.
\item The ideal in $o_{{\tilde{X}}}$ defined by ${\tilde{Y}}$ is invertible.
\item $i_\ast\circ R^i\beta_\ast=R^i\alpha_\ast\circ j_\ast$.
\item We have $
\cd \alpha^\ast\le 1$, $\cd \alpha_\ast\le 1$
\item $\alpha$ is proper
\item $o_{\tilde{X}}({\tilde{Y}})$ is relatively ample for $\alpha$.
\item $j$ makes $\tilde{Y}$ into a divisor in $\tilde{X}$ in the sense
  of enriched quasi-schemes (cfr \S\ref{ref:3.7b}).
\end{enumerate}
\end{theorems}
\begin{proof} These properties are either straightforward translations
  of the definitions or they follow from  properties of $\Dscr$
  which we have already proved. 
  \begin{enumerate}
\item As usual we separate two cases. If $\tau p=p$ then we 
let $\Escr_Y$ be as in the proof of Theorem \ref{ref:6.2.2a}. The
category equivalence between $\Gr(\Dscr_Y)$ and $\Gr(\Escr_Y)$ 
yields in that case what we want.

If  $p$ is smooth on $Y$ (in particular if $\tau p\neq p$) then 
$m_{Y,p}$ is invertible. Thus the same holds for $I_Y$ and we have
\[
\Dscr_Y=o_{Y}\oplus I_Y\oplus I_Y^{\otimes 2}\cdots
\]
It now follows from Propositions \ref{ref:3.11.4a} (with $\Sscr=0$) that
$\QGr(\Dscr_Y)\cong \Qch(Y)$. 

For use below we state the following 
formula for $\Nscr\in\Gr(\Dscr_Y)$.
\begin{equation}
\label{ref:6.21a}
\beta_\ast(\pi_Y\Nscr)=\dirlim \Nscr_n\otimes_{o_Y}
I^{\otimes{-n}}_{Y}
\end{equation}
We leave the easy proof to the reader.
\item This  follows from Theorem \ref{ref:6.2.2a}.2.
\item Using Theorem \ref{ref:6.2.2a}.1. one easily checks that
\begin{equation}
\label{ref:6.22a}
\pi \Mscr\otimes_{o_{\tilde{X}}} o_{\tilde{X}}(-{\tilde{Y}})= \pi(\Mscr \otimes_{\Dscr}
\Dscr(-\mu)((\mu-1)Y))
\end{equation}
What we want to prove now follows from the fact that
$\Dscr(-\mu)((\mu-1)Y)$ is an invertible graded bimodule over
$\Dscr$ (this is trivial if $\mu=1$ and if $\mu\neq 1$ it follows from
\eqref{ref:6.14a}). 
\item Define $\gamma=i\circ \beta$. Since $i_\ast$ is exact we have
  $R^i\gamma_\ast=i_\ast\circ R^i\beta_\ast$. Hence we have to show that
\begin{equation}
\label{ref:6.23a}
R^i\gamma_\ast=R^i\alpha_\ast\circ j_\ast
\end{equation}
 We have
  ${\tilde{Y}}=\Proj \Dscr'_{Y}$ where $\Dscr'_{Y}$ is as in
  \eqref{ref:6.18a}. \eqref{ref:6.19a} shows that
  $\Dscr'_{Y,p}$ satisfies the conditions for Proposition
  \ref{ref:3.9.12a}. Hence we can employ that proposition to obtain what we
  want.
\item This follows from Theorem \ref{ref:6.2.2a}.4, equation
  \eqref{ref:3.30a} and Proposition \ref{ref:6.1.2a}.
\item This follows from Theorem \ref{ref:6.2.2a}.3 together with
  part 
 of Proposition \ref{ref:3.9.7a}.
\item Since $\Dscr$ satisfies $\chi$ it follows from Proposition
\ref{ref:3.9.7a} that the canonical shift functor $\Mscr\mapsto \Mscr(1)$ is
relatively ample. Denote the corresponding bimodule by $o_{\tilde{X}}(1)$. If $p$ is
smooth on $Y$ then it follows from \eqref{ref:6.22a} that $o_{\tilde{X}}({\tilde{Y}})=o_{\tilde{X}}(1)$. So
in that case we are done. If $p$ is singular on $Y$ then by
\eqref{ref:6.14a}
 $-\otimes_{o_X} o_X(Y)$ defines an autoequivalence of $\Qch({\tilde{X}})$. Then \eqref{ref:6.22a} yields
\begin{equation}
\label{ref:6.24a}
o_{\tilde{X}}({\tilde{Y}})=o_{\tilde{X}}(\mu)\otimes_{o_X} o_X((1-\mu )Y)
\end{equation}
This yields what we want   since it is clear that 
$-\otimes_{o_X} o_X(Y)$ commutes with $R^i\alpha_\ast$ and
$L_i\alpha^\ast$. Note in passing that \eqref{ref:6.24a} also makes
sense in the case that $\tau p\neq p$ since then $\mu=1$.
\item Given 3. we only have to show that the map
  $\Oscr_{\tilde{X}}(-\tilde{Y})\r \Oscr_{\tilde{X}}$ is injective.
  Now $\Oscr_{\tilde{X}}(-\tilde{Y})$ is by definition
  $\alpha^\ast\Oscr_X\otimes_{o_{\tilde{X}}}
  o_{\tilde{X}}(-\tilde{Y})$ which according to 
  \eqref{ref:6.24a} is equal to $\pi(\Oscr_X\otimes_{o_X}
  \Dscr(-\mu)((\mu-1)Y))$. Hence by \eqref{ref:6.16a} it suffice to
  show that $\HTor_{o_X}^i(\Oscr_X,I^n_Y((\mu-1)Y))=0$ for $i>0$,
  $n\ge 0$. As in the commutative case one checks that
  $\HTor_{o_X}^i(\Oscr_X,I^n_Y((\mu-1)Y))=
  \HTor_{o_Y}^i(\Oscr_Y,I^n_Y((\mu-1)Y))$.  One now easily shows that
  $\HTor_{o_Y}^i(\Oscr_Y,I^n_Y((\mu-1)Y))=0$, for example using the
  fact that $I^n_Y\subset o_Y(nY)$ with the quotient being in
  $\Cscr_{f,p}$, together with the definition of $\HTor(-,-)$ (cfr.
  \eqref{ref:3.1b}).
\qed \end{enumerate}
\def\qed{}\end{proof}

\subsection{The normal bundle}
\label{ref:6.4b}
The main result of this section (equation \eqref{ref:6.26a}) will
be used in \S\ref{ref:6.5b}.

If $t:U\r V$ is a map of schemes, $\Nscr$ is a line bundle on $Y$ and
$\Fscr\subset \Nscr$ is a quasicoherent subsheaf then $t^{-1}(\Fscr)$
is defined as the image of $t^\ast(\Fscr)$ in $t^\ast(\Nscr)$. 

Now let us revert to the notations in use in the previous
sections. Recall that according to \eqref{ref:5.4a} we have $\Nscr_{Y/X}
=\Nscr_\tau$ where $\Nscr$ is an invertible sheaf on $Y$ and $\tau$ is
an automorphism of $Y$.

Define $\tau':{\tilde{Y}}\r {\tilde{Y}}$ as follows. If $p$ is a fixed point for
$\tau$ then $\tau$ extends in a natural way to the blowup of $Y$ in
$p$. We denote this extended map by $\tau'$. If $p$ is not a fixed
point then the map $\beta:{\tilde{Y}}\r Y$ is an isomorphism and we put
$\tau'=\beta^{-1}\tau\beta$. So in all cases we have
\begin{equation}
\label{ref:6.25a}
\beta\tau'=\tau\beta
\end{equation}

Our aim in this section is to prove the following formula
\begin{equation}
\label{ref:6.26a}
\Nscr_{{\tilde{Y}}/{\tilde{X}}}=\beta^{-1}(m_{Y,p}\Nscr)_{\tau'}
\end{equation}
In view of the commutative case this formula seems quite logical.
However there are some pitfalls. The main problem is that a priori ${\tilde{Y}}$
is only a quasi-scheme. Thus, although we can use the ordinary
definition of $\beta^{-1}(m_{Y,p}\Nscr)$, the fact that we can
consider the result as a bimodule, depends on the ``accidental'' event
that ${\tilde{Y}}$ is commutative. So to make sense of \eqref{ref:6.26a} we
have to bring in explicitly the identification of ${\tilde{Y}}$ with a
commutative scheme (which was given in the proof of Theorem
\ref{ref:6.3.1a}.1). Below we give the necessary computations. The
reader is advised to skim through the rest of this section.

It is most convenient to separate two cases, depending on whether $p$
is a fixed point or not. 

\noindent
{\bf $p$ is not a fixed point for $\tau$.} 
In this case $\mu=1$ and
$m_{Y,p}$ is invertible. Thus
$\beta^{-1}(m_{Y,p}\Nscr)=\beta^{\ast}(m_{Y,p}\Nscr)$. We have to
show that for $\Mscr\in \Qch({\tilde{Y}})$ we have
\begin{equation}
\label{ref:6.27a}
\Mscr\otimes_{o_{\tilde{Y}}}o_{\tilde{Y}}({\tilde{Y}})=\Mscr\otimes_{o_{\tilde{Y}}}(\beta^{\ast}
(m_{Y,p}\Nscr))_{\tau'})
\end{equation}
In this case the identification of ${\tilde{Y}}$ with a commutative scheme is
given by $\beta$. Thus if $\Tscr$ is a quasicoherent $\Oscr_Y$-module then
\[
\Mscr\otimes_{o_{\tilde{Y}}}\beta^{\ast}(\Tscr)=\beta^\ast(\beta_\ast\Mscr\otimes_{o_Y}
\Tscr)
\]
We use this in the computation below.
\begin{align*}
\Mscr\otimes_{o_{\tilde{Y}}}(\beta^{\ast}(m_{Y,p}\Nscr))_{\tau'}
&=\tau'_\ast(\Mscr\otimes_{o_{\tilde{Y}}}\beta^{\ast}(m_{Y,p}\Nscr))\qquad
\text{(See \eqref{ref:5.4a})}\\
&=(\beta^{-1})_\ast\tau_\ast\beta_\ast\beta^\ast(\beta_\ast\Mscr\otimes_{o_Y}
m_{Y,p}\Nscr)\\
&=(\beta^{-1})_\ast\tau_\ast(\beta_\ast\Mscr\otimes_{o_Y}m_{Y,p}\Nscr)
\end{align*}
Thus, using
\eqref{ref:6.24a}, \eqref{ref:6.27a} reduces to
\[
\beta_\ast(\Mscr(1))=\tau_\ast(\beta_\ast\Mscr\otimes_{o_Y}m_{Y,p}\Nscr)
\]
Now the righthand side of this equation is equal to 
\[
\beta_\ast\Mscr\otimes_{o_Y}m_{Y,p}\Nscr_\tau=
\beta_\ast\Mscr\otimes_{o_Y}I_Y
\]
So finally we have to show
\[
\beta_\ast(\Mscr(1))=\beta_\ast\Mscr\otimes_{o_Y}I_Y
\]
But this follows easily from \eqref{ref:6.21a}.

\noindent
{\bf $p$ is  a fixed point for $\tau$.} Now we identify ${\tilde{Y}}$ with the
ordinary commutative blowup of $Y$ at $p$. Denote the latter by
${\tilde{Y}}'$. We have ${\tilde{Y}}'=\Proj\Escr_Y$ where $\Escr_Y$ is as in the
proof of Theorem \ref{ref:6.2.2a}.
Let $\beta':{\tilde{Y}}'\r Y$ be the structure map. There is now a 
commutative diagram 
\[
\begin{CD}
{\tilde{Y}} @>\gamma>> {\tilde{Y}}'\\
@V\beta VV @V\beta' VV\\
Y @= Y
\end{CD}
\]
where $\gamma:{\tilde{Y}}\r {\tilde{Y}}'$ denotes the identification. Note however that
$\gamma$ does not commute with the canonical shift
functors. Translating \eqref{ref:6.17a} we find that
\begin{equation}
\label{ref:6.28a}
\gamma_\ast(\Mscr(1))=(\gamma_\ast\Mscr)(1)\otimes_{o_{Y}} o_Y(Y)
\end{equation}
If we now look back
at the definition of $\tau'$ then we see that it was in fact defined
on ${\tilde{Y}}'$ and then pulled back to ${\tilde{Y}}$ by $\gamma$. Thus 
\[
\tau'=\gamma^{-1}\tau''\gamma
\]
where $\tau''$ is the extension of $\tau$ to ${\tilde{Y}}'$.

We compute
\begin{align*}
\Mscr\otimes_{o_{\tilde{Y}}}\beta^{-1}(m^\mu_{Y,p}\Nscr)_{\tau'}
&=(\gamma^{-1})_\ast\tau''_\ast(\gamma_\ast\Mscr\otimes_{o_{{\tilde{Y}}'}}\beta^{'-1}
m^\mu_{Y,p}\Nscr)
\end{align*}
This equality is a formal computation, analogous to the one where
$p$ is not a fixed point.

Using \eqref{ref:6.24a} we now have to show
\begin{equation}
\label{ref:6.29a}
\gamma_\ast(\Mscr(\mu)\otimes_{o_Y} o_Y((1-\mu)Y))=
\tau''_\ast(\gamma_\ast\Mscr\otimes_{o_{{\tilde{Y}}'}}\beta^{'-1}
(m^\mu_{Y,p}\Nscr))
\end{equation}
From \eqref{ref:6.28a} it follows that the lefthand side of this equation is
equal to 
\begin{equation}
\label{ref:6.30a}
\gamma_\ast\Mscr(\mu)\otimes_{o_Y}o_Y(Y)=(\gamma_\ast\Mscr\otimes_{o_Y}
o_Y(Y)) (\mu)
\end{equation}
Now we translate everything to graded modules. Let $\gamma_\ast\Mscr$ be
represented by $\Pscr$. Then the righthand side of \eqref{ref:6.29a} is
represented by
\[
\tau_\ast(\Pscr\otimes_{\Escr_Y}(m^\mu_{Y,p}
\Escr_Y\otimes_{o_Y}\Nscr))
\]

\eqref{ref:6.30a} is represented by 
\[
\tau_\ast(\Pscr\otimes_{o_Y}\Nscr)(\mu)=\tau_\ast(\Pscr
\otimes_{\Escr_Y}(\Escr_Y(\mu)\otimes_{o_Y}\Nscr))
\]
Thus it is sufficient to show that up to right bounded bimodules
\[
m^\mu_{Y,p}
\Escr_Y=\Escr_Y(\mu)
\]
This is now clear.

\subsection{Birationality}
\label{ref:6.5b}
Here the notations are as in the previous section. We will show that $X$
and ${\tilde{X}}$ are isomorphic ``outside the $\tau$-orbit of $p$''. In particular
we may view $X$ and ${\tilde{X}}$ as being birational.

Define 
$\alpha^{-1}(\Cscr_p)$ as in \S\ref{ref:3.11b}. Thus the objects in
$\alpha^{-1}(\Cscr_p)$ are the objects in $\Qch({\tilde{X}})$ that are represented
by graded $\Dscr$-modules such that $\Mscr_n\in\Cscr_p$ for all $p$. 
$\alpha^{-1}(\Cscr_p)$ has the following slightly more intrinsic
description.
\begin{lemmas} 
\label{ref:6.5.1a}
$\alpha^{-1}(\Cscr_p)$ is the full subcategory of objects $\Mscr$ in
$\Qch({\tilde{X}})$ for which one  has $\alpha_\ast\Mscr(n{\tilde{Y}})\in\Cscr_p$ for all $n$.
\end{lemmas}
\begin{proof} 
From \eqref{ref:3.24a} together with \eqref{ref:6.24a} it easily follows
that $\alpha_\ast\Mscr(n{\tilde{Y}})\in\Cscr_p$ if
$\Mscr\in\alpha^{-1}(\Cscr_p)$. 

Let us now prove the converse inclusion.
Assume  that $\alpha_\ast\Mscr(n{\tilde{Y}})\in\Cscr_p$ for all $n$. We have to show that
$\Mscr\in\alpha^{-1}(\Cscr_p)$. Let $\Mscr=\pi\Nscr$ for
$\Nscr\in\Gr(\Dscr)$. It suffices to consider the case that
$\Nscr$ is noetherian. We will prove that $\Nscr_n\in\Cscr_p$ for $n\gg
0$. This implies $\Mscr\in\alpha^{-1}(\Cscr_p)$.

According to \eqref{ref:6.24a} 
\begin{equation}
\label{ref:6.31a}
\Mscr(n{\tilde{Y}})=\pi\Nscr(n\mu)\otimes_{o_X}o_X(n(1-\mu)Y)
\end{equation}
So if $\mu=1$ then this formula says that $\tilde{\Nscr}_n\in\Cscr_p$
for all $n$. Since $\Dscr$ satisfies $\chi$ (Theorem 
\ref{ref:6.2.2a}) it follows from lemma \ref{ref:3.9.3a} 
that $\tilde{\Nscr}_n=\Nscr_n$ for
$n\gg 0$. Hence we are done.

This reasoning still works if $\mu>1$ since $-\otimes_{o_X}o_X(Y)$
commutes with $\alpha_\ast$. Hence \eqref{ref:6.31a} yields that
$\tilde{\Nscr}_{n\mu}\in\Cscr_p$ for all $n$. Then according to lemma
\ref{ref:3.12.1a} we have  modulo $\Tors(\Dscr)$
\[
\Nscr=\tilde{\Nscr}=\tilde{\Nscr}^{(\mu)}\otimes_{\Dscr^{(\mu)}}\Dscr
\]
Let $\Nscr'$ be the module on the righthand side of this equation. It is
clear that $\Nscr'_n\in\Cscr_p$ for all $n$. Since $\Nscr$ and $\Nscr'$
are isomorphic in $\QGr(\Dscr)$ and $\Nscr$ is noetherian, it easily
follows that $\Nscr_n\in\Cscr_p$ for $n\gg 0$.
\end{proof}

\begin{propositions} 
\label{ref:6.5.2a}
We have 
\begin{enumerate}
\item $\alpha^\ast$ and $\alpha_\ast$ define inverse equivalences
  between $\Qch(X)/\Cscr_p$ and
  $\Qch({\tilde{X}})/\alpha^{-1}(\Cscr_p)$. 
\item The functor $\alpha_\ast$ sends $\alpha^{-1}(\Cscr_p)$ to
$\Cscr_p$ and the functor
$R^1\alpha_\ast$ sends $\Qch({\tilde{X}})$ to $\Cscr_p$. 
\item The functor $\alpha^\ast$ sends $\Cscr_p$ to
$\alpha^{-1}(\Cscr_p)$ and the functor $L_1\alpha^\ast$ sends $\Qch(X)$
to $\alpha^{-1}(\Cscr_p)$. 
\item The functors $R^i\alpha_\ast$ and $L_i\alpha^\ast$ preserve
coherent objects.
\end{enumerate}
\end{propositions}
\begin{proof}
\begin{enumerate}
\item
It is clear from Proposition \ref{ref:6.1.1a} and  Theorem
\ref{ref:5.5.10a}.1 that $\tilde{\Dscr}$ is strongly graded with respect
to $\Cscr_p$. What we have to prove now follows from Proposition
\ref{ref:3.11.4a}.
\item
The statement about $\alpha_\ast$ follows from \eqref{ref:3.24a}.
The statement about $R^1\alpha_\ast$ follows from applying $\alpha_\ast$ to
an injective resolution and using the fact that $\alpha_\ast$ is exact
modulo $\Cscr_p$ (by Proposition \ref{ref:3.11.4a}).
\item
This follows from the definition of $L_i\alpha^\ast$ together with
Propositions \ref{ref:6.1.1a} and \ref{ref:6.1.2a}.
\item This follows from Theorem \ref{ref:6.3.1a}.6 and lemma 
\ref{ref:3.1.17a}. 
\qed\end{enumerate}
\def\qed{}\end{proof}
\begin{corollarys} 
\label{ref:6.5.3a}
Assume that $q\in Y$ is not contained in the
  $\tau$-orbit of $p$. Let $q'$ be the unique point of ${\tilde{Y}}$ such that
  $\beta(q')=q$. Then $\alpha^\ast$ and $\alpha_\ast$ define
  inverse equivalences between $\Cscr_q$ and $\Cscr_{q'}$.
\end{corollarys}
\begin{proof} 
Using the foregoing proposition it is sufficient to prove the
following
\begin{enumerate}
\item $\Cscr_q\cap \Cscr_p=0$.
\item $\Cscr_{q'}\cap \alpha^{-1}(\Cscr_p)=0$.
\item $\alpha_\ast(\Cscr_{q'})\subset \Cscr_q$.
\item $\alpha^\ast(\Cscr_q)\subset \Cscr_{q'}$.
\end{enumerate}
1. is clear. To prove 3. it is sufficient to show that $\alpha_\ast
\Oscr_{\tau^{\prime n}q'}\in \Cscr_q$ where $\tau'$ is as in 
\S\ref{ref:6.4b}. This
follows from the fact that this is obviously true for $\beta$ by
\eqref{ref:6.25a}.

Now we prove 2. Suppose $\Mscr\in 
\Cscr_{q'}\cap \alpha^{-1}(\Cscr_p)$. Then $\alpha_\ast\Mscr(n{\tilde{Y}})\in
\Cscr_q\cap \Cscr_p=0$. From the relative ampleness of $o_{\tilde{X}}(n{\tilde{Y}})$ we
deduce that $\Mscr=0$.

Finally we prove 4. It is sufficient to show that
\begin{equation}
\label{ref:6.32a}
\alpha^\ast
\Oscr_{\tau^n q}=\beta^\ast \Oscr_{\tau^n q}
\end{equation}
First note that if $\Nscr\in\Cscr_{f,q}$ then $\hat{\Nscr}_p=0$. Hence
if $\Sscr\in\tilde{\Cscr}_{f,p}$ then by \eqref{ref:5.35a} we have
$\HTor_i^{o_X}(\Nscr,\Sscr)=0$. We deduce that 
\[
\Oscr_{\tau^n q}\otimes_{o_X}{I^n}=\Oscr_{\tau^n q}\otimes_{o_X}o_X(nY)=
\Oscr_{\tau^n q}\otimes_{o_Y}o_Y(nY)=\Oscr_{\tau^n q}\otimes_{o_Y}I^n_{Y}
\]
We obtain \eqref{ref:6.32a} by summing over all $n$ and applying $\pi$.
\end{proof}

\subsection{The exceptional curve}
\label{ref:6.6b}
In this section we will for simplicitly make use of the object
$\Oscr_X$ and consequently of the functor $(-)_{o_X}$ which goes from
bimodules on $X$ to objects in $\Mod(X)$ (cfr. \S\ref{ref:3.6b}).

From the compatibility of completion with $\Hom$ and tensor product we
deduce that the following functors
\[
\tilde{\Cscr}_{f,p} \xrightarrow{(-)_{o_X}} \Cscr_{f,p} \xrightarrow
{\Gamma(X,-)} \mod(k)
\]
are faithful and exact. On $o_p$ they act by 
$
o_p\mapsto \Oscr_p\mapsto k
$.

In the sequel we will write $\Oscr_L=\alpha^\ast O_{\tau p}$. 
Thus $\Oscr_L=\pi((\Dscr/m_{\tau p}\Dscr)_{o_X})$.

The following lemma gives a description of $\Dscr/m_{\tau p}\Dscr$
as $o_X-\Dscr$-bimodule.
\begin{lemmas}
There is an exact sequence of
$o_X-\tilde{\Dscr}$-bimodules.
\begin{equation}
\label{ref:6.33a}
0\r (o_X(-Y)\otimes_{o_X}\tilde{\Dscr})(1)\r \tilde{\Dscr}\r
\Dscr/m_{\tau   p}\Dscr\r 0
\end{equation}
\end{lemmas}
\begin{proof} 
From Proposition \ref{ref:6.1.3a} (and its proof) it follows easily that
\[
o_X(-Y)\otimes_{o_X} I^{n+1}=m_{\tau p} I^n
\]
This yields \eqref{ref:6.33a} by taking  direct sums over $n$.
\end{proof}

It was
observed by Smith and Zhang \cite{SmithZhang} that it is possible to 
develop a formalism  such that $\Oscr_L$ is really
the structure sheaf of a ``non-commutative curve'' $L$. We will say more
on this in \S\ref{ref:6.7b}.

In the special case that $p=\tau p$ one has that $m_{\tau p}\Dscr$ is
a twosided ideal. Then we can
simply define $L=\Proj \Dscr/m_p\Dscr$. In this case we will of course
use the notation $o_L$ for the algebra on $\tilde{X}$ corresponding to
the identity functor on $\Qch(L)$.

 We find
\begin{propositions} 
\label{ref:6.6.2a}
Assume that $p=\tau p$. Then
\begin{enumerate}
\item
$L\cong \PP^1$.
\item $L$ is embedded as a divisor in ${\tilde{X}}$.
\item There is a commutative diagram
\begin{equation}
\label{ref:6.34a}
\begin{CD}
L@>\gamma>> p\\
@Vu VV @V Vv V\\
{\tilde{X}} @>\alpha>> X
\end{CD}
\end{equation}
where $u,v$ are the inclusions and $\gamma$ is isomorphic to
the structure map $\PP^1\r \spec k$. 
\end{enumerate}
\end{propositions}
\begin{proof}
\begin{enumerate}
\item 
 Using Proposition
\ref{ref:5.6.2a} we see that $\Gr(\Dscr/m_p\Dscr)$ is
equivalent with $\Gr((\Dscr/m_p\Dscr)^\wedge)$ and a similar
result for $\QGr$. Using the
compatibility results for $\hat{(-)}$ we find 
\begin{equation}
\label{ref:6.35a}
(\Dscr/m_p\Dscr)^\wedge=R/m\oplus
(m/m^2)_{\phi^{-1}}\oplus(m^2/m^3)_{\phi^{-2}}\oplus\cdots
\end{equation} where $R$ is as in  Theorem \ref{ref:5.1.4a} and
$m$ is the maximal ideal of $R$. It is easily seen that the ring on
the right of \eqref{ref:6.35a} is a noetherian two generator quadratic
algebra with one relation.
It is well known that  the $\Proj$ of such a graded ring is $\PP^1$. 
\item From \eqref{ref:6.33a} it follows that 
up to
right bounded bimodules we have
\[
m_p\Dscr=\Dscr(1)(-Y)
\]
Hence, up to right bounded bimodules, $m_p\Dscr$ is an invertible
ideal and in particular $L$ is a divisor in $\tilde{X}$.
\item
This is a translation of the commutative diagram of $o_X$-algebras, given
by
\[
\begin{CD}
\Dscr/m_p\Dscr @<<< o_X/m_p\\
@AAA @AAA \\
\Dscr @<<< o_X
\end{CD}\qed
\]
\end{enumerate}
\def\qed{}\end{proof}
We are now in a position to prove the  following result.
\begin{theorems}
\label{ref:6.6.3a}
Assume $q\in {\tilde{Y}}$ if $\tau p\neq p$ and $q\in {\tilde{Y}}\cup L$ if $\tau p= p$. Then
$\Oscr_{q}$ has finite injective dimension in $\Qch({\tilde{X}})$.
\end{theorems}
\begin{proof} First note that $q\r \Oscr_q$ is a one-one
correspondence between the points on ${\tilde{Y}}$ and the simple objects in
$\Qch({\tilde{Y}})$. Similarly for $L$ and $\Qch(L)$. Therefore we interprete $L\cap {\tilde{Y}}$
set theoretically as those $q$ such that $\Oscr_q$ lies both in $\Qch({\tilde{Y}})$ and
in $\Qch(L)$. Put
\[
r=\begin{cases}
\beta(q)&\text{if $q\in {\tilde{Y}}$}\\
p&\text{if $q\in L$}
\end{cases}
\]
Using the diagrams \eqref{ref:6.34a} and \eqref{ref:6.20a} we see
that in all cases $\Oscr_r=\alpha_\ast \Oscr_q$. Thus $r$ is well defined.

Consider first the case $\tau p\neq p$. Then $\beta$ is an isomorphism.
If $r\not\in O_\tau(p)$ then according to Corollary \ref{ref:6.5.3a}
$\Cscr_r$ is equivalent with $\Cscr_q$. Since both these categories are
stable under injective hulls (Proposition \ref{ref:5.1.3a}) we are through in
this case.

However if $r\in O_\tau(p)$ then according to Theorem \ref{ref:5.1.4a}
$r$ is smooth on $Y$. Thus $q$ is also smooth on ${\tilde{Y}}$ hence we can apply
lemma \ref{ref:5.1.1a}.

Consider now the case $\tau p=p$. If $p\in L$ then we can apply lemma
\ref{ref:5.1.1a} again with $L$ as our curve. So assume now that $q\in
{\tilde{Y}}-{\tilde{Y}}\cap L$. We then claim that $r\neq p$. Suppose the contrary. Thus
there is an isomorphism
$\zeta:\Oscr_p\r\alpha_\ast \Oscr_q$. By adjointness we obtain a map
$\eta:\alpha^\ast \Oscr_p\r \Oscr_q$. This map must be non-zero since otherwise
$\zeta$ had to be zero also.
Since $\Oscr_q$ is simple we obtain that $\eta$ is surjective. 

Now the definition of $L$ implies that $\alpha^\ast \Oscr_p\in \Qch(L)$.
Hence since $\Qch(L)$ is closed in $\Qch({\tilde{X}})$ we find $\Oscr_q\in \Qch(L)$.
Contradiction.

Since $r\not\in O_\tau(p)=\{p\}$ we can now use the same reasoning as 
in the case $\tau p\neq p$ to conclude that $\Oscr_q$ has finite injective
dimension.
\end{proof}
\subsection{The structure of $\alpha^{-1}(\Cscr_p)$}
\label{ref:6.7b}
In this section we prove Proposition \ref{ref:6.7.1a} below.
It is a generalization of the main result of
\cite{SmithZhang} in our special case. Throughout $n=|O_\tau (p)|$. 
\begin{propositions}
\label{ref:6.7.1a} 
Let $L=\Qch(L)$ be the full abelian subcategory of $\Qch(\tilde{X})$
whose objects are
direct limits of subquotients of finite sums $\Oscr_L(m_1)\oplus\cdots
\oplus \Oscr_L(m_n)$ (note that this use of $L$ is consistent with the use in
\S \ref{ref:6.6b}, when $n=1$).

Let $S$ be one of the following graded rings.
\begin{enumerate}
\item
If $n=1$ then $S$ is the twist \cite{ATV2,Zhang} of $\gr_F R$ by the 
automorphism $\phi$, which was introduced in Theorem \ref{ref:5.1.4a}. Here $F$
denotes the $m$-adic filtration on $R$.
\item 
If $2\le n<\infty$ then $S=k[u,v]$ where $\deg u=1$, $\deg v=n$.
\item If $n=\infty$ then $S=k[x]$ with $\deg x=1$.
\end{enumerate}
Then
\[
L\cong 
\begin{cases}
\QGr S &\text{if $n<\infty$}\\
\Gr S &\text{if $n=\infty$}
\end{cases}
\]
This equivalence is compatible with the natural shift functors and sends
$\Oscr_L$ to $S$.
\end{propositions}
The proof of this result depends on $n$. If $n=\infty$ then the result
follows from \cite{SmithZhang}, so we treat this case first.

\medskip

\noindent
\textbf{The case {\mathversion{bold} $n=\infty$}} In this case we have
an exact sequence (by \eqref{ref:6.24a})
\[
0\r \Oscr_L(-1)\r \Oscr_L\r \Oscr_{L,\tilde{Y}}\r 0
\]
and
\begin{equation}
\label{ref:6.36a}
\Oscr_{L,\tilde{Y}}=j^\ast\alpha^\ast\Oscr_{\tau p}= \Oscr_{\tau
  \beta^{-1}(p)}
\end{equation}
(where we use $\alpha j=i\beta$). Proposition 
\ref{ref:6.7.1a}  follows
from \cite{SmithZhang} if we can show that $\Oscr_{\tau \beta^{-1}(p)}$ is the
only simple quotient of $\Oscr_L$.

Hence let $\Oscr_L\r\Sscr$ be such a simple quotient. Tensoring with
$\Oscr_{\tilde{Y}}$ we find that $\Sscr_Y$ is a quotient of $\Oscr_{\tau
\beta^{-1}(p)}$. Since $\Oscr_{\tau \beta^{-1}(p)}$ is simple, we either
have $\Sscr_{\tilde{Y}}=\Oscr_{\tau \beta^{-1}(p)}$ or $\Sscr_{\tilde{Y}}=0$. If we are in the
first case then $\Oscr_{\tau \beta^{-1}(p)}$ is a quotient of $\Sscr$ and
hence by the simplicity of $\Sscr$~: $\Sscr=\Oscr_{\tau \beta^{-1}(p)}$.
Hence  assume that we are in the second case. Let $\Tscr$ be some
graded quotient of $(\Dscr/m_{\tau p}\Dscr)_{o_X}$ such that $\Sscr=\pi
\Tscr$, Then $0=\Sscr_{\tilde{Y}}=\Sscr/\Sscr(-Y)=\Sscr/\Sscr(-1)$
implies that $\Tscr_m\cong \Tscr_{m+1}$ as $o_X$-modules for $m\gg 0$.

From the next lemma it follows that for $m\gg 0$, $\Tscr_m$
has a composition series starting with
$(\Oscr_{\tau p}(mY),\Oscr_{\tau p}((m-1)Y),\ldots)$, which is clearly
incompatible with the isomorphism $\Tscr_m\cong \Tscr_{m+1}$. Hence
we have  obtained a contradiction.
We conclude that $\Sscr=\Oscr_{\tau \beta^{-1} (p)}$ and so we can
  invoke the results of \cite{SmithZhang}.

The following lemma was used.
\begin{lemmas}
\label{ref:6.7.2a}
 Fix $t\in\NN$. Then  $(\Dscr/m_{\tau p}\Dscr)_{o_X,t}$
  is uniserial of length $t+1$ with composition series.
$(\Oscr_{\tau p}(tY),\ldots, \Oscr_{\tau p}(Y),\Oscr_{\tau p})$
(starting from the top).
\end{lemmas}
\begin{proof}
By definition $(\Dscr/m_{\tau p}\Dscr)_t$ is given by
\[
(m_p\cdots m_{\tau^{1-t}p})/(m_{\tau p} m_p\cdots m_{\tau^{1-t} p})
\]
We can compute this by completing (for example using Proposition
\ref{ref:5.2.2a}). We find 
\[
(\Dscr/m_{\tau
  p}\Dscr)^\wedge_t\cong(\ldots,0,R/m,\ldots,R/m,0,\ldots,\ldots)
\]
where the $R/m$ occur in positions $1,\ldots,1+t$. It is easy to see
that this is a uniserial right $C_p$-module with the correct composition
series.
\end{proof}

\medskip

\noindent
\textbf{The case {\mathversion{bold} $n=1$}}
This follows from the proof of Proposition \ref{ref:6.6.2a}.1.

\medskip

\noindent
\textbf{The case {\mathversion{bold} $2\le n<\infty$}}
From the viewpoint of  computations this is the most interesting case. 

According to Proposition \ref{ref:5.6.2a} we may clearly
assume that $X=\Spec C_p$, $Y=\Spec C_p/(N)$. Note that one has
$C_p/(N)=R/(U)\oplus\cdots\oplus R/(U)$.

 We should now analyze the blowup $\tilde{X}$ of
$X$ in the point defined by $m_0\subset C_p$ (which according to our current
conventions corresponds to $p$). However it turns out that it is
slightly more convenient to work with $m_{n-1}$. This does not alter
our results in any way since all  maximal ideals in $C_p$ are
conjugate under the automorphism induced by $N$.

By definition
\begin{align*}
\Dscr&=C_p\oplus m_{n-1} N^{-1} \oplus (m_{n-1}
N^{-1})^2\oplus\cdots\\
&=C_p\oplus m_{n-1} N^{-1} \oplus m_{n-1} m_{n-2} N^{-2}\oplus\cdots
\end{align*}
The $n$'th Veronese of $\Dscr$ is given by 
\[
\Dscr^{(n)}=C_0\oplus J N^{-n}\oplus J^2 N^{-2n}\oplus \cdots
\]
where $J=m_{n-1}m_{n-2}\cdots m_0$. Since $N^n J N^{-n} =J$ we find
that conjugation by $N^n$ induces an automorphism of
$\Dscr^{(n)}$. Twisting by this automorphism \cite{ATV2,Zhang}
yields that $\Gr(\Dscr^{(n)})$  is equivalent to $\Gr(\Uscr)$ where
$\Uscr$ is the ordinary Rees algebra of the ideal $J\subset C_p$.

Now before we continue, we remind the reader that
taking Veronese's and twisting is in general \emph{not} compatible with the
natural shift-functors. However here we have
$o_{\tilde{X}}(1)=o_{\tilde{X}}(\tilde{Y})$. That is, the shift
functor on $\tilde{X}$ is defined by a divisor. Hence to keep track of
this shift functor, we simply have to keep track of $\tilde{Y}$. As we
will see this is easy.

To make progress we have to compute  $J$ explicitly.
Using the material in Proposition \ref{ref:5.2.2a}
or directly we find that
\[
J=
\begin{pmatrix}
m          &(U)        &\cdots   &\cdots   &(U)      \\
\vdots     &\ddots     &\ddots   &         &\vdots   \\
\vdots     &           &\ddots   &\ddots   &\vdots   \\
\vdots     &           &         &\ddots   &(U)      \\
m          &\cdots     &\cdots   &\cdots   &m    
\end{pmatrix}
\]
Computing the powers of $J$ yields $J^n=m^{n-1} J$.

Let $T$ be the Reesring of $R$ associated to $m$. We find.
\[
\Uscr=
\begin{pmatrix}
T          &RU\oplus TU(-1)        &\cdots   &\cdots   &RU\oplus TU(-1)      \\
\vdots     &\ddots     &\ddots   &         &\vdots   \\
\vdots     &           &\ddots   &\ddots   &\vdots   \\
\vdots     &           &         &\ddots   &RU\oplus TU(-1)      \\
T          &\cdots     &\cdots   &\cdots   &T    
\end{pmatrix}
\]
In particular we have an inclusion 
\begin{equation}
\label{ref:6.37a}
\begin{pmatrix}
T          &TU(-1)        &\cdots   &\cdots   &TU(-1)      \\
\vdots     &\ddots     &\ddots   &         &\vdots   \\
\vdots     &           &\ddots   &\ddots   &\vdots   \\
\vdots     &           &         &\ddots   &TU(-1)      \\
T          &\cdots     &\cdots   &\cdots   &T    
\end{pmatrix}
\hookrightarrow \Uscr
\end{equation}
with rightbounded cokernel. Thus $\tilde{X}=\Proj\Dscr=\Proj
\Dscr^{(n)}\cong \Proj \Uscr$ is equal to the Proj of the lefthand
side of \eqref{ref:6.37a}.

There is a more elegant way to look at this. Let $X_c=\Spec R$,
$Y_c=\Spec R/(U)$, $\tilde{X}_c=\Proj T$. Let $T_Y$ be the Reesring of
$R/(U)$ and put $\tilde{Y}_c=\Proj Y_c$. Define $\Ascr$ to be the following
algebra on $\tilde{X}_c$~:
\[
\Ascr=
\begin{pmatrix}
o_{\tilde{X}_c}          &o_{\tilde{X}_c}(-\tilde{Y_c})        &\cdots   &\cdots   &o_{\tilde{X}_c}(-\tilde{Y_c})       \\
\vdots     &\ddots     &\ddots   &         &\vdots   \\
\vdots     &           &\ddots   &\ddots   &\vdots   \\
\vdots     &           &         &\ddots   & o_{\tilde{X}_c}(-\tilde{Y_c})  \\
o_{\tilde{X}_c}          &\cdots     &\cdots   &\cdots   &o_{\tilde{X}_c}    
\end{pmatrix}
\]
Then we obtain $\tilde{X}=\Spec \Ascr$ in $\QSch/X_c$.  A similar
computation yields that $\tilde{Y}=\Spec \Cscr$ where
$\Cscr=\diag(o_{\tilde{Y}_c},\ldots, o_{\tilde{Y}_c})$. Furthermore
$o_{\tilde{X}}(-\tilde{Y})$ is given by the twosided ideal 
\begin{equation}
\label{ref:6.38a}
\Jscr=\begin{pmatrix}
  o_{\tilde{X}_c} (-\tilde{Y_c}) &o_{\tilde{X}_c}(-\tilde{Y_c})
  &\cdots &\cdots &o_{\tilde{X}_c}(-\tilde{Y_c}) \\ \vdots &\ddots
  &\ddots & &\vdots \\ \vdots & &\ddots &\ddots &\vdots \\ \vdots & &
  &\ddots & o_{\tilde{X}_c}(-\tilde{Y_c}) \\ o_{\tilde{X}_c} &\cdots
  &\cdots &\cdots &o_{\tilde{X}_c}(-\tilde{Y_c})
\end{pmatrix}
\end{equation}
in $\Ascr$.

Now we compute the exceptional curve. By definition this will
correspond to $\pi(\Uscr/m_0\Uscr)$. A quick computation reveals that
\[
\Uscr/m_o\Uscr\cong (\gr R\,\, U\!\gr R(-1)\,\,\cdots\,\, U\!\gr R(-1))
\]
where $\gr R$ is the associated graded ring for the $m$-adic
filtration on $R$.

Let $\Oscr_{L_c}$ be the exceptional curve in $\tilde{X}_c$. Then we find
that the exceptional curve in $\tilde{X}$ is given by
\begin{equation}
\label{ref:6.39a}
\begin{split}
\Oscr_L&=(\Oscr_{L_c} \,\,\Oscr_{L_c}(-\tilde{Y}_c)\,\,\cdots\,\,
\Oscr_{L_c}(-\tilde{Y}_c))\\
&=(\Oscr_{L_c} \,\,\Oscr_{L_c}(-p'_c)\,\,\cdots\,\,
\Oscr_{L_c}(-p'_c))
\end{split}
\end{equation}
where $p'_c=L_c\cap \tilde{Y}_c$.

Tensoring \eqref{ref:6.39a} with positive and negative powers of
$\Jscr$ we find that $\Qch(L)$ is contained in the category of modules
over
\[
\Bscr=\begin{pmatrix}
  o_{L_c}  &o_{L_c}(-p'_c)
  &\cdots &\cdots &o_{L_c}(-p'_c) \\
 \vdots &\ddots
  &\ddots & &\vdots \\ \vdots & &\ddots &\ddots &\vdots \\ 
\vdots & &
  &\ddots & o_{L_c}(-p'_c) 
\\
 o_{L_c} &\cdots
  &\cdots &\cdots &o_{L_c}
\end{pmatrix}
\]

The inherited shift functor on $\Mod(\Bscr)$ is given by tensoring
with the inverse of 
\[
\Iscr=
\begin{pmatrix}
  o_{L_c}(-p'_c)  &o_{L_c}(-p'_c)
  &\cdots &\cdots &o_{L_c}(-p'_c) \\
 \vdots &\ddots
  &\ddots & &\vdots \\ \vdots & &\ddots &\ddots &\vdots \\ 
\vdots & &
  &\ddots & o_{L_c}(-p'_c) 
\\
 o_{L_c} &\cdots
  &\cdots &\cdots &o_{L_c}(-p'_c)
\end{pmatrix}
\]
To interprete this we have to remember that $L_c\cong \PP^1$. So we
can view $\Bscr$ as an ordinary sheaf of algebras on $\PP^1$ 
given by
\[
\Bscr=\begin{pmatrix}
 \Oscr_{\PP^1}  &\Oscr_{\PP^1}(-1)
  &\cdots &\cdots &\Oscr_{\PP^1}(-1) \\
 \vdots &\ddots
  &\ddots & &\vdots \\ \vdots & &\ddots &\ddots &\vdots \\ 
\vdots & &
  &\ddots & \Oscr_{\PP^1}(-1)
\\
\Oscr_{\PP^1} &\cdots
  &\cdots &\cdots &\Oscr_{\PP^1}
\end{pmatrix}
\]
and
$\Iscr$ as the corresponding sheaf of twosided ideals. With this new
point of view, $\Oscr_L$ is given by
\begin{equation}
\label{ref:6.40a}
\Oscr_L=(\Oscr_{\PP^1} \,\,\Oscr_{\PP^1}(-1)\,\,\cdots\,\,
\Oscr_{\PP^1}(-1))
\end{equation}
From this explicit interpretation it is now easy to verify that every
object in $\Mod(\Bscr)$ is a direct limit of subquotients of direct
sums of objects like \eqref{ref:6.40a}, tensored with powers of
$\Iscr$. Hence $L\cong \Spec \Bscr$. What remains to be shown is that
$\Spec \Bscr=\Proj S$.

One way to accomplish this is as follows. One considers the triple \cite{AZ}
$(\Mod(\Bscr),\Oscr_L,\Iscr^{-1})$ and one verifies that this triple is
ample. The corresponding graded ring 
\[
\Gamma_\ast(\Oscr_L)=\sum \Hom(\Oscr_L,\Oscr_L\otimes_{\Bscr}
\Iscr^{-n})
\]
 is equal to $S$. Hence according to \cite{AZ}~:
$\Mod(\Bscr)=\QGr(\Bscr)$. This finishes the proof of Proposition
\ref{ref:6.7.1a} in all cases.

\begin{lemmas}
\label{ref:6.7.3a}
If $\Mscr \in\alpha^{-1}(\Cscr_p)$ then there is a non-trivial map 
$\Oscr_L(t)\r \Mscr$ for some $t$.
\end{lemmas}
\begin{proof}
Assume that $\Mscr=\pi\Nscr$ where $\Nscr\in\Gr(\Dscr)$  is
torsionfree (for $\Tors(\Dscr)$). Let $u$ be such that $\Nscr_u\neq
0$. Since $\Nscr_u\in\Cscr_p$, $\Nscr_u$ will contain some
$\Oscr_{\tau^v p}$. Hence there is a non-trivial map
\[
(\Oscr_{\tau^v p}\otimes_{o_X} \Dscr)(-u)\r \Nscr
\]
Since $\Nscr$ is torsion free this yields a non-trivial map.
\begin{equation}
\label{ref:6.41a}
\alpha^\ast(\Oscr_{\tau^v p})(-u)\r \Mscr
\end{equation}
If it happens that $\tau^v p=\tau p$ then $\alpha^\ast(\Oscr_{\tau ^v
  p})=\Oscr_L$ and we are through.

Assume $\tau^v p\neq \tau p$. Then according to Proposition \ref{ref:8.3.2a}
\begin{equation}
\label{ref:6.42a}
\alpha^\ast(\Oscr_{\tau^v p})
 =\Oscr_{\tau^v \beta^{-1}(p)}
\end{equation}
There is a surjective map
\[
\Oscr_L=\alpha^\ast(\Oscr_{\tau p})\r \beta^\ast(\Oscr_{\tau
  p})=\Oscr_{\tau \beta^{-1}(p)} 
\]
Twisting yields a surjective map
\begin{equation}
\label{ref:6.43a}
\Oscr_L(v-1)\r \Oscr_{\tau \beta^{-1}(p)}(v-1)
=
\Oscr_{\tau^v \beta^{-1}(p)}
\end{equation}
Combining \eqref{ref:6.41a}\eqref{ref:6.42a}\eqref{ref:6.43a} yields what we want.
\end{proof}

\begin{corollarys}
\label{ref:6.7.4a}
Every  object in $\alpha^{-1}(\Cscr_p)\cap \coh(\tilde{X})$
is a finite extension of objects  in $\mod(L)$.
\end{corollarys}
\begin{proof}
This follows immediately from  lemma \ref{ref:6.7.3a}.
\end{proof}
\subsection{The strict transform}
\label{ref:6.8b}
This is a more specialized section. We introduce the notion of a
strict transform and its influence on the invariant $T_p(\Fscr)$
introduced in \S\ref{ref:5.7b}. \emph{As in
\S\ref{ref:5.7b} we assume that the $\tau$-orbit of $p$
has infinite order.}

 First note some lemmas.
\begin{lemmas}
\label{ref:6.8.1a}
  One has $\alpha^{-1}(\Cscr_p)\cap \mod(\tilde{X})\subset
  \trans_{\tilde{Y}}(\tilde{X})$.
\end{lemmas}
\begin{proof}
This follows for example from  Corollary \ref{ref:6.7.4a} together
with the fact that up to finite length modules every object in
$\mod(L)$ is isomorphic to $\Oscr_L^t$ for some $t$ (see
\cite{SmithZhang} or Proposition \ref{ref:6.7.1a}).
It then suffices to invoke \eqref{ref:6.36a}.
\end{proof}
\begin{lemmas}
The functors $\alpha^\ast$, $\alpha_\ast$, restrict to functors
between $\trans_Y(X)$ and $\trans_{\tilde{Y}}(\tilde{X})$. 
\end{lemmas}
\begin{proof}
That $\alpha^\ast$ preserves transversality  follows immediately
from $\alpha j=i\beta$. 

Let us now look at $\alpha_\ast$. Take an object $\Tscr\in
\trans_{\tilde{Y}}(\tilde{X})$. We have to show that
  $i^\ast\alpha_\ast\Tscr$ has finite length. For this it is
  sufficient to show that $\beta^\ast i^\ast\alpha_\ast
  \Tscr=j^\ast\alpha^\ast\alpha_\ast\Tscr$ has finite length. Now
  $\alpha^\ast\alpha_\ast\Tscr$ is isomorphic to $\Tscr$ modulo
  $\alpha^{-1}(\Cscr_p)\cap \mod(\tilde{X})$ (Theorem
  \ref{ref:6.5.2a}) so it is in $\trans_{\tilde{Y}}(\tilde{X})$
    by lemma \ref{ref:6.8.1a}. This proves what we want.
\end{proof}
If $\Fscr\in \Mod(X)$ then we define
\begin{align*}
\Fscr I^n&=\im (\Fscr\otimes I^{\otimes n}\r \Fscr(nY))\\
\Fscr\cdot \Dscr&=\oplus_n \Fscr I^n\\
\alpha^{-1}(\Fscr)&=\pi(\Fscr\cdot \Dscr)
\end{align*}
 We will use related notations such 
as $\Fscr\cdot \Dscr_Y$, $\beta^{-1}(\Fscr)$ when they apply.

The idea is that if $\Fscr\in \trans_Y(X)$ then $\alpha^{-1}(\Fscr)$
should correspond to the strict transform of $\Fscr$ on $\tilde{X}$.
Unfortunately, unlike in the commutative case $\alpha^{-1}(\Fscr)$
depends on $\Fscr$ itself and not only on  the image of  $\Fscr$
in $\mod(X)/\Cscr_{f,p}$. Therefore we modify this definition as follows.
\begin{definitions}
  Assume that $\Fscr\in \trans_Y(X)$ is $Y$-torsion free. Then the
  \emph{strict transform} $\alpha_s^{-1}(\Fscr)$ of $\Fscr$ is defined
  as $(\alpha^{-1}\circ N_p)(\Fscr)$ (see \S\ref{ref:5.7b} for the
  definition of $N_p$).
\end{definitions}
\begin{lemmas} 
\label{ref:6.8.4a}
$\alpha^{-1}(\Fscr)$ and $\alpha^{-1}_s(\Fscr)$ are
  isomorphic to $\alpha^\ast\Fscr$ modulo $\alpha^{-1}(\Cscr_p)$.
\end{lemmas}
\begin{proof} This follows from the definition of $\alpha^{-1}(\Fscr)$
  together with the fact that
  $\HTor_1^{o_X}(\Fscr,o_X(nY) /I^n)\in\Cscr_p$ by Theorem
  \ref{ref:5.5.10a}. 
\end{proof}
\begin{propositions}
\label{ref:6.8.5a}
Identify $\tilde{Y}$ and $Y$ via the map $\beta$.  With this
identification we have
$R=\hat{\Oscr}_{Y,p}=\hat{\Oscr}_{Y,\beta^{-1}(p)}$.  As usual let $m$
be the maximal ideal of $R$. Assume that $\Fscr\in\trans_Y(X)$ is
$Y$-torsion free.  Then we have the following
\begin{enumerate} 
\item $\alpha_s^{-1}(\Fscr)$ is $\beta^{-1}(p)$-normalized.
\item One has
\begin{equation}
\label{ref:6.44a}
T_{\beta^{-1}(p)}(\alpha^{-1}_s(\Fscr))=mT_p(\Fscr)
\end{equation}
\end{enumerate}
\end{propositions}
\begin{proof}
  Assume that $\Fscr$ is $p$-normalized.
  
  The first claim is that the exact sequence coming from the inclusion
  $o_X(-Y)\r o_X$
\[
0\r \Fscr((n-1)Y)\r \Fscr(nY)\r \Fscr_Y (nY)\r 0
\]
  restricts to an exact sequence
(for $n\ge 1$)
\[
0\r \Fscr I^{n-1} \r \Fscr I^n \r \Fscr_Y I_{Y}^n\r 0
\]
 Writing out everything this amounts to
checking the exactness and well-definedness of the complex
\[
0\r \Fscr m_p\cdots m_{\tau^{-n+2} p} (-Y)\r 
\Fscr m_p\cdots m_{\tau^{-n+1} p} \r
\Fscr_Y m_{Y,p}\cdots m_{Y,\tau^{-n+1} p} \r 0
\]
That this is indeed a complex is clear. To show it is exact it
suffices to show that
\begin{equation}
\label{ref:6.45a}
\begin{split}
\length(\Fscr(-Y)/ \Fscr m_p\cdots m_{\tau^{-n+2} p} (-Y))=&
\length(\Fscr/\Fscr m_p\cdots m_{\tau^{-n+1} p}
)\\
&\quad+
\length(\Fscr_Y /\Fscr_Y m_{Y,p}\cdots m_{Y,\tau^{-n+1} p})
\end{split}
\end{equation}
Since the finite length modules involved all lie in $\Cscr_p$ we can
check this by completing.   $\Fscr$ is $p$-normalized so we have $
\hat{\Fscr}_{Y,p}=(\cdots 0, T_p(\Fscr),0\cdots) $ with the non-zero
entry occuring in position $0$.  From the structure of
$\hat{\Fscr}_p$ which is given by lemma \ref{ref:5.7.4a}. 
we   compute that
$(\Fscr m_p\cdots m_{\tau^{-n+1} p})\hat{}_p=\hat{\Fscr}_p m_0 m_{-1}\cdots
m_{-n+1}
$
is equal to
\[
(\cdots T_p(\Fscr)  \ mT_p(\Fscr) \cdots mT_p(\Fscr) \ 0\cdots
0\cdots)
\]
where now the first $mT_p(\Fscr)$ occurs in location $-n+1$ and the
first $0$ in location~$1$.
In a similar way we find that  $(\Fscr_Y m_{Y,p}\cdots m_{Y,\tau^{-n+1}
  p})\hat{}_p$ is equal to
\[
(\cdots 0\ m T_p(\Fscr)\ 0\cdots)
\]
Now \eqref{ref:5.39a} immediately follows.

We conclude that we have a complex of graded $\Dscr$-modules
\begin{equation}
\label{ref:6.46a}
0\r \Fscr\cdot \Dscr (-1) \r  \Fscr\cdot \Dscr \r \Fscr_Y\cdot
\Dscr_Y\r 
\end{equation}
exact in degree $\ge 1$. A simple local computation shows that the
canonical map $\Fscr_YI_{Y}\otimes_{o_Y}
I^{n-1}_{Y}\r\Fscr_YI^n_{Y}$ is an isomorphism and hence, up to right
bounded objects, $\Fscr_Y\cdot\Dscr_Y=(\Fscr_Y I_Y\otimes\Dscr_Y)(-1)$.

Applying $\pi$ to  \eqref{ref:6.46a} we find the exact sequence 
\[
0\r \alpha^{-1}(\Fscr)(-\tilde{Y})\r \alpha^{-1}(\Fscr)\r
\beta^\ast(\Fscr_Y I_{Y})(-\tilde{Y})\r 0
\]
Now $\Fscr_Y I_Y=\Fscr_Y m_p (Y)$ and using that $\beta$ is an isomorphism
together with the fact that twisting by $Y$ and $\tilde{Y}$ on finite
length modules 
simply
amounts to applying $\tau$ 
we obtain
$\beta^\ast(\Fscr_Y I_Y)(-\tilde{Y})=\beta^\ast(\Fscr_Y m_p)$.
This yields 
$\alpha^{-1}(\Fscr)_{\tilde{Y}}=\beta^\ast(\Fscr_Y m_p)$ which
implies what we have to show. 
\end{proof}
\begin{propositions} 
\label{ref:6.8.6a} Assume that $\Fscr\in\trans_Y(X)$ is $Y$-torsion
free and that $q$ is not in the $\tau$-orbit of $p$ (but that $q$ also
has infinite $\tau$-orbit). Then
$T_{\beta^{-1}q}(\alpha_s^{-1}(\Fscr))=T_q(\Fscr)$.
\end{propositions}
\begin{proof} This is a local verification as in Proposition
  \ref{ref:6.8.5a}, but easier. 
\end{proof}
We would also like to understand $\alpha^{-1}(\Fscr)$ if $\Fscr$ is
not normalized. Let us denote the quotient functor $\mod(\tilde{X})\r
\mod(\tilde{X})/\Cscr_{f}$ by
$\eta$. The relevant result is the following.
\begin{propositions}
\label{ref:6.8.7a}
 Assume $\Fscr\in \trans_Y(X)$ is $Y$-torsion free.
Then
 $\eta\alpha^{-1}_s(\Fscr)$ is the minimal subobject of
$\eta\alpha^{-1}(\Fscr)$ such that the quotient lies in the image of
 $\alpha^{-1}(\Cscr_p)$. 
\end{propositions}
\begin{proof}
  First note that if $\Fscr'\subset \Fscr$ such that
  $\Fscr/\Fscr'\in\Cscr_p$ then
  $\alpha^{-1}(\Fscr)/\alpha^{-1}(\Fscr')\in \alpha^{-1}(\Cscr_p)$ by
  lemma \ref{ref:6.8.4a} and Proposition
  \ref{ref:6.5.2a}.

Let $\Fscr_n$ be the $p$-normalization of $\Fscr$ (see
\S\ref{ref:5.7b}). With a similar local
verification as in the proof of Proposition \ref{ref:6.8.5a} we find that the
inclusion $\Fscr_n(-mY)\subset \Fscr_n$, for $m\ge 0$ induces an
isomorphism $\alpha^{-1}(\Fscr_n(-mY))\r \alpha^{-1}(\Fscr_n)$ modulo
$\Cscr_{f}$. 

If we take $m$ large enough then we will have an inclusion
$\Fscr_n(-mY)\subset \Fscr$. Together with the result of the previous
paragraph this  yields an inclusion of $\alpha^{-1}(\Fscr_n)\subset
\alpha^{-1}(\Fscr)$, if we work modulo
$\Cscr_{f}$. Furthermore since $\Fscr_n$ is isomorphic to
$\Fscr$ modulo $\Cscr_p$, this  inclusion becomes an isomorphism when
viewed modulo $\alpha^{-1}(\Cscr_p)$.

Now let us assume that $\Fscr$ is $p$-normal and let $\Gscr$ be a
subobject of $\alpha^{-1}(\Fscr)$. Put $\Hscr=\alpha_\ast(\Gscr)$. If
we choose $m$ large enough then $\Fscr(-mY)\subset \Hscr$ and hence we
have inclusions $\alpha^{-1}\Fscr(-mY)\subset
\alpha^{-1}(\Hscr)\subset \Gscr \subset \alpha^{-1}(\Fscr)$.  Since as
above $\alpha^{-1} \Fscr(-mY)= \alpha^{-1}(\Fscr)$ modulo
$\Cscr_{f}$, we conclude that $\Gscr=\alpha^{-1}(\Fscr)$ modulo
$\Cscr_{f}$.
\end{proof}
If $\Fscr\in\trans_Y(X)$ then let us write $l_p(\Fscr)$ for the minimal
$n$ such that $m^n T_p(\Fscr)=0$.
We have the following result.
\begin{propositions} 
\label{ref:6.8.8a}
Let $\Fscr\in \trans_{\tilde{Y}}(\tilde{X})$ be
  $\tilde{Y}$-torsion free.
Then 
  $l_p(\alpha_\ast(\Fscr))\le l_{\beta^{-1}(p)}(\Fscr)+1$. 
\end{propositions}
\begin{proof}
  Let $\Gscr$ be the $p$-normalization of $\alpha_\ast(\Fscr)$. By
  Proposition \ref{ref:6.8.7a} we have modulo $\Cscr_f$ an inclusion
  $\alpha^{-1}(\Gscr)\subset \alpha^{-1}(\alpha_\ast(\Fscr))\subset
  \Fscr$. Hence we have $l_p(\alpha_\ast(\Fscr))=l_p(\Gscr)\le
  l_{\beta^{-1}(p)}(\alpha^{-1} (\Gscr))+1\le
  l_{\beta^{-1}(p)}(\Fscr)+1 $ (the first inequality follows from
  \eqref{ref:6.44a}).
\end{proof}
A similar verification yields
\begin{propositions}
\label{ref:6.8.9a}
Let $\Fscr\in \trans_{\tilde{Y}}(\tilde{X})$ be $\tilde{Y}$-torsion
free and assume that $q$ is not in the $\tau$-orbit of $p$ (but that
$q$ also has infinite $\tau$-orbit). Then
$T_q(\alpha_\ast(\Fscr))=T_{\beta^{-1}(q)}(\Fscr)$. 
\end{propositions}

\subsection{A result on $K_0$ of some categories}
\label{ref:6.9b}
In this section we use the results on strict transform to prove a
technical result which will be used later. \emph{We now assume that
  the $\tau$-orbit of every point on $Y$ has infinite order}. Thus in
particular $Y$ is smoooth.

For a collection of natural numbers $z=(z_o)_{o\in Y/\langle\tau\rangle}$
let us define $\trans_{Y,z}(X)$ as the full subcategory of $\trans_Y(X)$
consisting of objects $\Fscr$ such that $l_q(\Fscr)\le z_{\bar{q}}$
for $q\in Y$. We
define $M_z(X)$ as $\trans_{Y,z}(X)/\Cscr_f$. 

For $o\in Y/\langle\tau \rangle$ let us define $e_o$ by
$(\delta_{o,o'})_{o'\in Y/\langle\tau\rangle}$. Our aim is to prove the
following theorem
\begin{theorems}
\label{ref:6.9.1a}
 Let $\tilde{X}$ be the blowup of $X$ in $p$.
One has the following relation.
\[
K_0(M_z(\tilde{X}))\cong
\begin{cases}
 K_0(M_{z+e_{\bar{p}}}(X))\oplus\ZZ &\text{if $z_{\bar{p}}\ge 1 $}\\
 K_0(M_{z+e_{\bar{p}}}(X))&\text{if $z_{\bar{p}}=0$}
\end{cases}
\]
\end{theorems}
\begin{proof}
  By lemma \ref{ref:6.8.1a} $\alpha^{-1}(\Cscr_p)\cap \mod(\tilde{X})$
  is contained in $\trans_{\tilde{Y}}(\tilde{X})$.  Let us write
  $\Tscr=(\alpha^{-1}(\Cscr_p)\cap
  \mod(\tilde{X}))/\Cscr_{f,\beta^{-1}(p)}\subset
  \trans_{\tilde{Y}}(\tilde{X})/\Cscr_f$ and let us temporarily use the
  notation ${\bar M}_z(\tilde{X})=M_z(\tilde{X})/ (\Tscr\cap
  M_z(\tilde{X}))$.

We first show that
$\alpha^{-1}$, $\alpha_s^\ast$ define inverse equivalence 
  between $M_{z+e_{\bar{p}}}(X)$ and ${\bar M}_z(\tilde{X})$.
  We know already that $\alpha^\ast$ and $\alpha_\ast$ define inverse
  equivalences between $\mod(X)/\Cscr_{f,p}$ and
  $\mod(\tilde{X})/(\alpha^{-1}(\Cscr_p)\cap\mod(\tilde{X}))$. Since
  $\alpha_s^{-1}$ and $\alpha^\ast$ take the same values modulo
  $\alpha_s^{-1}(\Cscr_p)$,  it follows that $\alpha_s^{-1}$ and
  $\alpha_\ast$ also define inverse equivalences. Furthermore it
  follows from Propositions
  \ref{ref:6.8.5a},\ref{ref:6.8.6a},\ref{ref:6.8.8a},\ref{ref:6.8.9a} that
  $\alpha^{-1}$ and $\alpha_\ast$ restrict to equivalences between
  $M_{z+e_{\bar{p}}}(X)$ and ${\bar M}_z(\tilde{X})$.

We now compute the $K_0(\bar{M}_z(\tilde{X}))$ by the localization
sequence. 
\begin{equation}
\label{ref:6.47a}
K_0(\Tscr\cap
  M_z(\tilde{X}))\r K_0(M_z(\tilde{X}))\r 
  K_0(\bar{M}_z(\tilde{X}))\r 0
\end{equation}
If $z_{\bar{p}}=0$ then $\Tscr\cap
  M_z(\tilde{X})=0$ so in fact $M_z(\tilde{X})=\bar{M}_z(\tilde{X})$
  and we have the corresponding equality on $K_0$-groups. This proves
  what we need in the case $z_{\bar{p}}=0$.

Let us now consider the case $z_{\bar{p}}\ge 1$. In that case
$\mod(L)/\Cscr_{f,\beta^{-1}(p)}\subset  M_z(\tilde{X})$ and hence the
exact sequence \eqref{ref:6.47a} becomes
\[
K_0(\mod(L)/\Cscr_{f,\beta^{-1}(p)})\xrightarrow{\delta} K_0(M_z(\tilde{X}))\r 
  K_0(\bar{M}_z(\tilde{X}))\r 0
\]
where $\delta$ is the natural map. From Proposition
\ref{ref:6.7.1a} it easily follows that
$K_0(\mod(\tilde{X})/\Cscr_{f,\beta^{-1}(p)})=\ZZ$, so it remains to
show that $\delta$ is injective.
For $\Fscr\in \trans_{\tilde{Y}}(\tilde{X})$ let $t(\Fscr)$ be the
  degree of $\Div(\Fscr)$ (\S \ref{ref:5.7b}).

Clearly $t$ is additive on short exact sequences and
furthermore  it factors through $\trans_{\tilde{Y}}(\tilde{X})/
\Cscr_{f,\beta^{-1}(p)}$. The generator of
$K_0(\mod(L)/\Cscr_{f,\beta^{-1}(p)})$ is $[\bar{\Oscr}_L]$ and since
obviously $t(\bar{\Oscr}_L)=1$ we conclude that the image of
$[\bar{\Oscr}_L]$ cannot be zero under $\alpha$. Hence $\alpha$ must
be injective. 

We conclude that if $z_{\bar{p}}\ge 1$ then
$K_0(M_z(\tilde{X}))=K_0(\bar{M}_z(\tilde{X}))\oplus \ZZ$. This
finishes the proof.
\end{proof}

\section{Derived categories}
\label{ref:7a}

\subsection{Generalities}
\label{ref:7.1a}
If $Z$ is a quasi-scheme then in the sequel we use the notation
$D^\ast(Z)$ with $\ast=\phi,+,-,b$ for the standard derived categories of
$\Qch(Z)$. If $Z$ is  noetherian then we use $D_f^\ast(Z)$ for
the full subcategories of $D^\ast(Z)$ whose objects have coherent
cohomology.

If $\alpha:Y\r X$ is a map between quasi-schemes then since $\Qch(Y)$
has enough injectives it is trivial to define $R\alpha_\ast$. This is
unfortunately not the case for $L\alpha^\ast$.

We will say that $\Qch(X)$ has enough acyclic objects for
$\alpha^\ast$ if every object in $\Qch(X)$ is a quotient of an object
$\Uscr$ in $\Qch(X)$ such that $L_i\alpha^\ast \Uscr=0$ for $i>0$ (cfr
\S\ref{ref:3.10b}). The following lemma is easy (see \cite{RD}).
\begin{lemmas}
\label{ref:7.1.1a}
Assume that $\Qch(X)$ has enough acyclic objects for
$\alpha^\ast$. Then $L\alpha^\ast$ exists on $D^-(X)$ and can be
computed by acyclic resolutions. If $\alpha^\ast$ has finite
cohomological dimension then the same is true with $D(X)$ replacing $D^-(X)$.
 \end{lemmas}
Assume that $i:Y\r X$ makes $Y$ into a divisor in $X$ (\S\ref{ref:3.7b})
and  denote by $i^!$ the functor $\HHom(o_Y,-)$. 
From the resolution of $o_Y$ by two invertible bimodules
\begin{equation}
\label{ref:7.1b}
0\r o_X(-Y)\r o_X \r o_Y\r 0
\end{equation}
 it is clear that $i^!$ has
cohomological dimension one. We will need the following lemma.
\begin{lemmas}
Assume that $\Qch(X)$ has enough $Y$-torsion free objects. Then for
$\Fscr\in \coh(X)$ we have
\begin{equation}
\label{ref:7.2a}
Ri^!\Fscr=(Li^\ast\Fscr)(Y)[-1]
\end{equation}
\end{lemmas}
\begin{proof} For $\Fscr$ a torsion free object in $\Qch(X)$  this follows
  immediately by 
  applying $\Hom(-,\Fscr)$ to \eqref{ref:7.1b}. The general case follows
  from \cite[Prop. 7.4]{RD}.
\end{proof}

\subsection{Admissible compositions of maps between quasi-schemes}
In this section we point out a few problems with compositions of maps
between quasi-schemes which have no equivalent in the commutative
case. These problems will be solved in all concrete cases, but they
nevertheless represent a nuissance.

Assume that we have a composition of maps between quasi-schemes.
\begin{equation}
\label{ref:7.3a}
Z\xrightarrow{\beta} Y\xrightarrow{\alpha} X
\end{equation}
From the very definition of a morphism it is clear that
$\alpha_\ast\beta_\ast=(\alpha\beta)_\ast$. However in  contrast with
the commutative case there is no reason why the natural map
\begin{equation}
\label{ref:7.4a}
R(\alpha\beta)_\ast\rightarrow R\alpha_\ast\,R\beta_\ast
\end{equation}
should be an isomorphism 
(for trivial reasons, see the discussion after  lemma \ref{ref:7.2.7a} below).

This leads to the following definition.
\begin{definitions}
\label{ref:7.2.1a}
  A composition $(\alpha,\beta)$ as in \eqref{ref:7.3a} is
  \emph{admissible} if the identity \eqref{ref:7.4a} holds (as functors
  from $D^+(Z)$ to $D^+(X)$).
\end{definitions}
The following lemma is trivial.
\begin{lemmas}
\label{ref:7.2.2a}
The composition $(\alpha,\beta)$ is admissible if and only if for
every injective $E\in\Qch(Z)$ one has that $\alpha_\ast E$ is acyclic
for $\beta_\ast$.
\end{lemmas}
For simplicity we will use some variations on Definition
\ref{ref:7.2.1a}. If we work in $\Qsch/X$ (i.e. somewhat imprecisely~: ``if
$X$ is a fixed base quasi-schema''), then we will say that the map $\alpha$
is admissible if \eqref{ref:7.4a} holds.

Similarly if $\alpha:(Y,\Oscr_Y)\r (X,\Oscr_X)$ is a map of enriched
quasi-schemes then we will say that $\alpha$ is admissible if
\[
R\Gamma(X,-)\circ R\alpha_\ast=R\Gamma(Y,-)
\]
These conventions are to a certain extent compatible as seen by the
following lemma.
\begin{lemmas}
Assume that we have maps
\[
Z\xrightarrow{\beta} Y\xrightarrow{\alpha} \Spec R
\]
Put $\Oscr_Y=\alpha^\ast R$, $\Oscr_Z=\beta^\ast \Oscr_Y$. Then the
compositon $(\alpha,\beta)$ is admissible if and only if the induced
map of enriched quasi-schemes $\alpha:(Z,\Oscr_Z)\r (Y,\Oscr_Y)$ is
admissible.  
\end{lemmas}
Below we will give a few adhoc criteria which will allow us to show
that certain maps/compositions are admissible. They are all tautologies.
\begin{lemmas}
\label{ref:7.2.4a}
Assume that we have a composition as in \eqref{ref:7.3a}, but this time
assume that we are working  in $\Qsch/W$ for some base quasi-scheme
$W$. If the maps $\alpha$, $\beta$ and the composition
$(\alpha,\beta)$ are admissible then so is the map $\alpha\beta$.
\end{lemmas}
\begin{lemmas}
\label{ref:7.2.5a}
Assume that $\alpha:(Y,\Oscr_Y)\r (X,\Oscr_X)$ is a map of enriched
quasi-schemes. Then $\alpha$ is admissible if and only if
$L_i\alpha^\ast\Oscr_X=0$ for $i>0$.
\end{lemmas}
\begin{corollarys} 
\label{ref:7.2.6a}
Let $(X,\Oscr_X)$ be an enriched quasi-scheme.
Assume that $\Ascr\in \Alg(X)$. Then $\Spec \Ascr \r X$ is admissible
if and only if $\HTor^i_{o_X}(\Oscr_X,\Ascr)=0$ for $i>0$.
\end{corollarys}
\begin{corollarys}
\label{ref:7.2.7a}
Assume that $\alpha:(Y,\Oscr_Y)\r (X,\Oscr_X)$ is a map of enriched
quasi-schemes which embeds $Y$ as a divisor in $X$ (\S\ref{ref:3.7b}). Then
$\alpha$ is admisible.
\end{corollarys}
\begin{proof} This is a special case of the previous corollary, since
  $Y=\Spec o_Y$, where we consider $o_Y$ as an $o_X$ algebra.
\end{proof}
This  corollary indicates how to make a non-admissible map. Take
for example commutative schemes $Y,X$, $Y$ being a Cartier divisor in
$X$ and change the structure sheaf of $X$ into one which has $Y$-torsion.

Below quasi-schemes are often defined as $\Proj$'s of algebras. The
following lemma is the obvious analogue of Corollary \ref{ref:7.2.6a}.
\begin{lemmas}
\label{ref:7.2.8a}
  Let $(X,\Oscr_X)$ be an enriched quasi-scheme.  Assume that $\Ascr$
  is a noetherian graded algebra on $X$. Then $\Proj \Ascr \r X$ is
  admissible if $\HTor^{o_X}_i(\Oscr_X,\Ascr)=0$ is right bounded
  for $i>0$.
\end{lemmas}
\begin{lemmas}
\label{ref:7.2.9a}
Let $f:\Ascr\r\Bscr$ be as in Proposition \ref{ref:3.9.10a}. Let
$\bar{f}:\Proj\Bscr\r\Proj\Ascr$ be as in Proposition
\ref{ref:3.9.11a}. Then $\bar{f}$ defines an admissible morphism
of quasi-schemes in $\Qsch/X$.
\end{lemmas}
\begin{proof}
This follows from Proposition \ref{ref:3.9.12a}.
\end{proof}
The following  is also standard.
\begin{propositions} Assume that we have quasi-schemes and maps as in
  \eqref{ref:7.3a}. Assume that the composition $(\alpha,\beta)$ is
  admissible and that the maps $\alpha$, $\beta$, $\alpha\beta$
  satisfy the conditions of lemma \ref{ref:7.1.1a}. Then we have
  $L(\alpha\beta)^\ast= L\beta^\ast\,L\alpha^\ast$ (as functors from
  $D^-(X)$ to $D^-(Z)$).
\end{propositions}
\begin{proof}
It is easy to see that we have at least a natural transformation
\[
L\beta^\ast\,
  L\alpha^\ast\r L(\alpha\beta)^\ast
\]
In fact this can be deduced from the very definition of derived functors
(see \cite[\S5]{RD}).

Now it is clear that if $E$ runs trough the injectives in $\Qch(Z)$
then $\Hom_{D(Z)}(-,E)$ is a conservative system of functors on
$D(Z)$. Let $\Mscr\in D^-(X)$. Then we have by adjunction and admissibility
\begin{align*}
\Hom_{D(Z)}(L\beta^\ast\,
  L\alpha^\ast\Mscr,E)&=\Hom_{D(X)} (\Mscr,R\alpha_\ast\, R\beta_\ast E)\\
&=\Hom_{D(X)}(\Mscr,R(\alpha\beta)_\ast E)\\
&=\Hom_{D(Z)}(L(\alpha\beta)^\ast\Mscr,E)\qed
\end{align*}
\def\qed{}\end{proof}
\begin{remarks} It is easy to see that if we have an admissible
  composition as in \eqref{ref:7.3a} then similar results as the ones
  presented  above remain valid on
  unbounded derived categories provided suitable functors have
  finite cohomological dimension. In the sequel we will tacitly use such results,
  leaving the obvious proofs to the reader.
\end{remarks} 

We can now show that some commutative diagrams of quasi-schemes
we encountered previously  consist of admissible maps and admissible
compositions. 
\begin{lemmas} 
\label{ref:7.2.12a}
All maps and compositions of maps in   diagram \eqref{ref:6.20a}
  are admissible.
\end{lemmas}
\begin{proof}
$i$ and $j$ are admissible because they are divisors (see Theorem \ref{ref:6.3.1a}.8). $\beta$ is admissible since it is a
map between commutative schemes. To prove that $\alpha$ is admissible it
suffices by lemma \ref{ref:7.2.8a} to show that
$\HTor_i^{o_X}(\Oscr_X,\Dscr)=0$ for $i>0$. More precisely we have to
show that $\HTor_i^{o_X}(\Oscr_X,I^n_p)=0$ for $i>0$. Since
$I^n_p\subset o_X(nY)$ and the quotient is in $\tilde{C}_{f,p}$ this
follows from the compatibility of
$\HTor$ with completion (see Theorem \ref{ref:5.5.10a}).
\end{proof}

\begin{lemmas}
\label{ref:7.2.13a}
All maps and compositions in   diagram \eqref{ref:6.34a} are admissible.
\end{lemmas}
\begin{proof}
The admissibility of $(\alpha,u)$ follows from lemma \ref{ref:7.2.9a}. The
other compositions and maps we leave to the reader.
\end{proof}

\section{The derived category of a non-commutative blowup}
\label{ref:8a}
\subsection{The formalism of semi-orthogonal decompositions}
\label{ref:8.1a}
The material in this section is taken from \cite{Bondal2}. 

Let $\Ascr$ be a triangulated category and let $\Bscr$, $\Cscr$ be two
strict full triangulated subcategories of $\Ascr$.  $(\Bscr,\Cscr)$ is said to be a
\emph{semi-orthogonal pair} if $\Hom_\Ascr(B,C)=0$ for $B\in\Bscr$ and $C\in
\Cscr$. Define
\[
\Bscr^\perp=\{A\in\Ascr\mid \forall B\in\Bscr :\Hom_\Ascr(B,A)=0\}
\]
${}^\perp\Cscr$ is defined similarly.

The following result is a slight variation of the statement of
\cite[Lemma 3.1]{Bondal2}.
\begin{lemmas}
\label{ref:8.1.1a}
 The following statements are equivalent for a
  semi-orthogonal pair $(\Bscr,\Cscr)$.
\begin{enumerate}
\item
$\Bscr$ and $\Cscr$ generate $\Ascr$.
\item For every $A\in\Ascr$ there exists a distinguished triangle $B\r
  A\r C$ with $B\in\Bscr$ and $C\in\Cscr$.
\item $\Cscr=\Bscr^\perp$ and  the inclusion functor
  $i_\ast:\Bscr\r\Ascr$ has a right adjoint $i^!:\Ascr\r \Bscr$.
\item $\Bscr={}^\perp\Cscr$ and  the inclusion functor $j_\ast:\Cscr\r
  \Ascr$ has a left adjoint  $j^\ast:\Ascr\r \Cscr$
\end{enumerate}
If one of these conditions holds then the triangles in 2. are unique
up to unique isomorphism. They are necessarily of the form 
\[
i^! A\r A \r j^\ast A
\]
where the maps are obtained by adjointness from the identity maps
$i^!A\r i^! A$ and $j^\ast A \r j^\ast A$. In particular triangles as
in 2. are functorial.
\end{lemmas}
\begin{remarks} The notations $(i_*,i^!, j_*, j^*)$ are purely symbolic
  and shouldn't be interpreted as direct and inverse images. In fact
  in the main application  below (Theorem \ref{ref:8.4.1a}) $i_*$ will be
  given by an inverse image!
\end{remarks}
If a pair $(\Bscr,\Cscr)$ satisfies one of the conditions of the
previous lemma then we say that it is a \emph{semi-orthogonal
decomposition} of $\Ascr$.
For further reference we note the following diagram of arrows
\begin{equation}
\label{ref:8.1b}
\Cscr \begin{array}{c} j_\ast\\ \rightarrow\\ \leftarrow
  \\ j^\ast\end{array}
 \Ascr 
\begin{array}{c}\displaystyle  i^!\\ \rightarrow\\ 
  \leftarrow \\i_\ast
\end{array}
 \Bscr
\end{equation}
In the following lemma we give some relations between these arrows.
\begin{lemmas} One has~:
\begin{gather*}
 i^! i_\ast=\Id_\Bscr
\\
j^\ast j_\ast=\Id_\Cscr\\
j^\ast i_\ast=0\\
i^! j_\ast=0\\
\end{gather*}
\end{lemmas}  
In the sequel we will slightly extend the meaning of the notion of
semi-orthogonality. Assume that we have functors
\[
\Cscr\xrightarrow{j_\ast} \Ascr \xleftarrow{i_\ast} \Bscr
\]
which are fully faithful. Assume  furthermore that the essential images of
$\Bscr$ and $\Cscr$ in $\Ascr$ are  semi-orthogonal in
$\Ascr$. Then, if no confusion can arise, we wil also call $(\Bscr,\Cscr)$
a semi-orthogonal pair in $\Ascr$. Similarly for a semi-orthogonal
decomposition.

Semi-orthogonal decompositions can be constructed starting from a pair
of adjoint functors. For an arbitrary functor $F:\Ascr\r \Bscr$
between additive categories let us define $\ker F$ as the full subcategory
of $\Ascr$ whose objects satisfy $F(A)=0$.
\begin{lemmas}
\label{ref:8.1.4a}
 Assume that we have triangulated categories
  $\Ascr,\Bscr$ and a pair of adjoint functors $i_\ast:\Bscr\r \Ascr$,
  $i^!:\Ascr\r \Bscr$ such that $i^! i_\ast=\Id_\Bscr$. Then
  $i_\ast$ is an embedding of $\Bscr$ in $\Ascr$. The
  corresponding semi-orthogonal decomposition is given by 
$(\Bscr,\ker i^!)$.
\end{lemmas}

\subsection{Generalities}
\label{ref:8.2a}
In this section the notations and hypotheses of Section
\ref{ref:6a} will be in force. Our aim  is to give
a non-commutative version of a well-known theorem by Orlov
\cite{Orlov} which relates $D(\tilde{X})$ to $D(X)$.  To this end we
need the adjoint functors $R\alpha_\ast$ and $L\alpha^\ast$. In
particular we need that $X$ has enough acyclic objects for
$\alpha^\ast$ (by lemma \ref{ref:7.1.1a}).

Therefore  at this point we introduce an extra hypothesis
which will always hold in the applications. 

Let us call an object $\Lscr$ in $\coh(X)$ a \emph{line bundle} around
$Y$ if the map $\Lscr(-Y)\r \Lscr$ is injective and if
$\Lscr/\Lscr(-Y)$ is a line bundle on $Y$. Note that $\Oscr_X$ itself
is a line bundle on $Y$. We denote the additive category of objects
which are direct sums of line bundles on $Y$ by $\Vscr$.

{\def\thehypothesis{(**)}
\begin{hypothesis} Every object in $\Qch(X)$ is a quotient of an
  object in $\Vscr$.
\end{hypothesis}
}

 The following lemma is left to the reader.
\begin{lemmas} 
\label{ref:8.2.1a}
Assume $\Escr\in \Vscr$.
\begin{enumerate}
\item
$\HTor_i^{o_X}(\Escr,o_Y)=0$ for $i\neq 0$.
\item
$\HTor_i^{o_X}(\Escr,o_q)=0$ for $i\neq 0$.
\item
 $\HTor_i^{o_X}(\Escr,\Dscr)=0$ for $i\neq 0$.
\item
$\HTor_i^{o_X}(\Escr,\Dscr_Y)=0$ for $i\neq 0$.
\item
$L_i\alpha^\ast \Escr=0$ for  $i\neq 0$.
\end{enumerate}
\end{lemmas}
From this lemma together with Prop. \ref{ref:6.1.2a} it is easy
to see that $-\Lotimes_{o_X} \Dscr$, $-\Lotimes_{o_X} \Dscr_Y$,
$R\alpha^\ast$ etc\dots can be defined in the usual way \cite{RD} by
considering resolutions in $\Vscr$.

For further reference we recall the following.
\begin{lemmas}
\label{ref:8.2.2a} $R\alpha_\ast$ and $L\alpha^\ast$
 have finite cohomological dimension and commute with
direct sums.
\end{lemmas} 
\begin{proof} 
  The fact that $R\alpha_\ast$ and $L\alpha^\ast$ have finite
  cohomological dimension is proved in Theorem
  \ref{ref:6.3.1a}.

  Since $L\alpha^\ast$ is the left adjoint to $R\alpha^\ast$ it is
  clear that it is compatible with direct sums. Hence let us
  concentrate on $R\alpha_\ast$. From the discussion in Section
  \ref{ref:3.8b}  it follows that
  $R\alpha_\ast$ is equal to $R\omega(-)_0$ and furthermore for
  $\Mscr\in\Gr(\Dscr)$ there is a triangle
\begin{equation}
\label{ref:8.2b}
\Atriangle<1`-1`1; >[\Mscr`R\tau\Mscr`R\omega(\pi\Mscr);``]
\end{equation}
Thus it is sufficient to show that $R\tau$ commutes with direct
sums. According to Proposition \ref{ref:3.8.4a}
\[
\tau \Mscr=\injlim \HHom_{\Dscr}(\Dscr/\Dscr_{\ge n},\Mscr)
\]
According to Proposition \ref{ref:6.1.1a}, $I^n_p$ is a
coherent $o_X-o_X$ bimodule \S\ref{ref:3a} for all $n$. Hence this
holds also for $\Dscr/\Dscr_{\ge n}$. But then it is easy to see that
$\Dscr/\Dscr_{\ge n}$ is also coherent as $\Dscr-\Dscr$-bimodule. Thus
$\HHom_{\Dscr}(\Dscr/\Dscr_{\ge n},-)$ commutes with direct sums and
hence so does $\tau$. Now $\tau$ has finite cohomological dimension
(by \eqref{ref:8.2b} this is the same as $\alpha_\ast$ having finite
cohomological dimension).  It is now easy to see that the standard way
of defining $R\tau$ on $D(X)$ \cite{RD} is compatible with direct
sums.
\end{proof}
\subsection{Computation of some derived functors}
\label{ref:8.3a}

We need the following proposition.
\begin{propositions} 
\label{ref:8.3.1a}
Let $n\ge -1$.
\begin{enumerate} 
\item $R\alpha_\ast (\Oscr_L(n))= (\Dscr/m_{\tau p}\Dscr)_{n,o_X}$.
\item
 Assume that $\Escr\in \Vscr$. Then
$R\alpha_\ast(( \alpha^\ast
\Escr)(n))=\Escr\otimes_{o_X}\tilde{\Dscr}_n$ (cfr. \eqref{ref:6.15a}).
\end{enumerate}
\end{propositions}
\begin{proof}
  Since $\Dscr$ satisfies $\chi$ (Theorem
  \ref{ref:6.2.2a}) it is clear that the proposition is true for
  $n\gg 0$ (lemma \ref{ref:3.9.3a}). Our strategy will
  now be to use descending induction on $n$. To this end it is
  convenient to treat the cases $\tau p=p$ and $\tau p\neq p$
  separately. In the first case $L$ is divisor in $\tilde{X}$
  isomorphic to $\PP^1$ so we can use this. In the second case $p$ is
  smooth on $Y$ (Thm \ref{ref:5.1.4a}) and so 
  $\beta:\tilde{Y}\r Y$ is an isomorphism. It follows from this that a
  suitable analog of the proposition  is trivially true for $Y$. We can
  then lift this to $\tilde{X}$.
 
In the proof below we will write $\Escr_{\tau p}=\Escr\otimes_{o_X}
o_{\tau p}$ and similarly $\Escr_Y=\Escr\otimes_{o_X}o_Y$.

\begin{case} $\tau p=p$.
By admissibility we have
  $R\alpha_\ast
  \Oscr_L(n)=R\gamma_\ast \Oscr_L(n)$ (where $\gamma$ is as in diagram
  \eqref{ref:6.34a}). Thus 1. follows easily
  from the corresponding result for $\PP^1$.

For 2. we use the exact sequence obtained from
\eqref{ref:6.33a}
\begin{equation}
\label{ref:8.3b}
0\r \alpha^\ast(\Escr(-Y))(1)\r \alpha^\ast \Escr\r  \alpha^\ast(
\Escr_{\tau p})\r 0
\end{equation}
Since $\alpha^\ast(
\Escr_{\tau p})$ is a direct sum of copies of 
$\Oscr_L$ we deduce from 1. that 
\[
R^i\alpha_\ast (\alpha^\ast(
\Escr_{\tau p}))=
\begin{cases} 
\Escr\otimes_{o_X} (\Dscr/m_p\Dscr)_n &\text{if $i=0$}\\
0&\text{otherwise}
\end{cases}
\]
Assume that 2. is true for a
certain $n\ge 0$. Since the autoequivalence $-\otimes_{o_X} o_X(Y)$ is
compatible with $\alpha_\ast$ and $\alpha^\ast$ ($p$ is a fixed point)
2. will also be true for this $n$ if we
replace $\Escr$ by $\Escr(mY)$ for arbitrary $m$.

Tensoring  \eqref{ref:8.3b} with
$o_{\tilde X}(n-1)$
and applying the long exact sequence for $R\alpha_\ast$ yields
$R^i\alpha_\ast ((\alpha^\ast\Escr)(n-1))=0$ for $i>0$ and for $i=0$ we obtain
an exact sequence
\begin{equation}
\label{ref:8.4a}
0\r \Escr(-Y)\otimes_{o_X} \Dscr_n \r \alpha_\ast(\alpha^\ast (\Escr
)(n-1)) \r \Escr\otimes_{o_X} (\Dscr/m_p\Dscr)_n\r 0
\end{equation}
On the other hand if we tensor the exact sequence \eqref{ref:6.33a} on the left
with $\Escr$  and take the part of degree $n$ then we obtain an exact
sequence
\[
0\r \Escr(-Y)\otimes_{o_X} \Dscr_n\r \Escr \otimes_{o_X} \Dscr_{n-1}
\r \Escr\otimes_{o_X} (\Dscr/m_p\Dscr)_n\r 0
\]

Comparing with \eqref{ref:8.4a} using the five lemma completes the
proof of the proposition in the case $\tau p=p$.
\end{case}
\begin{case}$\tau p\neq p$. Now
  $p$ is smooth on $Y$ (Thm. \ref{ref:8.4.1a}) and the map
  $\beta:\tilde{Y}\r Y$ is an isomorphism.  
    From \eqref{ref:6.16a} we obtain an
  exact sequence
\begin{equation}
\label{ref:8.5a}
0\r \alpha^\ast(\Escr)(-1)\r \alpha^\ast(\Escr) \r \beta^\ast(\Escr_Y)\r 0
\end{equation}
By admissibility $R^i\alpha_\ast
((\beta^\ast\Escr_Y)(n))= R^i\beta_\ast( (\beta^\ast\Escr_Y)(n))$. For
$n\ge 0$ we have
\[
R^i\beta_\ast((\beta^\ast\Escr_Y)(n))=
\begin{cases}
\Escr_Y\otimes_{o_Y} \Dscr_{Y,n}& \text{if $i=0$}\\
0&\text{if $i>0$}
\end{cases}
\]
The case $i=0$ follows from \eqref{ref:6.21a}. The
case $i>0$ follows from the fact that $\beta$ is an isomorphism.

Assume that 2. is true for a certain $n$. Tensor \eqref{ref:8.5a} with
$o_{\tilde X}(n)$ and apply $R\alpha_\ast$. We obtain $R^i\alpha_\ast
((\alpha^\ast\Escr)(n-1))=0$ for $i\ge 2$ and an exact sequence
\begin{equation}
\label{ref:8.6a}
0\r \alpha_\ast((\alpha^\ast\Escr)(n-1))
\r
\Escr\otimes_{o_X}\Dscr_n
\r
\Escr\otimes_{o_X}\Dscr_{Y,n}
\r
R^1\alpha_\ast((\alpha^\ast\Escr)(n-1))
\r
0\end{equation}
Tensoring \eqref{ref:6.16a} on the left with $\Escr$ and taking the part of
degree $n$ yields an exact sequence
\[
0\r \Escr\otimes_{o_X}\Dscr_{n-1}\r
\Escr\otimes_{o_X} \Dscr_n\r
\Escr\otimes_{o_X} \Dscr_{Y,n}\r
\]
Comparing with \eqref{ref:8.6a} using the five lemma yields 2.
in this case.

To prove 1. we twist the exact sequence \eqref{ref:8.3b} by
$o_{\tilde{X}}(n)$ and put $\Escr=\Oscr_X$. Applying $R\alpha_\ast$ then
yields what we want.
\qed\end{case} \def\qed{}
\end{proof}

\begin{propositions}
\label{ref:8.3.2a}
Assume that $q\in Y$. If $q\neq \tau p$ then $L\alpha^\ast
\Oscr_q=\Oscr_{q'}$ with $q'$ such that $\beta(q')=q$ (such a $q'$ is
unique!).

For $q=\tau p$  we have
\[
L^i \alpha^\ast \Oscr_{\tau p}
=
\begin{cases}
\Oscr_L &\text{if $i=0$}\\
\Oscr_L(-1)&\text{if $i=1$}\\
0&\text{otherwise}
\end{cases}
\]
\end{propositions}

\begin{proof}
If $q\not \in O_\tau (p)$ then the result is an easy excercise, for
example using the fact that $I^n\subset o_X(nY)$ with quotient in
$\tilde{\Cscr}_{f,p}$ (see Prop. \ref{ref:6.1.1a}). 
Hence assume $q\in O_\tau(p)$.
We have by lemma \ref{ref:3.10.2a}
\[
L^i\alpha^\ast \Oscr_q=\pi \HTor^i_{o_X}(\Oscr_q,\Dscr)
\]
Now $\HTor^i_{o_X}(\Oscr_q,\Dscr)\in\Cscr_p(\Dscr)$ by Propositions
\ref{ref:6.1.1a} and \ref{ref:6.1.2a}. Since ``$\HTor$''
is compatible with completion (Theorem \ref{ref:5.5.10a})
it suffices to compute $\Tor^i_{C_p}(\hat{\Oscr}_q,\hat{\Dscr})$. This
is a mildly tedious calculation which we leave to the reader. One way
to proceed is to use Proposition \ref{ref:5.2.2a}
to construct a projective resolution of $\hat{\Oscr}_q$ of length 2
over $C_p$.
\end{proof}
\begin{lemmas}
\label{ref:8.3.3a}
\begin{enumerate} 
\item $\Hom_{o_{\tilde{X}}}(\Oscr_L,-)$ has finite cohomological dimension.
\item
The functor $\RHom_{o_{\tilde{X}}}(\Oscr_L,-)$ is defined on
$D(\tilde{X})$. For $\Mscr\in D(X)$ there is a triangle
\begin{equation}
\label{ref:8.7a}
\Atriangle<1`-1`1; >[R\alpha_\ast(\Mscr(- 1))(Y)`
\Oscr_{\tau p}\otimes_k\RHom_{D(\tilde{X})}(\Oscr_L,\Mscr)`R\alpha_\ast\Mscr;``]
\end{equation}
\item $\RHom(\Oscr_L,-)$ commutes with direct sums on $D(X)$.
\item $\RHom(\Oscr_L,\Oscr_L)=k$.
\end{enumerate}
\end{lemmas}

\begin{proof}
\begin{enumerate}
\item
Assume that $E$ is an injective object in $\Qch(\tilde{X})$. We have
\[
\Oscr_{\tau p}\otimes_k  \Hom_{o_{\tilde{X}}}(\Oscr_L,E)=
\HHom_{\Gr(\Dscr)} (\Dscr/m_{\tau
    p}\Dscr,\omega E)
\]
Now we use the exact sequence \eqref{ref:6.33a}. From exactness of $\pi$
it follows that $\omega E$ is injective. Furthermore assume that
$\Tscr$ is a right bounded $o_X-\Dscr$-module. Then adjointness yields
that $\Hom_{o_X}(-,\underline{\HHom}_{\Dscr}(\Tscr,\omega E))$
vanishes. Thus we deduce that $\underline{\HHom}_{\Dscr}(-,\omega E))$
is zero on right bounded $o_X-\Dscr$-bimodules. Hence we obtain from
\eqref{ref:6.33a} an exact sequence in $\Qch(X)$.
\[
0\r \HHom_{\Gr(\Dscr)}(\Dscr/m_{\tau
    p}\Dscr,\omega E)\r 
(\omega E)_0\r
(\omega E)_{-1}(Y)\r 
0
\]
Which can be rewritten as
\begin{equation}
\label{ref:8.8a}
0\r 
\Oscr_{\tau p}\otimes_k  \Hom_{o_{\tilde{X}}}(\Oscr_L,E)
\r 
\alpha_\ast E\r
\alpha_\ast (E(-1))(Y)\r 
0
\end{equation}
The fact that $\alpha_\ast$ has finite cohomogical dimension by
\ref{ref:6.3.1a} yields what we want.

\item The fact that $\RHHom_{o_{\tilde{X}}}(\Oscr_L,-)$ is defined on the
  whole of $D(X)$ follows from 1.  using \cite{RD}.  Now since $\alpha_\ast$
  also has finite cohomological dimension, every object in
  $\Qch(\tilde{X})$ has a resolution by objects which are acyclic for
  both $\alpha_\ast$  and $\Hom_{o_{\tilde{X}}}(\Oscr_L,-)$ (using the
  methods of \cite{RD}).
If $F$ is acyclic for
  $\Hom_{o_{\tilde{X}}}(\Oscr_L,-)$ then it is clear there will be a short
  exact sequence as in \eqref{ref:8.8a} with $E$ replaced by $F$. By
  considering resolutions with such acyclic objects one finds the
  triangle \eqref{ref:8.7a}.
\item This follows from 1. together with the construction of  
$\Hom_{o_{\tilde{X}}}(\Oscr_L,-)$ by acyclic resolutions. 
 \item This follows from substituting $\Mscr=\Oscr_L$ in the triangle
 \eqref{ref:8.7a} and using Proposition \ref{ref:8.3.1a}.
  \qed \end{enumerate}
\def\qed{}

\end{proof}

The following result is  well-known.
\begin{lemmas}
\label{ref:8.3.4a}
 Let $\Ascr$ be an abelian category.
Let $X\in D^*(\Ascr)$ with $*=+,-$ be an object such that 
$\Ext^i(H^j(X),H^{j+1-i}(X))=0$ for $i\ge 2$
and for all $j$. Then $X$ is isomorphic to the direct sum of its
homology groups. The same is true for $X\in D(\Ascr)$ if in addition
$\Ascr$ satisfies AB4,  has enough injectives  and $\Hom(H^j(X),-)$ has
finite cohomological dimension for all $j$.
\end{lemmas}
\begin{proof} Assume first that $X\in D^+(\Ascr)$. Write $H(X)=\oplus
H^i(X)[-i]$. Note that $H(X)$ is the category theoretic  direct sum of the
$H^i(X)[-i]$ in $D^+(\Ascr)$. We want to construct a quasi-isomorphism
$H(X)\r X$. To this end it is sufficient to construct maps $H^i(X)[-i]\r
X$ which induce isomorphisms on the $i$'t cohomology. Since $\tau_{\le i}
X\r X$ induces an isomorphism on $H^i$, it is clearly sufficient to show
that the canonical map $\tau_{\le i}X \r H^i(X)[-i]$ splits. From the
triangle 
\[
\tau_{\le i-1} X\r \tau_{\le i} X \r H^i(X)[-i]\r
\]
we find that we have to show that
\begin{equation}
\label{ref:8.9a}
\Hom(H^i(X)[-i],\tau_{\le i-1}
X[1])=0
\end{equation}
Now $\tau_{\le
 i-1} X$ is a bounded complex and hence \eqref{ref:8.9a} follows easily
 by induction from the hypotheses.
 
 The case $X\in D^-(\Ascr)$ is similar. Now assume $X\in D(\Ascr)$. In
 this case there
 are  two possible problems with the above reasoning.
 \begin{enumerate}
 \item $\oplus_i H^i(X)[-i]$ is perhaps no longer the direct sum of the 
 $H^i(X)[-i]$ in $D(\Ascr)$.
 \item $\tau_{\le i-1} X$  is now an unbounded complex so we can no
 longer verify \eqref{ref:8.9a} by induction.
 \end{enumerate}
 The first difficulty is resolved if we assume that $\Ascr$ satisfies
 AB4 \cite{Neeman}. For the second difficulty we have to show that
 $\Hom(H^i(X),-)$ is zero on $\Dscr(\Ascr)_{\le -N}$ for $N\gg 0$. Now
 according to \cite[Thm 5.1, Cor. 5.3]{RD}, if $\Ascr$ has enough
 injectives and $\Hom(H^i(X),-)$ has finite cohomological dimension
 then we can compute $\Hom(H^i(X),-)=H^0(\RHom(H^i(X),-))$ by acyclic
 resolutions. It follows easily that if $\cd \Hom(H^i(X),-)=t$ then an
 object in $\Dscr(\Ascr)_{\le -N}$  can be represented by an acyclic
 complex which is non-zero only in degree $\le -N+t$. Hence it
 suffices to take $N>t$.
 \end{proof}
 \begin{remarks} The reader may verify that the statement about
 $D(\Ascr)$ in the previous lemma is also true under the hypotheses AB4
 and AB4${}^\ast$. This is more elegant, but less useful in practice.
 \end{remarks}
 \begin{corollarys}
\label{ref:8.3.6a}
 Let $\Nscr$ be the additive subcategory of
 $\Qch(\tilde{X})$ whose objects are direct sums of copies of
 $\Oscr_L(-1)$. Then $\Nscr$ is a thick subcategory (closed under
 extensions).  Furthermore  the map
 \[
 D(k)\xrightarrow{-\otimes \Oscr_L(-1)} D_\Nscr(\Qch(X))
 \]
 is an equivalence.
 \end{corollarys}
\begin{proof}
This follows from lemma \ref{ref:8.3.4a} together with lemma \ref{ref:8.3.3a}.
\end{proof}

\subsection{The main theorem}
\label{ref:8.4b}

We prove the following Theorem \ref{ref:8.4.1a}. This is a
non-commutative version of a theorem by Orlov \cite{Orlov}. Our proof
is slightly different.
\begin{theorems}
\label{ref:8.4.1a}
  There is a semi-orthogonal decomposition of $D(\tilde{X})$ given by
  $(D(X),D(k))$. The diagram corresponding to
  \eqref{ref:8.1b} is as follows
\begin{equation}
\label{ref:8.10a}
  D(k)
\bfig
\putmorphism(0,100)(1,0)[``-\otimes_k \Oscr_L(-1)]{1100}{1}{a}
\putmorphism(0,-50)(1,0)[``F]{1100}{-1}{b}
\efig
D(\tilde{X})
\bfig
\putmorphism(0,100)(1,0)[``R\alpha_\ast]{300}{1}{a}
\putmorphism(0,-50)(1,0)[``L\alpha^\ast]{300}{-1}{b}
\efig
 D(X)
\end{equation}
where $F$ is the left adjoint to $-\otimes_k\Oscr_L(-1)$ (whose existence 
follows by lemma \ref{ref:8.1.1a}).

In \eqref{ref:8.10a} we may replace $D$ by $D^\ast$ where
$\ast=\emptyset,+,-,b$. Furthermore we may also replace $D^\ast$ everywhere by
$D_f^\ast$.
On $D_f^-(\tilde{X})$, $F$ is given by $\RHom(-,\Oscr_L(-1))^\ast$.
 \end{theorems}
\begin{proof}
\begin{step} 
\label{ref:1b}
$R\alpha_\ast L\alpha^\ast=\Id$. Since the objects in $\Vscr$ are
acyclic for $\alpha^\ast$, this follows from Proposition \ref{ref:8.3.1a}.2
on $D^-(X)$.
 
The general case is then routine. If $\Fscr\in D(X)$ then we have to
show that the adjunction mapping $\Fscr\r R\alpha_\ast
L\alpha^\ast\Fscr$ is a quasi-isomorphism, that is $H^i(\Fscr)\r
H^i(R\alpha_\ast L\alpha^\ast\Fscr)$ should be an isomorphism for all $i$. Since
$\alpha^\ast$, $\alpha_\ast$ have finite cohomological dimension, we
can test this (for a fixed $i$) by replacing $\Fscr$ by some $\tau_{\le
  N}\Fscr$ for $N\gg 0$. But then $\Fscr\in D^-(X)$ and this case was
already covered.
\end{step}
\begin{step}
\label{ref:2a} The composition of $=-\otimes_k \Oscr_L(-1)$ and 
$R\alpha_\ast$ is zero. This follows from the fact that
$R\alpha_\ast$ commutes with direct sums (lemma \ref{ref:8.2.2a})
together with $R\alpha_\ast\Oscr_L(-1)=0$, which was proved in
Proposition \ref{ref:8.3.1a}.
\end{step}

\begin{step} 
\label{ref:3b}
At this point $D(k)$ and $D(X)$ form a semi-orthogonal pair in
$D(\tilde{X})$. To show that they form a semi-orthogonal
decomposition we invoke lemma \ref{ref:8.1.4a}. Thus we have to prove
that $\ker R\alpha_\ast= D(k)$.  We prove first $\ker
R\alpha_\ast=D_\Nscr(\tilde{X})$ ($\Nscr$ as in Cor. \ref{ref:8.3.6a}).
Then by Corollary \ref{ref:8.3.6a} we find $D_\Nscr(\tilde{X})=D(k)$.
  
We claim first that
\[
\ker R\alpha_\ast\subset D_{\alpha^{-1}(\Cscr_p)}(\tilde{X})
\]
To prove this assume that $R\alpha_\ast\Mscr=0$. From the spectral
  sequence 
\[
R^p\alpha_\ast H^q(\Mscr)\r  R^{p+q}\alpha_\ast \Mscr
\]
we obtain short exact sequences 
\[
0\r R^1\alpha_\ast H^{i-1}(\Mscr)\r R^i\alpha_\ast\Mscr \r
R^0\alpha_\ast H^i(\Mscr)\r 0
\]
We conclude that for all $i$ one has $\alpha_\ast H^i(\Mscr)=0$, whence by
Proposition \ref{ref:6.5.2a} $\Mscr\in D_{\alpha^{-1}(\Cscr_p)}(\tilde{X})$. 
\end{step}
\begin{step}
\label{ref:4b}
 If $\phi:\Mscr\r\Nscr$ is a map in $D(\tilde{X})$
then we denote by $\cone(\phi)$ the cone of the triangle with base
$\phi$. $\cone(\phi)$ is unique up to non-unique isomorphism. 

We claim
\begin{equation}
\label{ref:8.11a}
\ker R\alpha_\ast=\{ \cone( L\alpha^\ast R\alpha_\ast\Mscr\r \Mscr
)\mid \Mscr\in D_{\alpha^{-1}(\Cscr_p)}(\tilde{X}) \}
\end{equation}
Since $R\alpha_\ast L\alpha^\ast=\Id$ it is clear that
$RHS\eqref{ref:8.11a}\subset LHS\eqref{ref:8.11a}$.  The opposite inclusion
is obvious since if $R\alpha_\ast\Mscr=0$ then according to Step
\ref{ref:3b}~: $\Mscr\in D_{\alpha^{-1}(\Cscr_p)}(\tilde{X})$. Hence 
$\Mscr=\cone( L\alpha^\ast R\alpha_\ast\Mscr\r\Mscr )\in
\RHS\eqref{ref:8.11a}$.
\end{step}

\begin{step}
\label{ref:5b}
  The formation of $\cone(L\alpha^\ast R\alpha_\ast \Mscr\r\Mscr)$ is
  functorial and compatible with shifts and  triangles.

According to lemma
\ref{ref:8.1.4a}, we have a semi-orthogonal decomposition
of $D(\tilde{X})$ given by $(\ker R\alpha_\ast,D(X))$.
 Since
in the triangle
\[
\Atriangle<1`-1`1; >[\Cscr`L\alpha^\ast R\alpha_\ast\Mscr`\Mscr;``]
\]
$L\alpha^\ast R\alpha_\ast\Mscr$ lies in the essential image  of
$D(X)$ and $\Cscr\in \ker R\alpha_\ast$ the good behaviour of
$\cone(L\alpha^\ast R\alpha_\ast\Mscr\r\Mscr)$ is a consequence of lemma
\ref{ref:8.1.1a}. 
\end{step}

\begin{step}
According to Step \ref{ref:4b} we have to show that if $\Mscr\in
D_{\alpha^{-1}(\Cscr_p)}(\tilde{X})$ then 
$\cone (L\alpha^\ast R\alpha_\ast\Mscr\r\Mscr)\in D_\Nscr(\tilde{X})$. That is
\begin{equation}
\label{ref:8.12a}
H^i(\cone (L\alpha^\ast R\alpha_\ast\Mscr\r\Mscr))\in\Nscr
\end{equation}
Now since $L\alpha^\ast$, $R\alpha_\ast$ have finite cohomological
dimension the truth of \eqref{ref:8.12a} will not be influenced by
the cohomology of $\Mscr$ in high degree. Hence to verify
\eqref{ref:8.12a} we may replace $\Mscr$ by a suitable $\tau_{\le
a}\tau_{\ge b}\Mscr$, that is, by a bounded complex. According to
Step \ref{ref:5b} and induction we may then assume that $\Mscr\in
\alpha^{-1}(\Cscr_p)$. 

Now   \eqref{ref:8.12a} is a statement about homology of complexes, and it is
easy to see that this homology is compatible with direct limits, when
$\Mscr\in\Qch(\tilde{X})$. Hence we may  assume that $\Mscr\in
\alpha^{-1}(\Cscr_p)\cap \coh(\tilde{X})$.

Now from lemma \ref{ref:8.4.2a} below it will follow that we
may in fact assume $\Mscr=\Oscr_L$, $\Oscr_L(-1)$ or $\Oscr_q$ with
$q\in \beta^{-1}(O_\tau(p))\setminus \beta^{-1}(\tau p)$. The last two
cases are trivial since by Propositions \ref{ref:8.3.1a} and
\ref{ref:8.3.2a} we find that $L\alpha^\ast R\alpha_\ast \Mscr=\Mscr$.

Hence assume $\Mscr=\Oscr_L$. Thus we have to compute
\begin{equation}
\label{ref:8.13a}
\cone(L\alpha^\ast R\alpha_\ast
\Oscr_L\r \Oscr_L)=\cone(L\alpha^\ast \Oscr_{\tau p}\r \Oscr_L)
\end{equation}
Since by Proposition \ref{ref:8.3.2a} and $\alpha^\ast\Oscr_{\tau
  p}=\Oscr_L$, and $\alpha^\ast$ has cohomological dimension one we
find that, up to shift, \eqref{ref:8.13a} is equal to $L_1\alpha^\ast
\Oscr_{\tau p}$, which according to Proposition \ref{ref:8.3.2a} is
equal to $\Oscr_L(-1)$.
\end{step}
\begin{step} Odd and ends. First of all since $R\alpha_\ast$ and
  $L\alpha^\ast$ have finite cohomological dimension and preserve
  coherent objects (by Theorem \ref{ref:6.3.1a}),
  it is clear that the above reasoning may be repeated for $D^\ast$
  and for $D_f^\ast$ with $\ast=\phi,+,-,b$.

  To show that $F$ on $D^-_f(\tilde{X})$ is equal to
  $\RHom(-,\Oscr_L(-1))^\ast$ we have to show that this functor is
  well-defined and is a left adjoint to $-\otimes \Oscr_L(-1)$. A
  quick verification shows that we have to check that
  $\coh(\tilde{X})$ contains enough objects $\Fscr$ such that the
  homology of $\RHom_{o_X}(\Fscr,\Oscr_L(-1))$ is finite dimensional
  and concentrated in degree zero.

For $\Fscr$ we take objects  of  the form $\alpha^\ast(\Escr)(-n-1)$ with
$n\in\NN$ and
$\Escr\in \Vscr\cap \coh(X)$. It is easy to see that there are enough of
those. We compute
\begin{align*}
\RHom_{o_{\tilde{X}}}(\alpha^\ast(\Escr)(-n-1),\Oscr_L(-1))
&=\RHom_{o_X}(\Escr,R\alpha_\ast(\Oscr_L(n)))\\
&=\RHom_{o_X}(\Escr,(\Dscr/m_{\tau p}\Dscr)_{n,o_X})
\end{align*}
It now follows for example from Proposition
\ref{ref:5.1.2a} that the homology of the last line of the
previous equation has the properties we want.  \qed\end{step}
\def\qed{}\end{proof}
\begin{lemmas}
\label{ref:8.4.2a}
$\alpha^{-1}(\Cscr_p)\cap \coh(\tilde{X})$ is the smallest subcategory
of $\Qch(\tilde{X})$ containing $\Oscr_L(-1)$, $\Oscr_L$ and $\Oscr_q$
with $q\in \beta^{-1}(O_\tau(p))\setminus \beta^{-1}(\tau p)$, and
which is closed under
\begin{enumerate}
\item
extensions;
\item
kernels of surjective maps;
\item
cokernels of injective maps.
\end{enumerate}
\end{lemmas}
\begin{proof}
By Corollary \ref{ref:6.7.4a} every object in $\alpha^{-1}(\Cscr_p)\cap
\coh(\tilde{X})$ is an extension of objects in $\Qch(L)\cap
\coh(\tilde{X})=\coh(L)$. Hence we have to verify the statement for
$\coh(L)$. But then we may as well verify the corresponding statement
for $S$. This  is a routine exercise, which we leave to the reader.
\end{proof}
\begin{corollarys}
\label{ref:8.4.3a}
If $\Qch(X)$ has finite injective dimension then so does
$\Qch(\tilde{X})$. 
\end{corollarys}
\begin{proof}
  Assume that $\Qch(X)$ has finite injective dimension.  We have to
  show that there exists an $N$ such that
  $\Hom_{D(\tilde{X})}(\Mscr,\Nscr)=0$ for all $\Mscr\in D^b_{\ge
    N}(\tilde{X})$ and $\Nscr\in D^b_{\le 0}(\tilde{X})$.

For an arbitrary morphism $\phi:\Mscr\r \Nscr$ we have a commutative
diagram of triangles
\[
\begin{CD}
  L\alpha^\ast R\alpha_\ast \Mscr @>>> \Mscr @>>> F\Mscr \otimes_k
  \Oscr_L(-1)@>>>\\ 
@V L\alpha^\ast R\alpha_\ast(\phi) VV @V\phi VV @V
  F\phi\otimes_k
  \Oscr_L(-1) VV\\
 L\alpha^\ast R\alpha_\ast \Nscr @>>> \Nscr @>>> F\Nscr\otimes_k
  \Oscr_L(-1)
  @>>>
\end{CD}
\]
Since $L\alpha^\ast$ is fully faithful and $\Qch(X)$ has finite
injective dimension we find that $L\alpha^\ast R\alpha_\ast(\phi)=0$ if
$N$ is large enough. Similarly $D^b(k)$ is a semi-simple category and
hence $F\phi=0$ if $N$ is large enough. This then implies that $\phi$
is in fact obtained from a map $\theta:F\Mscr\otimes_k
  \Oscr_L(-1)\r
L\alpha^\ast R\alpha_\ast \Nscr $. However by lemma \ref{ref:8.3.3a} we
find that $\Hom_{D(\tilde{X})}(\Oscr_L,-)$ has finite cohomological
dimension. Hence by making $N$ even larger if necessary we also find
that $\theta=0$. This finishes the proof.
\end{proof}

\section{Some results on graded algebras and their sections}
\label{ref:9a}
\subsection{Generalities}
\label{ref:9.1a}

In this section $\gamma:X\r \Spec k$ will be a noetherian quasi-scheme
over $k$ which is proper (see Def.\ \ref{ref:3.9.4a}). As
usual we put $\Oscr_X=\gamma^\ast k$.  We will also assume that
$\Oscr_X$ is noetherian.  By adjointness we have $\gamma_\ast
\Mscr=\Hom_{\Qch(X)}(\Oscr_X,\Mscr)=\Gamma(X,\Mscr)$. Similarly for
the derived functors. Hence the properness of $X$ simply means that if
$\Mscr\in \coh(X)$ then $H^i(X,\Mscr)$ is a finite dimensional
$k$-vector space for all $i$. The fact that $\Oscr_X$ is noetherian
implies that $H^i(X,-)$ commutes with direct limits.

Assume that $\Escr$ is a noetherian $\NN$-graded algebra on $X$ and
let $E=\Gamma(X,\Escr)$ (cfr \S\ref{ref:3.6b}). Then
$(\Escr)_{o_X}$ is an $E-\Escr$-bimodule, in the obvious sense. It is
clear that we have adjoint functors
\begin{equation}
\label{ref:9.1b}
\begin{split}
-\otimes_E \Escr_{o_X}:&\Gr(E)\r \Gr(\Escr) \\
\Gamma(X,-):&\Gr(\Escr)\r \Gr (E)
\end{split}
\end{equation}
Our aim will be to  relate $\QGr(\Escr)$ to $\QGr(E)$. To this end we
introduce the following definition.

\begin{definitions}
$\Escr$ is ample if for every  object $\Mscr$ in $\coh(X)$ we
have for $n\gg 0$  that $\Mscr\otimes_{o_X} \Escr_n$ is generated by
global sections and 
$H^i(X,\Mscr\otimes_{o_X} \Escr_n)=0$  for  $i>0$.
\end{definitions}
The following proposition was essentially proved in \cite{AVdB,VdB11} 
(see \cite[Thm 5.2]{VdB11}) under
more restrictive hypotheses. 
\begin{propositions} 
\label{ref:9.1.2a}
Let $\Escr$ be as above and assume that $\Gamma(X,-)$ has finite
cohomological dimension.
If $\Escr$ is ample then $E$ is noetherian. The functors in
\eqref{ref:9.1b} send $\gr(E)$ to $\gr(\Escr)$ and vice-versa.
Furthermore these functors define inverse equivalences between $\QGr(E)$
and $\QGr(\Escr)$.
\end{propositions}

Unfortunately ampleness of $\Escr$ is too strong for the
applications we have in mind. However we will show that under some extra
conditions the functors \eqref{ref:9.1b} may still be well behaved,
even if $\Escr$ is non-ample.

We introduce the following hypotheses.
{\def\thehypothesis{(***)}
\begin{hypothesis}
  $\Gamma(X,-)$ has finite cohomolical dimension and furthermore there
  is an injective graded $\Escr$-bimodule map $t:\Escr(-1)\r \Escr$
  such that $\Escr/\im t$ is ample. Furthermore the induced map
  $\Escr_{o_X}(-1)\r \Escr_{o_X}$ is injective.
\end{hypothesis}
}
Taking global sections yields a regular central element of $E$ in
degree one which we also denote by $t$.  We write
$\bar{\Escr}=\Escr/\im t$. If $\Nscr\in\Gr (\Escr)$ then $t$ defines a
map $\Nscr(-1)\r \Nscr$ in $\Gr(\Escr)$. We put $\bar{\Nscr}=\Nscr/\im
t$.
We say that $\Nscr$ is annihilated by $t$ if
$t:\Nscr(-1)\r \Nscr$ is the zero map. We say that $\Nscr$ is $t$-torsion
if $\Nscr$ is a union of objects which are each annihilated by some
$t^n$ ($n$ variable). We say that $\Nscr$ is $t$-torsion free if $t$ is
injective. Similar conventions apply to $\QGr(\Escr)$, $\Gr(E)$ and
$\Qgr(E)$.

 \begin{propositions} 
\label{ref:9.1.3a}
Let
 $\Nscr\in \gr(\Escr)$.  
 Then 
 \begin{enumerate}
\item Multiplication by $t$ on $H^i(X,\Nscr)$ is an automorphism  in high 
degree when $i>0$.
\item
The complex
\[
\Gamma(X,\Nscr(-1))\xrightarrow{t} \Gamma(X,\Nscr)\r
\Gamma(X,\bar{\Nscr})\r 0
\]
is exact in high degree.
\item $H^i(X,\Nscr)$ is finitely generated for all $i$.
\end{enumerate}
Furthermore we also have
\begin{enumerate}
\setcounter{enumi}{3}
\item $\Tor^E_i(-,\Escr)$ sends $\Tors(E)$ to $\Tors(\Escr)$ for all $i$.
 \item $E$ is noetherian.
\end{enumerate}
\end{propositions}
\begin{proof}
\begin{itemize}
\item[1.,2.]
Using the appropriate long exact sequences it follows that it is
sufficient to prove this  in the case that $\Nscr$ is
annihilated by $t$ and in the case that $\Nscr$ is $t$-torsion free.

In the first case we have that $\Nscr=\bar{\Nscr}$
and furthermore, for $i>0$, $H^i(X,\bar{\Nscr})$ is right bounded (by
ampleness of $\bar{\Escr}$). Hence 1.,2. are
 trivially true in this case.

Let us consider the second case. We apply the long exact sequence for
$\Gamma(X,-)$ to
\[
0\r \Nscr(-1)\r \Nscr \r \bar{\Nscr}\r 0
\]
Since $H^i(X,\bar{\Nscr})$ is right bounded, we immediately obtain the
statement about $H^i(X,\Nscr)$ for $i>1$. So let us assume $i\le 1$. For
$n\gg 0$ we have an exact sequence
\[
0\r \Gamma(X,\Nscr_{n-1})\xrightarrow{t} \Gamma(X,\Nscr_n)\r
\Gamma(\bar{\Nscr}_n)\r  H^1(X,\Nscr_{n-1})\xrightarrow{t}
H^1(X,\Nscr_n)\r 0
\]
So $\dim H^1(X,\Nscr_n)$ is descending and hence must eventually become
constant. Thus $t:H^1(X,\Nscr_{n-1})\r H^1(X,\Nscr_n)$ is an isomorphism
for $n\gg 0$. This proves 1.,2. for $i=0,1$.
\item[3.] In view of 1., the only non-trivial case is $i=0$. By
  Proposition \ref{ref:9.1.2a} $\Gamma(X,\bar{\Nscr})$ is finitely
  generated. Since according to 2., $\overline{\Gamma(X,\Nscr)}$ has
  finite colength in $\Gamma(X,\bar{\Nscr)}$ we find that
  $\overline{\Gamma(X,\Nscr)}$ is finitely generated. Hence by the
  graded version of Nakayama's lemma, $\Gamma(X,\Nscr)$ is also
  finitely generated.
\item[4.] Since $\Tor$ is compatible with direct limits 
we may prove the corresponding statement for ``$\tors$''. Furthermore,
using the appropriate long exact sequences, it follows that it is
sufficient to prove that $\Tor_i^E(V,\Escr)$ is right bounded for a
finite dimensional $E_0$-module $V$ (considered as $E$-module).

Since $V$ is annihilated by $t$, we have as usual
\begin{equation}
\label{ref:9.2a}
\Tor_i^E(V,\Escr)=\Tor_i^{\bar{E}}(V,\bar{\Escr})
\end{equation}
Put $\bar{E}'=\Gamma(X,\bar{\Escr})$. According to 2., applied to
$\Escr_{o_X}$ we find that the map $\bar{E}'\r \bar{E}$ is injective
and has finite dimensional cokernel.  Hence by \cite[Prop. 2.5]{AZ} and
Proposition \ref{ref:9.1.2a} one has
$\QGr(\bar{E}')=\QGr(\bar{E})=\QGr(\bar{\Escr})$.  The functor
realizing this equivalence is given by $-\otimes_{\bar{E}'}
\bar{\Escr}$. Hence this is in particular an exact functor. Now to
compute the righthand side of \eqref{ref:9.2a} we take a free resolution of
$V$. Since the grading on $V$ is right bounded this free resolution is
an exact sequence in $\QGr(\bar{E}')$. Hence it remains exact after
applying $-\otimes_{\bar{E}'} \bar{\Escr}$. In this way we obtain that
$\Tor_i^E(V,\Escr)\in \Tors(\Escr)$ for all $i$.
\item[5.]  By Proposition \ref{ref:9.1.2a} it follows that
  $\bar{E}'$ is noetherian. Since the map $\bar{E}\r \bar{E}'$ is
  injective and has finite dimensional cokernel, we deduce that
  $\bar{E}$ is also noetherian. Whence, by a Hilbert-type argument,
  $E$ is noetherian.\qed
\end{itemize}
\def\qed{}\end{proof}
For use below we write $W=\Proj \Escr$, $T=\Proj \bar{\Escr}$,
$V=\Proj E$, $S=\Proj \bar{E}$. It follows that $S,T$ are S
divisors in $W,V$.

As usual, we will denote the
quotient functors $\Gr(E)\r \QGr(E)$, $\Gr(\Escr)\r \QGr(\Escr)$ by $\pi$.
\begin{propositions} The functors
\begin{align*}
\delta^\ast:&\Qch(V)\r \Qch(W): \pi M\mapsto \pi(M\otimes_E \Escr)\\
\delta_\ast:&\Qch(W)\r \Qch(V):\pi N\mapsto
\pi\Gamma(X,N) \end{align*}
are well defined and  form an  adjoint pair.

  In
particular $(\delta^\ast,\delta_\ast)$ defines an morphism of quasi-schemes
\[
\delta:W\r V
\]
$\delta$ fits in a commutative diagram of quasi-schemes
\begin{equation}
\label{ref:9.3a}
\begin{CD}
T @>j>> W\\
@V\delta VV @V\delta VV\\
S @>i>> V
\end{CD}
\end{equation}
Here $i,j$ are the inclusion mappings. $\delta:T\r S$ is an
isomorphism.  If $\Mscr\in \Qch(T)$ then we have
\begin{equation}
\label{ref:9.4a}
\delta_\ast(\Mscr\otimes_{o_T} \Nscr_{T/W})=\delta_\ast(\Mscr)\otimes_{o_S}
\Nscr_{S/V}
\end{equation}
Thus the normal bundles on $T$ and $S$ correspond to each other under
the isomorphism $\delta$.
\end{propositions}
\begin{proof}
It is obvious that $\delta_\ast$ is well-defined. The fact that
$\delta^\ast$ is well-defined follows easily from Proposition
\ref{ref:9.1.3a}.4.

We now show that $\delta_\ast$, $\delta^\ast$ are adjoint functors. This
means that we have to construct a natural isomorphism between
\[
\Hom_{\QGr(E)}(\Mscr,\delta_\ast\Nscr)
\]
and
\[
\Hom_{\QGr(\Escr)}(\delta^\ast\Mscr,\Nscr)
\]
where $\Mscr=\pi M$, $\Nscr=\pi N$.  

By taking a presentation of $M$ we reduce to $M=E$. We have 
\begin{align*}
\Hom_{\QGr(E)}(E,\delta_\ast \Nscr)&=\dirlim_n \Hom_{\Gr(E)}
(E_{\ge n}, \Gamma(X,N))\\
&=\dirlim_n \Hom_{\Gr(\Escr)}(E_{\ge n}\otimes_E \Escr_{o_X},N)
\end{align*}
and we have to show that this is equal to 
\[
\dirlim_{n}\Hom_{\Gr(\Escr)}((\Escr_{o_X})_{\ge n},N)
\]
Put 
\begin{align*}
K_n&=\ker (E_{\ge n}\otimes_E \Escr_{o_X}\r (\Escr_{o_X})_{\ge n})\\
C_n&=\coker (E_{\ge n}\otimes_E \Escr_{o_X}\r (\Escr_{o_X})_{\ge n})
\end{align*}
It is now sufficient to show that $K_n,C_n$ are \emph{torsion} inverse
systems. That is for all $n$ there exists an $m$ such that the maps
$K_{n+m}\r K_n$, $C_{n+m}\r C_n$ are zero.

Tensoring the exact sequence 
\[
0\r E_{\ge n}\r E\r E/E_{\ge n}
\r 0
\]
with $\Escr_{o_X}$ and restricting to degrees $\ge n$ yields that 
\begin{align*}
K_n&=\Tor_1^E(E/E_{\ge n}, \Escr_{o_X})\\
C_n&=(E/E_{\ge n}\otimes_E \Escr_{o_X})_{\ge n}
\end{align*}
In particular it follows from Proposition \ref{ref:9.1.3a}.4 that $K_n$, $C_n$
are right bounded. Furthermore it is clear that $C_n$ is zero in degrees
$< n$, and by considering a minimal free resolution of $E/E_{\ge n}$ we
see that the same holds for $K_n$. It now follows clearly that $K_n$,
$C_n$ are torsion inverse systems.

That \eqref{ref:9.3a} is a commutative diagram is obvious and the fact
that $\delta:T\r S$ is an isomorphism is simply the ampleness of $\bar{\Escr}$.

The identity \eqref{ref:9.4a} is just a special case of the fact that
$\delta_\ast$ commutes with shift.
\end{proof}
From Proposition \ref{ref:9.1.3a} one easily deduces that the
shifts of $\pi E$ are acyclic for $\delta^\ast$ (in the sense of
\cite{RD}). There are clearly enough of those, so $L\delta^\ast$
exists.

It is easy to verify the following formulas
\begin{equation}
\label{ref:9.5a}
\begin{split}
R^i\delta_\ast \pi N&= \pi H^i(X,N)\\
L_i\delta^\ast \pi M&= \pi \Tor^E_i(M,\Escr)
\end{split}
\end{equation}
\begin{lemmas}
$R\delta_\ast$ sends $D^+_f(W)$ to $D^+_f(V)$ and
$L\delta^\ast$ sends $D^-_f(V)$ to $D^-_f(W)$.
\end{lemmas} 
\begin{proof} The statement about $L\delta^\ast$ follows from the fact
that $E$ and $\Escr$ are noetherian. The statement about $R\delta_\ast$
follows from Proposition \ref{ref:9.1.3a}.3.
\end{proof}
Thus in particular $\delta$ is \emph{proper} in the sense of
Definition \ref{ref:3.9.4a}.

The following is also easy.

\begin{propositions}
\label{ref:9.1.6a}
If $\Nscr\in \Tors_T(W)$ then $R^i\delta_\ast\Nscr=0$ for $i>0$. 
Similary if $\Mscr\in \Tors_S(V)$ then $L^i\delta^\ast \Mscr=0$ for
$i>0$. Finally $\delta_\ast$, $\delta^\ast$ define inverse equivalences 
between
$\Tors_T(W)$ and $\Tors_S(V)$ (see \S\ref{ref:3.7b} for notation).
\end{propositions}
\begin{proof} Since all functors involved commute with direct limits we
may assume that $\Mscr$, $\Nscr$ are coherent. The first two statements
are easily seen to be true for objects annihilated by $t$. The general
case follows from this by considering the appropriate long exact
sequences.

In particular we obtain that $\delta_\ast$ and $\delta^\ast$ are
exact.  We now have to show that the adjunction morphisms are
isomorphisms on coherent objects.  Again this is true for objects
annihilated by $t$ by the hypotheses on $\bar{\Escr}$. The general
case now follows by filtering objects in such a way that the
associated graded quotients are annihilated by $t$.
\end{proof}

The following proposition will be true in all our applications for
trivial reasons. However it may be interesting to note that one can
prove it in this generality.
\begin{propositions} 
$\Qch(W)$ has enough $T$-torsion free objects. In particular $Lj^\ast$
is defined.
\end{propositions}
\begin{proof}
  It suffices to show that every object in $\mod(W)$ is a quotient of
  a $T$-torsion free one. Let $\Mscr=\pi M$ where $M\in\gr \Escr$.
  With an argument as in the proof of lemma \ref{ref:3.7.3a} one
  finds a $M_0\subset M$ which is $t$-torsion free such that $M/M_0$ is
  $t$-torsion. Assume $t^n(M/M_0)=0$ and consider the following
  diagram with exact rows.
\begin{equation}
\label{ref:9.6a}
\begin{CD}
0 @>>> M_0 @>>> M @>>> M/M_0 @>>> 0\\
@. @VVV @VVV @|  \\
0 @>>> M_0/t^n M_0 @>>> M/t^n M @>>> M/M_0 @>>>0
\end{CD}
\end{equation}
According to Propositin \ref{ref:9.1.3a},
\begin{equation}
\label{ref:9.7a}
\Gamma(X,M)\r
\Gamma(X,M/t^n M)
\end{equation}
 is surjective in high degree. Hence there exists a
surjective map $E(-a)^b\r \Gamma(X,M)$ for certain $a,b$ such that the
composition with \eqref{ref:9.7a} has finite cokernel.

Tensoring with $\Escr_{o_X}$ and composing with $\Gamma(X,M)\otimes_E
\Escr_{o_X}\r M$ we obtain a map $\Escr_{o_X}(-a)^b\r M$ such that the
composition with $M\r M/t^n M$ has right bounded cokernel (using
Propositions \ref{ref:9.1.3a}.4 and
\ref{ref:9.1.6a}). Applying $\pi$ to the induced map $M_0\oplus
\Escr_{o_X}(-a)^b\r M$ together with \eqref{ref:9.6a} yields what we
want. 
\end{proof}

\begin{propositions} 
\label{ref:9.1.8a}
All maps and compositions  in \eqref{ref:9.3a} are
  admissible. Furthermore one has the following identities.
\begin{align}
i^\ast\delta_\ast&=\delta_\ast j^\ast\\
\label{ref:9.9a}
i^!\delta_\ast&=\delta_\ast j^!\\
\label{ref:9.10a}
Ri^!\,R\delta_\ast&= R\delta_\ast Rj^!\\
\label{ref:9.11a}
Li^\ast\, R\delta_\ast &= R\delta_\ast Lj^\ast
\end{align}
\end{propositions}
\begin{proof}
Checking admissibility is routine using the explicit formulas
\eqref{ref:9.5a}. We leave this to the reader.

To verify the compatibility $i^\ast\delta_\ast=\delta_\ast j^\ast$ we
may work in $\coh(X)$ (since everything is compatible with direct
limits). But then it is just a reformulation of Proposition
\ref{ref:9.1.3a}.2.

The other compatibility $i^!\delta_\ast=\delta_\ast j^!$ involves only
left exact functors, so it may be verified on injectives. Let $F\in
\Qch(W)$ be such an injective. We have an exact sequence
\[
0\r j^! F\r F \r F(T) \r 0
\]
Now $j^! F\in \Tors_T(W)$ hence $R^1\delta_\ast j^!
F=0$. 
Therefore we obtain an exact sequence
\begin{equation}
\label{ref:9.12a}
0\r \delta_\ast j^! F\r \delta_\ast F \r \delta_\ast (F(T)) \r 0
\end{equation}
Now $\delta_\ast$ is compatible with shift. That is $\delta_\ast
(F(T))=(\delta_\ast F)(S)$.
Since we also have an exact sequence
\begin{equation}
\label{ref:9.13a}
0\r i^!\delta_\ast F\r \delta_\ast F\r (\delta_\ast F)(S)
\end{equation}
we are through.

\eqref{ref:9.10a} follows from \eqref{ref:9.9a} provided we show for every
injective $E\in \Qch(X)$ that $j^!E$ is acyclic for $\delta_\ast$
and that $\delta_\ast E$ is acyclic for $i^!$. 

The statement about $j^!E$ is clear since $j^!E$ is supported on
$T$. The statement of $\delta_\ast E$ follows from the right exactness
of \eqref{ref:9.13a}. However \eqref{ref:9.13a} is the same as \eqref{ref:9.12a} and the
latter is right exact. This shows what we want.

\eqref{ref:9.11a} follows from \eqref{ref:9.10a} using \eqref{ref:7.2a}.
\end{proof}

In the next few paragraphs we investigate to what extent the map
$\delta:W\r V$ is an isomorphism.  

\begin{theorems}
\label{ref:9.1.9a}
\begin{enumerate}
\item
$\Iso_S(V)$ and $\Iso_T(W)$  are mapped  to each other under
$(\delta^\ast,\delta_\ast)$.
\item
$R^i\delta_\ast$ sends $\Qch(W)$ to $\Iso_S(V)$ for $i>0$.
\item
$L_i\delta^\ast$ sends $\Qch(V)$ to $\Iso_T(W)$ for $i>0$.
\item
$(\delta^\ast,\delta_\ast)$ define inverse equivalences between
$\Qch(W)/\Iso_T(W)$ and $\Qch(V)/\Iso_S(V)$.
\end{enumerate}
\end{theorems}
\begin{proof}
\begin{enumerate}
\item
This is clear by functoriality.
\item
This is precisely 1. of Proposition \ref{ref:9.1.3a}.
\item
We may assume that $\Mscr$ is coherent.
Furthermore by considering the appropriate long exact sequence it
suffices to consider the cases where $\Mscr$ is $t$-torsion and where
$\Mscr$ is $t$-torsion free.

If $\Mscr$ is $t$-torsion then $L_i\delta^\ast\Mscr=0$ for $i>0$ by
Proposition \ref{ref:9.1.6a} so the result is true. Hence assume that $\Mscr$ is
$t$-torsion free. Then the long exact sequence for $\delta^\ast$
applied to
\[
0\r \Mscr(-1)\r \Mscr\r \bar{\Mscr}\r 0
\]
yields that $t$ is an isomorphism between $L_i\delta^\ast\Mscr(-1)$ and $
L_i\delta^\ast\Mscr$ for $i>0$, proving 3. in this case also.
\item We have to show that the adjunction mappings are isomorphisms.
We already know that $\delta_\ast,\delta^\ast$ are exact functors between
$\Qch(W)/\Iso_T(W)$ and $\Qch(V)/\Iso_S(V)$ commuting with
direct limits. Hence it suffices to show that the adjunction mappings are
isomorphisms for    objects annihilated by $t$ and for objects
which are $t$-torsion free.

The case of objects annihilated by $t$ is already covered by Proposition
\ref{ref:9.1.6a}, so we concentrate on $t$-torsion free objects.

First assume that $\Nscr\in \coh(W)$ is a $t$-torsion free object. 
We have an exact sequence \begin{equation}
\label{ref:9.14a}
0\r \Nscr(-1)\r \Nscr \r \bar{\Nscr}\r 0
\end{equation}
Working modulo $\iso_W(T)$ this yields an exact sequence
\begin{equation}
\label{ref:9.15a}
0\r \delta^\ast\delta_\ast\Nscr(-1)\r \delta^\ast\delta_\ast\Nscr \r 
\delta^\ast\delta_\ast\bar{\Nscr}\r 0
\end{equation}
Now by  Proposition \ref{ref:9.1.6a} we already know that 
$ \delta^\ast\delta_\ast\bar{\Nscr}=\bar{\Nscr}$. Combining \eqref{ref:9.14a}
and \eqref{ref:9.15a} using the snake lemma now yields what we want.

The proof that the other adjunction mapping is an isomorphism is entirely
similar.\qed
\end{enumerate}
\def\qed{}\end{proof}
\subsection{The case of a blowing up}
\label{ref:9.2b}
Now we let $X,Y,p$ be as in \S\ref{ref:5a}.  Furthermore we also
assume that $X/k$ is proper.

Up to now we have developed everything without assuming projectivity.
 However in this section we throw in the towel and we introduce the
following hypothesis.
{\def\thehypothesis{(****)}
\begin{hypothesis}
$p$ is a smooth point on $Y$. $o_X(Y)$ is ample on $X$ and the
invertible $o_Y$-bimodule
$I_{Y}$ is ample on $Y$.
\end{hypothesis}
}
The fact that $X$ has an ample invertible bimodule means in particular
that Hypothesis (**) is satisfied. Furthermore it is also easy to
check that under the hypotheses one has $\cd \Gamma(X,-)\le 2$. This
is part of Hypothesis (***).

The smoothness of $Y$ in $p$ implies  that  $\tilde{Y}\r Y$
is an isomorphism  (see Theorem \ref{ref:6.3.1a}).
Also because of smoothness, we have by \eqref{ref:6.16a} an injective
map $t:\Dscr(-1)\r \Dscr$ such that $\bar{\Dscr}=\Dscr/\im t$ is equal
to $\Dscr_Y$ in degree $\ge 1$. By lemma \ref{ref:8.2.1a} it is
also clear that $t$ remains injective after applying the functor
$(-)_{o_X}$.

Since $I_{Y}$ is ample on $Y$ we find that $\Dscr_Y$ and hence
$\bar{\Dscr}$ is ample. Thus Hypothesis (***) holds and hence the
material in \S\ref{ref:9.1a} applies to the current situation. Let
$D$, $D_Y$ be the global sections of $\Dscr$ and $\Dscr_Y$
respectively and put $V=\Proj D$. We now have a commutative diagram of
$k$-quasi-schemes
\begin{equation}
\label{ref:9.16a}
\Atrianglepair<1`1`1`-1`1;>[Y`X`\tilde{X}`V;i`i`i`\alpha`\delta]
\end{equation}
where the arrows marked $i$ are divisors. All maps and compositions of
maps in this diagram are admissible by lemma \ref{ref:7.2.12a}, Proposition
\ref{ref:9.1.8a} and  lemma \ref{ref:7.2.4a}. Further properties of $\delta$
maybe deduced from \S\ref{ref:9.1a}. We define $\Oscr_M=\delta_\ast
\Oscr_L$ and we will refer to $\Oscr_M$ as the exceptional curve in
$V$.

Since we are in a specific situation here, we can of course prove
stronger results than those in \S\ref{ref:9.1a}. First of all note that we
have the following.
\begin{propositions}
\label{ref:9.2.1a}
$D$ satisfies $\chi$ \cite[Def. 3.7]{AZ}.
\end{propositions}
\begin{proof}
Denote the global sections of $\Dscr_Y$ by $D_Y$.
Since $I_{Y}$ is ample, $D_Y$ satisfies $\chi$
\cite{AZ}. Furthermore, by 
Proposition \ref{ref:9.1.3a}.2 the injective map $\bar{D}\r D_Y$ has finite
cokernel. Hence by \cite[Lemma 8.2]{AZ}, $\bar{D}$ will also satisfy
$\chi$. Then we can apply \cite[Thm 8.8]{AZ} to obtain that $D$ also
satisfies $\chi$.
\end{proof}
We also have the following important result.
\begin{propositions}
\label{ref:9.2.2a}
The cohomological dimension of $\delta_\ast$ is less than or equal to
one. Furthermore if $\Nscr\in\gr (\Dscr)$ then $R^1\delta_\ast\Nscr$ is a
finite extension of  quotients of $\Oscr_M$.
\end{propositions}
\begin{proof}
  Let $\Nscr\in \gr(\Dscr)$. By taking a resolution of $\Nscr$ we see
  that it is sufficient to treat the case $\Nscr=\Mscr\otimes_{o_X}
  \Dscr$ where $\Mscr\in\coh(X)$. Now define $\Escr=\oplus_{n\ge 0}
  o_X(nY)$. It is easy to see that $\Dscr$ embeds in $\Escr$ as graded
  $\Dscr$-bimodule. Put $\Nscr'=\Mscr\otimes_{o_X} \Escr$. Let $K,I,C$
  be respectively the kernel, image and cokernel of the induced map
  $\Nscr\r \Nscr'$. Thus we have exact sequences
\begin{gather*}
0\r K \r \Nscr\r I\r 0\\
0\r I\r \Nscr'\r C\r 0
\end{gather*}
It immediately follows from Theorem \ref{ref:5.5.10a}
that $K,C$ are in $\Cscr_p(\Dscr)$. Hence from \eqref{ref:9.5a} it
follows that $R^i\delta_\ast \pi K= R^i\delta_\ast \pi C=0$ for $i>0$.
Since $o_X(Y)$ is ample it  also follows from \eqref{ref:9.5a} that
$R^i\delta_\ast\pi\Nscr'=0$ for $i>0$. Plugging this in the above exact
sequences we find that $R^i\delta_\ast\pi\Nscr=0$ for $i>1$ and also
that $R^1\delta_\ast \pi\Nscr$ is a quotient of $\delta_\ast \pi C$.

Now $C$ itself is not necessarily coherent, but it is a union of $C_i$
which are coherent. Since $\delta_\ast$ is compatible with direct
limits we deduce that $R^1\delta_\ast \pi\Nscr$ is a quotient of some
$\delta_\ast \pi C_i$. From lemma \ref{ref:6.7.3a} it follows that
$\pi C_i$ has a finite filtration whose associated graded quotients
are quotients of $\Oscr_L(t)$. It then follows that $\delta_\ast C_i$ has
a finite filtration whose associated graded quotients are quotients of
$\Oscr_M(t)$. Since these quotients are in $\iso_Y(V)$ they are
invariant under shifting. Hence they are also quotients of $\Oscr_M$.
The result now follows for $R^1\delta_\ast \pi\Nscr$.
\end{proof}
The hypotheses for the following proposition will hold in our applications.
\begin{propositions}
\label{ref:9.2.3a}
  Assume that in addition to Hypotheses (****) we have
  $R\Gamma(X,\Oscr_X)=k$, $H^1(Y,I^n_{Y,p})=0$ for all $n$
  and $\Gamma(Y, \Oscr_Y)=k$.  Then
\begin{enumerate}
\item $R\Gamma(X,\Dscr_{o_X})=D$.
\item There is an exact sequence
\[
0\r D(-1)\xrightarrow{t} D \rightarrow D_Y\r 0
\]
\item $R\delta_\ast L\delta^\ast=\Id$ on $D^-(V)$.
\item $R\Gamma(\tilde{X},\Oscr_{\tilde{X}})=k$
\item $R\Gamma(V,\Oscr_V)=k$.
\end{enumerate}
\end{propositions}
\begin{proof}
\begin{itemize}
\item[1.,2.]
We have $\Dscr_0=\Oscr_X$ and furthermore by \eqref{ref:6.16a} there
are exact sequences for $n\ge 1$
\[
0\r \Dscr_{n-1}\r \Dscr_n \r I^n_{Y,p}\r 0
\]
Looking at the corresponding long exact sequence for $\Gamma(X,-)$
yields what we want by induction.
\item[3.] This follows from 1. since $L\delta^\ast$ can be computed by
  free resolutions. Note that $R\delta_\ast$ is defined on
  $D^-(\tilde{X})$ since according to Proposition \ref{ref:9.2.2a},
  $\delta_\ast$ has finite cohomological  dimension.
\item[4.] We have
\begin{align*}
R\Gamma(\tilde{X},\Oscr_{\tilde{X}})
&=R\Gamma(X, R\alpha_\ast \Oscr_{\tilde{X}})\qquad
\text{(admissibility)}\\
&=R\Gamma(X,\Oscr_X)\qquad (\text{Proposition \ref{ref:8.3.1a}.2}).\\
&=k \qquad\qquad (\text{by hypothesis})
\end{align*}
\item[5.] We have
\begin{align*}
R\Gamma(V,\Oscr_V)&=R\Gamma(V,R\delta_\ast \Oscr_{\tilde{X}})\qquad
\text{(by 1.)}\\
&=R\Gamma(\tilde{X},\Oscr_{\tilde{X}})\qquad\text{(admissibility)}\\
&=k \qquad\qquad (\text{by 4.})\qed
\end{align*}
\end{itemize}
\def\qed{}
\end{proof}

\section{Quantum plane geometry}
\label{ref:10a}
In this section we fix a ``triple'' $(Y,\sigma,\Lscr)$ with $Y$  
a smooth elliptic curve of degree three  in $\PP^2$,   $\sigma\in
\Aut(Y)$ a translation
and $\Lscr=i^\ast \Oscr_{\PP^2}(1)$ where $i:Y\r \PP^2$ denotes the
inclusion. 

We let $A=A(Y,\sigma,\Lscr)$ be the regular algebra associated to
this triple \cite{ATV1} and we put $X=\Proj A$. 
Since $A$ has the Hilbert series of a three dimensional polynomial
algebra it is customary to view $X$ as a non-commutative $\PP^2$.

$A$ contains a central element $g$ in degree three (determined up to a
scalar) such that $Y=\Proj A/gA$ \cite{ATV1,ATV2,AVdB}. Clearly
$Y\subset X$ is a divisor and so the theory of the
previous sections applies to it.

Since $X=\Proj A$, $\Qch(X)$ carries a canonical shift functor which
we denote by $\Mscr\mapsto \Mscr(1)$. The corresponding invertible bimodule is
denoted by $o_X(1)$. Thus  $o_X(Y)=o_X(1)^{\otimes 3}\overset{\text
{not.}}{=}
o_X(3)$. The $o_Y$ bimodule $o_Y(1)$ is defined similarly. By
construction we have $o_Y(1)=\Lscr_{\sigma}$. 

In addition we   use the functors
$\omega$, $\pi$,  $\tilde{-}$ which were defined in
\S\ref{ref:3.8b}.

\subsection{Multiplicities of some objects}
If $\Fscr\in \mod(X)$  then $\Gamma(X,\Fscr)$ is
finite dimensional. We denote by $\Ann_{\Oscr_X} \Fscr$ the kernel of
the  map $\Oscr_X\r \Gamma(X,\Fscr)^\ast \otimes_k\Fscr$
which is obtained from the map $\Gamma(X,\Fscr)\otimes_k \Oscr_X\r
\Fscr$. By using the exactness properties of the completion functor
 (\S\ref{ref:5a}) it follows easily that if $\Fscr\in \Cscr_f$
(cfr \S\ref{ref:5a}) then we have
\begin{equation}
\label{ref:10.1a}
\length (\Oscr_X/\Ann_{\Oscr_X}\Fscr)=\oplus_{p\in Y/\langle\tau\rangle}
\dim_k(C_p/\Ann_{C_p}(\hat{\Fscr}_p)) 
\end{equation}
Let $\Fscr\in \trans_Y(X)$ and let $p\in Y$. Then we define numbers
$t_{p,n}$ for $n\ge 1$ by $T_p(\Fscr)=\oplus_{n\ge 1}
\Oscr_{Y,p}/m^{t_{p,n}}_{Y,p}$ (note that to simplify the notation we
have dropped the completion sign). We consider $(t_{p,n})_n$ as a
partition of the length of $T_p(\Fscr)$ and we let $(r_{p,n}(\Mscr))_n$ stand for
the conjugate partition.  In the sequel we will consider the $(p,n)$
as the points which are infinitely near to $p$. The numbers $r_{p,n}(\Mscr)$
should be viewed as the multiplicities of those infinitely near
points.  We identify $(p,1)$ with $p$.

If $\Mscr\in \mod(X)$ then $\Mscr=\pi M$ for some finitely generated
$A$-module $M$. As usual \cite{ATV2} we have $\dim M_n=(e/d!)n^d+f(n)$ for
$n\gg 0$ where $d$ and $e$ are respectively the Gelfand-Kirillov
dimension and the multiplicity of $M$ and where $f(n)$ is a polynomial
of degree $<d$. We write
$d(\Mscr)=d-1$ and $e(\Mscr)=e$.

We conjecture the following (cfr \cite[Cor 3.7]{H}).
\begin{conjectures}
 Let $\Mscr\in \trans_Y(X)$ and assume that the
  image of $\Mscr$ in $\mod(X)/\Cscr_f$ is simple. Write $e=e(\Mscr)$ and
  $r_{p,n}=r_{p,n}(\Mscr)$ for $p\in Y$ and $n\ge 1$. Then the
  following inequality holds:
\begin{equation}
\label{ref:10.2a}
\sum_{p,n}
\frac{r_{p,n}(r_{p,n}-1)}{2}\le
\frac{(e-1)(e-2)}{2}
\end{equation}
\end{conjectures}
Note that the modern proof of the commutative analogue of this
conjecture 
uses resolution of singularities for a curve through blowing up the
ambient projective plane, together with the fact that an irreducible curve has
non-negative genus. Unfortunately we have not been able to generalize this
proof to the non-commutative case. Indeed it is rather the other way
round. We would like to use a result as \eqref{ref:10.2a} to
deduce properties of our non-commutative blowing up.

Imitating one of the classical proofs of \eqref{ref:10.2a} (see
\cite{Fulton}) as far as possible leads to Proposition
\ref{ref:10.1.2a} below.
\begin{propositions}
\label{ref:10.1.2a}
 Let $\Mscr\in \trans_Y(X)$ and assume that the
  image of $\Mscr$ in $\mod(X)/\Cscr_f$ is simple. Write $e=e(\Mscr)$ and
  $r_{p,n}=r_{p,n}(\Mscr)$ for $p\in Y$ and $n\ge 1$. Assume that
  $a_{p,n}$ are natural numbers such that 
\begin{equation}
\label{ref:10.3a}
\frac{f(f+3)}{2}\ge \sum_{\bar{p},n}\frac{a_{p,n}(a_{p,n}+1)}{2}
\end{equation}
for some integer $f<e$.
Then $ef\ge \sum_{\bar{p},n} r_{p,n} a_{p,n}$. Here the notation
$\bar{p}$ indicates that we take only one representant from each $\tau$-orbit.
\end{propositions}
\begin{proof}
The proof consists of several steps. 
\setcounter{step}{0}
\begin{step} 
\label{ref:1c}
Assume that $M\in \gr(A)$ contains no finite dimensional submodules and let
  $\Mscr=\pi M$. Then
\begin{equation}
\label{ref:10.4a}
\dim \{x\in A_n\mid M_0 x=0\}\ge \dim A_n-\dim
\Gamma(X,(\Oscr_X/\Ann\Mscr)(n))
\end{equation}
To prove this we note that the lefthand side of \eqref{ref:10.4a} is
the degree $n$ part of the kernel of the canonical map $A\r
M^\ast_0\otimes_k M$. If we denote this kernel by $L$ then we have to
prove
\[
\dim_k (A/L)_n\le \dim \omega(\Oscr_X/\Ann\Mscr)_n
\]
Now note that there is a canonical map $M_0\r \Gamma(X,\Mscr)$. If we
then look at the composition
\[
\Oscr_X\r \Gamma(X,\Mscr)^\ast \otimes_k\Mscr\r M^\ast_0\otimes_k
\Mscr
\]
we find $\Ann\Mscr\subset \pi L$.

Hence it suffices to prove that 
\[
\dim_k(A/L)_n\le \dim_k\omega(\Oscr_X/\pi L)_n=\dim_k (A/L)\tilde{}_n
\]
Now $A/L$ is a graded submodule of $M^\ast_0\otimes_k M$ and hence
$A/L$ contains no finite dimensional submodules. Therefore
$A/L\hookrightarrow (A/L)\tilde{}$. This proves what we want.
\end{step}
\begin{step}
\label{ref:2b}
Without loss of generality we may (and we will) assume that
$a_{p,n}=0$ if $r_{p,n}=0$. Furthermore by definition we have
$r_{p,1}\ge r_{p,2}\ge \cdots$ so by permuting the $a_{p,n}$ we may
also assume that $a_{p,1}\ge a_{p,2}\ge a_{p,3}\ge$, because in this
way \eqref{ref:10.3a} remains valid and $\sum_{p,n} r_{p,n}a_{p,n}$
does not decrease.
\end{step}
\begin{step}
\label{ref:3c}
Assume $\underline{a}_p=(a_{p,1},a_{p,2},\ldots)$ is a non-increasing set of
numbers, zero for large $n$.
We associate to such $\underline{a}_p$ an object in $\Cscr_{f,p}$.

Let $(b_{p,n})_{n\ge 1}$ be the partition which is conjugate to $a_{p,n}$.
We put $K_p(\underline{a}_p)=\prod_n
R/m^{b_{p,-n+1}}$. This definition is to be understood as defining a
  set of row vectors with the obvious right $C_p$-action.  Then we
  define $\Kscr_p(\underline{a}_p)$ as the object in $\Cscr_{f,p}$ such
  that $\hat{\Kscr}_p(\underline{a}_p)=K_p(\underline{a}_p)$.  
\end{step}
\begin{step} Below we need the following formula.
\begin{equation}
\label{ref:10.5a}
\dim (C_p/\Ann_{C_p}K_p(\underline{a}_p))=\sum_n
\frac{a_{p,n}(a_{p,n}+1)}{2}
\end{equation}
We leave the obvious verification to the reader. 
\end{step}
\begin{step}
\label{ref:5c}
  Now define $\Kscr=\oplus_p \Kscr_p(\underline{a}_p)$ and $K=\omega
  \Kscr$. 
For $f<e$ we are going to bound the dimension of 
\[
\{x\in A_f\mid K_{m} x=0\}
\]
for large $m$. This can be done by the method exhibited in Step
\ref{ref:1c}. We find that this dimension is bigger than
\[
\dim_k A_f-\dim_k \omega (\Oscr_X/\Ann \Kscr)_{m+f}
\]
Now
\[
\dim_k A_f=\frac{(f+1)(f+2)}{2}
\]
and 
\begin{align*}
\dim _k \omega (\Oscr_X/\Ann \Kscr)_{m+f}&=\length (O_X/\Ann \Kscr)\\
&=\sum_p\length (\Oscr_X/\Ann \Kscr_p(\underline{a}_p))\\
&=\sum_p \dim_k (C_p /\Ann_{C_p}(K_p(\underline{a}_p))\\
&=\sum_{p,n} \frac{a_{p,n}(a_{p,n}+1)}{2}
\end{align*}
where we have used
\eqref{ref:10.1a}\eqref{ref:10.5a}. So ultimately
we obtain for large $m$
\[
\dim_k\{x\in A_f\mid K_{m} x=0\}
\ge \frac{(f+1)(f+2)}{2}-\sum_{p,n} \frac{a_{p,n}(a_{p,n}+1)}{2}
\]
\end{step}
\begin{step}
\label{ref:6b}
Let $\Mscr$ be as in the statement of the proposition. Since
$(r_{p,n})_{p,n}$ is unaffected by $\Cscr_{f}$ we may and we will assume
that $\Mscr$ is $p$-normalized where $p$ runs through a fixed set of
representatives for the $\tau$-orbits in $Y$.  Then by lemma
\ref{ref:5.7.4a} we can recover the structure of $\hat{\Mscr}_p$ from
$T_p(\Mscr)$. Assume that $T_p(\Mscr)=\oplus_n R/m^{t_{p,n}}$. Then it
follows that $\oplus_n K_p(a_{p,1},\ldots,a_{p,t_{p,n}})$ is a
quotient of $\hat{\Mscr}_p$.

Hence $\oplus_{p,n} \Kscr_p(a_{p,1},\ldots,a_{p,t_{p,n}})$ is a quotient
of $\Mscr$. This yields that $L=\oplus_{p,n}
\omega \Kscr_p(a_{p,1},\ldots,a_{p,t_{p,n}})$ is a quotient of
$M=\omega\Mscr$ for large $m$.
\end{step}
\begin{step}
\label{ref:7b}
Since we have chosen $\Mscr$ to be normalized we have in particular
that $\Mscr$ contains no subobject supported on $Y$. Since the only
modules of GK-dimension one of $A$ are coming from objects supported on
$Y$ \cite{ATV2} it follows that $M=\omega \Mscr$ contains no
submodules of GK-dimension one. Since $\Mscr$ was also supposed to be
simple modulo $\Cscr_f$ it follows that $M$ is critical. I.e. $M$
contains no non-trivial submodules of multiplicity strictly smaller
than $e$.

If $0\neq x\in A_f$ then $e(A/xA)=f$. So it follows that
$\Hom_A(A/xA,M)=0$ which is the same as saying that multiplication by
$x$ is injective. Since $\GKdim M=0$ we have for large $m$: $\dim
M_m=em+s$ for some constant $s$.
It follows that  for large $m$ one has
\begin{equation}
\label{ref:10.6a}
\dim (M/Mx)_m=\dim M_m-\dim M_{m-f}=(em+s)-(e(m-f)+s)=ef
\end{equation}
Now by hypotheses
\[
\frac{f(f+3)}{2}+1=\frac{(f+1)(f+2)}{2}> 
\sum_{p,n} \frac{a_{p,n}(a_{p,n}+1)}{2} 
\]
Let $m$ be large. By virtue of Step \ref{ref:5c} there exists a
non-zero $x\in A_f$ such that $K_m x=0$. Now by construction every
indecomposable summand of $L$ 
 is  a quotient of $K$. Thus also $L_mx=0$.

Now by the Step \ref{ref:6b} there is a map $M\r L$, surjective in
high degree. Tensoring with $A/Ax$ yields a map $M/Mx\r L/Lx$ with the
same property.

We conclude
\[
ef=\dim_k(M/Mx)_m\ge \dim_k (L/xL)_m=\dim_k L_m=\sum_{p,n} 
\dim K_p(a_{p,1},\ldots, a_{p,t_{p,n}})
\]
Now an easy verification shows that 
\[
\sum_{n} 
\dim K_p(a_{p,1},\ldots, a_{p,t_{p,n}})=\sum_n a_{p,n} r_{p,n}
\]
which finishes the proof.
\qed\end{step}
\def\qed{}\end{proof}
The following proposition would be a trivial consequence of
\eqref{ref:10.2a}. With a lot more work we can also deduce it from 
Proposition \ref{ref:10.1.2a}.
\begin{propositions}
Assume that $\Mscr\in \trans_Y(X)$ represents a non-zero simple object
in $\mod(X)/\Cscr_f$. Let $r=(r_{p,n}(\Mscr))_{p,n}$, $e=e(\Mscr)$. If the number
of non-zero entries of $r$ is $\le 6$ then there are (up to
permutation) only two possibilities for $e,r$.
\begin{align*}
e=1,\qquad &r=(1,1,1)\\
e=2,\qquad &r=(1,1,1,1,1,1)
\end{align*}
\end{propositions}
\begin{proof}
  We will apply Proposition \ref{ref:10.1.2a}. Since this is
  obviously a numerical criterion we will write $r=(r_i)_{i\in I}$
  where $i$ runs through the pairs $(p,n)$ such that
  $r_{p,n}(\Mscr)\neq 0$. We will consider the entries of $r$ up to
  permutation. Put $m=|I|$.  We will apply Proposition
  \ref{ref:10.1.2a} with suitably chosen $(a_i)_{i\in I}$.
  
  First observe that Proposition \ref{ref:10.1.2a} remains true if
  we allow the $(a_i)_{i\in I}$ to be negative. Indeed if
\[
\frac{f(f+3)}{2}\ge \sum_{i\in I} \frac{b_i(b_i+1)}{2}
\]
with $b_i\in \ZZ$ then we define
\[
a_i=
\begin{cases}
b_i&\text{if $b_i\ge 0$}\\
-b_i-1&\text{if $b_i<0$}
\end{cases}
\]
Then $\sum_i b_i(b_i+1)/2=\sum_i a_i(a_i+1)/2$ and $a_i\ge b_i$. By 
Proposition \ref{ref:10.1.2a} we have $ef\ge \sum_i r_i a_i$ and thus $ef\ge \sum_i 
r_i b_i$.  

Now we can reformulate Proposition \ref{ref:10.1.2a} as follows:
the intersection of the closed ball
\[
B=\{(a_i)_{i\in I}\mid \sum \frac{a_i(a_i+1)}{2} \le
\frac{f(f+3)}{2}\}
\]
and the open half space 
\[
H=\{(a_i)_{i\in I}\mid \sum_i r_i a_i > ef\}
\]
contains no integral point. 

Now $B\cap H$ certainly contains an open ball with diameter equal to
the radius of $B$ minus the distance of $H$ to the center of $B$. 

The center of $B$ is $(-\frac{1}{2},\ldots,-\frac{1}{2})$ and the radius of $B$
is given by
\[
\sqrt{f(f+3)+\frac{m}{4}}
\]
The distance of $H$ to the center of $B$ is given by
\begin{equation}
\label{ref:10.7a}
\frac{ef+\frac{\sum_i r_i}{2}}{\|r\|}
\end{equation}
where $\|r\|=\sqrt{\sum_i r^2_i}$.

If we use the fact that $3e=\sum_i r_i$ then \eqref{ref:10.7a} becomes
equal to 
\[
\frac{e(f+\frac{3}{2})}{2}
\]
Thus $B\cap H$ contains an open ball of diameter
\[
\sqrt{f(f+3)+\frac{m}{4}}-
\frac{e(f+\frac{3}{2})}{\|r\|}
\]
Now it is easy to see that an open  ball of diameter $>\sqrt{m}$ must contain
a point with integral coordinates. Thus Proposition \ref{ref:10.1.2a} yields
\begin{equation}
\label{ref:10.8a}
\sqrt{f(f+3)+\frac{m}{4}}-
\frac{e(f+\frac{3}{2})}{\|r\|}\le\sqrt{m}
\end{equation}
From $\sum_i r_i=3e$ we obtain the estimate $\|r\|\ge
3e/\sqrt{m}$. Combining this with \eqref{ref:10.8a} yields 
\[
\sqrt{f(f+3)+\frac{m}{4}}\le\left( 
   \frac{f}{3}+\frac{3}{2}\right)\sqrt{m} 
\]
We now take $f=e-1$. Then
\[
(e-1)(e+2)+\frac{m}{4}\le
\left( \frac{e}{3}+\frac{7}{6}\right)^2m
\]
which yields 
\[
\frac{(e-1)(e+2)}{
\left(
 \frac{e}{3}+\frac{7}{6}
\right)^2
-
\frac{1}{4}}
\le m
\]
If we combine this with the hypotheses $m\le 6$ then we obtain $e\le
13$.

Explicit enumeration of all possibilities for $e,r$ (using a
computer) yields that $B\cap H$ always contains an integral point for
$e\ge 3$. So the remaining possibilities are $e=1,2$.

If $e=2$ then if follows easily from Proposition \ref{ref:10.1.2a}
that $r=(1,1,1,1,1,1)$. However if $e=1$ then Proposition
\ref{ref:10.1.2a} gives no information whatsoever. But is
$e(\Mscr)=1$ then $\Mscr=\pi M$ where $M$ is a so-called
``line-module'' \cite{ATV2}. We will analyze this case directly.

According to Proposition \ref{ref:5.7.2a} we may change
$\Mscr$ in such a way that $\Mscr$ is $p$ normal, where $p$ runs
through a set of representatives of the $\tau$-orbits on $Y$.  Then we
can read of $r$ from the decomposition of $\Mscr/\Mscr(-Y)$ into
uniserial $\Oscr_Y$-modules.

As indicated above $\Mscr=A/xA$ where $x\in A_1$. But then
$\Mscr/\Mscr(-Y)$ is given by the zeroes of $x\in \Gamma(Y,\Lscr)$.
Since $\Lscr$ defines an embedding of $Y$ in $\PP^2$ this can be
interpreted as the scheme theoretic intersection of $Y$ and $V(x)$.
This is a union of uniserial schemes, such that every point in $Y$
occurs at most once in the reduced locus. Hence we obtain that
$r=(1,1,1)$.
\end{proof}
\subsection{Classification of lines and conics}
The cases $e(\Mscr)=1$ and $e(\Mscr)=2$ represent lines and
conics. The first type of object has been studied in \cite{ATV2} and
the second type of object has been studied in  \cite{Ajitabh}.

If $\Mscr\in \trans_Y(X)$ then let us denote by
$\overline{\Div}(\Mscr)$ the image of $\Div(\Mscr)$ in
$\NN(Y/\langle\tau\rangle)$ under the map which sends $p\in Y$ to its
$\tau$-orbit.

Let us call $M\in \Gr(A)$ \emph{standard} if $M_0\neq 0$ and $M_i=0$ if
$i<0$. Let us call $M\in \gr(A)$ Cohen-Macaulay of dimension $d$ if
$M$ has GK-dimension $d$ and if $\Ext^s_A(M,A)=0$  unless $s=3-d$.
The following is easy to see and will be used below (cfr \cite{Aj1}).
\begin{lemmas}
  Assume that $M\in \gr(A)$ has GK-dimension two and that $M$ contains
  no submodules of GK-dimension one. Then $\tilde{M}$
   is Cohen-Macaulay.
\end{lemmas}
According to \cite{Ajitabh} and \cite{ATV2} there are the following
three classes of critical Cohen-Macaulay modules of GK-dimension two
and multiplicity $1$ or $2$.
\begin{enumerate}
\item Lines: modules of the form $A/xA$, $x\in A_1-\{0\}$.
\item Conics of the first kind: critical graded $A$-modules of the form
  $A/xA$, $x\in A_2-\{0\}$.
\item Conics of the second kind: critical graded $A$-modules with
  minimal resolution
\[ 0\r A(-1)^2\r A^2\r M\r 0\]
\end{enumerate}

Shifts of such modules will be called shifted lines and conics (of the
first and the second kind). They exhaust all critical Cohen-Macaulay
modules of GK-dimension two and multiplicity $\le 2$ \cite{Ajitabh}.

Below $L$ will be an effective divisor on $Y$ such that
$\Lscr=\Oscr_Y(L)$.  The following results from \cite{Ajitabh}
describe the divisor classes of modules of multiplicity $\le 2$.
\begin{propositions}
\label{ref:10.2.2a}
\begin{enumerate} 
\item If $M=N(s)$ where $N$ is a line module then $\Div\pi M$ is of the form
  $\sigma^s L$ for some $s$.
\item
 If $M=Q(s)$ where $Q$ is a conic of the first kind then $\Div
\pi M\sim \sigma^{s} L+ \sigma^{s-1}L$. 
\item
If $M=Q(s)$ where $Q$ is a 
conic of the second kind then $\Div \pi M\sim\sigma^{s} L+ \sigma^s L$.
\end{enumerate}
\end{propositions}
If $D\in \NN(Y/\langle\tau\rangle)$ and $H$ is a divisor  on $E$
then we  say that $D$ is compatible with $H$ if there are
$(p_i)_{i=1,\ldots,l}$ with $l=\deg H$ such that $D=\sum_{i=1}^{l}
\overline{p}_i$ and  $H\sim \sum_{i=1}^l p_i$. 

From the previous proposition we deduce the following:
\begin{lemmas}
\label{ref:10.2.3a}
Assume that $M$ is a Cohen-Macaulay module of GK-dimension 2. If
$e(M)=1$ then $\overline{\Div}(\pi M)$ is compatible with $L$. If $e(M)=2$ then
$\overline{\Div}(\pi M)$ is compatible with $2L$. 
\end{lemmas}
\begin{proof}
This follows from the fact that the image of $\sigma^s L$ in
$\NN(Y/\langle\tau\rangle)$ is compatible with $L$. This follows in
turn from the fact that $\tau=3\sigma$ and $\deg L=3$.
\end{proof}
We also deduce the following lemma.
\begin{lemmas}
\label{ref:10.2.4a}
Assume that $M\in\gr(A)$ is a Cohen-Macaulay module of GK-dimension 2.
Then the fact whether it is a shifted conic of the first or second
kind can be recognized from   $\Div \pi M$.
\end{lemmas}
On the other hand we can't recognize the kind of a conic from its image
in $\trans_Y(X)$. In fact we have the following easy lemma.
\begin{lemmas}
\label{ref:10.2.5a}
Assume that $\Nscr\in \trans_Y(X)$ is such that $e(\Nscr)=2$. Then
$\Nscr$ is equivalent modulo $\Cscr_f$ to an object of the form
$\pi(A/xA)(m)$.
\end{lemmas}
\begin{proof} We may assume that $\Nscr$ is $Y$-torsion free. Let
  $N=\omega(\Nscr)$. If $N$ is a shifted conic of the 
  first kind then we are done. Assume that this is not the case. Let
  $p$ be a point in the support of $\Div(\Nscr)$ and let $\Nscr'$ be
  the kernel of the associated map $\Nscr\r \Oscr_p$. Put
  $N'=\omega\Nscr'$. Then by the formula \eqref{ref:5.39a} together
  with the discussion preceding lemma \ref{ref:10.2.4a} we deduce
  that $N'$ is of the first kind.
\end{proof}
\begin{theorems}
\label{ref:10.2.6a}
 The map $\overline{\Div}$ defines a bijection between
\begin{enumerate} 
\item  simple objects in $\trans_Y(X)/\Cscr_f$ of multiplicity $1$ and
  elements of $\NN(Y/\langle \tau\rangle)$ compatible with $L$;
\item simple objects in $\trans_Y(X)/\Cscr_f$ of multiplicity $2$ and
  elements of $\NN(Y/\langle\tau\rangle)$ compatible with $2L$ that
  are not the sum of two elements of $\NN(Y/\langle\tau\rangle)$
  compatible with $L$.
\end{enumerate}
\end{theorems}
For the proof we need the following easy lemma.
\begin{lemmas} 
\label{ref:10.2.7a}
Assume that $\Mscr\in \trans_Y(X)$ is $Y$-torsion free and assume that
$\Div(\Mscr)=\sum^l_{i=1} p_i$ where the $(q_i)_i$ are in the
$\tau$-orbit of the $(p_i)_i$. Then there exists $\Nscr\in
\trans_Y(X)$ which is $Y$-torsion free and which is equivalent with
$\Mscr$ modulo $\Cscr_f$ such that $\Div(N)=\sum_i q_i$.
\end{lemmas}
\begin{proof} To prove this we first replace $\Mscr$ by $\Mscr(mY)$,
  $m$ large and then apply the formula \eqref{ref:5.39a}.
\end{proof}
\begin{proof}[Proof of Theorem \ref{ref:10.2.6a}] 
  Note that by lemma \ref{ref:10.2.3a} if $\Mscr\in \trans_Y(X)$ then
  $\overline{\Div}(\Mscr)$ is compatible with $L$ or $2L$ if
  $e(\Mscr)$ is $1$ or $2$.

We first prove that
  $\overline{\Div}$ is a bijection between
simple objects in $\trans_Y(X)/\Cscr_f$ of multiplicity $1$ and
  elements of $\NN(Y/\langle \tau\rangle)$ compatible with
  $L$. Surjectivity is obvious so we consider injectivity.

Assume that
  $\Nscr_{1,2}\in \trans_Y(X)$
  have multiplicity $2$ and  $\overline{\Div}(
  \Nscr_1)=\overline{\Div}(\Nscr_2)$. According to lemma
  \ref{ref:10.2.7a} and the previous discussion we can assume that 
  $\Nscr_1=\pi (A/xA)(m)$, $\Nscr_2=(A/yA)(n)$ in such a way that
  $\Div(\Nscr_1)=\Div(\Nscr_2)$. Now consider $x,y$ as global sections
  of $\Lscr$ on $Y$. Then $\Div \Nscr_1=\sigma^m\operatorname{div}(x)$, $\Div
  \Nscr_2=\sigma^n \operatorname{div}(y)$. Since $\sigma$ has infinite order this
  implies $m=n$ and $x=y$ (up to a scalar).
  
Now we prove the second part of the theorem.
   We consider injectivity first. Assume that $\Nscr_{1,2}\in
  \trans_Y(X)$ have multiplicity $=2$ and are such that
  $\overline{\Div}(\Nscr_1)=\overline{\Div}(\Nscr_2)$. Using lemmas
  \ref{ref:10.2.5a},\ref{ref:10.2.7a} 
we may assume that $\Nscr_1=\pi(A/xA)(m)$ and
  $\Div(\Nscr_1)=\Div(\Nscr_2)$, where $\Nscr_2$ is in addition
   $Y$-torsion free. But then by lemma
  \ref{ref:10.2.4a} we find that $\Nscr_2$ is obtained from a
  shifted conic of the first kind, whence
  $\Nscr_2=\pi((A/yA)(n))$. Now exactly as in the case $e=1$ this
  implies $x=y$, $m=n$.
  
Now we have to describe the image of $\overline{\Div}$. 
First let
$D\in \NN(Y/\langle\tau\rangle)$ be compatible with $2L$ and not be
the sum of two elements compatible with $L$. Then according to lemma
\ref{ref:10.2.3a} if $\overline{\Div}(A/xA)=D$ then $A/xA$ is
critical.

Assume on the other hand that $\Mscr\in\trans_Y(X)$ is such that
$e(\Mscr)=2$ and $\overline{\Div}(\Mscr)=D_1+D_2$, $D_{1,2}$ compatible
with $L$. We claim that $\Mscr$ is not simple modulo
$\Cscr_f$. Without loss of generality we may assume that $\Mscr$ is
$Y$-torsion free.

Using lemma \ref{ref:10.2.7a} we may assume that $\Div \Mscr=E_1+E_2$
with $E_1\sim L$, $E_2\sim \sigma^{-1}L$. 
 Put $M=\omega\Mscr$.   

It follows from 
 Proposition \ref{ref:10.2.2a} that $M=A/xA$ with $x\in A_2-\{0\}$. Let $N$ be
 the line module corresponding to $E_1$. Since the divisor of $x$
 contains $E_1$ it follows that $A/xA$ maps surjectively to
 $N/Ng$.  Since $x\in A_2$ this implies
 that $M$ maps surjectively to $N$. Hence $A/xA$ is not critical.
\end{proof}
\section{Blowing up $n$ points in the elliptic
quantumplane}
\label{ref:11a} 
\subsection{Derived categories}
\label{ref:11.1a}In this section our notations and conventions will be as in Section
\S\ref{ref:10a} except that  we
do \emph{not} assume that $\sigma$ has infinite order.

We choose $n$   points $p_1,\ldots,p_n\in Y$ where $n\le 8$ (at this
point not necessarily distinct). We
define quasi-schemes $X_j$, $j=1,\ldots,n+1$, $\tilde{X_j}$,
$j=1,\ldots,n$ 
containing $Y$ as a divisor. We do this as follows~:
\begin{enumerate}
\item $X_1=X$.
\item $\tilde{X_j}$ is the blowup of $X_j$ in the point $p_j$.
\item $X_{j+1}$ is derived from the pair $(X_j,p_j)$ in the same way
  as $V$ was derived from $(X,p)$ in \S\ref{ref:9.2b}.
\end{enumerate}
Of course we are only allowed to make this construction if Hypothesis
(****) of \S \ref{ref:9.2b} holds for the triples $(Y,X_j,p_j)$. Let us
check this now by induction.

Smoothness of the $p_j$ is by hypotheses.

Ampleness of $o_{X_j}(Y)$ for $j\ge 2$ follows from the fact that
$X_j$ is the $\Proj$ of a graded algebra which satisfies $\chi$
(Proposition \ref{ref:9.2.1a}) and $Y$ is defined by a central element
of degree one. This argument also works for $X_1$ if we note that
$X_1=\Proj A^{(3)}$ where $A^{(3)}$ denotes the $3$-Veronese of $A$.

Hence the only thing that remains to be checked is that $I_{Y,p_j}$ is
ample (we have now included the point $p_j$ in the prior notation
$I_Y$). We check ampleness by induction. Assume that $I_{Y,p_{t}}$ is
ample for $t<i$.  By definition we have $I_{Y,p_j}=m_{Y,p_j}
\Nscr_{Y/X_j}$.  By \eqref{ref:9.4a} we have that
$\Nscr_{Y/X_j}=\Nscr_{Y/\tilde{X}_{j-1}}$ and by \eqref{ref:6.26a} we
have $\Nscr_{Y/\tilde{X}_{j-1}}=m_{p_{j-1}} \Nscr_{Y/X_{j-1}}$.

Since $\Nscr_{Y/X}$ is obtained by the functor $-\otimes Ag^{-1}$ and
$Ag^{-1}=A(3)$  one  deduces from \cite{AVdB} that
\[
\Nscr_{Y/X}=(\Lscr_\sigma)^{\otimes 3}
=(\Lscr\otimes_{\Oscr_Y} \sigma^\ast(\Lscr)\otimes_{\Oscr_Y}
\sigma^{2\ast}(\Lscr))_{\sigma^3}
\]
Thus $\Nscr_{Y/X}=\Nscr_\tau$ with $\tau=\sigma^3$ and
$\Nscr=(\Lscr\otimes_{\Oscr_Y} \sigma^\ast(\Lscr)\otimes_{\Oscr_Y}
\sigma^{2\ast}(\Lscr))_{\sigma^3}$. In particular $\deg\Nscr=9$.

Define
\begin{equation}
\label{ref:11.1b}
\Nscr_j=(m_{Y,p_1}\cdots
m_{Y,p_{j-1}}\Nscr)_\tau
\end{equation}
By the above discussion we obtain
$\Nscr_{Y/X_{j}}=\Nscr_{Y/\tilde{X}_{j-1}}=(\Nscr_j)_\tau$ and
$I_{Y,p_j}=(\Nscr_{j+1})_\tau$.

Since $\deg\Nscr_{j}=9-j+1$ we obtain that $\Nscr_{j+1}$ has positive
degree if $j\le 9$. This is true by the restriction $n\ge 8$. Hence
$I_{Y,p_j}$ is ample.

With our current notations the following diagram replaces 
\eqref{ref:9.16a}
\begin{equation}
\label{ref:11.2a}
\Atrianglepair<1`1`1`-1`1;>[Y`X_j`\tilde{X_j}`X_{j+1};i`i`i`\alpha_j`\delta_j]
\end{equation}
We will denote the exceptional curves in $\tilde{X}_j$  by $\Oscr_{L_j}$ and
their direct images in $X_{j+1}$ by $\Oscr_{M_{j+1}}$.
\begin{lemmas}
\label{ref:11.1.1a}
The hypotheses for Proposition
  \ref{ref:9.2.3a} hold for $X_j,Y,p_j$. In particular we have~:
\begin{enumerate}
\item
  $R\Gamma(X_j,\Oscr_{X_j})=R\Gamma(\tilde{X}_j,\Oscr_{\tilde{X}_{j}})=k$. 
\item $R\delta_{{j},\ast} \,L\delta^\ast_{j}$ is the identity on
  $D^{-}(X_{{j}+1})$.
\end{enumerate}
\end{lemmas}
\begin{proof} The fact that  the hypotheses for Proposition
  \ref{ref:9.2.3a} hold is  verified by induction, starting from the
  easy fact that
  $R\Gamma(X,\Oscr_X)=k$.
\end{proof}

We borrow the following definition from the commutative case.
\begin{definitions}  
\label{ref:11.1.2a}
Assume
  that $(q_{j})_{{j}=1,\ldots,n}$, $n\le 8$ are points in $\PP^2$.
The points $(q_{j})_{j}$ are in general position if and
  only if
\begin{itemize}
\item All points are different.
\item
No three  points lie on a line.
\item
No six  points lie on a conic.
\item If $n=8$ then not all points lie on a singular cubic divisor
 in such a way that one of the points lies on the
  singularity.
\end{itemize}

\end{definitions}
The following will be our main theorem.
\begin{theorems}
\label{ref:11.1.3a}
  Assume that the points $(\sigma p_{j})_{j}$ are in general position
  (with respect to the embedding of $Y\subset \PP^2$, fixed in the
  beginning of this section). Then the following holds for all ${j}$~:
\begin{enumerate}
\item $\delta_{{j}\ast},\delta^\ast_{j}$ have cohomological dimension $\le 1$.
\item $L\delta_{j}^\ast$, $R\delta_{{j}\ast}$ are inverse equivalences
 between $D(\tilde{X}_{j})$ and $D(X_{{j}+1})$.
  \item  $\Qch(X_{j})$, $\Qch(\tilde{X}_{j})$ have finite injective dimension.
\end{enumerate}
\end{theorems}
\begin{proof}  Without loss of generality we may  assume that ${j}=n$.
  1.,2. of the theorem will be proved
  directly but for 3. we will need induction. That is, we will assume
  that  $\Qch(X_{j})$ has finite injective dimension for
  ${j}<n+1$. Note that this is clearly true if ${j}=1$ since $X_1$ is the
  $\Proj$ of a graded algebra of finite global dimenion. By Corollary
\eqref{ref:8.4.3a} we find that  $\Qch(\tilde{X}_{j})$ also has
finite injective dimension for ${j}<n+1$.

 We
  already know by lemma \ref{ref:11.1.1a} that $R\delta_{n,\ast}\,
  L\delta^\ast_n$ is the identity on $D^{-}(X_{n+1})$. We will start
  by showing that  $L\delta^\ast_n\, R\delta_{n,\ast}$ is the identity
  on $D^{-}_f(\tilde{X}_{n})$. By lemma \ref{ref:8.1.4a} this means that we
  have to show that $\ker R\delta_{n,\ast}=0$. Hence assume that
  $\Mscr\in D^-_f(\tilde{X}_n)$ is such that
  $R\delta_{n,\ast}\Mscr=0$. We will construct
  $\Qscr_{j}\in D^-_f(X_{j})$ for ${j}=1,\ldots, n$ in such a way that the
  following holds~:
\begin{itemize}
\item[(a)]
$R\Gamma(X_{j},\Qscr_{j})=0$ (note that $R\Gamma(X_{j},-)$ is well defined
on $D^-(X_{j})$ because 
$\Qch(X_{j})$ has finite injective dimension by the induction hypotheses).
\item[(b)]
There are triangles
\begin{align}
\label{ref:11.3a}
&\Atriangle<1`-1`1;>[\Uscr_{n}`L\alpha^\ast_{n}\Qscr_n`\Mscr;``]&\\
\label{ref:11.4a}
&\Atriangle<1`-1`1;>[\Uscr_{{j}}`L\alpha^\ast_{{j}}\Qscr_{j}`L\delta^\ast_{j}\Qscr_{{j}+1};``]
\qquad \text{for ${j}=1,\ldots n-1$}
\end{align}
where the $\Uscr_{j}$ are direct sums of shifts in the derived category
of copies of $\Oscr_{L_{j}}(-Y)$.
\end{itemize}
To do this we define
\begin{align*}
\Qscr_n&=R\alpha_{n,\ast}\Mscr\\
\Qscr_{j}&=R\alpha_{{j},\ast}L\delta^\ast_{{j}} \Qscr_{{j}+1}
\end{align*}
By adjointness we have maps
\begin{align*}
L\alpha^\ast_n\Qscr_n&\r \Mscr\\
L\alpha^\ast_{j} \Qscr_{j}&\r L\delta^\ast_{{j}} \Qscr_{{j}+1}
\end{align*}
which become the identity after applying $R\alpha_{{j},\ast}$,
${j}=1,\ldots,n$. The existence of the triangles
\eqref{ref:11.3a}\eqref{ref:11.4a} now follows from Theorem
\ref{ref:8.4.1a}. 
By admissibility we also have
\[
R\Gamma(X_n,\Qscr_n)=R\Gamma(\tilde{X}_n,\Mscr)=R\Gamma(X_{n+1},
R\delta_{\ast,n}\Mscr )=0
\]
and
\[
R\Gamma(X_{j},\Qscr_{j})=R\Gamma(\tilde{X}_{j}, L\delta^\ast_{{j}}
\Qscr_{{j}+1})=
R\Gamma(X_{{j}+1},\Qscr_{{j}+1})
\]
(for the last equality we have used lemma \ref{ref:11.1.1a}.2).
Hence by induction $R\Gamma(X_{j},\Qscr_{j})=0$. This finishes the proof
of (a) and (b) above.

To continue we use restriction to $Y$. Applying $Li^\ast$ to
\eqref{ref:11.3a} and \eqref{ref:11.4a} yields triangles

\begin{align} \label{ref:11.5a}
&\Atriangle<1`-1`1;>[Li^\ast\Uscr_{n}`Li^\ast\Qscr_n`Li^\ast\Mscr;``]\\
&\Atriangle<1`-1`1;>[Li^\ast\Uscr_{{j}}`Li^\ast\Qscr_{j}`Li^\ast\Qscr_{{j}+1};``]
\qquad \text{for ${j}=1,\ldots n-1$} \label{ref:11.6a}
\end{align}
where the $Li^\ast\Uscr_{j}$ are direct sums of shifts in the derived category
of $\Oscr_{\tau p_{j}}(-Y)=\Oscr_{p_{j}}$.

Now by hypotheses $R\delta_{\ast,n} \Mscr=0$. Thus by  \eqref{ref:9.11a}
we obtain
$
0= Li^\ast R\delta_{\ast,n} \Mscr=Li^\ast\Mscr
$. Hence $Li^\ast \Qscr_n=Li^\ast\Uscr_n$.

Taking into account that $\RHom_{D(Y)}(\Oscr_{p_{j}},\Oscr_{p_t})=0$ if
${j}\neq t$ we deduce that the map $Li^\ast\Qscr_n\r Li^\ast
\Uscr_{n-1}$, must be the zero map in \eqref{ref:11.6a}. Hence
 $Li^\ast \Qscr_{n-1}=Li^\ast \Uscr_{n-1}[-1]\oplus
Li^\ast\Qscr_n=Li^\ast\Uscr_{n-1}[-1]\oplus Li^\ast\Uscr_n$.
Continuing  yields
\begin{equation}
\label{ref:11.7a}
Li^\ast \Qscr_1=Li^\ast \Uscr_1[-n+1]\oplus\cdots \oplus
Li^\ast\Uscr_n
\end{equation}
In particular $Li^\ast\Qscr_1$ is a direct sum of shifts in the
derived category of copies of $\Oscr_{p_{j}}$. 

Now $\Qscr_1$ lives on $X_1=X$ and this is a very well understood
quasi-scheme. In fact it follows from lemma \ref{ref:11.1.4a} below that
necessarily $\Qscr_1=0$. But then from \eqref{ref:11.7a} we immediately
deduce that $\Uscr_{j}=0$.

From the triangle \eqref{ref:11.4a} we then find that
$L\delta_1^\ast\Qscr_2=0$. Applying $R\delta_{\ast,1}$, using lemma
\ref{ref:11.1.1a} yields that $\Qscr_2=0$. Continuing by induction we
eventually find that $\Qscr_n=0$, and hence by \eqref{ref:11.3a} we finally
obtain $\Mscr=0$.

 At this point we have partially proved 2. Let us now prove
 that $\delta^\ast_n$ has cohomological dimension $\le 1$. 
Let  $\Mscr\in \Qch(X_n)$ and put $\Nscr=L\delta_n^\ast\Mscr$. From
the spectral sequence
 
\[
R^p\delta_{n,\ast} H^q(\Nscr)\r  R^{p+q}\delta_{n,\ast} \Nscr
\]
we obtain exact sequences
\begin{equation}
\label{ref:11.8a}
0\r R^1\delta_{n,\ast} H^{i-1}(\Nscr)\r R^i\delta_{n,\ast}\Nscr \r
R^0\delta_{n,\ast} H^i(\Nscr)\r 0
\end{equation}
Since
$R^i\delta_{n,\ast}\Nscr=R^i\delta_{n,\ast}L\delta^\ast_n\Mscr=H^i(\Mscr)$
we deduce that $R^i\delta_{n,\ast}\Nscr=0$ for $i<0$. By
\eqref{ref:11.8a} this yields that $R^1\delta_{n,\ast} H^{i}(\Nscr)=0$ for
$i<-1$ and $R^0\delta_{n,\ast} H^i(\Nscr)$ for $i<0$. In particular
$R\delta_{n,\ast} H^i(\Nscr)=0$ for $ i<-1$. By the part of 2. already
proved, this implies $H^i(\Nscr)=$ for $i<-1$. Hence $\cd
\delta^\ast_n\le 1$. This finishes the proof of 1.

Since it now follows that $L\delta^\ast_n$ is defined on the unbounded
derived category, we can  prove 2. completely. We have to show that
the adjunction mappings are isomorphisms on the unbounded derived
category. Take for example $\Mscr\in D(\tilde{X}_n)$.  We have to show
that $\cone(L\delta^\ast_n R\delta_{n,\ast} \Mscr\r \Mscr)$ is zero. In
the same way as in Step 6 of the proof of Theorem \ref{ref:8.4.1a}
we reduce to the case $\Mscr\in \coh(\tilde{X}_n)$. But then $\Mscr\in
D^-_f(\tilde{X}_n)$ and this case was already handled.

The proof that the other adjunction morphism is an isomorphism is
entirely similar.

By the induction hypotheses (and the discussion in the first paragraph
of this proof) we know that $\Qch(\tilde{X}_n)$ has finite injective
dimension. By an argument similar to the proof of
  Corollary
\ref{ref:8.4.3a} we then find that  $\Qch({X}_{n+1})$ has
finite injective dimension (using 1. and 2.).
 This finishes the proof.
\end{proof}

\begin{lemmas}
\label{ref:11.1.4a} Let $Y,X$ and $(p_{j})$ be as above. Assume  that the
$\sigma p_{j}$ are in general position. Let 
$\Tscr\in D_f^-(X)$ 
be an object such that
\begin{enumerate}
\item
 $R\Gamma(X,\Tscr)=0$
\item
The homology of $Li^\ast\Tscr$ is a direct sum of copies of
$\Oscr_{p_{j}}$.  
\end{enumerate}
Then $\Tscr=0$.
\end{lemmas}
\begin{proof}

Define  $\Escr=\Oscr_X(-2)\oplus \Oscr_X(-1)\oplus \Oscr_X$ and 
$H=\End_{o_X}(\Escr)$. Thus
\[
H=\begin{pmatrix}
k& 0 & 0\\
A_1 & k & 0\\
A_2 & A_1 & k
\end{pmatrix}
\]
By a  standard generalization of  \cite{Beilinson} it follows that
the functor 
\begin{equation}
\label{ref:11.9a}
F:D^-_f(X)\r D^-_f(H): \Tscr=\RHom(\Escr,\Tscr)
\end{equation}
is an equivalence.

Let $e_i$, $i=1,2,3$ be the diagonal idempotents in $H$ and let $P_i=e_iH$
be the corresponding projectives. Under $F$ we have the following
correspondences
\[
\Oscr_X(-2)\leftrightarrow P_1\qquad  \Oscr_X(-1)\leftrightarrow P_2
\qquad \Oscr_X\leftrightarrow P_3
\]
From  \eqref{ref:11.9a} we deduce that 
\[
R\Gamma(X,\Tscr)=0\iff (F\Tscr)e_3=0
\]
Put $H'=(1-e_3)H(1-e_3)=\bigl(\begin{smallmatrix} k&0\\ A_3 &
  k\end{smallmatrix} \bigr)$.

Right modules over $H$ can be written in block form as row vectors $(M_1,
M_2, M_3)$. Similarly $H'$-modules can be written as
$(M_1,M_2)$. Sending 
$(M_1,M_2)\mapsto (M_1,M_2,0)$ defines an exact
functor $I:\Mod(H')\r \Mod(H)$ which extends to an exact functor
$I:D^-_f(H')\r D^-_f(H)$. It is easy to see that this functor defines an
equivalence of $D^-_f(H')$ with the full subcategory of $D^-_f(H)$
consisting of objects $T$ such that $Te_3=0$.
Hence we find in particular that 
$\ker R\Gamma(X,-)$ is equivalent to $ D^-_f(H')$.

Now $H'$ is hereditary, and hence by lemma
\ref{ref:8.3.4a} an element of $D^-_f(H')$ is the sum of its
homology.
Assume that $T\in\mod(H')$. Then $T$ has a minimal resolution of length
$2$
\[
0\r W\otimes_k P_1\r U\otimes_k P_2 \oplus V\otimes_k P_1 \r T\r 0
\]
The complex $W\otimes_k P_1\r U\otimes_k P_2 \oplus V\otimes_k P_1$
splits as a direct sum of $W\otimes_k P_1\r U\otimes_k P_2$ and
$V\otimes_k P_1$.

Now we view $T$ as a $H$-module via the functor $I$ defined above. We
find that  $F^{-1}T$ is a sum of
\begin{equation}
\label{ref:11.10a}
W\otimes_k \Oscr_X(-2)\r U\otimes_k \Oscr_X(-1)
\end{equation}
and $V\otimes_k \Oscr_X(-2)$. We conclude that 
 $\Tscr$ is a direct sum of shifts
in the derived category of 
complexes of the form \eqref{ref:11.10a} and of complexes of the form
$\Oscr_X(-2)$. 

By hypotheses we know in addition that the homology of $Li^\ast\Tscr$
is given by direct sums of copies of $\Oscr_{p_{j}}$. This yields  that
$Li^\ast\Tscr$ must be a direct sum of complexes
\[
W\otimes_k \Oscr_Y(-2)\xrightarrow{\phi} U\otimes_k \Oscr_Y(-1)
\]
where $\phi$ is injective and $\coker \phi$ is a direct sum of copies of
$\Oscr_{p_{j}}$. In particular $\dim W=\dim U$.

For convenience we tensor \eqref{ref:11.10a} with $o_Y(2)$. 
Taking into account that
$\Oscr_Y(1)=\sigma_\ast(\Lscr)$
we obtain a
 complex
 \[
 W\otimes_k \Oscr_Y\rightarrow U\otimes_k \sigma_\ast(\Lscr)
 \]
 whose cokernel which is a direct sum of copies of
 $\Oscr_{p_{j}}(2)=\Oscr_{\sigma^{2} p_{j}}$. 

Now we invoke Lemma \ref{ref:11.1.6a} below with $\Mscr=\sigma_\ast(\Lscr)$ and
$q_{j}=\sigma^2 p_{j}$.
 By hypotheses the $q_{j}$ are in general
position with respect to the embedding defined by $\Mscr$. We obtain
$W=U=0$ and the proof is done.
\end{proof}
The above proof was based on lemma \ref{ref:11.1.6a} below. To prove this
lemma we need some notions of the theory of commutative blowing up.

Assume that $Z$ is a smooth (commutative!) surface let $Y\subset Z$ be a
 divisor. Let $p\in Y$ be a smooth point and let
$\psi:\tilde{Z}\r Z$
and $\tilde{Y}$  be respectively
the  blowup of $Z$ in $p$ and the strict transform of
$Y$. As usual $\tilde{Z}=\Proj (\oplus_n m_p^n)$, where $m_p$ is the
maximal ideal of $\Oscr_Z$ defined by $p$. Let $L\subset \tilde{Z}$ be the
exceptional curve.

Let $\Pscr$ be a one-dimensonal coherent Cohen-Macaulay module on $Z$.
We define the \emph{strict transform} $\psi^{-1}(\Pscr)$ of $\Pscr$ as
$\pi (\oplus_n m_p^n\Pscr)$. The following lemma is easily proved.
\begin{lemmas} 
\label{ref:11.1.5a}
Let the notation be as above. Assume that $Y$
  is not contained in the support of $\Pscr$ and furthermore that the
  $p$-primary component of $i^\ast \Pscr$ is of the form
  $\Oscr_p^{\oplus u}$. Then  $\Supp(\psi^{-1}\Pscr)\cap
  \tilde{Y}\cap L=\emptyset$. 
\end{lemmas}
We use this lemma to prove the following result.

\begin{lemmas} 
\label{ref:11.1.6a}
Let $Y$ be embedded as a cubic divisor in $\PP^2$
  through a very ample line bundle $\Mscr$ of degree three. Assume that
  $(q_{j})_{{j}=1,\ldots,n}\in Y$, $n\le 8$,
  are smooth points in general position. Then it is impossible to have a
  map $\Oscr_Y^u\xrightarrow{\phi} \Mscr^u$ whose cokernel is a direct
  sum of copies of $\Oscr_{q_{j}}$, unless $u=0$.
\end{lemmas}
\begin{proof}
Let $t:Y\r \PP^2$ be the embedding. The map $\phi$ can be uniquely
lifted to a map $\Oscr_{\PP^2}^u\xrightarrow{\mu} \Oscr_{\PP^2}^u$
such that $\phi=t^\ast \mu$. Let $\Pscr=\coker  \mu$. Then $\Pscr$ is
a  one dimensional Cohen-Macaulay module on $\PP^2$ such that
$t^\ast\Pscr$ is a direct sum of copies of $\Oscr_{q_{j}}$. We claim that this is
impossible unless $\Pscr=0$. 

To do this we perform the blowup of $\PP^2$ in $q_1,\ldots,q_n$. Let
$\tilde{\Pscr}$ and $\tilde{Y}$ be the (iterated) strict transforms of
$\Pscr$ and $Y$. Then according to lemma \ref{ref:11.1.5a} we have $\Supp
\tilde{\Pscr}\cap \tilde{Y}=\emptyset$. On the other hand it follows
from the Nakai criterion that if the $q_{j}$ are in general position
then $\tilde{Y}$ is ample. This is clearly a
contradiction.
\end{proof}

\subsection{Exceptional simple objects}
If $Y\subset X$ is a commutative curve contained as a divisor in a
quasi-scheme $X$ then we will call a simple object in $\mod(X)$
\emph{exceptional} if it is not of the form $\Oscr_p$ for $p\in Y$.
One of the aims of these notes is to count the exceptional simple
objects in the quasi-schemes $X_n$ which were introduced in the
previous section. So far we have only been able to do this under some
additional hypotheses as can be seen from our main result below
(Theorem \ref{ref:11.2.1a}). In this section we use the same notations
and hypothese as
in sections \S\ref{ref:10a} and
\S\ref{ref:11.1a}. We assume in addition that $\tau$ has
infinite order.

Throughout we will choose a group law on $Y$ in such a way that if
$p,q,r\in Y$ lie on a line then $p+q+r=0$. Furthermore we choose
a fixed set $F\subset E$ of representatives for the $\tau$-orbits on
$Y$. We partially order $\NN F$ by putting $y\le z$ if $y_p\le z_p$
for all $p\in F$. If $z\in\NN F$ then we write $|z|=\sum_{p\in F} z_p$,
$N(z)=\sum_{p\in F} z_p p$.

For $n\in \NN$ we define
\[
H_n=\{y\in \NN F\mid |y|=n,N(y)\in \ZZ\tau\}
\]
and for  $z\in \NN F$ we also define
\begin{gather*}
  A_z=\{y\in H_3\mid y\le z\}\\
  B_z=\{y\in H_6\mid y\le z, \text{$y$ is not the sum of two elements of
    $H_3$}\}
\end{gather*}
We now have the following theorem.
\begin{theorems}
\label{ref:11.2.1a}
Let $p_1,\ldots,p_n\in Y$, $n\le 6$ be such that $(\sigma p_i)_i$ are
in general position (Def. \ref{ref:11.1.2a}).  For $p\in F$ let
$z_p$ be the cardinality of the intersection of $\{p_1,\ldots,p_n\}$
with the $\tau$-orbit of $p$. Let $O$ be the number of non-zero $z_p$.
Then the number of non-isomorphic exceptional simple objects in
$\mod(X_{n+1})$ is equal to $ n+|A_z|+|B_z|-O $.
\end{theorems}
The proof of this theorem will follow rather easily from our previous
results.  We start with the following lemma.
\begin{lemmas}
\label{ref:11.2.2a} Assume that $\Ascr$, $\Bscr$ are two abelian categories
and $\Cscr\subset\Ascr$, $\Dscr\subset \Bscr$ are two abelian
subcategories closed under subquotients. Assume that there are inverse equivalence $F$, $G$
between $D^b(\Ascr)$ and $D^b(\Bscr)$. Assume furthermore that for all
$i$, $H^iF$ sends $\Cscr$ to $\Dscr$ and $H^iG$ sends $\Dscr$ to
$\Cscr$. Then the maps
\begin{equation}
\label{ref:11.11a}
\begin{gathered}
\bar{F}:\Cscr\r\Dscr:[C]\mapsto \sum (-)^i [H^i(FC)]\\
\bar{G}:\Dscr\r \Cscr:[D]\mapsto \sum (-)^i [H^i(GD)]
\end{gathered}
\end{equation}
define isomorphisms between the Grothendieck groups of $\Cscr$ and
$\Dscr$. 
\end{lemmas}
\begin{proof}  Let $\bar{\Cscr}$ and $\bar{\Dscr}$ be the closures of
  $\Cscr$ and $\Dscr$ under extensions. Then clearly $F$ and $G$
  define inverse equivalences between $D^b_{\bar{\Cscr}}(\Ascr)$ and
  $D^b_{\bar{\Dscr}}(\Bscr)$. Hence we obtain (using standard
  isomorphisms for Grothendieck groups)
\[
K_0(\Cscr)
\cong K_0(\bar{\Cscr})\cong K_0(D^b_{\bar{\Cscr}}(\Ascr))
\cong K_0(D^b_{\bar{\Dscr}}(\Bscr))\cong K_0(\bar{\Dscr})=K_0(\Dscr)
\]
The composition of these isomorphisms (and their inverses) is given by
\eqref{ref:11.11a}.
\end{proof}
We deduce the following result.
\begin{lemmas} Let $p_1,\ldots,p_n\in Y$, $n\le 8$ be such that
  $(\sigma p_i)_i$ are in general position. Let $z\in \NN(Y/\langle
  \tau\rangle)$. Then $K_0(M_z(\tilde{X_i}))$ and $K_0(M_z(X_{i+1}))$
  are isomorphic for $i\le n$ (see \S\ref{ref:6.9b} for
  notations). 
\end{lemmas}
\begin{proof}
  By Theorem \ref{ref:11.1.3a} we know that $L\delta_{i}^\ast$
  and $R\delta_{i,\ast}$ are inverse equivalences between
  $D^b(\tilde{X}_{i})$ and $D^b(X_{{i+1}})$.
  
  Now using the explicit construction of the functors $T_p$ in
  Proposition \ref{ref:5.7.2a} it is easy to verify that
  for $\Mscr\in \trans_Y(X_{{i}+1})$ we have
  $T_p(\Mscr)=T_p(\delta_{i}^\ast\Mscr)$. Similarly using Proposition
  \ref{ref:7.2.4a} we have for $\Nscr\in \trans_Y(\tilde{X})$ the
  identity $T_p(\delta_{{i},\ast}\Nscr)=T_p(\Nscr)$.  Finally using
  Theorem \ref{ref:9.1.9a} we find that the higher derived functors of
  $\delta^\ast_{i}$ and $\delta_{{i},\ast}$ map $\mod(X_{{i}+1})$ to
  $\iso_Y(\tilde{X}_{{i}})$ and $\mod(\tilde{X}_{i})$ to
  $\iso_Y(X_{i+1})$.

It follows that we can apply lemma \ref{ref:11.2.2a} to obtain an
isomorphism between $M_z(\tilde{X}_{i})$ and $M_z(X_{{i}+1})$.
\end{proof}
We now obtain.
\begin{lemmas}
\label{ref:11.2.4a}
  Let $p_1,\ldots,p_n\in Y$, $n\le 8$ be such that $(\sigma p_i)_i$
  are in general position. Let $z$ and $O$ be as in the statement of
  Theorem \ref{ref:11.2.1a} (where we identify $\NN F$ with
  $\NN(Y/\langle\tau\rangle)$). Then the number of non-isomorphic
  exceptional simple objects in $\mod(X_{n+1})$ is equal to
$
n+\rk K_0(M_z(X_0))-O
$.
\end{lemmas}
\begin{proof} By iterating Theorem \ref{ref:6.9.1a} and lemma
  \ref{ref:11.2.4a} we easily find
\[
\rk K_0(M_0(X_{n+1}))=n+\rk K_0(M_z(X_1))-O
\]
It is now easy to see that $M_0(X_{n+1})$ is equivalent to $\iso_Y(X_{n+1})$.
Clearly the exceptional simple objects in $\mod(X_{n+1})$ coincide
with the simple objects in $\iso_Y(X_{n+1})$. It now suffices to check
that $\iso_Y(X_{n+1})$ is a finite length category.

Note that by construction $X_{n+1}=\Proj A_{n+1}$ where
  $A_{n+1}$ is a noetherian graded ring containing a regular central
  element $t$ in degree one (see  \S\ref{ref:9.1a} and in particular
  Prop. \ref{ref:9.1.3a}). If $M$ is a graded $A_{n+1}$-module and
  $\Mscr=\pi M$ then multiplication by $t$ corresponds to the map
  $\Mscr(-Y)\r \Mscr$.
  It easily
  follows that if $M$ is a graded $A_{n+1}$-module of Gelfand-Kirillov
  dimension $>1$ then $\pi M$ is not in $\iso_Y(X_{n+1})$. Hence 
  objects in $\iso_Y(X_{n+1})$ come from graded $A$-modules of
  Gelfand-Kirillov dimension one. By considering multiplicity one
  finds that the latter form a finite length category when viewed in
  $\Proj A_{n+1}$.
\end{proof}
\begin{proof}[Proof of Theorem \ref{ref:11.2.1a}]
By lemma \ref{ref:11.2.4a} we have to compute $K_0(M_z(X_0))$ in the
special case that $|z|\le 6$ (by our assumption that $n\le 6$). 

Remember that $X_0=\Proj A$ where $A$ is a three dimensional elliptic
Artin-Schelter regular algebra. By considering the
action of the central element in degree three it is easily seen that
every object in $\trans_Y(X_0)$ is of the form $\pi M$ with $\GKdim
M\le 2$. Since the graded $A$-modules with $\GKdim \le 1$ correspond
to objects in $\Cscr_f$ \cite{ATV1} it follows by considering
multiplicity that $\trans_Y(X_0)/\Cscr_f$ is a finite length category.
The same holds for $M_z(X_0)$ so a basis for $K_0(M_z(X_0))$ is given
by the isomorphism classes of simple objects. These simple objects
have been classified in Theorem \ref{ref:10.2.6a}. There are in 1-1
correspondence with the following two sets.
\begin{align*}
  A'_z&=\{D\in \NN(Y/\langle\tau\rangle)\mid  D\le z, \text{$D$
    is compatible
    with $L$}\}\\
    B'_z&=\{D\in\NN(Y/\langle\tau\rangle)\mid  D\le z,\text{$D$
      is compatible with $2L$}\\
&\qquad \text{ and $D$ is not a sum $D_1+D_2$ with
    $D_i$ compatible with $L$}\}
\end{align*}
where $L$ represents the divisor of a line in $\PP^2$. Thus
$\rk K_0(M_z(X_{n+1}))=|A'_z|+|B'_z|$. It is now easy to see that
$A_z$ is in bijection with $A'_z$ and similarly $B_z$ is in bijection
with $B'_z$. This finishes the proof.
\end{proof}
\section{Non-commutative cubic surfaces}
\label{ref:12a}
 In this section we
recycle notations and assumptions from  section
\S\ref{ref:11a}.  We will
however assume in addition 
 that $n=6$. Thus will fix six points $(p_{j})_{j}$ on
$Y$ and our aim will be to study $X_7$. Since we will make no use
of Theorem \ref{ref:11.1.3a}, we will not assume that the
points $(\sigma p_{j})$ are in general position. Furthermore we will
also not assume that $\tau$ has infinite order.

 By construction $X_7=\Proj F$ for a certain graded
$k$-algebra $F$. Since the hypotheses for Proposition \ref{ref:9.2.3a}
hold on $X_6$ (lemma \ref{ref:11.1.1a}) we find by that proposition that $F$
contains a regular central element $t$ in degree one such that
$\bar{F}=F/tF$ is the twisted homogeneous coordinate ring \cite{AVdB}
associated to the triple $(Y,\Nscr_7,\tau)$. Here $\Nscr_7$ is a
line bundle of degree $9-6=3$, defined by \eqref{ref:11.1b}.

From these data we can compute the Hilbert-series of $F$. We find
$H(F,s)=(1-s^3)/(1-s)^4$. This suggests that $X_7$ should be viewed
as a non-commutative cubic surface. 
In this section we substantiate this
intuition by showing that there exists a  4-dimensional
Artin-Schelter regular algebra $P$ \cite{AS}, containing a normal
element $C$ in degree three such that $F=P/(C)$.  We can then view
$\Proj P$ as a quantum $\PP^4$ which contains $X_7$ as a cubic divisor.

\begin{remarks} Of course the commutative analogue of this is
  well-known. If one blows up 6 points in general position in $\PP^2$
  then one obtains a cubic surface in $\PP^4$ \cite{H}. However the
  reader may wonder why, in the non-commutative case, we don't need
  that our points are in general position. The explanation is of
  course that $X_7$ is not a straight blowing up of $X_1$, but is
  constructed by repeatedly applying the constructions $X_{j}\mapsto
  \tilde{X}_{{j}}\mapsto X_{{j}+1}$.  In the commutative case, if the points
  are in general position, then $\delta_{j}$ is an isomorphism between
  $\tilde{X}_{j}$ and $X_{{j}+1}$ and so in that case $X_7$ is indeed a straight
  blowing up of $X_1$. This will in general not be true in the
  non-commutative case, except in a derived sense. See Theorem
  \ref{ref:11.1.3a}.
\end{remarks}

The construction of $P$ is easy. It follows from \cite{ATV1}
that $\bar{F}$ has a (minimal) presentation
\[
\bar{F}=k[x_1,x_2,x_3]/(r_1,r_2,r_3,C_3)
\]
where $\deg x_1=1$, $\deg r_i=2$, $\deg C_3=3$.  One deduces that 
 $F$ has a presentation
\[
F=k[x_1,x_2,x_3,t]/(r'_1,r'_2,r_3',C_3',[t,x_1],[t,x_2],[t,x_3])
\]
where   $r'_i$, $C'_3$ are homogeneous liftings of $r_i$, $C_3$. 

We now put
\[
P=k[x_1,x_2,x_3,t]/(r'_1,r'_2,r'_3,[t,x_1],[t,x_2],[t,x_3]) 
\]
and we will show that $P$ is Artin-Schelter regular. To this end we
make use of the fact that by \cite{ATV1} one knows that
$\bar{P}=P/tP=k[x_1,x_2,x_3]/(r_1,r_2,r_3)$ is a three-generator
three-dimensional Artin-Schelter regular algebra \cite{AS}. We can then use
the criterion \cite[Cor. 2.7]{LSV}.

To state this criterion we let $R_P$ and  $R_{\bar{P}}$ stand for the
relations of degree two in $P$ and $\bar{P}$. Suppose we have the
following
\begin{enumerate}
\item Left and right multiplication by $t$ is injective on $P_1$.
\item The image of $(P_1\otimes R_P)\cap (R_P\otimes P_1)$ under the
  natural map $P^{\otimes 3}\r \bar{P}^{\otimes 3}$ is   
$(\bar{P}_1\otimes R_{\bar{P}})\cap (R_{\bar{P}}\otimes \bar{P}_1)$.
\end{enumerate}
Then according to \cite[Cor. 2.7]{LSV}, $P$ will be Artin-Schelter
regular with Hilbert series $1/(1-s)^4$.

We now verify conditions 1.,2.  Condition 1. follows from the
observations that  $P_{\le 2}=F_{\le 2}$
and that $F$ is a domain  since $\bar{F}$ is a domain.

Hence we concentrate on condition 2. To simplify the notations we put
$V=\sum_i kx_i$, $R_2=\sum_i kr_i\subset V^{\otimes 2}$,
$R_3=kC_3\subset  V^{\otimes 3}$, $W=V\oplus kt$,
$S_2=(\sum_i kr'_i)+ (\sum_ik[t,x_i])\subset W^{\otimes 2}$,
$S_3=kC_3'\subset W^{\otimes 3}$.  With these notations
\begin{align*}
\bar{F}&=k[V]/(R_2\oplus R_3)\\
F&=k[W]/(S_2\oplus S_3)\\
P&=k[W]/(S_2)\\
\bar{P}&=k[V]/(R_2)
\end{align*}
Note that $R_P=S_2$,
$R_{\bar{P}}=R_2$.

 We first claim that the
following complexes are exact in degrees $\le 3$. 
\begin{equation}
\label{ref:12.1a}
0\r (V\otimes R_2\cap R_2\otimes V)\otimes \bar{F}\r
(R_2\oplus R_3)\otimes \bar{F} \r V\otimes \bar{F}\r \bar{F}\r k \r0
\end{equation}
\begin{equation}
\label{ref:12.2a}
0\r (W\otimes S_2\cap S_2\otimes W)\otimes F\r
(S_2\oplus S_3)\otimes {F} \r W\otimes {F}\r {F}\r k \r0
\end{equation}
Let us first consider \eqref{ref:12.1a}. The only place where exactness
is non-obvious is
at $(R_2\oplus R_3)\otimes \bar{F}$. Hence it is sufficient
to show that the alternating sum of the Hilbert series of the terms in
\eqref{ref:12.1a} is zero in
degrees $\le 3$. This easily follows from the fact 
\[
\dim (V\otimes R_2\cap R_2\otimes V)=1
\]
which is true because  $\bar{P}$ is Koszul.

We use the same method to check the exactness of
\eqref{ref:12.2a}. This time we need 
\[
\dim(W\otimes S_2\cap S_2\otimes W)=4
\]
Now we know that the dimensions of 
\begin{align*}
F_2&=W^{\otimes 2}/ (S_2)\\
F_3&=W^{\otimes 3} /(S_3+W\otimes S_2+ S_2\otimes W)
\end{align*}
are equal to $10$ and $19$ respectively. This yields that
\begin{align*}
\dim S_2&=16-10=6\\
\dim (S_3+W\otimes S_2+ S_2\otimes W)&=64-19=45
\end{align*}
Since by \cite{ATV1} we have $R_3\cap (V\otimes R_2+R_2\otimes V)=0$, it
also follows that
$
S_3\cap(W\otimes S_2+S_2\otimes W)=0
$.
Thus 
\[
\dim (W\otimes S_2+S_2\otimes W)=44
\]
Hence we obtain that
\[
\dim(W\otimes S_2\cap S_2\otimes W)=4\times 6+6\times 4-44=4
\]
This proves what we want.

Now we tensor \eqref{ref:12.2a} with $\bar{F}$ and we combine the result
with \eqref{ref:12.1a} to form the following commutative diagram.

\begin{tiny}
\[
\begin{CD}
 @. @. W\otimes \bar{F} @>\gamma>> k\otimes \bar{F} @. @. @.\\
@. @. @V\alpha VV @V\beta VV @. @. @.\\
0 @>>> (W\otimes S_2\cap S_2\cap W)\otimes \bar{F} @>>> (S_2\oplus
S_3)\otimes \bar{F} @>>> W\otimes \bar{F} @>>> \bar{F} @>>> k@>>> 0\\
@. @VVV @VVV @VVV @VVV @.\\
0 @>>> (V\otimes R_2\cap R_2\cap V)\otimes \bar{F}  @>>> (R_2\oplus
R_3)\otimes \bar{F} @>>> V\otimes \bar{F} @>>> \bar{F} @>>> k @>>> 0\\
@. @. @VVV @VVV @. @. @.\\
 @. @. 0 @. 0 @. @. @.\\
\end{CD}
\]
\end{tiny}

$\alpha$, $\beta$ and $\gamma$ are defined by
\begin{align*}
\alpha(w\otimes 1)&=(tw-wt)\otimes 1\\
\beta(1\otimes 1)&=t\otimes 1\\
\gamma(w\otimes 1)&=1\otimes \bar{w}
\end{align*}
The homology of the middle complex is given by $\Ext^i_F(k,\bar{F})$
in degrees $\le 3$. In particular this complex is exact at $(S_2\oplus
S_3)\otimes \bar{F}$. 

It is also clear that $\gamma$ is surjective in degrees $\ge 1$. A
trivial diagram chase now shows that
\[
W\otimes S_2\cap S_2\otimes W\r V\otimes R_2\cap R_2\otimes V
\]
is surjective. This completes the proof of conditions 1. and 2. above.

So at this point we know that $P$ is Artin-Schelter regular with
Hilbert series $1/(1-s)^4$. We still have to show that $C'_3$ is a
regular normalizing element in $P$.

By looking at Hilbert series it is clear that $t$ is a regular central
element in $P$.  Hence since $\bar{P}$ is a domain by \cite{ATV1}, the
same holds for $P$.  In particular $C'_3$ is regular in $P$.  Looking
at Hilbert series of $P$ and $F$ reveals that the twosided ideal
$(C_3')$ in $P$ must be free of rank one on the left and on the right.
Hence $(C_3')=C_3'P=PC_3'$ and thus $C_3'$ is normalizing. This
completes the proof.

\appendix
\section{Two-categories}
\label{ref:Aa}
A 2-category is a category where the homsets themselves
are categories. The objects of such a category are called 0-cells, the
arrows are called 1-cells and the arrows between arrow are called
2-cells. Such 2-cells are drawn as follows
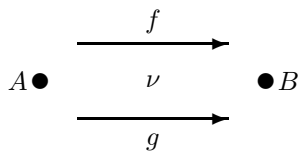
\setlength{\unitlength}{1mm}
\thicklines
\[
\begin{picture}(40,20)
\put(2,10){\makebox(0,0){$A$}}
\put(38,10){\makebox(0,0){$B$}}
\put(20,2){\makebox(0,0){$g$}}
\put(20,18){\makebox(0,0){$f$}}
\put(20,10){\makebox(0,0){$\nu$}}
\put(10,5){\vector(1,0){20}}
\put(10,15){\vector(1,0){20}}
\put(5,10){\circle*{2}}
\put(35,10){\circle*{2}}
\end{picture}
\]
As usual arrows can be composed, and so can 2-cells. It turns out that
2-cells even have two compositions. Vertical ones
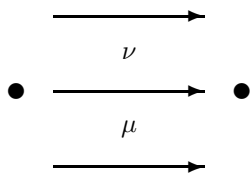
\[
\begin{picture}(30,20)
\put(5,0){\vector(1,0){20}}
\put(5,10){\vector(1,0){20}}
\put(5,20){\vector(1,0){20}}
\put(0,10){\circle*{2}}
\put(30,10){\circle*{2}}
\put(15,5){\makebox(0,0){$\mu$}}
\put(15,15){\makebox(0,0){$\nu$}}
\end{picture}
\]
denoted by $\mu\cdot \nu$, which come from the composition in
$\Hom(A,B)$ and horizontal ones
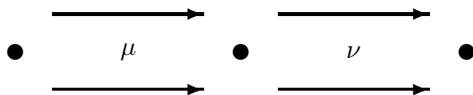
\[
\begin{picture}(60,10)
\put(0,5){\circle*{2}}
\put(30,5){\circle*{2}}
\put(60,5){\circle*{2}}
\put(15,5){\makebox(0,0){$\mu$}}
\put(45,5){\makebox(0,0){$\nu$}}
\put(5,0){\vector(1,0){20}}
\put(35,0){\vector(1,0){20}}
\put(5,10){\vector(1,0){20}}
\put(35,10){\vector(1,0){20}}
\end{picture}
\]
denoted by $\mu\nu$ which come from the fact that the pairing
$\Hom(B,C)\times \Hom(A,B)\r \Hom(A,C)$ has to be a bifunctor. Between
those two compositions there is a natural compatibility. Assume that
one has the following diagram
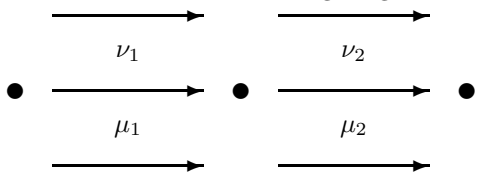
\[
\begin{picture}(60,20)
\put(5,0){\vector(1,0){20}}
\put(5,10){\vector(1,0){20}}
\put(5,20){\vector(1,0){20}}
\put(0,10){\circle*{2}}
\put(30,10){\circle*{2}}
\put(15,5){\makebox(0,0){$\mu_1$}}
\put(15,15){\makebox(0,0){$\nu_1$}}
\put(35,0){\vector(1,0){20}}
\put(35,10){\vector(1,0){20}}
\put(35,20){\vector(1,0){20}}
\put(60,10){\circle*{2}}
\put(45,5){\makebox(0,0){$\mu_2$}}
\put(45,15){\makebox(0,0){$\nu_2$}}
\end{picture}
\]
then one has
$\mu_1\mu_2\cdot\nu_1\nu_2=(\mu_1\cdot\nu_1)(\mu_2\cdot\nu_2)$.  Of
course since a set is a category with only the identity arrows we can
consider every category trivially as a 2-category.

The archetypical example of a 2-category is ``$\mathbf{Cat}$,'' the
category of all categories. In this case the objects are the
categories (living in some universe), the arrows are the functors and
the 2-cells are the natural transformations. It is therefore not
surprising that the standard properties of categories and functors can
be mimicked inside a 2-category.

For example an arrow $A\xrightarrow{f}B$ is a left adjoint of an arrow
$B\xrightarrow{g}A$ if there is a unit $\eta:\Id_{A}\r gf$ and a
counit $\epsilon:fg\r \Id_B$ satisfying the standard associativity
conditions.  As usual $g$ is determined by $f$ up to unique
isomorphism. If in this situation the $\eta$ and $\nu$ are
isomorphisms then we call $f,g$ inverse equivalences and we say that $A$
and $B$ are equivalent.

In a 2-category it is natural not only to consider ordinary
commutative diagrams (so-called ``strict'' commutative diagrams) but
also pseudo-commutative diagrams. These  diagrams are commutative
up to \emph{explicit} isomorphism. For example the notation
\begin{equation}
\label{ref:A.1a}
\begin{picture}(40,40)
\put(2,20){\makebox(0,0){$h$}}
\put(20,2){\makebox(0,0){$j$}}
\put(38,20){\makebox(0,0){$f$}}
\put(20,38){\makebox(0,0){$g$}}
\put(20,20){\makebox(0,0){$\nu$}}
\put(10,5){\vector(1,0){20}}
\put(5,30){\vector(0,-1){20}}
\put(10,35){\vector(1,0){20}}
\put(35,30){\vector(0,-1){20}}
\put(5,5){\circle*{2}}
\put(35,5){\circle*{2}}
\put(5,35){\circle*{2}}
\put(35,35){\circle*{2}}
\put(12.5,32.5){\vector(-1,-1){5}}
\end{picture}
\end{equation}
means that there is an isomorphism $fg\xrightarrow{\nu}jh$. Often $\nu$
is  clear from the context. In such a case we
sometimes tacitly ignore $\nu$ and treat \eqref{ref:A.1a} as a real
commutative diagram.

The naturality of such pseudo-commutative diagrams is reflected in the
definition of a pseudo-functor \cite{KS} between 2-categories, which
we give below. Assume that $\Cscr$, $\Dscr$ are 2-categories. A
pseudo-functor $T:\Cscr\r \Dscr$ associates to every object of $\Cscr$
an object of $\Dscr$, to every arrow $f:A\r B$ of $\Cscr$ an arrow
$T(f):T(A)\r T(B)$ of $\Dscr$ and to every 2-cell $\nu:f\r g$ a two-cell
$T(\nu):T(f)\r T(g)$. If $T$ were an ordinary functor then we would
require that for compositions of arrows $fg$ one has $T(fg)=T(f)T(g)$.
However for a pseudo-functor we only require the existence of
isomorphisms $\eta_{f,g}:T(f)T(g)\r T(fg)$  which we consider as being part of
the description of $T$.  The data describing $T$ has to satisfy a
list of compatibilities which may be summarized by saying that every
diagram that can commute must commute.

One can go on and define natural transformations between
pseudo-functors and even natural transformations between natural
transformations (``modifications'').  In this way pseudo-functors
between 2-categories form themselves a 2-category and the category of
all 2-categories is a 3-category!

We will call a pseudo-functor $S:\Cscr\r \Dscr$ an \emph{equivalence} if
for every object $D$ in $\Dscr$ there exists an object $C$ in
  $\Cscr$ such that $D$ is equivalent to $S(C)$ (essential
  surjectivity) and for all objects $A,B$ in $\Cscr$ the canonical map
\[
\Hom_\Cscr(A,B)\r \Hom_\Dscr(SA,SB)
\]
is an equivalence of categories.
As usual $\Cscr$ and $\Dscr$ are said to be equivalent if there exists
an equivalence $S:\Cscr\r \Dscr$. One verifies that such an $S$ has a
quasi-inverse (in an appropriate sense) and hence ``equivalence of two-categories
is symmetric''.

An example of a pseudo-functor between 2-categories is given by
adjunction. Assume that $\Cscr$ is a 2-category in which every arrow $f$
possesses a right adjoint $Rf$. Then $R$ defines a
pseudo-functor $R:\Cscr\r \Cscr^{\text{opp}}$. If every arrow $f$ has
also a left adjoint $Lf$ then $R$, $L$ are inverse equivalences of
2-categories between $\Cscr$ and $\Cscr^{\text{opp}}$. 

One more bit of notation. If $\Cscr$ is a 2-category and $A$ is an
object of $\Cscr$ then the relative category $\Cscr/A$ is the
2-category of pairs $(B,f)$ where $B$ is an object of $\Cscr$ and
$f:B\r A$ is an arrow. An arrow $(B,f)\r (C,g)$ in $\Cscr/A$ is given
by an arrow $h:B\r C$ together with an isomorphism $\mu:f\r gh$. A
2-cell $(h,\mu)\r (h',\mu')$ is given by an isomorphism $\nu:h\r h'$
such that $(\Id_g\nu)\cdot \mu=\mu'$.

\def\pretend#1\haswidth#2{\setbox0\hbox{#2}\hbox to \wd0{#1\hss}}
\section{Summary of notations}
\begin{tabbing}
\pretend{\bf
Symbol}\haswidth{$\cohBIMOD(o_X-o_X)$\hskip
0.2cm}\=\pretend{\bf Meaning}\haswidth{locally closed subvarieties
that form a stratification of $PX_B$\hskip 0.2cm}\=\bf
Section\=\\[\smallskipamount] $k$\> an algebraically closed
field\>\strut\hfill\S\ref{ref:1a}\>\\
$\Inj(\Dscr)$\>injective objects in $\Dscr$\>\strut\hfill\S\ref{ref:3.1a}\>\\
$\Lscr(\Dscr,\Cscr)$\>left exact functors from $\Dscr$ to
$\Cscr$\>\strut\hfill\S\ref{ref:3.1a}\>\\ $\BIMOD(\Cscr-\Dscr)$\>the
opposite categorie of
$\Lscr(\Dscr,\Cscr)$\>\strut\hfill\S\ref{ref:3.1a}\>\\
$\Bimod(\Cscr-\Dscr)$\>objects in $\BIMOD(\Cscr-\Dscr)$ that have a
left adjoint\>\strut\hfill\S\ref{ref:3.1a}\>\\ $\otimes$\>composition of
bimodules\>\strut\hfill\S\ref{ref:3.1a}\>\\ $\MOD(\Cscr)$\>the category
$\BIMOD(\Ab-\Cscr)$\>\strut\hfill\S\ref{ref:3.1a}\>\\
$\HHom_\Cscr(\Mscr,-)$\> the left exact functor represented by the
bimodule $\Mscr$\>\strut\hfill\S\ref{ref:3.1a}\>\\ $\HExt$\>the derived functor
of $\HHom$\>\strut\hfill\S\ref{ref:3.1a}\>\\ $\HTor$\>a kind of derived
functor of ``$\otimes$'' \>\strut\hfill\S\ref{ref:3.1a}\>\\
$\ALG(\Dscr)$\>the algebra objects in
$\BIMOD(\Dscr-\Dscr)$\>\strut\hfill\S\ref{ref:3.1a}\>\\
$\Alg(\Dscr)$\>the algebra objects in
$\Bimod(\Dscr-\Dscr)$\>\strut\hfill\S\ref{ref:3.1a}\>\\
$\Mod(\Ascr)$\>the module category of the algebra
$\Ascr$\>\strut\hfill\S\ref{ref:3.1a}\>\\ ''$\invlim$''\>the virtual
inverse limit\>\strut\hfill\S\ref{ref:3.1a}\>\\
$\BIGR(\Cscr-\Dscr)$\>graded ``weak'' $\Cscr$-$\Dscr$ bimodules
\>\strut\hfill\S\ref{ref:3.2b}\>\\ $\Bigr(\Cscr-\Dscr)$\>graded
$\Cscr$-$\Dscr$ bimodules\>\strut\hfill\S\ref{ref:3.2b}\>\\
$\GRALG(\Dscr)$\>the algebra objects in
$\BIGR(\Dscr-\Dscr)$\>\strut\hfill\S\ref{ref:3.2b}\>\\ $\Gralg$\>the
algebra objects in
$\Bigr(\Dscr-\Dscr)$\>\strut\hfill\S\ref{ref:3.2b}\>\\
$\Gr(\Ascr)$\>the graded modules over the algebra
$\Ascr$\>\strut\hfill\S\ref{ref:3.2b}\>\\
$\Sscr(\Ascr)$\>a certain Serre subcategory of $\Gr(\Ascr)$
\>\strut\hfill\S\ref{ref:3.3b}\>\\
$i_\ast$, $i^\ast$, $i^!$\>functors associated to an inclusion $i$ of
quasi-schemes\>\strut\hfill\S\ref{ref:3.4b}\>\\
$\Mod(X)$\>the category of objects associated to a quasi-scheme $X$\>
\strut\hfill\S\ref{ref:3.6b}\>\\
$\Oscr_X$\>a distinguished object for an ``enriched'' quasi-scheme.
\>\strut\hfill\S\ref{ref:3.6b}\>\\
$\Gamma(X,-)$\>notation for the functor $\Hom_X(\Oscr_X,-)$
\>\strut\hfill\S\ref{ref:3.6b}\>\\ 
$o_X$\>the bimodule on $X$ represented by the identity functor\>
\strut\hfill\S\ref{ref:3.6b}\>\\
$\Bimod(X)$\>notation for
$\Bimod(\Mod(X)-\Mod(X))$\>\strut\hfill\S\ref{ref:3.6b}\>\\ 
$\Alg(X)$\>notation for $\Alg(\Mod(X)$\>\strut\hfill\S\ref{ref:3.6b}\>\\
$\Sch$\>the category of quasi-compact, quasi-separated
schemes\>\strut\hfill\S\ref{ref:3.6b}\>\\ 
$\QSch$\>the category of quasi-schemes\>\strut\hfill\S\ref{ref:3.6b}\>\\ 
$\QSch/X$\>the category of quasi-schemes over $X$\>\strut\hfill\S\ref{ref:3.6b}\>\\ 
$\Spec \Ascr$\>a quasi-scheme with module category $\Mod(\Ascr)$\>\strut\hfill\S\ref{ref:3.6b}\>\\ 
$o_X(-Y)$\>a subbimodule of $o_X$ associated to a divisor
$Y\subset X$\>\strut\hfill\S\ref{ref:3.7b}\>\\
$o_X(nY)$\>notation for $o_X(-Y)^{\otimes
-n}$\>\strut\hfill\S\ref{ref:3.7b}\>\\
$\Mscr(nY)$\>notation for $\Mscr\otimes o_X(nY)$\>\strut\hfill\S\ref{ref:3.7b}\>\\
$\Nscr_{Y/X}$\>the ``normal bundle'' of $Y$ in $X$\>\strut\hfill\S\ref{ref:3.7b}\>\\
$\Tors_Y(X)$, 
$\Iso_Y(X)$\>certain categories associated to  $Y\subset X$\>
\strut\hfill\S\ref{ref:3.7b}\>\\ 
$\Tors(\Ascr)$\>torsion modules over a graded algebra $\Ascr$\>\strut\hfill\S\ref{ref:3.8b}\>\\ 
$\QGr(\Ascr)$\>the category
$\Gr(\Ascr)/\Tors(\Ascr)$\>\strut\hfill\S\ref{ref:3.8b}\>\\ 
$\tau$\>the ``torsion functor''\>\strut\hfill\S\ref{ref:3.8b}\>\\ 
$\pi$\>the quotient functor $\Gr(\Ascr)\r\QGr(\Ascr)$\>\strut\hfill\S\ref{ref:3.8b}\>\\
$\omega$\>the right adjoint to $\pi$\>\strut\hfill\S\ref{ref:3.8b}\>\\
$(\tilde{-})$\>the composition $\omega\pi$\>\strut\hfill\S\ref{ref:3.8b}\>\\
$\Proj \Ascr$\>A quasi-scheme whose category is $\QGr(\Ascr)$\>\strut\hfill\S\ref{ref:3.8b}\>\\
$\Pqsch/X$\>``projective'' quasi-schemes over $X$\>\strut\hfill\S\ref{ref:3.8b}\>\\ 
$Q\Sscr(\Ascr)$\>the image of $\Sscr(\Ascr)$ in $\QGr(\Ascr)$\>\strut\hfill\S\ref{ref:3.11b}\>\\ 
$\alpha^{-1}(\Sscr)$\>an alternative notation for
$Q\Sscr(\Ascr)$\>\strut\hfill\S\ref{ref:3.11b}\>\\  
$\Ascr^{(n)}$\>the $n$'th Veronese of $\Ascr$\>\strut\hfill\S\ref{ref:3.12b}\>\\ 
$L_i\alpha^\ast$\>a kind of derived functor to $\alpha^\ast$\>\strut\hfill\S\ref{ref:3.10b}\>\\ 
$\PC(A)$\>the category of pseudo-compact $A$-modules\>\strut\hfill\S\ref{ref:4a}\>\\ 
$\Top(A)$\>the category of topological $A$-modules\>\strut\hfill\S\ref{ref:4a}\>\\ 
$\Dis(A)$\>the category of discrete $A$-modules\>\strut\hfill\S\ref{ref:4a}\>\\ 
$\PCFin(A)$\>the category of pseudo-compact finite length modules\>\strut\hfill\S\ref{ref:4a}\>\\ 
$\PC(A-B)$\>the category of pseudo-compact
$A$-$B$-bimodules\>\strut\hfill\S\ref{ref:4a}\>\\ 
$\ctimes$\>completed tensor product\>\strut\hfill\S\ref{ref:4a}\>\\ 
$X$\>usually a fixed quasi-scheme\>\strut\hfill\S\ref{ref:5.1a}\>\\ 
$Y$\>usually a fixed commutative curve which is a divisor in
$X$\>\strut\hfill\S\ref{ref:5.1a}\>\\ 
$p$\>usually a fixed point on $Y$
\>\strut\hfill\S\ref{ref:5.1a}\>\\
$\tau$\>an automorphism of $Y$ associated to
$\Nscr_{Y/X}$\>\strut\hfill\S\ref{ref:5.1a}\>\\ 
$\Nscr_\tau$\>the twisting of the line bundle $\Nscr$ by the
automorphism $\tau$
\>\strut\hfill\S\ref{ref:5.1a}\>\\
$O_\tau(p)$\>the $\tau$-orbit of $p$
\>\strut\hfill\S\ref{ref:5.1a}\>\\
$\Oscr_p$\>the object in $\Mod(X)$ corresponding to $p\in
Y$\>\strut\hfill\S\ref{ref:5.1a}\>\\
$\Cscr_f$\>finite length objects supported on $Y$ \>\strut\hfill\S\ref{ref:5.1a}\>\\
$\Cscr$\>direct limits of objects in $\Cscr_f$
\>\strut\hfill\S\ref{ref:5.1a}\>\\
$\Cscr_p$\>objects in $\Cscr$ supported on the $\tau$-orbit of $p$
\>\strut\hfill\S\ref{ref:5.1a}\>\\
$\Cscr_{f,p}$\>the finite length objects in $\Cscr_p$
\>\strut\hfill\S\ref{ref:5.1a}\>\\
$C_p$\>the pseudo-compact ring associated to $\Cscr_p$
\>\strut\hfill\S\ref{ref:5.1a}\>\\
$N$\>a canonical normal element in $C_p$
\>\strut\hfill\S\ref{ref:5.1a}\>\\
$\hat{(-)}_p$\>the completion functor
\>\strut\hfill\S\ref{ref:5.1a},\S\ref{ref:5.3a},\S\ref{ref:5.4b}\>\\
$R$\>the completion of $\Oscr_Y$ at
$p$\>\strut\hfill\S\ref{ref:5.1a}\>\\ 
$m$\>the maximal ideal of $R$
\>\strut\hfill\S\ref{ref:5.1a}\>\\
$U$\>the generator of the maximal ideal of $R$
\>\strut\hfill\S\ref{ref:5.1a}\>\\
$m_i$\>the $i$'th maximal ideal of $C_p$
\>\strut\hfill\S\ref{ref:5.1a}\>\\
$S_i$\>the $i$'th simple $C_p$-module
\>\strut\hfill\S\ref{ref:5.1a}\>\\
$e_i$\>the $i$'th diagonal primitive idempotent in $C_p$
\>\strut\hfill\S\ref{ref:5.1a}\>\\
$P_i$\>the $i$'th pseudo-compact projective of $C_p$
\>\strut\hfill\S\ref{ref:5.1a}\>\\
$o_q$\>the bimodule on $X$ corresponding to $q\in Y$
\>\strut\hfill\S\ref{ref:5.1a}\>\\
$\cohBIMOD(o_X-o_X)$\>coherent bimodules on $X$
\>\strut\hfill\S\ref{ref:5.5b}\>\\
$\tilde{\Cscr}_{f,p}$\>finite extensions of $o_{\tau^i p}$
\>\strut\hfill\S\ref{ref:5.5b}\>\\
$\trans_Y(X)$\>objects ``transversal'' to $Y$
\>\strut\hfill\S\ref{ref:5.5b}\>\\
$\Div(\Fscr)$\>the divisor on $Y$ associated to an object in $\mod(X)$
\>\strut\hfill\S\ref{ref:5.5b}\>\\
$\Fscr_Y$\>notation for $\Fscr/\Fscr(-Y)$
\>\strut\hfill\S\ref{ref:5.7b}\>\\
$T_p(\Fscr)$\>an invariant for $\Fscr\in\trans_Y(X)$
\>\strut\hfill\S\ref{ref:5.7b}\>\\
$N_p(\Fscr)$\>the $p$-normalization of $\Fscr$
\>\strut\hfill\S\ref{ref:5.7b}\>\\
$m_q$\>the ideal in $o_X$ associated to $o_q$
\>\strut\hfill\S\ref{ref:6.1a}\>\\
$m_{Y,q}$\>the ideal in $o_Y$ associated to $o_q$
\>\strut\hfill\S\ref{ref:6.1a}\>\\
$I$\>notation for $m_p(Y)$
\>\strut\hfill\S\ref{ref:6.1a}\>\\
$I_Y$\>notation for $m_{Y,p}(Y)$
\>\strut\hfill\S\ref{ref:6.1a}\>\\
$\mu$\>the multiplicity of $p$ on $Y$
\>\strut\hfill\S\ref{ref:6.1a}\>\\
$\Dscr$\>the Rees algebra for $I$
\>\strut\hfill\S\ref{ref:6.2b}\>\\
$\Dscr_Y$\>the Rees algebra for $I_Y$
\>\strut\hfill\S\ref{ref:6.2b}\>\\
$\tilde{X}$\>the blowup of $X$ in $p$
\>\strut\hfill\S\ref{ref:6.3b}\>\\
$\tilde{Y}$\>the strict transform of $Y$ in $\tilde{X}$
\>\strut\hfill\S\ref{ref:6.3b}\>\\ 
$\alpha$, $\beta$, $i$, $j$\>maps in a diagram relating
$X$,$\tilde{X}$, $Y$,$\tilde{Y}$ 
\>\strut\hfill\S\ref{ref:6.3b}\>\\ 
$\tau'$\>the lifting of $\tau$ to $\tilde{Y}$
\>\strut\hfill\S\ref{ref:6.3b}\>\\ 
$\alpha^{-1}(\Cscr_p)$\>analog for the objects supported on the
exceptional curve
\>\strut\hfill\S\ref{ref:6.5b}\>\\ 
$L$\>the exceptional curve 
\>\strut\hfill\S\ref{ref:6.6b}\>\\ 
$\Mod(L)$\>the objects on $L$
\>\strut\hfill\S\ref{ref:6.6b}\>\\ 
$\Oscr_L$\>the ``structure sheaf'' on $L$
\>\strut\hfill\S\ref{ref:6.6b}\>\\ 
$\alpha^{-1}$\>the non-normalized strict transform
\>\strut\hfill\S\ref{ref:6.8b}\>\\
$\alpha_s^{-1}$\>the normalized strict transform
\>\strut\hfill\S\ref{ref:6.8b}\>\\
$l_p(\Fscr)$\>the Loewy length of $T_p(\Fscr)$
\>\strut\hfill\S\ref{ref:6.8b}\>\\
$M_z(X)$\>A certain subcategory of $\trans_Y(X)/\Cscr_f$
\>\strut\hfill\S\ref{ref:6.9b}\>\\
$D^\ast(Z)$, $D^\ast_f(Z)$\>derived categories of a quasi-scheme $Z$
\>\strut\hfill\S\ref{ref:10a}\>\\
$(Y,\Lscr,\sigma)$\>a ``triple'' in the sense of \cite{ATV1}
\>\strut\hfill\S\ref{ref:10a}\>\\
$A$\>the Artin-Schelter regular algebra associated to
$(Y,\Lscr,\sigma)$
\>\strut\hfill\S\ref{ref:10a}\>\\
$g$\>the canonical central element in $A$
\>\strut\hfill\S\ref{ref:10a}\>\\
$o_X(1)$\>the invertible bimodule corresponding to the shift on $\Gr(A)$
\>\strut\hfill\S\ref{ref:10a}\>\\
$r_{p,n}(\Mscr)$\>the multiplicities of $\Mscr$ at points infinitely
near to $p$
\>\strut\hfill\S\ref{ref:10a}\>\\
$e(\Mscr)$\>the multiplicity of $\Mscr$
\>\strut\hfill\S\ref{ref:10a}\>\\
$\overline{\Div}(\Mscr)$\>a variant on $\Div(\Mscr)$
\>\strut\hfill\S\ref{ref:10a}\>
\end{tabbing}

\ifx\undefined\bysame
\newcommand{\bysame}{\leavevmode\hbox to3em{\hrulefill}\,}
\fi

\end{document}